\documentclass[11pt]{article}

\usepackage[margin=1.3in]{geometry}
\usepackage{mathtools}
\usepackage{graphics}
\usepackage{epsfig}
\usepackage{multirow}
\usepackage{array}
\usepackage{extarrows}
\usepackage{graphicx}
\usepackage{xcolor}
\mathtoolsset{showonlyrefs}
\usepackage{mathabx}

\usepackage{amssymb}

\usepackage[symbol]{footmisc}

\usepackage{comment}
\usepackage[titletoc,title]{appendix}
\numberwithin{equation}{section}
\usepackage{amsthm}
\usepackage{enumerate}
\usepackage{amsmath}

\usepackage{latexsym}
\usepackage{xcolor}
\usepackage{pgf}
\usepackage[numbers]{natbib}
\usepackage{algorithm, algpseudocode}
\usepackage{hyperref}
\hypersetup{
    colorlinks,%
    citecolor=black,%
    filecolor=black,%
    linkcolor=black,%
    urlcolor=black
}
\usepackage{subfig}
\numberwithin{equation}{section}

\newtheorem{theorem}{ \noindent T{\footnotesize HEOREM}}
\newtheorem{prop}{ \noindent P{\footnotesize ROPOSITION}}[section]
\newtheorem{lemma}{ \noindent L{\footnotesize EMMA}}[section]
\newtheorem{coro}{ \noindent C{\footnotesize OROLLARY}}
\newtheorem{assumption}{ \noindent A{\footnotesize SSUMPTION}}[section]
\newtheorem{remark}{ \noindent R{\footnotesize EMARK}}[section]

\newcommand{\E}{\mathbb{E}}
\newcommand{\Prob}{\mathbb{P}}

\newcommand{\R}{\mathbb{R}}
\newcommand{\C}{\mathbb{C}}

\renewcommand\Re{\operatorname{Re}}
\renewcommand\Im{\operatorname{Im}}

\newcommand{\T}{\top}

\newcommand{\diag}{\mathrm{diag}}
\newcommand{\osc}{\mathrm{osc}} 
\newcommand{\off}{\mathrm{Off}} 
\newcommand{\db}{\mathrm{db}} 
\newcommand{\orc}{\mathrm{orc}} 

\newcommand{\GBias}{\operatorname{GBias}}
\newcommand{\blockO}[2]{\mathbb O(#1)\oplus\mathbb O(#2)}

\DeclareMathOperator{\Span}{Span}

\DeclareMathOperator{\orth}{orth}

\title{Geometric bias in eigenspace perturbation under random heterogeneous noise}
\author{Fengkai Liu\thanks{Department of Mathematics, Hong Kong University of Science and Technology. Email: fliuar@connect.ust.hk} 
\and Ke Wang\thanks{Department of Mathematics, Hong Kong University of Science and Technology. Email: kewang@ust.hk}
\and Wanjie Wang\thanks{Department of Statistics and Data Science, National University of Singapore. Email: wanjie.wang@nus.edu.sg}
}
\date{\today}

\begin{document}

\maketitle

\begin{abstract}
    Spectral methods rely fundamentally on the stability of principal eigenspaces under random perturbations. Classically, this stability is quantified by the Davis-Kahan and Wedin theorems, which bound the eigenspace error using the operator norm of the noise and the relevant spectral gaps. While these worst-case bounds are sharp for arbitrary deterministic perturbations, they can be wasteful in the low-rank signal-plus-random-noise setting, as they fail to capture the fine-grained interaction between the signal geometry and the noise distribution. In this paper, we study the spectral perturbation of signal-plus-noise matrices corrupted by sparse, random noise with an arbitrary, inhomogeneous variance profile. We demonstrate that under heterogeneous noise variances, the empirical eigenvectors suffer a systematic, deterministic geometric bias that is entirely invisible to classical perturbation bounds. By leveraging the Quadratic Vector Equation (QVE) and establishing fine-grained isotropic local laws, we derive near-optimal, non-asymptotic perturbation bounds for the leading eigenspaces in the operator and $2\to\infty$ norms.  The bounds separate the usual signal-to-noise contribution, stochastic fluctuations, and structured geometric bias terms determined by the alignment between the signal eigenspaces and the row-wise variance profile. We further develop refined rowwise bounds that adapt to the variance-weighted leverage of the signal space, leading to sharper guarantees in delocalized regimes.

As applications, we establish strong consistency
of adjacency spectral clustering for degree-corrected stochastic block models with heterogeneous degree parameters and unbalanced
communities, recovering the logarithmic expected-degree scale in the regular balanced case. We also study spectral embedding for generalized
random dot product graphs, showing that the full signal embedding admits sharp rowwise control, whereas spectral truncation can retain a
systematic geometric bias determined by the omitted signal directions and the heterogeneous variance profile.
\end{abstract}

\tableofcontents

\section{Introduction and main results}
\label{sec:intro}

Spectral methods provide a simple and powerful way to extract low-dimensional
structure from high-dimensional data.  They underlie principal component
methods, low-rank matrix denoising and approximation, spectral clustering and
community detection, latent-position estimation in networks, and many other
procedures in modern statistics and data science.  In these problems, the
statistical object of interest is often encoded by a small number of dominant
eigenvalues and eigenvectors of a population matrix, while the available data
produce only a noisy version of that matrix.  Understanding the stability of
the leading spectral structure is therefore a basic step in understanding the
performance of the resulting statistical procedure.

We study the signal-plus-noise model
\begin{equation}
    \widetilde A=A+E,
    \label{eq:model}
\end{equation}
where $A\in\R^{n\times n}$ is a deterministic symmetric matrix of low rank and
$E$ is a symmetric random perturbation with independent centered entries up to
symmetry.  Our main tasks are to understand the locations of the empirical
outlier eigenvalues and the perturbation of the corresponding signal
eigenspaces, both in the usual operator-norm geometry and in the
$2\to\infty$ norm that measures the largest rowwise error.  We also study a
closely related object, the weighted spectral embeddings,  that arises naturally in applications.

Classical perturbation results such as Weyl's inequality and the
Davis--Kahan--Wedin $\sin\Theta$ theorem control these quantities through the
operator norm of the perturbation and the relevant spectral gaps.  Such bounds
are sharp for arbitrary deterministic perturbations, but they can be
conservative when the perturbation is random.  A substantial literature has
therefore developed stochastic perturbation bounds that exploit independence,
concentration, leave-one-out arguments, and resolvent methods; see
Section~\ref{sec:related} for a detailed discussion.  The present paper
addresses a feature that becomes essential when the noise is heterogeneous:
the variance profile does not merely change the size of the random
fluctuations, but can also change the geometry of the population signal space.

More specifically, heterogeneous noise produces a systematic
second-order effect governed by its row-variance profile.
As a consequence, two population signal directions that are orthogonal in the
usual Euclidean geometry need not remain orthogonal in the variance-weighted
geometry induced by the perturbation.  This creates a systematic mixing of
signal directions, which we call the \emph{geometric bias}.  The phenomenon is
invisible to a perturbation bound based only on $\|E\|$.  It also disappears in
the homogeneous case, where all row variances are equal, but can be of leading
order for a genuinely inhomogeneous variance profile.

Our first contribution is a non-asymptotic perturbation theory that
retains the geometry induced by the heterogeneous variance profile.
Using the Quadratic Vector Equation (QVE) together with isotropic and
compressed local laws, we separate the eigenspace error into a
deterministic geometric bias, stochastic mixing within the signal
space, and leakage into the population null space.  For truncated
eigenspaces, the omitted signal directions interact with the variance
profile and can produce a systematic error amplified by a small
within-signal spectral gap.

Our second contribution is a row-adaptive refinement of the
$2\to\infty$ bound. A leave-one-out local-law argument replaces the
global null-space fluctuation scale by a variance-weighted leverage
scale tailored to the signal eigenvectors, while leaving the
signal-space mixing terms unchanged. In the full signal eigenspace,
the geometric-bias and within-signal mixing terms disappear. The same
refinement also yields sharper bounds for weighted spectral embeddings,
which are used in our network applications.

We illustrate the two sides of the theory through two statistical
applications: for the degree-corrected stochastic block model, where the full
signal eigenspace is retained, the row-adaptive bound yields strong consistency
under heterogeneous degree parameters and unbalanced communities, reducing to the logarithmic expected-degree scale in the regular balanced case; for
generalized random dot product graphs, we instead focus on the choice of
embedding dimension and show that truncation introduces the geometric bias
predicted by our theory through the interaction of the omitted signal
directions, signed spectral gaps, and heterogeneous variance profile.

\bigskip
\noindent\textbf{Roadmap.}
The remainder of this section states the precise top-$k$ eigenspace
perturbation theorem and its row-adaptive refinements and consequences.
Section~\ref{sec:application} develops the two statistical applications.
Section~\ref{sec:geometric-bias} explains the geometric-bias mechanism through
a rank-two model and the QVE expansion, and includes numerical experiments and
an oracle debiasing procedure.  Section~\ref{sec:related} discusses the related
literature.  The technical arguments are developed in
Sections~\ref{sec:qve-preliminaries}--\ref{sec:outlier-ev}: we first derive
the QVE expansion responsible for the variance-profile bias, then establish
the isotropic, compressed, and row-adaptive local-law inputs and the outlier
eigenvalue locations, and finally combine these estimates with blockwise
eigenvector equations.  Additional proofs and extensions are collected in the
appendices.

\bigskip
\noindent\textbf{Notation.}
For a vector $v = (v_1, \dots, v_n)^\T \in \mathbb{R}^n$, we let $\|v\|$ denote its Euclidean $\ell_2$ norm and $\|v\|_\infty := \max_i |v_i|$ denote its maximum absolute entry. The oscillation of a vector $v$ is denoted by $\osc(v) := \max_i v_i - \min_i v_i$. For a matrix $B$, we define its off-diagonal part as $\off(B) := B - \diag(B)$, where $\diag(B)$ is the diagonal matrix formed by the diagonal entries of $B$. We use various matrix norms throughout the paper. Let $\|B\|$ denote the operator norm, and let $\|B\|_F$ denote the Frobenius norm. The $2\to\infty$ norm of a matrix is defined as $\|B\|_{2,\infty} := \max_i \|e_i^\T B\|$, which is the maximum Euclidean norm of the rows of $B$. The infinity operator norm is denoted by $\|B\|_{\infty \to \infty} := \max_i \sum_j |B_{ij}|$. We let $\mathbb{O}(k)$ denote the set of $k \times k$ orthogonal matrices. For two non-negative sequences $X_n$ and $Y_n$, we write $X_n \lesssim Y_n$ (or $Y_n \gtrsim X_n$) if there exists an absolute constant $c > 0$ such that $X_n \le c Y_n$ for all sufficiently large $n$. We write $X_n \asymp Y_n$ if both $X_n \lesssim Y_n$ and $X_n \gtrsim Y_n$ hold. Throughout the paper, $C, C'$ denote positive absolute constants whose exact values may change from line to line.

\subsection{Top-\texorpdfstring{$k$}{k} eigenspace perturbation}
\label{subsec:main-top-k}

We now introduce the precise model notation and assumptions.  Write the
rank-$r$ eigendecomposition of the signal as
\[
    A=U\Lambda U^\T=\sum_{i=1}^r\lambda_i u_i u_i^\T,
\]
and order the nonzero eigenvalues algebraically.  Let
\[
    r_+:=\#\{j:\lambda_j>0\},  \qquad r_-:=r-r_+,
\]
so that
\[
    \lambda_1\ge\cdots\ge\lambda_{r_+}>0>
    \lambda_{r_++1}\ge\cdots\ge\lambda_r.
\]
The ordered eigenvalues and their corresponding eigenvectors of $\widetilde{A}$ are denoted as $\widetilde{\lambda}_i,\widetilde{u}_i$ for $1\le i \le n$. 
In the main text we focus on the right-edge outliers and the top-$k$ ($1\le k\le r_+$) positive signal eigenspace and its perturbed counterpart
\begin{equation}
    U_k:=(u_1,\ldots,u_k),
    \qquad
    \widetilde U_k:=(\widetilde u_1,\ldots,\widetilde u_k).
    \label{eq:Uk}
\end{equation}
We write
$\Lambda_k:=\diag(\lambda_1,\ldots,\lambda_k)$ and define
$\widetilde\Lambda_k$ analogously.  For an index set
$J\subseteq[r]$, $U_J$ denotes the submatrix of $U$ formed by the columns
indexed by $J$.  The negative-spike case follows by applying the same argument
to $-A$; see Remark~\ref{rem:general-space}.

Throughout the main results, fix $D>0$.  We impose the following assumptions on random noise $E$: 
\begin{assumption}[Structural assumption on the noise]\label{assump:noise}Let $E$ be an $n\times n$ symmetric random matrix with independent centered entries (up to symmetry) and variance profile $\Sigma = (\sigma_{ij}^2)_{1\le i,j\le n}$ where $\mathbb{E}(E_{ij}^2) = \sigma_{ij}^2.$ Assume there exists a parameter $K\ge 1$ such that the sub-Gaussian norm\footnote{The sub-Gaussian norm of a random variable $X$ is defined as $\|X\|_{\psi_2} :=\inf\{t>0: \E(\exp(X^2/t^2))\le 2 \}.$ In particular, $\Prob(|X|\ge t) \le 2 \exp(-t^2/\|X\|_{\psi_2}^2)$. See \cite{Vbook} for details.} of each entry satisfies:
$$\|E_{ij}\|_{\psi_2}\le  K \sigma_{ij} .$$
Equivalently, we can write $E_{ij}=\sigma_{ij}\xi_{ij}$ where the random variables $\xi_{ij}$ have mean 0, variance 1, and satisfy
$ \max_{i,j} \|\xi_{ij}\|_{\psi_2} \le K.$

Let $\sigma_{\max}:=\max_{i,j} \sigma_{ij}>0$. We further assume that there exists a parameter $\beta\ge 1$ such that
\begin{equation}\label{eq:condition-3rd4th}
    \E |E_{ij}|^3\le
\beta K\sigma_{\max}\,\sigma_{ij}^2,\qquad \E |E_{ij}|^4 \le \beta^2 K^2\sigma_{\max}^2\,\sigma_{ij}^2.
\end{equation}
\end{assumption}
\begin{remark}[The role of $\beta$]
\label{rem:beta}
The parameter $\beta$ separates the third- and fourth-moment controls
from the sub-Gaussian norm $K$ and is useful in sparse regimes.
Under the sub-Gaussian assumption alone, one may always take
$ \beta\lesssim \sqrt{\log(eK)},$ 
but this can introduce an artificial logarithmic loss. For example,
for centered Bernoulli noise with edge probability $p$, one has
$K\lesssim p^{-1/2}$, whereas the third and fourth moments satisfy
the sharper bounds in \eqref{eq:condition-3rd4th} with $\beta=O(1)$.
Thus keeping $\beta$ explicit allows our results to retain the correct
logarithmic scale in sparse network applications.
\end{remark}

\medskip
Let us denote $R_i:=\sum_{j=1}^n \sigma_{ij}^2$. Define the diagonal matrix \[\mathcal R :=\diag(R_1,\dots,R_n)\]and the \emph{geometric bias matrix}
\begin{align}\label{def:V}
    \mathcal{V}:= U^\T \mathcal{R} U = (u_i^\T \mathcal{R} u_j)_{1\le i,j \le r}.
\end{align}
To capture the heterogeneity of the noise profile, we define the \emph{variance oscillation} parameter:
\begin{equation*}
    \osc(\mathcal R)=R_{\max}-R_{\min},
\end{equation*}
where 
\[
R_{\max}:=\max_i R_i\qquad \text{and}\qquad R_{\min}:=\min_i  R_i.
\] 

Assumption~\ref{assump:noise} is imposed throughout. In addition, we
work under one of the following two auxiliary regimes.

\begin{assumption}[Auxiliary regularity conditions]
\label{assump:extra}
Fix $D>0$. We assume that one of the
following conditions holds.
\begin{enumerate}
    \item[(i)] \emph{General sub-Gaussian regime:}
    \[
        \sqrt{R_{\max}}
        \ge C_0(D+8)K\sigma_{\max}\log n.
    \]

    \item[(ii)] \emph{Bounded-entry refinement:} If
    there exists an absolute constant $c_0\ge1$ such that $\max_{i,j}|E_{ij}| \le c_0K\sigma_{\max}$ almost surely, then we only require
    \[
        \sqrt{R_{\max}}   \ge C_0(D+8)K\sigma_{\max}\sqrt{\log n}.
    \]
\end{enumerate}
Here $C_0>0$ is a sufficiently large absolute constant.
\end{assumption}
The second regime is a refinement of the first noise model rather than
a separate distributional assumption: in both cases
Assumption~\ref{assump:noise} holds, while the  bounded-entry condition in (ii) allows the variance condition to be weakened by a
factor $\sqrt{\log n}$. Note that the bounded-entry condition implies \eqref{eq:condition-3rd4th} with
$\beta\le c_0$.

\begin{remark}[On the Assumption \ref{assump:extra}]\label{rem:auxiliary-assumptions} The regularity condition prevents any single maximal entry from dominating its row sum.   Under this condition, our  estimate (see Lemma \ref{lem:bdE}) yields $\|E\|\le 3\sqrt{R_{\max}}$ with high probability, thereby guaranteeing a  macroscopic signal-to-noise separation ($\lambda_k\ge 3\|E\|$) for any supercritical spike $\lambda_k\ge 9\sqrt{R_{\max}}$. Both assumptions naturally accommodate sparse regimes. For example, for centered Bernoulli noise with parameter $p$,  Assumption~\ref{assump:extra} requires  $p\gtrsim \log n/n$.
\end{remark}

For technical convenience, throughout the main results we impose
\begin{equation}
    \lambda_1\le R_{\max}^3,
    \qquad
    R_{\max}\le n^{100},
    \qquad
    \sigma_{\max}^{-1}\le n^{100}.
    \label{eq:technical-scale-conditions}
\end{equation}
These conditions only restrict the range in which we carry out the
local-law analysis and are not essential to the perturbation
phenomenon. The polynomial bounds are extremely mild; the exponent
$100$ is arbitrary and may be replaced by any sufficiently large fixed
constant. In particular, they are automatic in the network models
considered below, where $R_{\max}\le n$.

The upper bound $\lambda_1\le R_{\max}^3$ merely excludes an
ultra-strong signal regime. When $|z|\gg\|E\|$, the resolvent admits a
direct Neumann-series expansion, and the corresponding perturbation
bounds follow by simpler arguments. We therefore exclude this regime
to avoid a separate treatment; see \cite{Wang24} for related arguments
without such an upper restriction.

\medskip

Define the \emph{spectral gap} 
\begin{align}\label{def:spectral-gap}
    \delta_k : = \begin{cases}
        \lambda_k - \lambda_{k+1}, & \text{ if } k< r_+;\\
        \lambda_k - 3\sqrt{R_{\max}}, & \text{ if } k = r_+.
    \end{cases}
\end{align}
With the convention that $\delta_0:=+\infty$.

We define the regime-dependent fluctuation scale
\begin{equation}
\label{def:M}
M_D:=
\begin{cases}
C'(D+6)^{3/2}
r\beta K\sigma_{\max}(\log n)^2,
& \text{in the sparse sub-Gaussian regime},\\[1mm]
C'(D+6)
r c_0 K\sigma_{\max}\log n,
& \text{in the sparse bounded-entry regime},
\end{cases}
\end{equation}
where $C'>0$ is a sufficiently large absolute constant.

Furthermore, for fixed $k\in [r_+]$, define index sets 
\begin{align*}
    &\mathcal{J}:=\{k+1,\dots,r_{+}\}, \quad \mathcal{I}:=[r]\setminus \mathcal{J},\\
    &\mathcal{N}:=\{r_{+}+1,\dots,r\}, \quad \mathcal{K}:=[r]\setminus \mathcal{N},
\end{align*} 
and denote 
\begin{align}\label{def:bias-k}
\mathcal{B}_k:=\frac{10\sqrt{k}}{\delta_k \lambda_k}\|\mathcal{V}_{\mathcal{J}\mathcal{I}}\| +\frac{4 \sqrt{k}}{\lambda_k^2}\|\mathcal{V}_{\mathcal{N}\mathcal{K}}\| + \frac{80\sqrt{k}}{\delta_k \lambda_k}\cdot \frac{R_{\max}}{\lambda_k^2}\osc(\mathcal R),
\end{align}
where $\mathcal{V}_{\mathcal{J}\mathcal{I}} = U_{\mathcal J}^\T \mathcal{R} {U}_{\mathcal I}$ and $\mathcal{V}_{\mathcal{N}\mathcal{K}} = U_{\mathcal N}^\T \mathcal{R} {U}_{\mathcal K}$. If $k=r_+$ (thus $\mathcal J=\emptyset$), then $\|\mathcal{V}_{\mathcal{J}\mathcal{I}}\|=0$.  

We also set
\[\rho_k\equiv \rho_{k,D} := 10 M_D + 15\frac{\osc(\mathcal{R})}{\lambda_k}.\]

\begin{theorem}[Top-$k$ eigenspace perturbations]\label{thm:top-k-eigenspace} Fix $D>0$. Under Assumptions \ref{assump:noise} and  \ref{assump:extra}, for any $k\in [r_+]$ such that $\lambda_k \ge 9 \sqrt{R_{\max}} + 100\rho_{k,D}$ and 
\begin{align}\label{def:gap-assump-main}
\delta_k \ge \sqrt{k}\left(65 \rho_{k,D} + 20 \frac{ \osc(\mathcal R)}{\lambda_k}\right),
\end{align}
we have the following bounds with probability at least $1-n^{-D}$:

\textup{(i)} For the operator norm bound, we have
\begin{align*}
    \|\sin\angle(U_k, \widetilde{U}_k)\|
    \le \mathcal{B}_k + 10\sqrt{k}\frac{M_D}{\delta_k}+2\frac{\|E\|}{\lambda_k}.
\end{align*}

\textup{(ii)} For the $2\to\infty$ norm bound, we first have 
\begin{align}
   \| \widetilde{U}_k-P_{U_k} \widetilde{U}_k\|_{2,\infty} \le \|U\|_{2,\infty} \left(\mathcal{B}_k +10\sqrt{k}\frac{M_D}{\delta_k}+3\sqrt{k} \frac{\osc{(\mathcal R)}}{\lambda_k^2} \right) + 3 \sqrt{k} \frac{M_D}{\lambda_k}.
\end{align}  
Consequently,
\begin{align}
\min_{O\in \mathbb{O}(k)} \| \widetilde{U}_k - U_k O\|_{2,\infty} &\le \| \widetilde{U}_k-P_{U_k} \widetilde{U}_k\|_{2,\infty} + \|U_k\|_{2,\infty} \|\sin\angle(U_k, \widetilde{U}_k)\|^2,\\
&\le 4\|U\|_{2,\infty} \left( \mathcal{B}_k +20\sqrt{k}\frac{M_D}{\delta_k}+ 3\sqrt{k} \frac{\osc{(\mathcal R)}}{\lambda_k^2} +27 \frac{R_{\max}}{\lambda_k^2}\right) + 3 \sqrt{k} \frac{M_D}{\lambda_k}.
\end{align}
\end{theorem}

Theorem~\ref{thm:top-k-eigenspace} separates three distinct
mechanisms in the eigenspace perturbation. The first is the
geometric-bias term $\mathcal B_k$, which comes from mixing between
the target eigenspace and the non-target signal directions. Its leading
components are governed by
\[
    \mathcal V_{\mathcal J\mathcal I}
    =U_{\mathcal J}^\T\mathcal R\,U_{\mathcal I},
    \qquad
    \mathcal V_{\mathcal N\mathcal K}
    =U_{\mathcal N}^\T\mathcal R\,U_{\mathcal K},
\]
and quantify how the signal eigenspaces interact with the
row-variance profile $\mathcal R=\diag(R_i)$. These couplings vanish
when $\mathcal R=cI$. Indeed, by orthogonality, \[ \|\mathcal V_{\mathcal J\mathcal I}\| = \|U_{\mathcal J}^\T (\mathcal R-cI)\,U_{\mathcal I}\|\le \frac{1}{2} \osc(\mathcal R)\]
by choosing $c=\frac{1}{2}(R_{\max} + R_{\min})$. Likewise, $\|\mathcal V_{\mathcal N\mathcal K}\| \le \frac{1}{2} \osc(\mathcal R)$.  

The term $M_D/\delta_k$ is the stochastic contribution associated
with mixing between the target eigenspace and the remaining signal
directions. It therefore appears in both the principal-angle and
$2\to\infty$ norm bounds. By contrast, the standalone term $M_D/\lambda_k$
in the $2\to\infty$ norm bound comes from the component of the empirical
eigenspace in the null space of $A$. For the principal-angle loss, this rowwise
null-space leakage can be controlled at the coarser operator-norm scale
$\|E\|/\lambda_k$.

\begin{remark}[Bounds for general eigenspaces]
\label{rem:general-space}
Theorem~\ref{thm:top-k-eigenspace} is stated for the top positive
eigenspace only to simplify notation. Analogous bounds hold for bottom
and intermediate signal eigenspaces by symmetry and elementary
projector decompositions. See Appendix~\ref{app:general-simple-facts}
for details.
\end{remark}

\begin{remark}[Extension to rectangular matrices] Our eigenspace perturbation results extend to singular subspace perturbation for rectangular matrices via the standard Hermitian dilation technique. See Appendix \ref{app:rectangular} for a brief discussion.
\end{remark}

\subsection{Refined perturbation bounds and consequences}
Although Theorem~\ref{thm:top-k-eigenspace} focuses on the
operator-norm principal angle and the $2\to\infty$ loss, its proof
actually controls the empirical eigenspace in both the non-target
signal directions and the null space of $A$. Following the framework
of \cite{Wang24}, these estimates can be extended to other error
metrics, including Frobenius and unitarily invariant norms, linear
forms, max norms, and weighted rowwise losses. We only record the
extensions needed for the applications below.

In our heterogeneous setting, the same decomposition also allows us
to refine one term in the $2\to\infty$ bound of
Theorem~\ref{thm:top-k-eigenspace}. The term $M_D/\delta_k$ comes from
stochastic mixing with the non-target signal directions and is
unchanged by the refinement below. In contrast, the standalone term
$M_D/\lambda_k$ arises from the null-space action and can be sharpened
using the row-adaptive local law from
Section~\ref{subsec:local-law-inputs}, at the cost of a slightly
stronger signal-strength condition. This refinement is important in
the applications, where it allows us to reach the logarithmic-degree
regime.

Before stating the refined bounds, we record a stronger variance condition required by the leave-one-out argument.
\begin{assumption}[Row-adaptive variance condition]
\label{assump:row-adaptive}
For fixed $D>0$, we assume
\begin{equation}
    \sqrt{R_{\max}}\ge
    \begin{cases}
        C_1 K\sigma_{\max}(D+r)\log n,
        &\text{in the sparse sub-Gaussian regime},\\[1mm]
        C_1 K\sigma_{\max}\sqrt{(D+r)\log n},
        &\text{in the sparse bounded-entry regime},
    \end{cases}
    \label{eq:variance-dominance-ra}
\end{equation}
where $C_1>0$ is a sufficiently large absolute constant.
\end{assumption}


We further define the regime-dependent entry scale
\begin{equation}\label{eq:L-ra}
    L_{D,r}^{\rm ra}:=
    \begin{cases}
        K\sigma_{\max} \sqrt{(D+r)\log n},
        &\text{in the sparse sub-Gaussian regime},\\
        K\sigma_{\max},
        &\text{in the sparse bounded-entry regime}.
    \end{cases}
\end{equation}
In the sub-Gaussian regime, $L_{D,r}^{\rm ra}$ is the truncation level used in the leave-one-out argument, whereas in the bounded-entry regime the deterministic entry bound is used
directly.

Denote 
\begin{equation}
    \mathfrak R_D(U)  :=\sqrt{(D+r)\log n}\,\Gamma(U)  +
    \left(  {(D+r)\log n}\, L_{D,r}^{\rm ra}+M_D  \right)\|U\|_{2,\infty},
    \label{eq:signal-row-scale}
\end{equation}
where
\[
    \Gamma(U):= \max_i \Big(\sum_j\sigma_{ij}^2\|e_j^\T U\|^2 \Big)^{1/2}.
\]
\begin{prop}[Refined rowwise top-$k$ perturbation]
\label{prop:refined-top-k}
Fix $D>0$ and $1\le k<r_+$. Suppose
Assumptions~\ref{assump:noise} and \ref{assump:row-adaptive} hold.
Assume 
\begin{align*}
\lambda_k \ge C\sqrt{(D+r)R_{\max}\log n}\quad \text{and} \quad \delta_k \ge C\sqrt{k}\left(M_D +  \frac{ \osc(\mathcal R)}{\lambda_k}\right).
\end{align*} \label{eq:refined-top-k-row}
Then, with probability at least $1-n^{-D}$,
\begin{align}
    \|(I-P_{U_k})\widetilde U_k\|_{2,\infty}
    \le{} \|U\|_{2,\infty} \left(  \mathcal B_k
        +10\sqrt{k}\frac{M_D}{\delta_k}
        +3\sqrt{k}\frac{\osc(\mathcal R)}{\lambda_k^2} \right)+
    C'\sqrt{k}\frac{\mathfrak R_D(U)}{\lambda_k}.
 \end{align}
Consequently,
\begin{align*}
\min_{O\in \mathbb{O}(k)} \| \widetilde{U}_k - U_k O\|_{2,\infty} 
&\le C'\|U\|_{2,\infty} \left( \mathcal B_k
        +\sqrt{k}\frac{M_D}{\delta_k}
        +\sqrt{k}\frac{\osc(\mathcal R)}{\lambda_k^2}
        +\frac{R_{\max}}{\lambda_k^2}  \right)+ C'\sqrt{k}\frac{\mathfrak R_D(U)}{\lambda_k}.
\end{align*}
\end{prop}
Compared with Theorem~\ref{thm:top-k-eigenspace}, the only change is that the 
null-space term of order $M_D/\lambda_k$ is replaced by $\mathfrak R_D(U)/\lambda_k$. The latter can be far smaller
for delocalized signal spaces: unlike the global scale $M_D$,
$\mathfrak R_D(U)$ depends on the variance-weighted row leverage
$\Gamma(U)$, with its remaining terms scaled by $\|U\|_{2,\infty}$. For instance, for a fixed-rank incoherent signal
with $\|U\|_{2,\infty}\asymp n^{-1/2}$ and $R_i\asymp d$, in the
sparse bounded-entry regime,
\[
    \mathfrak R_D(U)
    \lesssim
    \sqrt{\frac{d\log n}{n}}
    +\frac{\log n}{\sqrt n},
\]
whereas $M_D$ is of order $\log n$. Thus the refinement can be much sharper in the delocalized regimes relevant to the network applications below.

Next, we present a clean refinement for the full signal eigenspace, where the signal-space mixing terms disappear. This is a direct consequence of Proposition~\ref{prop:refined-top-k}.
\begin{coro}[Full signal eigenspace]
\label{cor:full-signal}
Suppose $A\succeq0$ has rank $r \ge 1$. Fix $D>0$. Suppose
Assumptions~\ref{assump:noise} and \ref{assump:row-adaptive} hold.
Assume $$\lambda_r \ge C\left(  \sqrt{(D+r)R_{\max}\log n}+M_D \right).$$ 
Then, with probability at least $1-n^{-D}$,
\begin{align*}
    \min_{O\in\mathbb O(r)} \|\widetilde U_r-UO\|_{2,\infty}
    \le  C'\sqrt r\, \frac{\mathfrak R_D(U)}{\lambda_r} + C'\|U\|_{2,\infty}
    \frac{\sqrt r\,\osc(\mathcal R)+R_{\max}}{\lambda_{r}^2}.
\end{align*}
\end{coro}

The same row-adaptive refinement carries over to the weighted spectral embedding given below. For $k\in[r_+]$, define the population and empirical rank-$k$
spectral truncations
\[
A_k:=U_k\Lambda_k U_k^\T,
\qquad
\widetilde A_k:= \widetilde U_k\widetilde\Lambda_k\widetilde U_k^\T.
\]
We write
\begin{align}\label{def:rank-k}
A=A_k+A_{>k}, \qquad A_{>k} := U_{>k}\Lambda_{>k}U_{>k}^\T = \sum_{s=k+1}^r\lambda_su_su_s^\T,
\end{align}
where $U_{>k}:=(u_{k+1},\ldots,u_r)$ and
$\Lambda_{>k}:=\diag(\lambda_{k+1},\ldots,\lambda_r)$, with the
convention $A_{>r}=0$. We also introduce the weighted spectral
embeddings
\[
Y_k:=U_k\Lambda_k^{1/2}, \qquad \widetilde Y_k := \widetilde U_k\widetilde\Lambda_k^{1/2}.
\]
For later use, set
\[
\tau_k := \left(\sum_{s=1}^k\lambda_s\right)^{1/2}
\]
and
\[
\eta_k := \|U_k^\T EU_k\| + \bigl(\|\Lambda_{>k}\|+\|E\|\bigr)\left(\mathcal B_k + 10\sqrt{k}\frac{M_D}{\delta_k} + 2\frac{\|E\|}{\lambda_k} \right)^2.
\]
\begin{prop}[Refined weighted spectral embedding]
\label{prop:refined-weighted-embedding}
Under the assumptions of Proposition~\ref{prop:refined-top-k}, with
probability at least $1-n^{-D}$,
\begin{align}
    \min_{O\in\mathbb O(k)}
    \|\widetilde Y_kO-Y_k\|_{2,\infty}
    \le{}&
    \|U\|_{2,\infty}
    \left(
        \sqrt{\lambda_k}\,\mathcal B_k\mathbf 1_{\{k<r\}}
        +
        10\frac{M_D\tau_k}{\delta_k}\mathbf 1_{\{k<r\}}
        +
        3\sqrt{k}
        \frac{\osc(\mathcal R)}{\lambda_k^{3/2}}
    \right)
    \notag\\
    &+
    C\sqrt{k}
    \frac{\mathfrak R_D(U)}{\sqrt{\lambda_k}}
    +
    \frac{\|U_k\|_{2,\infty}}{\sqrt{\lambda_k}}\eta_k .
    \label{eq:refined-weighted-embedding}
\end{align}
\end{prop}
For the full signal embedding, the corresponding simplification is obtained directly by setting the non-target signal space to be empty; we use this form in the proof of the network application below. The proofs of the above results are given in Appendix~\ref{app:proof-embedding-approx}.

\section{Applications}\label{sec:application}
\subsection{Strong consistency for the degree-corrected stochastic block model} \label{subsec:dcsbm}

The degree-corrected stochastic block model (DCSBM) extends the stochastic block
model with vertex-specific degree parameters \cite{KarrerNewman2011}. Let
$r\ge2$ be the number of communities, let $g:[n]\to[r]$ be the community
assignment, and let $Z\in\{0,1\}^{n\times r}$ be the corresponding membership
matrix. We allow $r=r_n$. 

Let $ \Theta:=\diag(\theta_1,\ldots,\theta_n)$ with $\theta_i>0$ contain the vertex-specific degree parameters, 
and let $B\in\R^{r\times r}$ be a symmetric block-connectivity matrix. Conditional
on the deterministic probability matrix
\begin{equation*}
    P:=\Theta ZBZ^\top\Theta,
\end{equation*}
the observed adjacency matrix
$\widetilde A=(\widetilde A_{ij})_{i,j=1}^n$ satisfies $\widetilde A_{ij}\sim\mathrm{Bernoulli}(P_{ij}) $ for $i\le j$ independently up to symmetry.  Throughout, every community is nonempty, $0\le P_{ij}\le1$, and $B$ has full rank $r$.    
Write
\[
    \widetilde A=P+E,
    \qquad
    P=U\Lambda U^\top,
    \qquad
    E:=\widetilde A-P.
\]

We next introduce the quantities that describe community imbalance and degree
heterogeneity. For $a\in[r]$, set the community size and community weighted mass parameters as
\begin{equation*}
    n_a:=\#\{i:g(i)=a\},
    \qquad
    s_a:=\sum_{i:g(i)=a}\theta_i^2,
    \qquad
    s_-:=\min_{a\in[r]}s_a.
\end{equation*}
Define the diagonal matrix 
\begin{equation*}
    S:=Z^\top \Theta^2 Z=\diag(s_1,\ldots,s_r)
\end{equation*} 
and the population row scales
\begin{equation*}
    \mu_i:=\frac{\theta_i}{\sqrt{s_{g(i)}}},
    \qquad
    \mu_-:=\min_{i\in[n]}\mu_i,
    \qquad
    \mu_+:=\max_{i\in[n]}\mu_i.
\end{equation*}
Let
\begin{equation*}
    d_i:=\sum_{j=1}^nP_{ij},
    \qquad
    d_{\max}:=\max_{i\in[n]}d_i,
    \qquad
    p_{\max}:=\max_{i,j}P_{ij}.
\end{equation*}
In our setting, 
\begin{equation*}
    R_i:=\sum_{j=1}^nP_{ij}(1-P_{ij}),
    \qquad
    R_{\max}:=\max_iR_i,
    \qquad
    \mathcal R:=\diag(R_1,\ldots,R_n),
\end{equation*}
and define the variance-weighted leverage scale
\begin{equation*}
    \Gamma_{\mathrm{dc}}  :=  \max_{i\in[n]}
    \left(  \sum_{j=1}^n  P_{ij}(1-P_{ij})\mu_j^2 \right)^{1/2}.
\end{equation*}
Finally, denote the extreme singular values of $B$ as $\sigma_{\min}(B)=\min_{1\le j\le r}|\lambda_j(B)|$ and $\sigma_{\max}(B)=\|B\|$. 
Note that the nonzero eigenvalues of $P$ coincide with those of
$S^{1/2}BS^{1/2}$. To see this, set $X_0:=\Theta ZS^{-1/2}$. Thus $X_0^\top X_0=I_r$ and $P=X_0(S^{1/2}BS^{1/2})X_0^\top$. Hence, the smallest population signal of $P$ is denoted by 
\begin{equation}
    \tau_n := \min_{1\le j\le r}|\lambda_j(P)|=\sigma_{\min}\!\left(S^{1/2}BS^{1/2}\right) .
    \label{eq:dcsbm-exact-signal}
\end{equation}
In particular, we have
\begin{equation}
\tau_n \ge \sigma_{\min}(B) s_-.
    \label{eq:dcsbm-compare}
\end{equation}
\begin{assumption}[Heterogeneous DCSBM]
\label{assum:DCSBM}
The matrix $B=B_n$ is symmetric, entrywise nonnegative, and nonsingular.  Moreover, the graph is sparse in the sense that $p_{\max}=o(1)$.
\end{assumption}

Assumption~\ref{assum:DCSBM} places no lower bound on the community
proportions and no comparability condition on the degree parameters. It also
allows $\sigma_{\min}(B)$ and $\sigma_{\max}(B)$ to vary with $n$ and $r$. In particular, $\sigma_{\max}(B)$ may grow with $r$. 

Assumption~\ref{assum:DCSBM} allows the number of communities $r_n$ to increase with $n$. As $n$ increases, it is natural to have more communities. However, the discussion about theoretical guarantee of community detection with an increasing $r_n$ is quite limited. On SBM, \cite{banks2016information} provides the average degree requirement for nontrivial community recovery, and \cite{Zhang2015Minimax} establishes a minimax theory for both strong consistency and weak consistency. However, both works do not involve degree heterogeneity and unbalanced community sizes.

\paragraph{Spectral clustering algorithm.} Our goal is to recover the unknown labels $g(1),\ldots,g(n)$ from $\widetilde A$. We use
the same raw-adjacency eigenspace and spherical normalization as in Lei and
Rinaldo~\cite{LeiRinaldo2015}. To accommodate highly unbalanced communities,
Algorithm~\ref{alg:dcsbm} combines this embedding with the farthest-first
traversal of Gonzalez~\cite{Gonzalez1985}, followed by nearest-center
assignment. Unlike a global $r$-means objective, the final step depends only on
geometric separation and does not weight a community in proportion to its size.

\begin{algorithm}
    \caption{Farthest-first spectral clustering for the DCSBM}
    \label{alg:dcsbm}
    \begin{algorithmic}
        \State \textbf{Input:} adjacency matrix
        $\widetilde A\in\R^{n\times n}$ and number of communities $r$.
        \State \textbf{Output:} estimated labels
        $\widetilde{\mathbf g}\in[r]^n$.
        \State \emph{Step 1.} Compute the $r$ eigenvectors of $\widetilde A$
        corresponding to its $r$ eigenvalues with largest absolute values, and
        collect them as $\widetilde U_r\in\R^{n\times r}$.
        \State \emph{Step 2.} Normalize the rows:
        $\widetilde z_i:=e_i^\top\widetilde U_r/\|e_i^\top\widetilde U_r\|$, $i\in[n]$.
        \State \emph{Step 3.} Choose any $i_1\in[n]$. For
        $\ell=2,\ldots,r$, choose $i_\ell\in\operatorname*{argmax}_{i\in[n]}
            \min_{1\le s<\ell}\|\widetilde z_i-\widetilde z_{i_s}\|.$
        \State \emph{Step 4.} Assign each vertex to its nearest selected center:
       $
            \widetilde g(i)\in \operatorname*{argmin}_{1\le\ell\le r}
            \|\widetilde z_i-\widetilde z_{i_\ell}\|.
        $
    \end{algorithmic}
\end{algorithm}

To see why the algorithm works, the population eigenvector matrix satisfies
\[
    e_i^\top U=\mu_i e_{g(i)}^\top V
\]
for some $V\in\mathbb O(r)$; see
\eqref{eq:dcsbm-population-rows} in Appendix~\ref{app:dcsbm-proof}.
Row normalization therefore removes the vertex-specific factor $\mu_i$ and maps
all vertices in community $a$ to the same population center $e_a^\top V$. These
$r$ centers have pairwise distance $\sqrt2$ (see
Lemma~\ref{lem:dcsbm-population-geometry}).

A uniform rowwise perturbation of size at most a sufficiently small multiple of
$\mu_-$ preserves this separation after row normalization. The farthest-first
step then selects one representative from each community, regardless of the
community sizes, and nearest-center assignment recovers every label exactly.

For fixed $D>0$, set
\begin{equation*}
    t_{D,r}:=(D+r)\log n
\end{equation*}
and define the error scale
\begin{equation}
    \varepsilon_{D,n}(\tau_n) :=
    \sqrt r\, \frac{  \sqrt{t_{D,r}}\,\Gamma_{\mathrm{dc}} +t_{D,r}\mu_+ }{\tau_n}
    +  (\sqrt r+1)\mu_+\frac{d_{\max}}{\tau_n^2}.
    \label{eq:dcsbm-intrinsic-error-tau}
\end{equation}
Since $\tau_n \ge \sigma_{\min}(B) s_-$ as in \eqref{eq:dcsbm-compare}, we get $\varepsilon_{D,n}(\tau_n) \le \varepsilon_{D,n}(\sigma_{\min}(B) s_-)$.
\begin{theorem}[Strong consistency for heterogeneous DCSBMs]
\label{thm:dcsbm-general}
Suppose Assumption~\ref{assum:DCSBM} holds. For every fixed $D>0$, there
exist constants $C_{1,D},C_{2,D},c_D>0$, independent of $B$, such that, if
\begin{equation}
    d_{\max}\ge C_{1,D}t_{D,r}
    \qquad\text{and}\qquad
    \tau_n\ge C_{1,D}\sqrt{t_{D,r}d_{\max}},
    \label{eq:dcsbm-general-size}
\end{equation}
then, with probability at least $1-n^{-D}$,
\begin{equation}
    \min_{O\in\mathbb O(r)}  \|\widetilde U_rO-U\|_{2,\infty}
    \le C_{2,D}\cdot \varepsilon_{D,n}(\tau_n).
    \label{eq:dcsbm-general-rowwise}
\end{equation}
If, in addition,
\begin{equation}
    \varepsilon_{D,n}(\tau_n) \le c_D\mu_-,
    \label{eq:dcsbm-general-recovery}
\end{equation}
then Algorithm~\ref{alg:dcsbm} exactly recovers the communities up to a
permutation of their labels.
\end{theorem}
\begin{remark}
Since $\tau_n \ge \sigma_{\min}(B) s_-$, the signal condition in \eqref{eq:dcsbm-general-size} may be replaced by 
\begin{equation*}
    \sigma_{\min}(B) s_- \ge C_{1,D} \sqrt{t_{D,r}d_{\max}},
\end{equation*}
which is stronger but directly verifiable.  Accordingly, the right-hand side of
\eqref{eq:dcsbm-general-rowwise} may be replaced by
$C_{2,D}\cdot \varepsilon_{D,n}(\sigma_{\min}(B) s_-)$, and exact recovery follows whenever $\varepsilon_{D,n}(\sigma_{\min}(B) s_-) \le c_D\mu_-$.
\end{remark}

Theorem~\ref{thm:dcsbm-general} does not impose uniform bounds on the
singular values of $B$, which is useful when the number of communities
$r=r_n$ grows and the spectral scale of $B$ may vary with $r$. The dependence
on $r$ is retained explicitly through $t_{D,r}=(D+r)\log n$ and the quantities
$\tau_n$, $d_{\max}$, and $\Gamma_{\mathrm{dc}}$. Thus the admissible growth of $r$ depends on the accompanying signal and sparsity scales. For example, if $r\asymp\log n$, then $t_{D,r}\asymp(\log n)^2$ and the signal condition becomes
\[
d_{\max} \gtrsim_D \log^2n, \qquad     \tau_n \gtrsim_D \sqrt{d_{\max}}\,\log n+(\log n)^2.
\]
When $r$ is fixed, \cite{SuWangZhang2020} establishes strong consistency when the expected degree exceeds $C\log n$. When $r$ increases, the required degree also grows with $r$, reflecting the need for a denser network to achieve exact community recovery.
We next consider fixed $r$, where these quantities admit a simpler
interpretation in terms of community imbalance and the degree scale.

The following corollary gives an explicit singular-value formulation when the
degree parameters remain comparable but the community sizes may be highly
unbalanced. Assume that $r$ is fixed and that, for constants
$c_\theta,C_\theta>0$ and a scale $\theta_0=\theta_0(n)>0$,
\begin{equation}
    c_\theta\theta_0
    \le\theta_i\le
    C_\theta\theta_0
    \qquad\text{for all }i\in[n].
    \label{eq:dcsbm-degree-regular}
\end{equation}
We also assume $\theta_0^2 \sigma_{\max}(B)=o(1)$. Since \[ p_{\max} = \max_{i,j} P_{ij} =\max_{i,j} \theta_i \theta_j B_{g(i),g(j)} \le C_\theta^2\theta_0^2 \max_{a,b} B_{ab}\le C_\theta^2\theta_0^2\sigma_{\max}(B) ,\]  this implies that the condition $p_{\max} = o(1)$ in Assumption~\ref{assum:DCSBM} is satisfied. Let
\begin{equation*}
    m_n:=\min_{a\in[r]}n_a,
    \qquad
    N_n:=\max_{a\in[r]}n_a,
    \qquad
    \kappa_n:=\frac{N_n}{m_n},
    \qquad
    d_-:=m_n\theta_0^2.
\end{equation*}
Thus, $\mu_-\asymp 1/\sqrt{N_n}$ and $s_-\asymp
m_n\theta_0^2=d_-$. 
\begin{coro}[Unbalanced communities with a varying block matrix]
\label{cor:dcsbm-unbalanced}
Suppose Assumption~\ref{assum:DCSBM} and
\eqref{eq:dcsbm-degree-regular} hold, and suppose $r$ is fixed. For every
fixed $D>0$, there exist constants $C_{1,D},C_{2,D}>0$, depending only on
$D,r,c_\theta$, and $C_\theta$, such that, if
\begin{equation}
    d_- = m_n\theta_0^2 \ge
    C_{1,D}\frac{\sigma_{\max}(B)}{\sigma_{\min}(B)^2}
    \left( \kappa_n\log n+\kappa_n^{3/2} \right),
    \label{eq:dcsbm-unbalanced-degree}
\end{equation}
then, with probability at least $1-n^{-D}$,
\begin{align}
    \min_{O\in\mathbb O(r)} \|\widetilde U_r O-U\|_{2,\infty}
    \le \frac{C_{2,D}}{\sqrt{N_n}}
    \Bigg( \frac{\sqrt{\sigma_{\max}(B)}}{\sigma_{\min}(B)} \sqrt{\frac{\kappa_n\log n}{d_-}}
        + \frac{\sigma_{\max}(B)}{\sigma_{\min}(B)^2} \frac{\kappa_n^{3/2}}{d_-}
    \Bigg).
    \label{eq:dcsbm-unbalanced-rowwise}
\end{align}
Moreover, $C_{1,D}$ can be chosen sufficiently large so that
Algorithm~\ref{alg:dcsbm} exactly recovers the communities up to a permutation
of their labels.
\end{coro}
The corollary is particularly clear in the regular case. If the
communities are balanced, so that $m_n\asymp N_n\asymp n$ and $\kappa_n\asymp1$, and $\sigma_{\max}(B)\asymp \sigma_{\min}(B) \asymp1$, then Corollary~\ref{cor:dcsbm-unbalanced} reduces to the logarithmic degree condition
\[
    n\theta_0^2\gtrsim_D\log n.
\]
Thus the spectral algorithm achieves exact recovery down to the optimal
logarithmic expected-degree regime (see \cite{Abbe2017}). More generally, if
$B=c_nB_0$ for a fixed well-conditioned matrix $B_0$, then
$\sigma_{\max}(B)\asymp \sigma_{\min}(B) c_n$, and the condition becomes
\[
    c_n n\theta_0^2\gtrsim_D\log n,
\]
again corresponding to logarithmic expected degree. In this case, the leading
term in \eqref{eq:dcsbm-unbalanced-rowwise} is
\[
    n^{-1/2}
    \sqrt{\frac{\log n}{c_n n\theta_0^2}}.
\]
For unbalanced communities, condition~\eqref{eq:dcsbm-unbalanced-degree}
further quantifies the losses due to the weakest block direction, through
$\sigma_{\min}(B)$, and the increased noise associated with $\sigma_{\max}(B)$ and community imbalance. The proofs of Theorem~\ref{thm:dcsbm-general} and Corollary~\ref{cor:dcsbm-unbalanced} are deferred to Appendix~\ref{app:dcsbm-proof}.

\begin{figure}[tbp]
\centering
\begin{minipage}{0.48\textwidth}
    \centering
    \includegraphics[width=1\linewidth]{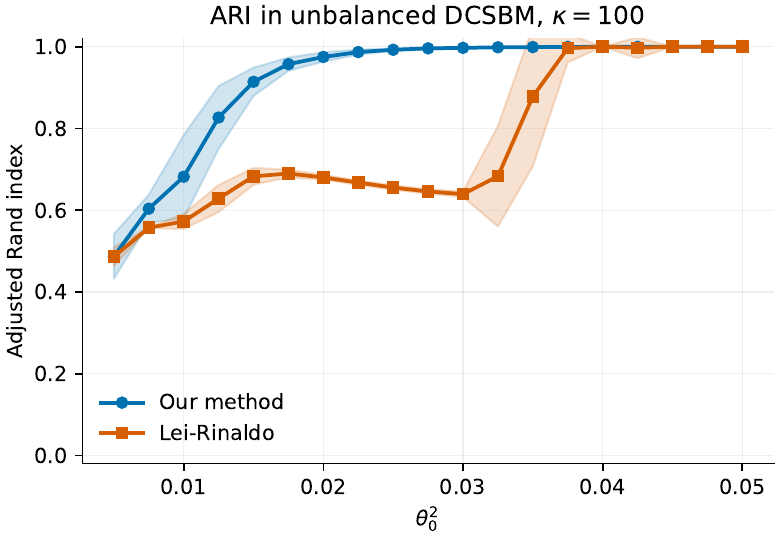}
    \smallskip
    \\{\small (a) Fixed $\kappa=100$, varying $\theta_0^2$.}
\end{minipage}
\hfill
\begin{minipage}{0.48\textwidth}
    \centering
    \includegraphics[width=1\linewidth]{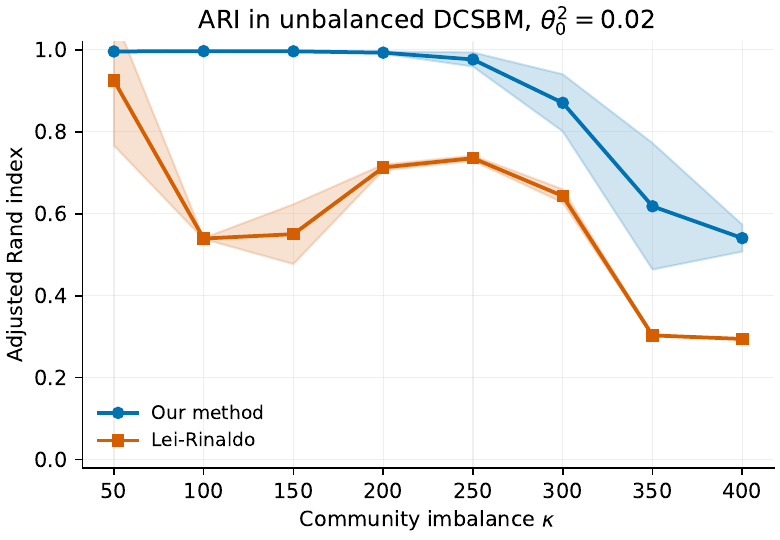}
    \smallskip
    \\{\small (b) Fixed $\theta_0^2=0.02$, varying $\kappa$.}
\end{minipage}
\caption{ARI under unbalanced DCSBMs. The curves report the mean over 100 Monte Carlo repetitions, with shaded regions representing one standard deviation. The proposed method is more robust to community-size imbalance than classical spectral clustering.}
\label{fig:unbalanced-dcsbm}
\end{figure}

\paragraph{Numerical simulations.} We consider the DCSBM
\[
P=\theta_0^2\Theta Z BZ^\top\Theta,
\qquad \theta_i\overset{\mathrm{iid}}{\sim}\mathrm{Unif}(0.8,1.2),
\]
where $\theta_0^2$ denotes the overall network density. We consider $r = 5$ communities with a severe imbalance level. For each setting, we report mean Adjusted Rand Index (ARI) \cite{hubert1985comparing} with one-standard-deviation bands over 100 repetitions. A larger ARI indicates a better community detection result. 
We apply our Algorithm \ref{alg:dcsbm} and the classical spectral community detection algorithm in \cite{LeiRinaldo2015}, where the only difference is at the $k$-means step. To avoid extreme values in empirical analysis, we take the center of 10 farthest data points in Step 3. The results are shown in Figure \ref{fig:unbalanced-dcsbm}.

In the first experiment, the community sizes are \((8000,3000,800,300,80)\), so \(\kappa=100\), with \(\operatorname{diag}(B)=(1,1,3,8,20)\) and off-diagonals of $B$ as 0.2. The diagonal entries mean that smaller communities usually have denser connections. As \(\theta_0^2\) increases, our proposed method recovers the communities substantially earlier than the standard spectral community detection method in \cite{LeiRinaldo2015}, where the main difference is the $k$-means step. In sufficiently dense networks, both methods perform well.

In the second experiment, \(n=12000\), \(\theta_0^2=0.02\), and the community sizes follow
$n_k=n\frac{\kappa^{-(k-1)/4}}{\sum_{j=1}^5\kappa^{-(j-1)/4}}$, $k = 1, 2, 3, 4, 5$, so that $n_1/n_5 = \kappa$. We take \(\kappa=50,\ldots,400\). For each $\kappa$, we set the diagonals of $B$ matrix as 
\[
B_{kk}=C_{k}\frac{n_1}{n_k},
\qquad
C=\left(1,\frac14,\frac14,\frac14,\frac14\right).
\]
Hence, as $\kappa$ increases, the smallest community $n_5$ enjoys that the same order of $n_5 B_{55}$. 
Our proposed method remains nearly perfect through \(\kappa=250\) and consistently outperforms the classical spectral community detection method as imbalance increases. When $\kappa$ exceeds 250, $B_{55}$ cannot increase since the probability must be cut at 1. Hence, the performance declines due to the shrunk signal strength. 

\paragraph{Relation to earlier DCSBM results.}
Lei and Rinaldo~\cite{LeiRinaldo2015} established weak consistency for
spherical adjacency spectral clustering under the DCSBM. Related spectral
procedures include regularized clustering~\cite{QinRohe2013}, SCORE
normalization~\cite{Jin15}, and normalized-adjacency methods
\cite{GulikersLelargeMassoulie2017}. Lei and Zhu~\cite{LeiZhu2017}
subsequently obtained exact recovery at logarithmic degree by combining a
consistent initializer with sample splitting and cross-clustering; their main
presentation assumes community sizes proportional to $n$. Su, Wang, and
Zhang~\cite{SuWangZhang2020} establish strong consistency for regularized
spectral clustering under a different normalized-Laplacian framework for balanced communities. When $r\to\infty$, existing theory is much more limited and is primarily available for balanced SBMs; see, for example, \cite{banks2016information,Zhang2015Minimax}.
Theorem~\ref{thm:dcsbm-general} instead analyzes the raw adjacency eigenspace directly and allows the number of communities $r=r_n$ to increase with $n$. Its guarantees retain arbitrary imbalance and degree heterogeneity through $s_-$, $\mu_\pm$, and $\Gamma_{\mathrm{dc}}$. Corollary~\ref{cor:dcsbm-unbalanced} provides an explicit imbalance-dependent guarantee without sample splitting or a
local refinement step.

\subsection{Generalized random dot product graphs}
\label{sec:grdpg}

The generalized random dot product graph (GRDPG) provides a broad
latent-position framework for network models, including stochastic
block models and their variants; see
\cite{RubinDelanchyEtAl2022,AthreyaEtAl2018}. In particular,
degree heterogeneity can be incorporated through vertex-specific
scalings of the latent positions, so the DCSBM considered in
Section~\ref{subsec:dcsbm} may be viewed as a structured special case.
Since the preceding application already treats heterogeneous degree
parameters and unbalanced communities, here we impose a regular latent
design and focus on a different question: the effect of the embedding
dimension.

When the full signal dimension is retained, adjacency spectral
embedding estimates the latent positions up to the natural
non-identifiability of the GRDPG. In practice, however, one may embed
into a smaller dimension $k<r$, either because the signal rank is
unknown or because a lower-dimensional representation is desired.
Our perturbation results show that such truncation can introduce a
systematic geometric bias when the variance profile is heterogeneous.

Consider a GRDPG with deterministic latent
positions $y_1,\ldots,y_n\in\R^r$. Let
\[
    J:=\diag(I_p,-I_q),
    \qquad p+q=r,
\]
and write $Y=(y_1,\ldots,y_n)^\T\in\R^{n\times r}$. Conditional on the
latent positions, the edges are independent with
\[
    \widetilde A_{ij}\sim\operatorname{Bernoulli}(P_{ij}),
    \qquad
    P_{ij}=\rho_n y_i^\T Jy_j,
    \qquad i<j,
\]
where $\rho_n$ controls the network density. We allow independent
self-loops for theoretical simplicity. The population adjacency matrix is
\[
    P:=\E\widetilde A=\rho_nYJY^\T.
\]
The matrix $J$ specifies the signature of the latent inner product,
allowing $P$ to have both positive and negative eigenvalues; the ordinary
random dot product graph corresponds to $q=0$. The primary inferential
goal is to recover the latent positions $Y$ from $\widetilde A$.

Define
\begin{equation}
    d_n:=n\rho_n,
    \qquad
    \Delta_n:=\frac1nY^\T Y,
    \qquad
    \bar y_n:=\frac1n\sum_{i=1}^ny_i.
    \label{eq:grdpg-scales}
\end{equation}
We work first under a regular-design condition.

\begin{assumption}[Regular GRDPG]
\label{ass:grdpg}
There are fixed constants $c,C>0$ such that
\begin{equation}
    \max_i\|{y_i}\|\le C,
    \qquad
    cI_r\preceq\Delta_n\preceq CI_r,
    \qquad
    c\le y_i^\T J\bar y_n\le C,
    \quad i\in[n],
    \label{eq:grdpg-design}
\end{equation}
and
\begin{equation}
    0\le P_{ij}=\rho_n y_i^\T Jy_j\le1-c,
    \qquad i,j\in[n].
    \label{eq:grdpg-probabilities}
\end{equation}
\end{assumption}
In Assumption~\ref{ass:grdpg}, the bound $\max_i\|y_i\|\le C$ keeps
individual latent positions bounded, while $cI_r\preceq\Delta_n\preceq
CI_r$ ensures they span all $r$ directions at a nondegenerate scale.
Since $\sum_{j=1}^n P_{ij}=d_n\,y_i^\T J\bar y_n$, the quantity
$y_i^\T J\bar y_n$ is proportional to the expected degree of vertex $i$,
so the bounds $c\le y_i^\T J\bar y_n\le C$ force all expected degrees to
be of order $d_n=n\rho_n$. Finally, \eqref{eq:grdpg-probabilities}
keeps all edge probabilities valid and uniformly below one. Together,
these conditions make the nonzero eigenvalues of $P$ of order $d_n$
while allowing the model to remain sparse as $\rho_n\to0$.

Let
\[
    P=U\Lambda U^\T,\qquad \Lambda= \diag(\lambda_1, \ldots, \lambda_r)
\]
be its rank-$r$ eigendecomposition,  with the positive eigenvalues listed first and the negative eigenvalues last, each ordered by decreasing magnitude. Set
\[
    X:=U |\Lambda|^{1/2},
    \qquad
|\Lambda|:=\diag(|\lambda_1|,\ldots,|\lambda_r|).
\]
Then $XJX^\T=P$.  Since $P=\rho_nYJY^\T$, there exists $Q_0\in\mathbb{O}(p,q):=\{Q\in\R^{r\times r}:Q^\T JQ=J\}$ such that
\begin{equation}
    XQ_0=\sqrt{\rho_n}\,Y.
    \label{eq:grdpg-factor-equivalence}
\end{equation}
Thus $X$ represents the latent positions $Y$ up to a scaling and an indefinite orthogonal transformation.

For the observed adjacency matrix $\widetilde A$, select its $k$  eigenvalues of largest magnitude and reorder the selected eigenpairs so that the positive eigenvalues precede the negative ones, with decreasing magnitude within each sign class. Denote the resulting eigenvalue matrix by
\[
    \widetilde\Lambda_k := \diag(\widetilde\lambda_1,\ldots,\widetilde\lambda_k).
\]
Denote $ |\widetilde\Lambda_k| = \diag(|\widetilde\lambda_1|,\ldots,|\widetilde\lambda_k|).$ The $k$-dimensional adjacency spectral embedding is
\[
    \widetilde X_k := \widetilde U_k|\widetilde\Lambda_k|^{1/2}.
\]
When $k=r$, the selected eigenpairs are the $r$ empirical eigenvalues
of largest magnitude, ordered by the same sign convention as the
population eigenvalues. In this case,
$\widetilde X_r$ estimates $X$ up to a
block-orthogonal transformation in $\blockO{p}{q}$, and hence the scaled
latent positions $\sqrt{\rho_n}Y$ up to the intrinsic
$\mathbb O(p,q)$ nonidentifiability.

In practice, however, the latent dimension $r$ is typically unknown, so the embedding dimension $k$ must be selected from the data. The statistical behaviour depends qualitatively on whether $k=r$, $k<r$, or $k>r$. The choice of the embedding dimension is itself a model-selection problem. Existing methods include network cross-validation
\cite{ChenLei2018} and spectral procedures that retain eigenvalues above an estimated noise floor \cite{HongCape2026}. Our purpose here is complementary. Conditional on a chosen embedding dimension, we quantify the rowwise error of the resulting embedding and determine how the heterogeneous Bernoulli variances distort its population target.

\paragraph{Full embedding.} We first consider the correctly specified case $k=r$, where the adjacency spectral embedding targets the full latent position matrix, omitting no signal directions. The next theorem gives a finite-sample $2\to\infty$ error bound for this embedding and, consequently, a uniform latent-position recovery guarantee down to the logarithmic expected-degree regime. The logarithmic-degree threshold follows from Proposition~\ref{prop:refined-weighted-embedding}; we defer the proof to Appendix~\ref{app:grdpg-proof}.

\begin{theorem}[Full GRDPG spectral embedding]
\label{thm:grdpg-full}
Suppose Assumption~\ref{ass:grdpg} holds.  For every fixed $D>0$, there are constants
$C_{1,D},C_{2,D}>0$ such that, if
\begin{equation}
    d_n\ge C_{1,D}\log n,
    \label{eq:grdpg-log-degree}
\end{equation}
then, with probability at least $1-n^{-D}$, the $r$ empirical eigenvalues of largest magnitude consist of $p$ positive and $q$ negative signal outliers, and
\begin{align}
    \min_{W\in \blockO{p}{q}} \|{\widetilde X_r W-X}\|_{2,\infty} &\le C_{2,D}\sqrt{\frac{\log n}{n}},
    \label{eq:grdpg-full-spectral}\\
    \inf_{Q\in\mathbb O(p,q)}
    \|{\rho_n^{-1/2}\widetilde X_r Q-Y}\|_{2,\infty} &\le C_{2,D}\sqrt{\frac{\log n}{d_n}}.
    \label{eq:grdpg-full-latent}
\end{align}
Consequently, the rescaled adjacency spectral embedding is uniformly consistent whenever
$d_n/\log n\to\infty$.
\end{theorem}

Theorem~\ref{thm:grdpg-full} strengthens existing uniform recovery
guarantees for the GRDPG. Rubin-Delanchy, Cape, Tang, and Priebe \cite{RubinDelanchyEtAl2022}
established uniform consistency and asymptotic normality under an
i.i.d.\ latent-position model and a polylogarithmic expected-degree
condition, i.e., $d_n=\omega\!\left(\log^{4c}n\right)$ for $c > 1$. In contrast,
Theorem~\ref{thm:grdpg-full} is finite-sample, allows a deterministic
latent design, and is valid down to the logarithmic degree regime
$d_n\gtrsim\log n$. Moreover, the bound $\sqrt{{\log n}/{n}}$ in \eqref{eq:grdpg-full-spectral} has the same order as the minimax lower bound of Yan and Levin
\cite{YanLevin2023} in the bounded-condition-number regime. Thus the usual adjacency spectral embedding attains the optimal rowwise scale, up to constants, at essentially the weakest degree scale at which uniform recovery can be expected. 

\paragraph{Truncated embeddings and geometric bias.}
The situation is qualitatively different when $k<r$. Under-embedding necessarily loses information about the omitted signal directions, but it can also distort the retained coordinates themselves. Our interest is in this second effect: heterogeneous noise variances induce a systematic geometric bias through the interaction between retained and omitted signal directions. This is complementary to the recent work of Taing and Levin \cite{TaingLevin2026}, which studies dimension misspecification relative to the full latent-position target. Here we instead compare $\widetilde X_k$ with the truncated population spectral factor $X_k$ and quantify the additional bias within the retained signal space.

Write
\begin{equation}
    P
    =
    U_+\Lambda_+U_+^\T
    -
    U_-\Lambda_-U_-^\T,
    \label{eq:grdpg-signed-decomposition}
\end{equation}
where
\[
    \Lambda_+
    =
    \diag(\lambda_1^+,\ldots,\lambda_p^+),
    \qquad
    \Lambda_-
    =
    \diag(\lambda_1^-,\ldots,\lambda_q^-),
\]
and both sequences are decreasing and positive.  

Fix $0\le k_+\le p$ and $0\le k_-\le q$ with
$k:=k_++k_-<r$, and define
\[
    U_k := (U_{+,k_+},U_{-,k_-}),
    \qquad
    \Lambda_k := \diag(\Lambda_{+,k_+},-\Lambda_{-,k_-}),
    \qquad
    X_k:=U_k|\Lambda_k|^{1/2},
\]
where
\[
    \Lambda_{+,k_+}
    :=
    \diag(\lambda_1^+,\ldots,\lambda_{k_+}^+),
    \qquad
    \Lambda_{-,k_-}
    :=
    \diag(\lambda_1^-,\ldots,\lambda_{k_-}^-),
\]
and empty blocks are omitted. The empirical embedding $\widetilde X_k$
is defined analogously from the selected positive and negative empirical
outliers.

For the centered adjacency noise $E:=\widetilde A-P$, we have $R_i:=\sum_{j=1}^nP_{ij}(1-P_{ij})$ and $\mathcal R:=\diag(R_1,\ldots,R_n).$ Define the within-sign boundary gaps
\[
    \delta_{+,n}
    :=
    \lambda_{k_+}^+-\lambda_{k_++1}^+,
    \qquad
    \delta_{-,n}
    :=
    \lambda_{k_-}^--\lambda_{k_-+1}^-,
\]
whenever $0<k_+<p$ and $0<k_-<q$, respectively. We adopt the convention
that $\delta_{+,n}=\infty $ if $k_+\in\{0,p\}$ and  $\delta_{-,n}=\infty$ if $k_-\in\{0,q\}$.

For $k_+>0$, set
\[
    \lambda_{+,*}:=\lambda_{k_+}^+,
    \qquad
    \mathcal V_{+,1}:= \big\|  U_{+,>k_+}^\T \mathcal R (U_{+,k_+},U_-) \big\|,
    \qquad
    \mathcal V_{\pm} := \|U_-^\T\mathcal R U_+\|,
\]
where $\mathcal V_{+,1}:=0$ if $k_+=p$, and define
\[
    \mathcal B_{k_+}^+
    :=
    \frac{10\sqrt{k_+}}
         {\delta_{+,n}\lambda_{+,*}}
    \left(
        \mathcal V_{+,1}
        +
        8\frac{R_{\max}}{\lambda_{+,*}^2}
        \osc(\mathcal R)
    \right)
    +
    4\sqrt{k_+}
    \frac{\mathcal V_{\pm}}{\lambda_{+,*}^2}.
\]
Set $\mathcal B_0^+:=0$. Define
$\lambda_{-,*}$ and $\mathcal B_{k_-}^-$ analogously by interchanging
the positive and negative signal blocks, and set $\mathcal B_0^-:=0$.

Finally, define the geometric-bias contribution to the truncated embedding by
\begin{equation}
    \GBias_{k_+,k_-}
    :=
    \|U\|_{2,\infty}
    \left(
        \sqrt{\lambda_{+,*}}\,\mathcal B_{k_+}^+
        +
        \sqrt{\lambda_{-,*}}\,\mathcal B_{k_-}^-
    \right),
    \label{eq:grdpg-GBias}
\end{equation}
where a term is omitted when the corresponding retained sign block is empty.

The leading components of
\eqref{eq:grdpg-GBias} are
\begin{equation}
    \|U\|_{2,\infty} \frac{ \sqrt{k_+} \big\|U_{+,>k_+}^\T
            \mathcal R  U_{+,k_+}\big\| }{ \delta_{+,n}\sqrt{\lambda_{+,*}} },
    \qquad
    \|U\|_{2,\infty}
    \frac{ \sqrt{k_-}  \big\| U_{-,>k_-}^\T \mathcal RU_{-,k_-} \big\| }{ \delta_{-,n}\sqrt{\lambda_{-,*}} }.
    \label{eq:grdpg-leading-bias}
\end{equation}
They measure the failure of the retained and omitted signal directions
to remain orthogonal under the variance-weighted inner product induced
by $\mathcal R$.  Now we present the truncated embedding result. Denote 
\[
    \delta_n:=\min\{\delta_{+,n},\delta_{-,n}\},
\]

\begin{prop}[Truncated GRDPG embedding]
\label{prop:grdpg-truncated}
Suppose Assumption~\ref{ass:grdpg} holds. For every fixed $D>0$,there exist constants $C_{1,D},C_{2,D}>0$ such that , if
\[
    d_n\ge C_{1,D}\log n,
    \quad\text{and} \quad
    \delta_n\ge C_{1,D}\log n,
\]
then, with probability at least $1-n^{-D}$,
\begin{equation}
    \min_{W\in\blockO{k_+}{k_-}}
    \|{\widetilde X_kW-X_k}\|_{2,\infty}
    \le
    \GBias_{k_+,k_-}
    +
    C_{2,D}\left(
        \sqrt{\frac{\log n}{n}}
        +
        \sqrt{\frac{d_n}{n}}\,
        \frac{\log n}{\delta_n}
    \right).
    \label{eq:grdpg-truncated-main}
\end{equation}
\end{prop}

Proposition~\ref{prop:grdpg-truncated} separates the systematic geometric bias from the stochastic rowwise error. The leading bias is governed by the variance-weighted couplings \[ U_{+,>k_+}^\T\mathcal R U_{+,1:k_+} \quad\text{and}\quad U_{-,>k_-}^\T \mathcal R U_{-,1:k_-}, \] and is amplified when the truncation boundary lies near another same-sign signal eigenvalue. Thus networks with comparable signal eigenvalues and noise levels may nevertheless exhibit different truncated-embedding errors because their variance profiles interact differently with the signal eigenspaces. In particular, the leading bias vanishes when $\mathcal R=cI$. This provides information complementary to methods that select eigenvalues above an estimated noise floor \cite{HongCape2026}: such methods determine which directions are distinguishable from noise, whereas our bound also quantifies the stability of truncating within the signal spectrum. In the well-separated regime $\delta_n\ge c d_n$, the geometric bias is of smaller order, and Proposition~\ref{prop:grdpg-truncated} simplifies to 
\begin{equation} 
\min_{W\in\blockO{k_+}{k_-}} \|{\widetilde X_k W-X_k}\|_{2,\infty} \lesssim \sqrt{\frac{\log n}{n}}. 
\label{eq:grdpg-truncated-macroscopic} 
\end{equation} 
A short verification of \eqref{eq:grdpg-truncated-macroscopic} is given in the appendix. The over-embedded case $k>r$ involves noise-bulk eigenvectors rather than signal outliers and is outside the scope of the present perturbation theory.

\section{Geometric bias, sharpness, and debiasing}\label{sec:geometric-bias}
\subsection{Interpretation of the geometric bias}
The present variance-profile model introduces a new deterministic geometric bias term
\[
\mathcal B_k=
\frac{10\sqrt{k}}{\delta_k \lambda_k}\left( \|\mathcal{V}_{\mathcal{J}\mathcal{I}}\|+8\frac{R_{\max}}{\lambda_k^2}\osc(\mathcal R)\right)+4\sqrt{k}\frac{\|\mathcal{V}_{\mathcal{N}\mathcal{K}}\|}{\lambda_k^2},
\]
which captures the interaction between the signal eigenspaces and the inhomogeneous row-variance profile $\mathcal R=\diag(R_i)$. In the homogeneous noise setting where $\mathcal R = c I$, the geometric bias term $\mathcal B_k$ vanishes; our current bounds in Theorem \ref{thm:top-k-eigenspace} match the two-term bounds studied in \cite{Wang24}, up to logarithmic factors and rank factors. 

The two quantities $\mathcal{V}_{\mathcal{J}\mathcal{I}}, \mathcal{V}_{\mathcal{N}\mathcal{K}}$ play different roles in the bias term. The term $\mathcal V_{\mathcal J \mathcal I}$ captures how the target top-$k$ cluster interacts with the lower positive cluster $\{k+1,\dots,r_+ \}$, appearing at the gap-dependent scale $(\delta_k\lambda_k)^{-1}$. By contrast, $\mathcal V_{\mathcal N \mathcal K}$ captures interactions between negative and nonnegative signal directions, but at the weaker scale $\lambda_k^{-2}$. In practice, the contribution from $\mathcal V_{\mathcal N \mathcal K}$ is often negligible. For instance, when all spikes are nonnegative, we have $\mathcal N=\emptyset$ and this term vanishes. Even when negative spikes exist, the term involving $\mathcal V_{\mathcal N\mathcal K}$ is usually dominated by $\mathcal V_{\mathcal J \mathcal I}$ as it has a stronger suppression factor $\lambda_k^{-2}$.

We believe that the leading geometric bias terms $$\frac{\|\mathcal V_{\mathcal J\mathcal I}\|}{\delta_k\lambda_k}, \quad \frac{\|\mathcal V_{\mathcal N\mathcal K}\|}{\lambda_k^2}$$ are intrinsic to the variance-profile model. 

\paragraph{A rank-two model illustration.} To explain where these terms come from, consider a simple rank-2 model $A= \lambda_1 u_1u_1^\T+\lambda_2 u_2u_2^\T$ with $\lambda_1 > \lambda_2>0$ and the leading empirical eigenvector $\widetilde{u}_1$. Let $Q$ be the orthogonal projection onto the null space of $A$. From the decomposition
$\widetilde{u}_1 = u_1 u_1^\T \widetilde{u}_1 + u_2 u_2^\T \widetilde{u}_1 + Q \widetilde{u}_1,$
taking the $\ell_2$ norm on both sides and rearranging terms, we have 
\[\sin^2\angle(u_1,\widetilde{u}_1)= (u_2^\T \widetilde{u}_1)^2 + \|Q \widetilde{u}_1\|^2.\] The term $u_2^\T \widetilde{u}_1$ captures the interaction between signal directions. In our proofs, this quantity is controlled by projecting the eigenvector equation
$\widetilde u_1=G(\widetilde\lambda_1)A\widetilde u_1$
onto $u_2$, where $G(z):=(zI-E)^{-1}$ is the resolvent of $E$.  To precisely estimate this term, from the eigenvector equation $(A+E) \widetilde{u}_1 = \widetilde{\lambda}_1 \widetilde{u}_1$, we solve for 
\[ \widetilde{u}_1 = (\widetilde{\lambda}_1 -E)^{-1} A \widetilde{u}_1 = G(\widetilde{\lambda}_1) A \widetilde{u}_1.\]
A key technical input in our proof is to show that the random resolvent $G$ can be precisely approximated by its QVE approximation, a deterministic matrix $\Phi$ (defined \eqref{def:QVE}). This $\Phi$ isolates the deterministic contribution of the variance profile. Indeed, from \[u_2^\T\widetilde u_1\approx u_2^\T\Phi(\widetilde\lambda_1)A\widetilde u_1 = \lambda_1(u_2^\T \Phi(\widetilde\lambda_1) u_1) u_1^\T \widetilde u_1 + \lambda_2(u_2^\T \Phi(\widetilde\lambda_1) u_2) u_2^\T \widetilde u_1 ,\] plugging in the expansion of $\Phi\approx z^{-1}I + z^{-3}\mathcal R$ in \eqref{eq:Phi-expansion}, the leading off-diagonal term is
\[ u_2^\T \Phi(\widetilde\lambda_1)u_1 \approx \frac{u_2^\T\mathcal R u_1}{{\widetilde\lambda_1}^3}.\]
 Substituting this back and approximating $\widetilde\lambda_1 \approx \lambda_1$ and $u_1^\T \widetilde u_1 \approx 1$, we obtain 
\[ u_2^\T \widetilde u_1 \approx \frac{1}{\lambda_1^2} (u_2^\T\mathcal R u_1) + \frac{\lambda_2}{\lambda_1} u_2^\T \widetilde u_1.\]
Rearranging this equation to solve $u_2^\T \widetilde u_1$ yields
\[ |u_2^\T \widetilde u_1| \approx \frac{|u_2^\T\mathcal R u_1|}{\delta_1 \lambda_1},\]
which mirrors the $\frac{\|\mathcal V_{\mathcal J\mathcal I}\|}{\delta_k\lambda_k}$ term in our theorem. If instead $\lambda_1>0>\lambda_2$, the same calculation reveals 
\[ |u_2^\T \widetilde u_1 | \approx \frac{|u_2^\T\mathcal R u_1|}{\lambda_1(\lambda_1 - \lambda_2)}= \frac{|u_2^\T\mathcal R u_1|}{\lambda_1(\lambda_1 + |\lambda_2|)}.\]
We bound $\frac{|u_2^\T\mathcal R u_1|}{\lambda_1(\lambda_1 + |\lambda_2|)}\le \frac{|u_2^\T\mathcal R u_1|}{\lambda_1^2}$ in our perturbation results for simplicity. 

The preceding calculation explains the leading bias terms $\frac{\|\mathcal V_{\mathcal J\mathcal I}\|}{\delta_k\lambda_k}, \frac{\|\mathcal V_{\mathcal N\mathcal K}\|}{\lambda_k^2}$.  By contrast, the higher-order term in $\mathcal B_k$ involving
\[\frac{R_{\max}}{\delta_k\lambda_k^3}\osc(\mathcal R)\]
comes from bounding the diagonal remainder $\varepsilon(z)$ in the second-order  expansion of $\Phi(z)$ in \eqref{def:QVE}:
\[\Phi(z)= \frac{1}{z}I_n + \frac{1}{z^3} \mathcal{R} + \varepsilon(z). \]
We control this remainder using a uniform oscillation estimate, which may not be sharp. A higher-order expansion of $\Phi(z)$ would reveal additional structured variance-profile terms beyond $\mathcal R$ and replace the uniform $\osc(\mathcal R)$ bound with finer geometric couplings, yielding a sharper higher-order term. However, isolating these higher-order structures requires substantially more delicate analysis to control the deterministic remainders and extract the finer stochastic fluctuations. This is beyond the scope of the current paper.

\subsection{Sharpness and numerical experiments}\label{sec:optimality}


Before presenting the numerical experiments, we briefly discuss the sharpness and necessity of the different terms in our perturbation bounds. In the homogeneous-noise setting studied in \cite{Wang24}, eigenspace perturbation is governed by two mechanisms:
stochastic mixing between distinct signal directions, of order $\sqrt{k}/\delta_k$ times the relevant fluctuation scale, and leakage into the noise subspace, of order $\|E\|/\lambda_k$. As discussed in
\cite[Section~4.1]{Wang24}, these contributions are near-optimal in the i.i.d.\ Gaussian noise case.

Under a heterogeneous variance profile, Theorem~\ref{thm:top-k-eigenspace} contains an additional contribution, namely the geometric-bias term $\mathcal B_k$. Its leading components are determined by the variance-weighted couplings between the target and non-target signal eigenspaces, and vanish in the homogeneous case $\mathcal R=cI$. Thus the first two mechanisms are inherited from the
homogeneous theory, whereas $\mathcal B_k$ is specific to variance
heterogeneity.

The rank-two example above shows explicitly that this geometric contribution appears at the predicted scale and therefore cannot in general be removed from the perturbation bound. The simulations below are designed to isolate this effect and illustrate that the deterministic bias persists under random perturbations and is visible at finite sample sizes.

\paragraph{Numerical simulations.} We consider a rank-2 model $A= \lambda_1 u_1u_1^\T +\lambda_2 u_2u_2^\T$ with spikes $\lambda_1 > \lambda_2 > 0$ and eigenvectors $u_1 = \frac{1}{\sqrt 2}(e_1 + e_2)$ and $u_2 = \frac{1}{\sqrt 2}(e_1 - e_2)$. The noise variance profile has row sums $\mathcal R_{\Delta} = \text{diag}(\rho+\Delta, \rho-\Delta,\rho,\dots,\rho)$, where $\rho>0$ is the baseline variance and $0\le \Delta\le 0.8\rho$ is the heterogeneity parameter. This gives the geometric coupling $\mathcal V_{12}=u_2^\T \mathcal R_{\Delta} u_1 =\Delta$.

\begin{figure}[tbp]
\centering
\begin{minipage}{0.48\textwidth}
    \centering
    \includegraphics[width=\linewidth]{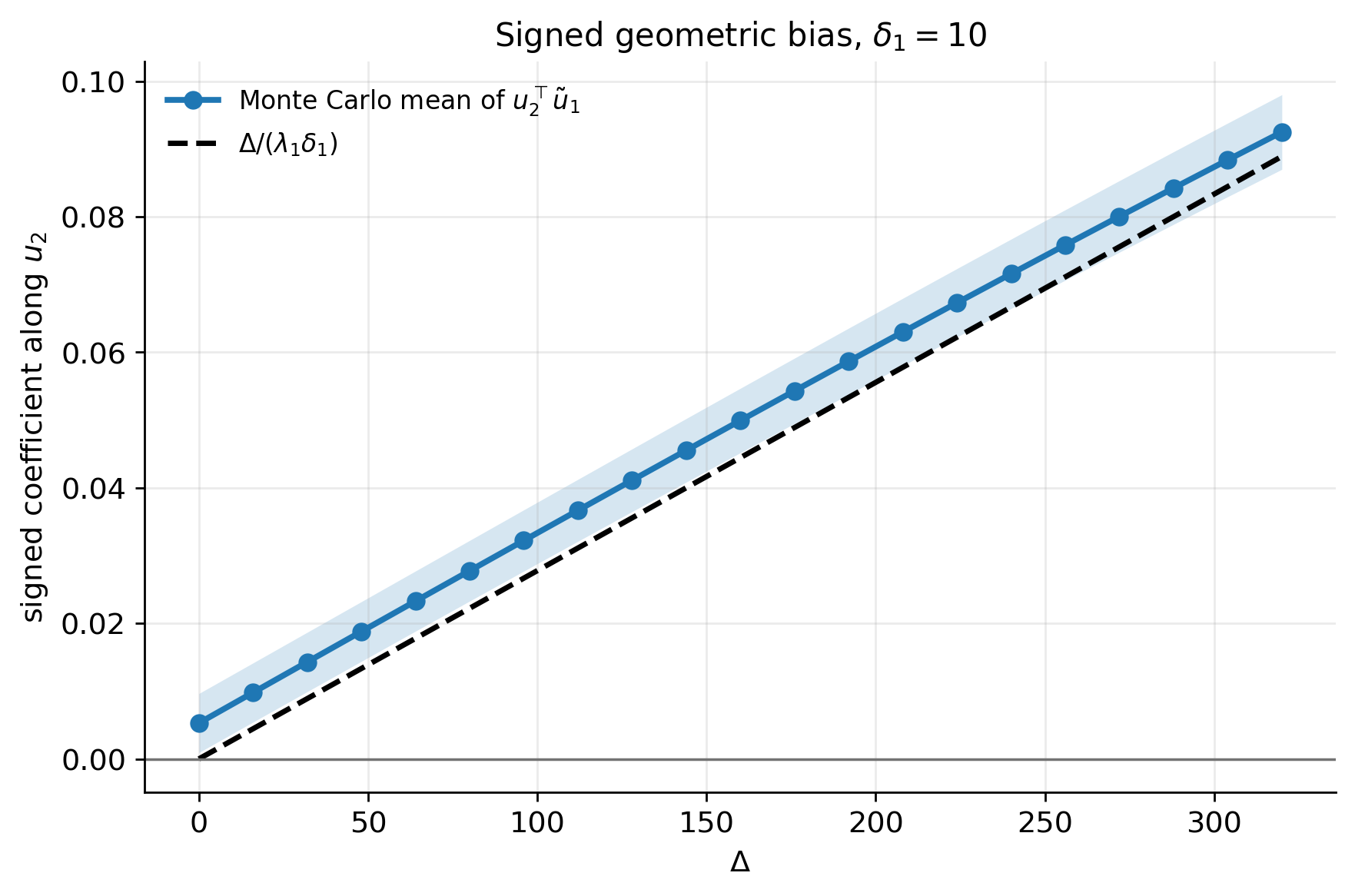}
    \smallskip
    \\{\small (a) $\delta_1=10$.}
\end{minipage}
\hfill
\begin{minipage}{0.48\textwidth}
    \centering
    \includegraphics[width=\linewidth]{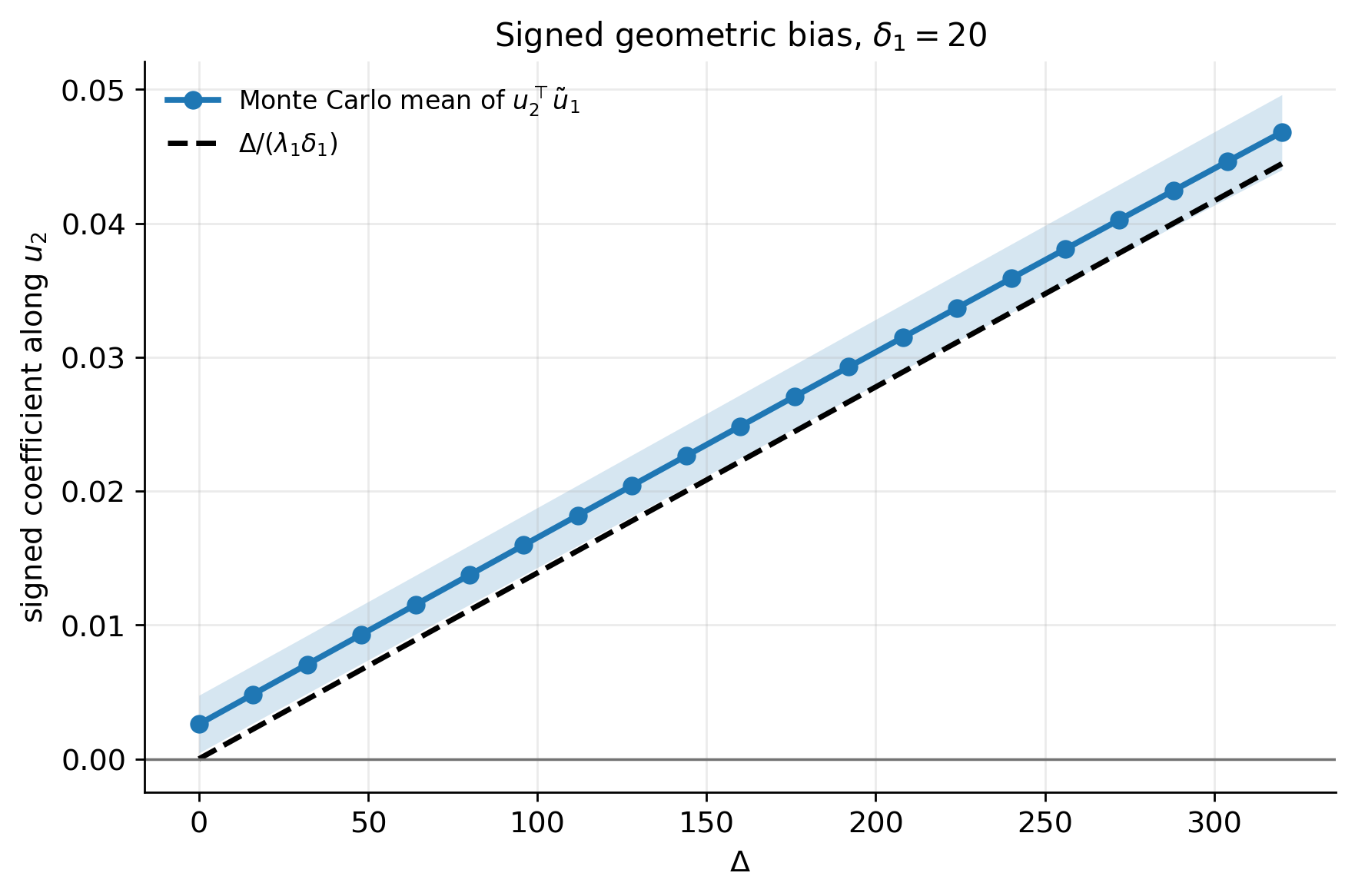}
    \smallskip
    \\{\small (b) $\delta_1=20$.}
\end{minipage}
\caption{Geometric bias under varying spectral gaps. Setup: dimension $n=1200$, baseline noise $\rho=400$, top spike $\lambda_1=360$, heterogeneity parameter $\Delta \in \{0,16,32,\ldots,320\}$, with 300 Monte Carlo samples. Each panel shows the empirical alignment $\mathbb{E}[u_2^\T \widetilde u_1]$ versus $\Delta$, compared with the theoretical prediction $\Delta/(\lambda_1\delta_1)$.}
\label{fig:geometric-bias-gap-scan}
\end{figure}

To show that geometric bias is deterministic rather than random, we compute the Monte Carlo average of $u_2^\T \widetilde u_1$, where in each trial we align the empirical eigenvector so that $u_1^\T \widetilde u_1 \ge 0$. Since the noise is zero-mean, averaging removes random fluctuations while preserving the deterministic bias. Figure \ref{fig:geometric-bias-gap-scan} compares these empirical averages with our theoretical prediction $\frac{\Delta}{\delta_1 \lambda_1}$ for two different spectral gaps $\delta_1 = \lambda_1 - \lambda_2$.

\subsection{Oracle de-biasing}
The goal of this subsection is modest. We show, at the oracle level, that the leading $\mathcal V_{\mathcal J\mathcal I}$-bias in $\mathcal B_k$ can be removed. As discussed earlier, $\mathcal V_{\mathcal J\mathcal I}$ causes the main deviation from previous homogeneous bounds. When $k=r_+$, $\mathcal J=\emptyset$  and this bias vanishes. For intermediate eigenspaces ($k<r_+$), correction may help when the deterministic bias dominates the stochastic fluctuations, i.e.,
\[ \frac{\|\mathcal V_{\mathcal J\mathcal I}\|}{\lambda_k}\gg M.\]

For $s\in [k]$, consider the block matrices $U_{\mathcal{J}}^\T \Phi(\widetilde{\lambda}_s) U_{\mathcal I}$ and $ U_{\mathcal{J}}^\T \Phi(\widetilde{\lambda}_s) U_{\mathcal J}$. 
Define 
\begin{align}\label{def:Ts}
    \mathcal T_s:= (I- U_{\mathcal{J}}^\T \Phi(\widetilde{\lambda}_s) U_{\mathcal J}\Lambda_{\mathcal J} )^{-1} U_{\mathcal{J}}^\T \Phi(\widetilde{\lambda}_s) U_{\mathcal I} \Lambda_{\mathcal I}.
\end{align}
The next result guarantees that $\mathcal T_s$ is well-defined on the event considered in Theorem~\ref{thm:top-k-eigenspace}. Its proof is provided in Appendix \ref{app:J-block-stability}.

\begin{lemma}\label{lem:J-block-stability}Assume $\mathcal J=\{k+1,\ldots,r_+\} \neq \emptyset$. Under the assumptions and with the same probability as in Theorem~\ref{thm:top-k-eigenspace}, for every $s\in [k]$,
\[
\left\| \left(I-U_{\mathcal J}^\T\Phi(\widetilde\lambda_s)U_{\mathcal J}\Lambda_{\mathcal J}\right)^{-1}\right\| \le 60\frac{\widetilde\lambda_s}{\delta_k}.
\]
\end{lemma}
Define the oracle corrected vectors
\begin{align*}
    w_s^{\orc}:= P_U\widetilde u_s - U_{\mathcal J}{\mathcal T}_s U_{\mathcal I}^\T \widetilde u_s,
\qquad s=1,\ldots,k.
\end{align*}
Let
\begin{equation}\label{eq:def-W-orc}
W_k^{\orc}:=(w_1^{\orc},\ldots,w_k^{\orc}),
\qquad
U_k^{\orc}:=\orth(W_k^{\orc}).
\end{equation}
We denote by $\orth(W)$ the orthonormalization of the matrix $W$, which consists of orthonormal columns spanning the same subspace as $W$. 

Define the oracle error term
\begin{align*}
    \mathcal{E}_k^{\orc}:= 4\sqrt{k} \frac{\|\mathcal{V}_{\mathcal{N}\mathcal{K}}\|}{\lambda_k^2} + 24 \sqrt{k} \frac{R_{\max}}{{\lambda}_k^4} \osc(\mathcal R) + 124\sqrt{k} \frac{M_D}{\delta_k}.
\end{align*}
\begin{prop}[Oracle blockwise de-biasing]\label{prop:oracle-block-debias}Assume that the conditions of Theorem \ref{thm:top-k-eigenspace} hold. The following holds with the same probability as Theorem \ref{thm:top-k-eigenspace}:
\begin{align}
    \|\sin\angle(U_k^{\orc},U_k)\| &\le 2\mathcal{E}_k^{\orc},\\
    \min_{O\in\mathbb O(k)}\|U_k^{\orc}-U_kO\|_{2,\infty}&\le 3 \|U\|_{2,\infty} \mathcal{E}_k^{\orc}.
\end{align}
\end{prop}
The proof of Proposition \ref{prop:oracle-block-debias} is deferred to Appendix \ref{app:oracle-block-debias}. The key point is that the oracle correction subtracts the deterministic contribution from the $\mathcal J$-block predicted by the QVE. This removes the leading $\mathcal V_{\mathcal J\mathcal I}$-term in $\mathcal B_k$. Additionally, the projection $P_U$ removes the null-space component, which explains why the signal-noise term 
$\|E\|/\lambda_k$ in
Theorem~\ref{thm:top-k-eigenspace} is absent from Proposition~\ref{prop:oracle-block-debias}.

The oracle estimator uses the unknown population eigenspace. A natural plug-in version replaces the population blocks by empirical outlier blocks. We provide a plug-in de-biasing implementation and explain its limitations in the Appendix \ref{app:plug-in}. Developing optimal data-driven bias correction methods is left for future work.

\section{Related literature}\label{sec:related}
This section highlights perturbation results relevant to heterogeneous random
noise. A broader review of matrix perturbation, rowwise eigenspace analysis,
and statistical applications can be found in~\cite{Wang24}. Here we focus on
work most closely related to variance heterogeneity and random matrix
perturbation with non-identically distributed noise.

\paragraph{Perturbation bounds beyond worst-case theory.}
Classical perturbation results, including Weyl's inequality~\cite{weyl1912}
and the Davis-Kahan-Wedin sin~$\Theta$ theorem~\cite{DK,Wedin}, control
eigenvalue and eigenspace errors through the operator norm of the perturbation
and the relevant spectral separation. While sharp for arbitrary deterministic perturbations, these bounds can be conservative for structured random perturbations. Modern refinements include population-gap
variants, Schatten and unitarily invariant norm bounds, perturbation
projection error bounds, and other deterministic or stochastic improvements;
see, for example,
\cite{YWS15,VL13,CZ18,LHZ21,ZZ22,MR3310977,MR3256861,MR4260218,Zhong17,
MR2846302,OVW2,BV23}. 

Extensive work has also addressed entrywise and rowwise perturbation bounds. Representative works include
\cite{FWZ17,AFWZ20,AFW22,CLCPC21,CTP19,CFMW19,Lei19,ZB18,BV23,YCF21,YW24,YL24}.
Perturbation of linear and bilinear forms of eigenvectors and singular vectors
has been studied in
\cite{MR3565274,MR3960915,EBW18,LCPC21,CWC21,Ag23,ALP22}. These results play an
important role in high-dimensional statistics, including matrix completion,
spectral clustering, ranking, community detection, submatrix localization, and
Gaussian mixture models; see, for example,
\cite{MR2565240,Cands2010MatrixCW,KMO09,AKS98,ACV14,CLR17,CX16,Mc01,
VLBB08,Vu18,LZZ21,BKRSW11,KBR11,BI13,BIS15,MW15,DGGR14,FR10,HWX16,MRZ15,
BBH18,DHB23}. For a statistical survey of spectral methods and perturbation
bounds, see~\cite{CCFM21}. 

In the homogeneous Gaussian-noise setting,
\cite{OVW22,Wang24} obtained sharp singular-subspace perturbation estimates
showing that the error is governed by a signal-to-noise term and a stochastic
signal-space fluctuation term. Recent work of Tran and Vu
\cite{TV24-expansion,TV24-norm} develops combinatorial contour-expansion and
relative-norm perturbation bounds that exploit structural interaction between
the signal and the perturbation. Our work shares the goal of moving beyond worst-case perturbation theory, but through a different mechanism: 
we use variance-profile local laws and QVE expansions to identify a
deterministic geometric bias induced by heterogeneous noise.

\paragraph{Heterogeneous noise in statistical spectral problems.}
There is a substantial statistical literature on spectral estimation under
heteroskedastic or heterogeneous noise. Zhang, Cai and Wu~\cite{ZCW22}
introduced HeteroPCA, based on diagonal deletion and imputation, to correct
the bias caused by heterogeneous diagonal noise in spiked covariance models.
Yan, Chen and Fan~\cite{YCF21} developed inference for heteroskedastic PCA with
missing data, obtaining distributional guarantees for the estimated principal
subspace and entrywise inference for the spiked covariance matrix. Related
work on high-dimensional heteroskedastic PCA and weighted PCA includes
\cite{HBF18,HGBF21,HYFB23}. Zhang and Mondelli~\cite{ZM24} recently
studied rank-one matrix denoising with doubly heteroscedastic noise and
derived asymptotic fundamental limits and optimal spectral methods.

These works are close in motivation but differ from ours in both model and
mechanism. Much of this literature concerns sample covariance, missing-data,
or rectangular denoising models, where the main issue is often diagonal-noise
bias correction, whitening, optimal shrinkage, or statistical inference. By
contrast, we study a symmetric additive signal-plus-noise model with
arbitrary entrywise variance profile and possible sparsity. The bias identified
here is not merely a diagonal variance bias; it is a geometric bias governed by
$\mathcal V=U^\T \mathcal R U$, which records how the row-wise variance
profile changes the geometry of the signal eigenspaces.

We also cite related work on rowwise, entrywise, and inferential perturbation
for completeness, including
\cite{FWZ17,CTP19,AFWZ20,AFW22,CLCPC21,Lei19,BV23,YCF21,YW24,ALP22,Ag23,
CWC21,CC25}; these works concern fine-grained perturbation and
distributional inference in related spectral models, but they do not address
the variance-profile geometric bias studied here.

\paragraph{Related literature in random matrix theory.}
Our work is also connected to the random matrix literature on finite-rank
deformations and outlier eigenvectors. The BBP transition was introduced by
Baik, Ben Arous and P\'ech\'e~\cite{BBP05}, with subsequent developments for
spiked covariance models, matrix denoising, and finite-rank deformations in
\cite{BS06,DP07,BY08,BGN,BGN1,BGM11,BEKYY2016,CM2017,CD18,BDWW20,BDW21,
BW21,Ben20,FFHL19,OVW}. The variance-profile setting is closely related to the
theory of general Wigner-type matrices and the quadratic vector equation.
Ajanki, Erd\H{o}s and Kr\"uger~\cite{AEK19} developed the
QVE/Dyson-equation framework and proved local laws and universality for
matrices with general variance profiles. Our use of the QVE is more elementary
and takes place in the large-$z$ outlier domain, where the solution is stable
and admits a direct expansion.

Recent asymptotic work has begun to study outliers, detection, and inference
limits in spiked models with variance-profile, inhomogeneous, or structured
noise. Bhattacharya, Chakrabarty and Hazra~\cite{BCH26} analyze outlier
eigenvalues and eigenvectors of generalized Wigner matrices with finite-rank
deterministic perturbations. Guionnet, Ko, Krzakala and Zdeborov\'a~\cite{GKKZ25}
study low-rank matrix estimation with inhomogeneous output channels and
detection thresholds. Bao, Cheong, Lee and Li~\cite{BCLL25} study signal
detection from spiked noise via asymmetrization. We also mention the physics
work of Ferreira and Metz~\cite{FM26} on the BBP transition in an inhomogeneous
rank-one spiked Wigner model. These results are complementary to ours: they
are mainly asymptotic and focus on limiting distributions, information limits,
detection thresholds, or BBP transitions, whereas our results are finite-sample
perturbation bounds for deterministic low-rank signals under sparse
independent variance-profile noise.


\section{Deterministic QVE estimates and preliminary bounds}\label{sec:qve-preliminaries}
Consider the resolvent (Green function) of $E$:
\begin{equation*} 
G(z) = (zI - E)^{-1}, \qquad z\in\mathbb{C}^+. 
\end{equation*} 
For $z\in \C^{+}$,  define a deterministic diagonal matrix \[\Phi\equiv\Phi(z) =\mathrm{diag}(\phi_1(z),\ldots,\phi_n(z))\] where $\phi(z)=(\phi_1(z),\ldots,\phi_n(z))$ is the solution to the Quadratic Vector Equations (QVE) given by: \begin{align}\label{def:QVE}
\frac{1}{\phi_i(z)} = z-\sum_{j=1}^n \sigma_{ij}^2 \phi_j(z),\qquad i \in [n],
\end{align}
with $\Im \phi_i(z)<0$ for $z\in \C^{+}$. Each $\phi_i(z)$ is analytic
on $\C^+$. It is shown in \cite[Theorem 2.1]{AEK19} that such a solution exists, is unique, and extends continuously to the real axis outside the spectrum (so in particular $\phi_i(z)>0$ for $z$ sufficiently large). 

We collect deterministic estimates for the QVE solution $\Phi(z)$ and elementary bounds for resolvents $G(z)$. These estimates are used in two places: to identify the deterministic locations of outlier eigenvalues, and to extract the geometric bias from the expansion of $\Phi(z)$.

Recall that
\[
    R_i=\sum_{j=1}^n\sigma_{ij}^2, \qquad R_{\max}=\max_i R_i,  \qquad \mathcal R=\operatorname{diag}(R_1,\ldots,R_n).
\]

The next result provides basic bounds on $\phi_i(z)$ and the second-order expansion of these $\phi_i$'s for large $z$. The proof is provided in Appendix \ref{sec:proof-lem:bdphi}.
\begin{lemma}\label{lem:bdphi}Assume $z\in \C$ satisfying $|z| \ge  \sqrt{6 R_{\max}}$. 

\noindent\textup{(i)} For every $i\in [n]$, we have
    \[\frac{1}{2|z|}\le |\phi_i(z)| \le \frac{3}{2|z|}\] and
    \[\frac{3}{4}|z| \le \big|z-\sum_j \sigma_{ij}^2 \phi_j(z)\big| \le \frac{5}{4}|z|.\]
Furthermore, we have
\begin{align}\label{eq:phi-diff-bd}
  \|\Phi'(z)\| =\max_i  |\phi_i'(z)| \le \frac{4}{|z|^2}
\quad\text{and} \quad
\|\Phi''(z)\| =\max_i|\phi_i''(z)| \le \frac{28}{|z|^3}.
\end{align}
    
\noindent\textup{(ii)} For $i\in [n]$,  we have the second-order expansion:
\begin{align}\label{eq:phi-expand}
    \phi_i(z) = \frac{1}{z} + \frac{R_i}{z^3} + \varepsilon_i(z) \quad \text{with} \quad |\varepsilon_i(z)|\le \frac{45}{8}\frac{R_{\max}^2}{|z|^5}. 
\end{align}
In particular, 
\begin{align}\label{eq:Phi-expansion}
    \Phi(z) = \frac{1}{z}I_n + \frac{1}{z^3} \mathcal{R} + \varepsilon(z),
\end{align}
where 
\[ \varepsilon(z)= \diag(\varepsilon_1(z),\cdots,\varepsilon_n(z))  \quad \text{and} \quad \|\varepsilon(z)\| \le \frac{45}{8}\frac{R_{\max}^2}{|z|^5}. \]
Furthermore, we have
\begin{align}\label{eq:Phi-prime-expansion}
    \Phi'(z)= -\frac{1}{z^2}I_n-\frac{3}{z^4}\mathcal R+\varepsilon'(z),
\end{align}
where
\[\varepsilon'(z)= \diag(\varepsilon_1'(z),\cdots,\varepsilon_n'(z))  \quad \text{and} \quad \|\varepsilon'(z)\|\le \frac{243}{5}\frac{R_{\max}^2}{|z|^6}.
\]
\end{lemma}

The entrywise bound in Lemma~\ref{lem:bdphi} is not sufficient for the
eigenspace perturbation theorem. For the error term $\varepsilon(z) =\diag(\varepsilon_1(z),\dots,\varepsilon_n(z))$ in the second-order expansion of $\Phi(z)$ in \eqref{eq:phi-expand}, we control its \emph{oscillation}. For $z\in \C$ satisfying $|z| \ge \sqrt{6 R_{\max}}$, denote
\begin{equation}
    \osc(\varepsilon(z)):= \max_{i,j} |\varepsilon_i(z)-\varepsilon_j(z)|.
\end{equation}
If $z$ is a real number, then $\osc(\varepsilon(z))=\max_i \varepsilon_i(z) -\min_i \varepsilon_i(z).$ Recall that \[\osc(\mathcal R)=R_{\max}-R_{\min}.\] The bound on $\osc(\varepsilon(z))$ is given below; its proof is deferred to Appendix \ref{app:proof-lem-osc}.
\begin{lemma}\label{lem:osc} For any complex number $z$ with $|z| \ge \sqrt{6 R_{\max}}$, we have
\begin{align}
    \osc(\varepsilon(z))\le \frac{15 R_{\max}}{|z|^5} \osc(\mathcal R).
\end{align}
\end{lemma}

The same QVE expansion also determines the typical locations around
which the outlier eigenvalues concentrate. We prove the next result in Appendix \ref{app:proof-prop-vartheta}.
\begin{prop}[Deterministic outlier location]\label{prop:vartheta}
Let $j\le r_+$ and assume $\lambda_j \ge 6\sqrt{R_{\max}}$. Define
\[
\widehat\alpha_j(z):=1-\lambda_j u_j^\T \Phi(z)u_j \quad\text{for}\quad  z\ge 3\sqrt{R_{\max}}.
\]
Then $\widehat\alpha_j$ has a unique real zero $\vartheta_j \in [3\sqrt{R_{\max}},\infty)$.
Moreover, this zero is simple and satisfies
\begin{align}\label{eq:det-outlier}
    \left| \vartheta_j - \lambda_j - \frac{\mathcal V_{jj}}{\lambda_j}\right|\le 145\frac{R_{\max}}{\lambda_j^3} \osc(\mathcal R).
\end{align}
\end{prop}
Since $\widehat\alpha_j(\lambda_j)$ is increasing and $\widehat\alpha_j(\lambda_j) \le 0$, we also have 
\begin{align}\label{eq:compare-outlier}
\vartheta_j \ge \lambda_j.    
\end{align}

The following result provides crude bounds on the operator norms of Green functions $G(z)=(zI-E)^{-1}$ and $G^{(k)} =(zI-E^{(k)})^{-1}$, where $E^{(k)}$ is obtained from $E$ with the $k$-th row and column replaced by zero. 
\begin{lemma}[Lemma 16 from \cite{OVW22}]\label{lem:bdG} For $z\in \C$ with $|z|\ge 2 \|E\|$, we have
\begin{align}\label{eq:crudebdG}
    \|G(z)\|\le \frac{2}{|z|},\quad \|G^{(k)}\|\le \frac{2}{|z|}.
\end{align}
\end{lemma}

\begin{lemma}[Bound on the noise matrix]\label{lem:bdE} Let $E=(E_{ij})_{n\times n}$ be a symmetric matrix whose entries are centered and independent (up to symmetry) random variables. Assume the entries $E_{ij}$'s are sub-Gaussian with $\|E_{ij}\|_{\psi_2}\le  K \sigma_{ij}.$ Denote $R_{\max}=\max_i \sum_j \sigma_{ij}^2$. 
Then for any $D>0$, we have
\[\Prob\left(\|E\|\le 2.9\sqrt{R_{\max}} + c K \sigma_{\max} {(D+2)\log n}\right) \ge 1-3n^{-D}\]
for an absolute constant $c>0$.
\end{lemma}
The result is a direct consequence of the sharp non-asymptotic norm bounds established in Bandeira and van Handel \cite{BvH16}. We provide a short derivation in Appendix \ref{app:lem-bdE}. 


\section{Proof strategy and technical inputs}\label{sec:technical}
This section states the probabilistic inputs used in the proof of
Theorem~\ref{thm:top-k-eigenspace} and explains how they combine with the
deterministic QVE estimates from Section~\ref{sec:qve-preliminaries}. The two
main probabilistic ingredients are the isotropic local laws for the
resolvent $G(z)$ and the location theorem for the outlier eigenvalues of
$A+E$.

\subsection{Isotropic local laws}
\label{subsec:local-law-inputs}
Consider $G(z)=(zI-E)^{-1}$ for $z\in \C^+$. Recall the QVE approximation $\Phi(z)$ defined in \eqref{def:QVE}. Our first main technical step is to establish an isotropic local law for $G(z)$. Denote 
\begin{align}
    \mathfrak m_D:=\begin{cases}
C_{\rm gen}'(D+6)^{3/2} \beta K\sigma_{\max}(\log n)^2,
& \text{in the sparse sub-Gaussian regime},\\[1mm]
C_{\rm bd}'(D+6)c_0 K\sigma_{\max}\log n,
& \text{in the sparse bounded-entry regime}.
\end{cases}
\end{align}
Here, $C_{\rm gen}'$ and $C_{\rm bd}'$ are sufficiently large absolute constants.
\begin{theorem}[Isotropic local law]\label{thm:local} Fix $D>0$. Suppose Assumptions~\ref{assump:noise} and \ref{assump:extra} hold. For any deterministic unit vectors $v,w\in\R^n$ and every fixed $z\in \C$ with $|z|\ge 6\sqrt{R_{\max}}$, we have
\[\left|v^\T (G(z)-\Phi(z))w \right|\le \frac{\mathfrak m_D}{|z|^2}\]
with probability at least $1-n^{-D-4}$.
\end{theorem}
Theorem \ref{thm:local} is proved in Section \ref{sec:local-law-proof}. We also need a version that is uniform over the spectral domain where the
outlier eigenvalues are located. Define
\begin{align}\label{def:S-domain}
    \mathcal S_{\rm out}:=
    \left\{ z\in\mathbb C: 6\sqrt{R_{\max}}\le |z|\le 2R_{\max}^3 \right\}.
\end{align}
This result is a direct consequence of Theorem \ref{thm:local}; its proof is provided in Appendix~\ref{app:local-law-uniform}. Recall that $A=U \Lambda U^\T$.
\begin{coro}[Uniform compressed local law]\label{lem:local-law-uniform} Fix $D>0$. Suppose the assumptions of Theorem~\ref{thm:local} hold.  Then, with probability at least $1-n^{-D}$,
\begin{align}
    &\sup_{z\in \mathcal{S}_{\rm out}} |z|^2 \| U^\T (G(z)-\Phi(z))U\| \le C' r\cdot{\mathfrak m_D} = M_D,\label{eq:uniform-compressed-local-law}\\
    & \sup_{z\in \mathcal{S}_{\rm out}} |z|^2 \max_{1\le i \le n}\| e_i^\T (G(z)-\Phi(z))U\| \le  M_D. \label{eq:uniform-row-signal-local-law}
\end{align}
Consequently, 
\begin{equation}\label{eq:uniform-row-signal-local-law-Q}
    \sup_{z\in\mathcal S_{\rm out}} |z|^2 \max_{1\le i\le n}
    \left\|  e_i^\T Q (G(z)-\Phi(z)) U \right\| \le 2 M_D,
\end{equation}
where $Q=I-UU^\T$  and $M_D$ is defined in \eqref{def:M}.
\end{coro}

We next establish a row-adaptive local law that retains the rowwise
variation of the deterministic test matrix. This refinement will be
important in the applications, where it allows us to reach the sharp
degree regime. For flexibility, we formulate the result for a
deterministic matrix $B\in\R^{n\times r}.$ The proof combines
a leave-one-out argument with the isotropic local law in
Theorem~\ref{thm:local}. The leave-one-out step requires a slightly
stronger variance-dominance condition. 

Fix $D>0$. For simplicity, define
\begin{equation}\label{eq:t-Dm}
    t_{D,r}:=C(D+r)\log n,
\end{equation}
where $C>0$ is a sufficiently large absolute constant. For the deterministic matrix $B\in\R^{n\times r}$, define
\begin{equation}\label{eq:gamma-def}
    \gamma_i(B)^2 := \sum_{j=1}^n \sigma_{ij}^2\|e_j^\T B\|^2,
    \qquad
    \Gamma(B):=\max_{1\le i\le n}\gamma_i(B).
\end{equation}
Note that $\gamma_i(B) =\left(\E\|e_i^\T EB\|^2\right)^{1/2}$
is the  standard-deviation scale for the $i$th row of $EB$,
while $\Gamma(B)$ is the largest such rowwise fluctuation scale.

Finally, define the row-adaptive spectral domain
\begin{equation} \label{eq:S-ra}
    \mathcal S_{\rm ra} :=
    \left\{ z\in\C: C_{\rm ra}\sqrt{t_{D,r}R_{\max}} \le |z|\le 2R_{\max}^3 \right\},
\end{equation}
where $C_{\rm ra}>0$ is a sufficiently large absolute constant.

We are now ready to state the row-adaptive local law. Recall that $L_{D,r}$ is defined in  \eqref{eq:L-ra}.
\begin{coro}[Uniform row-adaptive local law]\label{coro:row-adaptive}
Fix $D>0$. Suppose Assumption~\ref{assump:noise} and
Assumption~\ref{assump:row-adaptive} hold. Let
$B\in\R^{n\times r}$ be deterministic with $\|B\|\le n^{100}$.
Then, with probability at least $1-n^{-D}$,
\begin{equation}\label{eq:row-adaptive-main}
    \sup_{z\in\mathcal S_{\rm ra}}  |z|^2\|(G(z)-\Phi(z))B\|_{2,\infty}
    \le C\left[ \sqrt{t_{D,r}}\,\Gamma(B) + \left( t_{D,r}L_{D,r}^{\rm ra}  +\mathfrak m_D \right) \|B\|_{2,\infty} \right].
\end{equation}
Here $C>0$ is an absolute constant. 
\end{coro}

The proof of Corollary~\ref{coro:row-adaptive} is given in Appendix~\ref{app:local-law-uniform}.

\subsection{Outlier eigenvalue location}
\label{subsec:outlier-location}
We next state the outlier-location results used in the proof of the eigenspace perturbation results.  For any $\lambda_j\ge 6\sqrt{R_{\max}}$, Proposition \ref{prop:vartheta} guarantees that there exists a unique solution $\vartheta_j$ to the equation $\lambda_j^{-1}- u_j^\T \Phi(z)u_j=0$ for $z\ge 3\sqrt{R_{\max}}$. Besides, 
\[
\vartheta_j \approx \lambda_j + \frac{\mathcal V_{jj}}{\lambda_j}.
\]
For a simple supercritical eigenvalue $\lambda_k$, we will show that $\vartheta_k$ is the typical location of $\widetilde{\lambda}_k$.

Recall that
\[\rho_k\equiv\rho_{k,D} := 10 M_D + 15\frac{\osc(\mathcal{R})}{\lambda_k}.\]
Denote the index set for the supercritical eigenvalues as
\begin{align}\label{def:O+}
    \mathcal{O}_{+}:=\{l \in [r_+] : \lambda_l \ge 6\sqrt{R_{\max}}\}.
\end{align}
\begin{theorem}[Individual outlier eigenvalue location]\label{thm:outlier-ev} Fix $D>0$ and $k\in [r_+]$. If $\lambda_k \ge 6 \sqrt{R_{\max}} + 100\rho_k$ is a simple eigenvalue and the gap condition holds:
\begin{align*}
\min_{l <k } (\vartheta_l -\vartheta_k) &\ge 50 \rho_k \quad \text{and} \quad \min_{l\in\mathcal{O}_+,\, l >k } (\vartheta_k -\vartheta_l) \ge 50 \rho_k
\end{align*}
with the convention that the minimum over an empty set is $+\infty$,  then,  with probability at least $1-n^{-D}$,
\[ |\widetilde{\lambda}_k- \vartheta_k| \le \rho_{k}.\]
\end{theorem}
The proof of Theorem \ref{thm:outlier-ev} is given in Section \ref{sec:outlier-ev}. Furthermore, we provide the following result on the outlier eigenvalue clusters, which will be used in the proof of the top-$k$ eigenspace perturbations. 
\begin{theorem}[Lower edge of the top-$k$ outlier clusters]\label{thm:outlier-cluster} Fix $D>0$ and $k\in [r_+]$. Define 
\begin{align}
    \underline{\vartheta}_k:=\min_{1\le t \le k} \vartheta_t.
\end{align}
If $\lambda_k \ge 6 \sqrt{R_{\max}} + 100\rho_k$ and the gap condition holds:
\begin{align}\label{def:gap-cluster}
    \min_{l\in\mathcal{O}_+,\, l >k} (\underline{\vartheta}_k - \vartheta_l) \ge 50\rho_k
\end{align}
with the convention that the minimum over an empty set is $+\infty$, then, with probability at least $1-n^{-D}$,
\[
\widetilde\lambda_k \ge \underline{\vartheta}_k-\rho_k.
\]
In particular, for every $s\in [k]$,
\[
\widetilde\lambda_s\ge \underline{\vartheta}_k-\rho_k.
\]    
\end{theorem}
\begin{remark}
    Theorem \ref{thm:outlier-cluster} holds under a stronger but more conventional assumption that the $k$-th deterministic root lower-bounds the preceding roots (i.e., $\vartheta_l \ge \vartheta_k$ for all $l < k$). Under that assumption, $\underline{\vartheta}_k = \vartheta_k$. However, we introduce $\underline{\vartheta}_k$ to avoid the internal ordering of these $\vartheta_s$'s for $s\le k$. Indeed, the deterministic roots can physically cross and permute, even when the underlying $\lambda_s$'s are strictly ordered. For instance, even if $\lambda_l > \lambda_{s}$ for $l<s\in [k]$, as given in Proposition \ref{prop:vartheta}, $\vartheta_l \approx \lambda_l + \frac{\mathcal{V}_{ll}}{\lambda_l}$ could be dominated by $\vartheta_{s} \approx \lambda_{s} + \frac{\mathcal{V}_{ss}}{\lambda_s}$. This is because the shift $\mathcal{V}_{ll} = u_l^\T \mathcal{R} u_l$ depends on the alignment of the $l$-th population eigenvector with the variance profile $\mathcal{R}$. In the presence of anisotropic noise, it is entirely possible for a lower-ranked spike's variance projection (e.g., $\mathcal{V}_{ss}$) to be significantly larger than that of a higher-ranked spike (e.g., $\mathcal{V}_{ll}$). 
\end{remark}
The proof of Theorem \ref{thm:outlier-cluster} parallels that of Theorem \ref{thm:outlier-ev} and can be found in Section \ref{sec:outlier-ev}.

\subsection{Proof roadmap}
\label{subsec:proof-roadmap}
We briefly explain how the technical inputs above are used to prove
Theorem~\ref{thm:top-k-eigenspace}. The proof is carried out on a
high-probability event where $ \|E\|\le 3\sqrt{R_{\max}}$ and the isotropic local laws Theorem \ref{thm:local} and Corollary \ref{lem:local-law-uniform} hold.

The signal condition $\lambda_k\ge 9\sqrt{R_{\max}}$ 
ensures that the outlier eigenvalues satisfy
\begin{align}\label{eq:0605-outlier-crude}
   \widetilde\lambda_s \ge  \lambda_k-\|E\| \ge 6\sqrt{R_{\max}}
\end{align}
for all $s\in [k]$ by Weyl's inequalty. Thus the uniform local law Corollary \ref{lem:local-law-uniform} applies at the random eigenvalues $z=\widetilde\lambda_s$.

The deterministic part of the proof starts from the eigenvector equation
\[
    \widetilde u_s=G(\widetilde\lambda_s)A\widetilde u_s,
    \qquad s\in [k].
\]
Projecting this identity onto the non-target signal directions and the null space gives separate bounds for how much $\widetilde u_s$ extends outside the target space. The null-space projection produces the usual signal-to-noise contribution 
$\|E\|/\lambda_k$. For the signal-space terms, we decompose
\[
    G(z)=\Phi(z)+\Xi(z),\qquad \Xi(z):=G(z)-\Phi(z).
\]
The compressed local law controls the contribution of $\Xi(z)$, while the
large-$z$ expansion of $\Phi(z)$ identifies the deterministic
variance-profile contribution. This yields the stochastic term involving
$M_D$ and the geometric bias term $\mathcal B_k$. The final operator-norm
and $2\to\infty$ estimates follow from these blockwise bounds through deterministic reductions.

\section{Proof of the main perturbation bounds: Theorem \ref{thm:top-k-eigenspace}}\label{sec:top-k-eigenspace}
Set 
\begin{align}
    \Xi(z):=G(z)-\Phi(z).
\end{align}

Fix $D>0$. Define the events
\[
    \mathcal E_{\rm op}:=\{\|E\|\le 3\sqrt{R_{\max}}\}
\]
and 
\[
    \mathcal E_{\rm loc}(D):= \left\{ \sup_{z\in\mathcal S}  |z|^2\|U^\T \Xi(z) U\| \le M_D \right\}\cap \left\{ \sup_{z\in\mathcal S_{\rm out}} |z|^2 \max_{1\le i\le n}\left\|  e_i^\T Q \Xi(z) U \right\| \le 2 M_D \right\},
\]
where $\mathcal S_{\rm out}:=\left\{ z\in\mathbb C: 6\sqrt{ R_{\max}}\le |z|\le 2R_{\max}^3 \right\}$ and $Q=I-U U^\T$.
Set
\[
    \mathcal E_D:=\mathcal E_{\rm op}\cap \mathcal E_{\rm loc}(D).
\]
By Lemma~\ref{lem:bdE} and Corollary~\ref{lem:local-law-uniform}, we have
\[
    \mathbb P(\mathcal E_D)\ge 1-n^{-D}.
\]
All estimates below are deterministic on $\mathcal E_D$, except when the
outlier-location theorem is invoked. In that step we intersect with the
corresponding outlier-location event, which still has probability at least
$1-n^{-D}$ after adjusting constants.

On the event $\mathcal E_D$, note that our signal assumption $9\sqrt{R_{\max}} \le \lambda_k \le R_{\max}^3$ ensures $\widetilde{\lambda}_s \in \mathcal S_{\rm out}$ for all $s\le k$ by Weyl's inequality (see \eqref{eq:0605-outlier-crude}). 

For simplicity, we write 
\[M\equiv M_D.\]

We first present deterministic bounds on the difference between the eigenspaces $U_k$ and $\widetilde{U}_k$. Let $J=[r]\setminus[k]$ denote the complement index set and $U_J = (u_{k+1},\dots,u_r)$. Let $P_W$ be the orthogonal projection onto the subspace $W$.
\begin{prop}\label{prop:spacepert} For $1\le k \le r+$, assume $\lambda_k \ge 2\|E\|$. Then the following deterministic bounds hold.

\textup{(i)} For the operator norm bound, we have
\begin{align}\label{eq:sin-det}
    \|\sin\angle(U_k, \widetilde{U}_k)\| \le \|U_J^\T \widetilde{U}_k\|_F + 2\frac{\|E\|}{\lambda_k}.
\end{align}

\textup{(ii)} For the $2\to\infty$ norm bound, we have
\begin{align}\label{eq:row-det}
\min_{O\in \mathbb{O}(k)} \| \widetilde{U}_k - U_k O\|_{2,\infty} \le \| \widetilde{U}_k-P_{U_k} \widetilde{U}_k\|_{2,\infty} + \|U_k\|_{2,\infty} \|\sin\angle(U_k, \widetilde{U}_k)\|^2,
\end{align}
where \begin{align}\label{eq:row-proj-det}
   \| \widetilde{U}_k-P_{U_k} \widetilde{U}_k\|_{2,\infty} &\le \|U\|_{2,\infty} \cdot \|U_J^\T \widetilde{U}_k\|_F + \frac{8}{3}\sqrt{k}  \|U\|_{2,\infty} \frac{\osc(\mathcal{R})}{\widetilde{\lambda}_k^2}\nonumber\\
   &+ 2 \sqrt{\sum_{s=1}^k \widetilde{\lambda}_s^2  \max_{1\le i\le n}\left\|  e_i^\T Q \Xi(\widetilde{\lambda}_s) U \right\|^2}.
\end{align}
\end{prop}
 Given these deterministic bounds, the main technical contribution of our perturbation analysis is the estimation of \[\|U_J^\T \widetilde{U}_k\|_F = \sqrt{\sum_{s=1}^k \|U_J^\T \widetilde{u}_s\|^2}.\] We first present a deterministic bound for each $\|U_J^\T \widetilde{u}_s\|$. The geometric bias matrix $\mathcal{V} = U^\T \mathcal{R} U$ is defined in \eqref{def:V}. For any index sets $\mathtt I, \mathtt  J\in [r]$ with $U_{\mathtt I}=(u_i)_{i\in \mathtt  I}$ and $U_{\mathtt J} = (u_j)_{j\in \mathtt J}$, we denote 
\[ \mathcal{V}_{\mathtt I \mathtt J}:= U_{\mathtt I}^\T \mathcal{R} \,U_{\mathtt J}.\] For a matrix $B$, we denote by $\off(B)$ the matrix obtained by setting all diagonal entries to zero.
\begin{prop}\label{prop:detbd-vector} Define index sets 
\begin{align*}
    &\mathcal{J}:=\{k+1,\dots,r_{+}\}, \quad \mathcal{I}:=[r]\setminus \mathcal{J},\\
    &\mathcal{N}:=\{r_{+}+1,\dots,r\}, \quad \mathcal{K}:=[r]\setminus \mathcal{N}.
\end{align*} 
Fix $s\in [k]$. Define the gap
\[ \Delta_{\mathcal J}^{\Phi}(\widetilde{\lambda}_s) :=\widetilde{\lambda}_s\min_{j\in \mathcal{J}}|1-\lambda_j u_j^\T\Phi(\widetilde{\lambda}_s)u_j|\]
with the convention that the minimum over an empty set is $+\infty$. 
If the following gap condition holds:
\begin{align}\label{eq:gap-det-J}
    \Delta_{\mathcal J}^{\Phi}(\widetilde{\lambda}_s) \ge \frac{2\lambda_{k+1}}{\widetilde{\lambda}_s^2}\|\off(\mathcal{V}_{\mathcal J \mathcal{J}})\| + \frac{13\lambda_{k+1}R_{\max}}{\widetilde{\lambda}_s^4}\osc(\mathcal{R}),
\end{align}
then we have 
\begin{align}\label{eq:0506-individual}
   \|U_J^\T \widetilde{u}_s\| \le  &\frac{4}{\Delta_{\mathcal J}^{\Phi}(\widetilde{\lambda}_s)}\left(\frac{\|\mathcal{V}_{\mathcal{J}\mathcal{I}}\|}{\widetilde{\lambda}_s} + \frac{13}{2}\frac{R_{\max}}{\widetilde{\lambda}_s^3}\osc(\mathcal{R}) + \widetilde{\lambda}_s^2 \|U^\T \Xi(\widetilde{\lambda}_s)U\| \right) \\
   &+ 2\sqrt{3} \left( \frac{\|\mathcal{V}_{\mathcal{N}\mathcal{K}}\|}{\widetilde{\lambda}_s^2} + \frac{13}{2}\frac{R_{\max}}{\widetilde{\lambda}_s^4} \osc(\mathcal R) + \widetilde{\lambda}_s\|U^\T \Xi(\widetilde{\lambda}_s)U\|\right).
\end{align}  
\end{prop}
Note that when $\mathcal J= \emptyset$, we have $\Delta_{\mathcal J}^{\Phi}(\widetilde{\lambda}_s)=\infty$, so the corresponding term on the right-hand side of \eqref{eq:0506-individual} vanishes.  The proofs of Propositions \ref{prop:spacepert} and \ref{prop:detbd-vector} are deferred to Appendix \ref{app:perturb-prepare}. 
On the event $\mathcal E_D$, the condition $\lambda_k\ge9\sqrt{R_{\max}}$ implies $\lambda_k\ge3\|E\|\ge2\|E\|,$ so Proposition~\ref{prop:spacepert} applies.

Now we continue with the estimation of $\|U_J^\T \widetilde{u}_s\|$ and split the discussion into two cases: $k=r_+$ or $k<r_+$. 

\medskip
\noindent{\underline{Case 1.} Assume $k= r_+$}. Then $\mathcal J = \emptyset$ and $\Delta_{\mathcal J}^{\Phi}(\widetilde{\lambda}_s)=+\infty$. In this case, $$J=\mathcal{N}=\{r_{+}+1,\dots,r\}.$$ We just bound 
\[\|U_J^\T \widetilde{u}_s\| \le 2\sqrt{3} \left( \frac{\|\mathcal{V}_{\mathcal{N}\mathcal{K}}\|}{\widetilde{\lambda}_s^2} + \frac{13}{2}\frac{R_{\max}}{\widetilde{\lambda}_s^4} \osc(\mathcal R) + \widetilde{\lambda}_s\|U^\T \Xi(\widetilde{\lambda}_s)U\|\right).\]
By Theorem \ref{thm:outlier-cluster} (note that the gap condition \eqref{def:gap-cluster} holds trivially), for every $s\in [k]$,
\[
\widetilde\lambda_s\ge \underline{\vartheta}_k-\rho_k \ge \lambda_k -\rho_k \ge \frac{99}{100}\lambda_k.
\] 
Thus,
\begin{align}\label{eq:0601-lateruse}
    \|U_J^\T \widetilde{u}_s\|=\|U_{\mathcal N}^\T \widetilde{u}_s\| < 4 \frac{\|\mathcal{V}_{\mathcal{N}\mathcal{K}}\|}{\lambda_k^2} +24\frac{R_{\max}}{{\lambda}_k^4} \osc(\mathcal R)+ 4 \frac{M}{\lambda_k}
\end{align}
and 
\begin{align}\label{eq:0422-edge}
    \|U_J^\T \widetilde{U}_k\|_F = \sqrt{\sum_{s=1}^k \|U_J^\T \widetilde{u}_s\|^2} \le \sqrt{r_+}\left( 4 \frac{\|\mathcal{V}_{\mathcal{N}\mathcal{K}}\|}{\lambda_{r_+}^2} +24\frac{R_{\max}}{{\lambda}_{r_+}^4} \osc(\mathcal R)+ 4 \frac{M}{\lambda_{r_+}}\right).
\end{align}

\medskip
\noindent{\underline{Case 2.} Assume $k< r_+$}. Under our gap assumption
\begin{align}\label{eq:0421-gap}
\delta_k \ge 65 \rho_k + 20 \frac{ \osc(\mathcal R)}{\lambda_k} \ge 65 \rho_k + 20 \frac{\lambda_{k+1}}{\lambda_k^2} \osc(\mathcal R),
\end{align}
we first show that \eqref{def:gap-cluster} holds, that is, 
\[
\min_{l\in\mathcal{O}_+,\, l >k} (\underline{\vartheta}_k - \vartheta_l) \ge 50\rho_k \quad \text{where} \quad \underline{\vartheta}_k=\min_{1\le t \le k} \vartheta_t.
\]
If $\mathcal{O}_+\setminus[k] = \emptyset$, then the inequality holds because the left-hand side is $+\infty$ by our notation. Assume $\mathcal{O}_+\setminus[k]  \neq \emptyset$.

Suppose \[\underline{\vartheta}_k = \vartheta_{t}\] for some $1\le t\le k$. Observe from \eqref{eq:compare-outlier} that 
\begin{align}\label{eq:0420-simplefact}
    \vartheta_{t}\ge \lambda_t \ge \lambda_k > \lambda_{k+1} \ge \lambda_j
\end{align}
for any $j \in \mathcal{J}=\{k+1,\dots,r_{+}\}$.

Define two auxiliary functions
\[
\mathcal{L}(x):= x + \frac{R_{\min}}{x}-145\frac{R_{\max}}{x^3}\osc(\mathcal R),\quad \mathcal{U}(x):= x + \frac{R_{\max}}{x}+145\frac{R_{\max}}{x^3}\osc(\mathcal R)
\]
for $x\ge 6\sqrt{R_{\max}}$. Note that both $\mathcal{L}(x)$ and $\mathcal{U}(x)$ are increasing functions. Note that $\lambda_t\ge 9\sqrt{R_{\max}}$. Thus,
\[
\mathcal{L}(\lambda_t) \le \vartheta_t \le \mathcal{U}(\lambda_t).
\]
For any $l\in \mathcal{O}_+$ satisfying $l>k$, since $\lambda_l \le \lambda_{k+1}$, by the definition of $\vartheta_l$ for and Proposition \ref{prop:vartheta}, we have 
\[
\mathcal{L}(\lambda_l) \le \vartheta_l\le \mathcal{U}(\lambda_l) \le \mathcal{U}(\lambda_{k+1}).
\]
In particular, for any $l \in \mathcal{O}_+\setminus[k]$ (note that $t\in \mathcal{O}_+$), we have
\begin{align*}
    \vartheta_t - \vartheta_l &\ge \mathcal{L}(\lambda_t) - \mathcal{U}(\lambda_{k+1})\\
    &= (\lambda_t -\lambda_{k+1}) +\frac{R_{\min}}{\lambda_t} - \frac{R_{\max}}{\lambda_{k+1}} - 145\frac{R_{\max}}{\lambda_t^3}\osc(\mathcal R)- 145\frac{R_{\max}}{\lambda_{k+1}^3}\osc(\mathcal R).
\end{align*}
Rewriting
\[
\frac{R_{\min}}{\lambda_t} - \frac{R_{\max}}{\lambda_{k+1}}=\frac{\lambda_{k+1}-\lambda_t}{\lambda_t \lambda_{k+1}}R_{\min} + \frac{R_{\min}-R_{\max}}{\lambda_{k+1}}
\]
and using the assumption $\lambda_t,\lambda_{k+1} \ge 9\sqrt{R_{\max}}\ge 6\sqrt{R_{\min}}$, we further get
\begin{align*}
    \vartheta_t - \vartheta_l &\ge \left(1-\frac{R_{\min}}{\lambda_t\lambda_{k+1}} \right)(\lambda_t -\lambda_{k+1}) + \frac{R_{\min}-R_{\max}}{\lambda_{k+1}}- \frac{145}{36}\frac{\osc(\mathcal R)}{\lambda_t}- \frac{145}{36}\frac{\osc(\mathcal R)}{\lambda_{k+1}}\\
    &\ge \frac{35}{36} \delta_k - \left(1+\frac{145}{36}\right)\frac{\osc(\mathcal R)}{\lambda_{k+1}}- \frac{145}{36}\frac{\osc(\mathcal R)}{\lambda_t}\\
    &\ge \frac{35}{36} \delta_k - \frac{181}{36}\frac{\osc(\mathcal R)}{\lambda_{k+1}}- \frac{145}{36}\frac{\osc(\mathcal R)}{\lambda_k},
\end{align*}
where we used $\lambda_t \ge \lambda_k$ in the last inequality. To proceed, note that 
\[ 
\frac{\osc(\mathcal R)}{\lambda_{k+1}} = \frac{\osc(\mathcal R)}{\lambda_{k}} + \delta_k\frac{\osc(\mathcal R)}{\lambda_k\lambda_{k+1}} \le \frac{\osc(\mathcal R)}{\lambda_{k}} + \frac{1}{36}\delta_k.
\]
We get
\begin{align*}
    \vartheta_t - \vartheta_l &\ge \left(\frac{35}{36} - \frac{181}{36^2} \right)\delta_k - \frac{326}{36} \frac{\osc(\mathcal R)}{\lambda_{k}} \ge 50 \rho_k
\end{align*}
by our assumption \eqref{eq:0421-gap} on $\delta_k$. Hence, the conclusion of Theorem \ref{thm:outlier-cluster} holds. In particular, for every $s\in [k]$,
\begin{align}\label{eq:0421-LB-edge}
  \widetilde\lambda_s\ge \underline{\vartheta}_k-\rho_k = {\vartheta}_t-\rho_k.  
\end{align}

Recall $\Delta_{\mathcal J}^{\Phi}$ defined in Proposition \ref{prop:detbd-vector}. Next, we show that for any $s\in [k]$,
\begin{align}\label{eq:0601-Delta-lwbd}
    \Delta_{\mathcal J}^{\Phi}(\widetilde{\lambda}_s)=\widetilde{\lambda}_s\min_{j\in \mathcal{J}}|1-\lambda_j u_j^\T\Phi(\widetilde{\lambda}_s)u_j| \ge \frac{1}{2}\delta_k.
\end{align}
Denote $\alpha_j(z)= \lambda_j^{-1} - u_j^\T \Phi(z) u_j$ for each $j \in [r]$. Then 
\[
\Delta_{\mathcal J}^{\Phi}(\widetilde{\lambda}_s)=\min_{j\in \mathcal{J}}\big|\widetilde{\lambda}_s\lambda_j\alpha_j(\widetilde{\lambda}_s)\big|.
\]
For each $j \in \mathcal{J}$, by the expansion of $\Phi(z)$ in \eqref{eq:Phi-expansion}, we have 
\begin{align*}
    \alpha_j(z) = \frac{1}{\lambda_j} - u_j^\T \Phi(z) u_j = \frac{1}{\lambda_j}-\frac{1}{z} - \frac{\mathcal{V}_{jj}}{z^3} - \varepsilon_{jj}(z), 
\end{align*}
where $\varepsilon_{jj}(z):=u_j^\T \varepsilon(z) u_j$. 
Similarly,
\begin{align*}
    \alpha_t(z) = \frac{1}{\lambda_t}-\frac{1}{z} - \frac{\mathcal{V}_{tt}}{z^3} - \varepsilon_{tt}(z).
\end{align*}
From the difference $\alpha_j(z)-\alpha_t(z)$, we find
\begin{align*}
    \alpha_j(z) = \left(\frac{1}{\lambda_j} - \frac{1}{\lambda_t}\right) + \alpha_t(z) -\frac{\mathcal{V}_{jj}-\mathcal{V}_{tt}}{z^3} -(\varepsilon_{jj}(z)-\varepsilon_{tt}(z)).
\end{align*}
Therefore,
\begin{align*}
  \widetilde{\lambda}_s\lambda_j \alpha_j(\widetilde{\lambda}_s) &=  \widetilde{\lambda}_s\left(1- \frac{\lambda_j}{\lambda_t}\right) + \widetilde{\lambda}_s\lambda_j\alpha_t(\widetilde{\lambda}_s) -\frac{\lambda_j}{\widetilde{\lambda}_s^2}\left( \mathcal{V}_{jj}-\mathcal{V}_{tt}\right)-\widetilde{\lambda}_s\lambda_j(\varepsilon_{jj}(\widetilde{\lambda}_s)-\varepsilon_{tt}(\widetilde{\lambda}_s))\\
  &:= T_1 + T_2 + T_3 +T_4.
\end{align*}
It follows that
\begin{align*}
  \left| \widetilde{\lambda}_s\lambda_j \alpha_j(\widetilde{\lambda}_s) \right|\ge  T_1 + T_2 - |T_3| -|T_4|.
\end{align*}
We estimate $T_i$'s respectively. For $T_1$, by \eqref{eq:0420-simplefact}, we have
\begin{align*}
    T_1 = \widetilde{\lambda}_s\left(1- \frac{\lambda_j}{\lambda_t}\right) \ge \widetilde{\lambda}_s\left(1- \frac{\lambda_{k+1}}{\lambda_k}\right) = \frac{\widetilde{\lambda}_s}{\lambda_k} \delta_k.
\end{align*}
Then, by \eqref{eq:0421-LB-edge}, we get
\begin{align*}
    T_1\ge (\vartheta_t -\rho_k) \frac{\delta_k}{\lambda_k} = \frac{\vartheta_t}{\lambda_k}\delta_k - \rho_k  \frac{\delta_k}{\lambda_k} \ge \delta_k - \rho_k
\end{align*}
since $\frac{\vartheta_t}{\lambda_k}\ge 1$ by \eqref{eq:0420-simplefact} and $\frac{\delta_k}{\lambda_k}\le 1$.

Next, we estimate $T_2= \widetilde{\lambda}_s\lambda_j\alpha_t(\widetilde{\lambda}_s)$. By the same estimates as \eqref{eq:04112026-1} and \eqref{eq:04112026-2}, we have
\[\frac{5}{2x^2} \le \alpha_t'(x) \le \frac{2}{x^2}\quad \text{for } x \ge 3\sqrt{R_{\max}}.\]
In particular, $\alpha_t(x)$ is increasing for $x \ge 3\sqrt{R_{\max}}$. By our definition, $\alpha_t(\vartheta_t)$=0. Since $\widetilde{\lambda}_s \ge \vartheta_t - \rho_k$ by \eqref{eq:0421-LB-edge}, 
\[ \alpha_t(\widetilde{\lambda}_s) \ge \alpha_t (\vartheta_t - \rho_k)=\alpha_t (\vartheta_t - \rho_k)- \alpha_t(\vartheta_t)=-\rho_k \alpha'_t(\zeta), \]
where the last equation is due to the mean value theorem and $\zeta$ is a value between $\vartheta_t-\rho_k$ and $\vartheta_t$.
Then, by \eqref{eq:0420-simplefact},
\[
\widetilde{\lambda}_s\lambda_j \alpha_t(\widetilde{\lambda}_s) \ge -\widetilde{\lambda}_s\lambda_j\frac{2}{\zeta^2} \rho_k \ge - 2\frac{\widetilde{\lambda}_s\lambda_{k+1}}{(\vartheta_t-\rho_k)^2}\rho_k.
\]
Since $\vartheta_t-\rho_k \ge \frac{99}{100}\vartheta_t$ and $\widetilde{\lambda}_s \le \frac{3}{2}\lambda_s \le \frac{3}{2}\vartheta_s$ by Weyl's inequality and \eqref{eq:compare-outlier}, we have
\begin{align*}
    T_2 \ge -2 (\frac{100}{99})^2 \frac{\frac{3}{2}\vartheta_s \lambda_{k+1}}{\vartheta_t^2} \rho_k > - 4\rho_k.
\end{align*}
Finally, for $T_3$ and $T_4$, by Lemma \ref{lem:osc}, the bounds $\widetilde{\lambda}_s \le \frac{3}{2}\lambda_s$ and \[\widetilde{\lambda}_s \ge \vartheta_t -\rho_k \ge \frac{99}{100}\vartheta_t \ge \frac{99}{100} \lambda_k,\] we have
\begin{align*}
    |T_3| + |T_4| &\le \frac{\lambda_j}{\widetilde{\lambda}_s^2} \osc(\mathcal R) + \widetilde{\lambda}_s \lambda_j \cdot \osc(\varepsilon(\widetilde{\lambda}_s))\\
    &\le 4 \frac{\lambda_{k+1}}{\lambda_k^2} \osc(\mathcal R) + \frac{15\cdot 2^4}{36}\frac{\lambda_{k+1}}{\lambda_k^2} \osc(\mathcal R)=\frac{32}{3}\frac{\lambda_{k+1}}{\lambda_k^2} \osc(\mathcal R).
\end{align*}

Combining the above estimates, we arrive at
\begin{align*}
    \widetilde{\lambda}_s\lambda_j \alpha_j(\widetilde{\lambda}_s) \ge \delta_k -5\rho_k - \frac{32}{3}\frac{\lambda_{k+1}}{\lambda_k^2} \osc(\mathcal R) \ge \frac{1}{2}\delta_k
\end{align*}
by our gap assumption \eqref{eq:0421-gap}. Hence, 
\[
\Delta_{\mathcal J}^{\Phi}(\widetilde{\lambda}_s)=\min_{j\in \mathcal{J}}\big|\widetilde{\lambda}_s\lambda_j\alpha_j(\widetilde{\lambda}_s)\big| \ge \frac{1}{2}\delta_k.
\]
Furthermore, \eqref{eq:gap-det-J} in Proposition \ref{prop:detbd-vector} holds. To see this, note that the right-hand side of \eqref{prop:detbd-vector}
\begin{align*}
&\frac{2\lambda_{k+1}}{\widetilde{\lambda}_s^2}\|\off(\mathcal{V}_{\mathcal J \mathcal{J}})\| + \frac{13\lambda_{k+1}R_{\max}}{\widetilde{\lambda}_s^4}\osc(\mathcal{R})\\
&\le 2 \left(\frac{100}{99}\right)^4 \frac{ \lambda_{k+1}}{\lambda_k^2}\osc(\mathcal R) + \frac{13}{36}\left(\frac{100}{99}\right)^4 \frac{ \lambda_{k+1}}{\lambda_k^2}\osc(\mathcal R)\\
&< 3 \frac{ \lambda_{k+1}}{\lambda_k^2}\osc(\mathcal R) < \frac{1}{2}\delta_k \le \Delta_{\mathcal J}^{\Phi}(\widetilde{\lambda}_s).
\end{align*}
We continue from the bound established in Proposition \ref{prop:detbd-vector}. On the event $\mathcal{E}_D$, substituting the uniform local law bound $\|U^\T \Xi(\widetilde{\lambda}_s)U\| \le M/\widetilde{\lambda}_s^2$ and using  $\widetilde{\lambda}_s\ge \frac{99}{100}\lambda_k$, 
\begin{align}
   \|U_J^\T \widetilde{u}_s\| \le  &\frac{4}{\Delta_{\mathcal J}^{\Phi}(\widetilde{\lambda}_s)}\left(\frac{\|\mathcal{V}_{\mathcal{J}\mathcal{I}}\|}{\widetilde{\lambda}_s} + \frac{13}{2}\frac{R_{\max}}{\widetilde{\lambda}_s^3}\osc(\mathcal{R}) + \widetilde{\lambda}_s^2 \|U^\T \Xi(\widetilde{\lambda}_s)U\| \right) \nonumber \\
   &+ 2\sqrt{3} \left( \frac{\|\mathcal{V}_{\mathcal{N}\mathcal{K}}\|}{\widetilde{\lambda}_s^2} + \frac{13}{2}\frac{R_{\max}}{\widetilde{\lambda}_s^4} \osc(\mathcal R) + \widetilde{\lambda}_s\|U^\T \Xi(\widetilde{\lambda}_s)U\|\right)\label{eq:0625-columnbd}\\
   &\le \frac{8}{\delta_k}\left(\frac{100}{99} \frac{\|\mathcal{V}_{\mathcal{J}\mathcal{I}}\|}{\lambda_k} + \frac{13}{2}\left(\frac{100}{99}\right)^3 \frac{\osc(\mathcal R)R_{\max}}{\lambda_k^3} + M  \right)\nonumber\\
   & + 2\sqrt{3} \left( \left(\frac{100}{99}\right)^2 \frac{\|\mathcal{V}_{\mathcal{N}\mathcal{K}}\|}{{\lambda}_k^2} + \frac{13}{2}\left(\frac{100}{99}\right)\frac{ \osc(\mathcal R)R_{\max}}{{\lambda}_k^2} + \frac{100}{99}\frac{M}{\lambda_k}\right).\nonumber
\end{align}  
Now the bound above simplifies to 
\begin{align*}
   \|U_J^\T \widetilde{u}_s\| < \frac{10}{\delta_k \lambda_k}\left( \|\mathcal{V}_{\mathcal{J}\mathcal{I}}\|+8\frac{R_{\max}}{\lambda_k^2}\osc(\mathcal R)\right) +3 \frac{\|\mathcal{V}_{\mathcal{N}\mathcal{K}}\|}{\lambda_k^2} +10\frac{M}{\delta_k}.
\end{align*} 

Finally, we establish the bound 
\begin{align}\label{eq:0421-Ubd}
\|U_J^\T \widetilde{U}_k\|_F &= \sqrt{\sum_{s=1}^k \|U_J^\T \widetilde{u}_s\|^2}\\
&\le \frac{10\sqrt{k}}{\delta_k \lambda_k}\left( \|\mathcal{V}_{\mathcal{J}\mathcal{I}}\|+8\frac{R_{\max}}{\lambda_k^2}\osc(\mathcal R)\right) +3 \sqrt{k}\frac{\|\mathcal{V}_{\mathcal{N}\mathcal{K}}\|}{\lambda_k^2} +10\sqrt{k}\frac{M}{\delta_k}.
\end{align}

Combining \eqref{eq:0422-edge}, for any $k\in[r_+]$, we always have
\begin{align}\label{eq:0424-biasB}
    \|U_J^\T \widetilde{U}_k\|_F \le \mathcal{B}_k+10\sqrt{k}\frac{M}{\delta_k},
\end{align}
where 
\[\mathcal{B}_k:=\frac{10\sqrt{k}}{\delta_k \lambda_k}\left( \|\mathcal{V}_{\mathcal{J}\mathcal{I}}\|+8\frac{R_{\max}}{\lambda_k^2}\osc(\mathcal R)\right) +4 \sqrt{k}\frac{\|\mathcal{V}_{\mathcal{N}\mathcal{K}}\|}{\lambda_k^2} .\]

Consequently, by \eqref{eq:sin-det},
\begin{align*}
    \|\sin\angle(U_k, \widetilde{U}_k)\|
    \le \mathcal{B}_k + 10\sqrt{k}\frac{M}{\delta_k}+2\frac{\|E\|}{\lambda_k}.
\end{align*}

From \eqref{eq:row-proj-det}, using \eqref{eq:0424-biasB}, on the event $\mathcal E_D$,  we also have 
\begin{align*}
   &\| \widetilde{U}_k-P_{U_k} \widetilde{U}_k\|_{2,\infty} \\
   &\le \|U\|_{2,\infty} \cdot \|U_J^\T \widetilde{U}_k\|_F + \frac{8}{3}\sqrt{k} \|U\|_{2,\infty} \frac{\osc{(\mathcal R)}}{\widetilde{\lambda}_k^2}+ 2 \sqrt{\sum_{s=1}^k \widetilde{\lambda}_s^2  \max_{1\le i\le n}\left\|  e_i^\T Q \Xi(\widetilde{\lambda}_s) U \right\|^2}\\
   &\le \|U\|_{2,\infty} \left(\mathcal{B}_k + 10\sqrt{k}\frac{M}{\delta_k}+3\sqrt{k} \frac{\osc{(\mathcal R)}}{\lambda_k^2} \right) + 2M\sqrt{\sum_{s=1}^k \frac{1}{\widetilde{\lambda}_s^2}}\\
   &\le \|U\|_{2,\infty} \left(\mathcal{B}_k + 10\sqrt{k}\frac{M}{\delta_k}+3\sqrt{k} \frac{\osc{(\mathcal R)}}{\lambda_k^2} \right) + 3 \sqrt{k} \frac{M}{\lambda_k}
\end{align*}
using $\widetilde{\lambda}_k \ge \frac{99}{100}\lambda_k$. Now, from \eqref{eq:row-det}, we obtain that 
\begin{align*}
    &\min_{O\in \mathbb{O}(k)} \| \widetilde{U}_k - U_k O\|_{2,\infty}\\
    &\le \| \widetilde{U}_k-P_{U_k} \widetilde{U}_k\|_{2,\infty} + \|U_k\|_{2,\infty} \|\sin\angle(U_k, \widetilde{U}_k)\|^2\\
    &\le 4\|U\|_{2,\infty} \left( \mathcal{B}_k +20\sqrt{k}\frac{M}{\delta_k}+ 3\sqrt{k} \frac{\osc{(\mathcal R)}}{\lambda_k^2} +27 \frac{R_{\max}}{\lambda_k^2}\right) + 3 \sqrt{k} \frac{M}{\lambda_k}
\end{align*}
by our assumptions on $\lambda_k$ and $\delta_k$.

\section{Proof of the isotropic local law: Theorem \ref{thm:local}}\label{sec:local-law-proof}
Let $z\in \C$ satisfy $|z|\ge 6\sqrt{R_{\max}}$. Recall that 
\[ G(z)=(zI-E)^{-1}.\]
Denote 
\begin{equation}\label{def:H}
    H\equiv H(z)=\mathrm{diag} (h_1(z),\dots,h_n(z)),\qquad h_i\equiv h_i(z) = \sum_{j=1}^n \sigma_{ij}^2 \phi_{j}(z).
\end{equation}
Then the QVE can be written as \[\Phi(z)=(zI-H)^{-1}.\]

Since $(zI-E)G=I$ and $(zI-H)\Phi=I$, subtracting these identities, we have
\[z(G-\Phi) = EG -H\Phi =(E-H)G + H(G-\Phi).\]
Rearranging the terms, we further obtain $(zI-H)(G-\Phi) = (E-H)G$. Hence, for deterministic unit vectors $v,w\in\mathbb R^n$,  
\begin{equation}\label{eq:local-law-basic-decomp}
v^\T(G-\Phi)w = v^\T(z-H)^{-1} (E-H)Gw.
\end{equation}
Set \[\widehat{v}^\T:= v^\T(z-H)^{-1}.\] Then $$v^\T(G-\Phi)w =\widehat{v}^\T (E-H)Gw.$$ 
We also define
\begin{equation}\label{def:L}
    L\equiv L(z) = \mathrm{diag}(l_1(z),\dots,l_n(z)),\qquad l_i\equiv l_i(z) = \sum_{j=1}^n \sigma_{ij}^2 G_{jj}(z).
\end{equation}
Then \eqref{eq:local-law-basic-decomp} gives the decomposition
\begin{equation}\label{eq:local-law-two-terms}
    v^\T (G-\Phi)w =
    \widehat v^\T (E-L)Gw + \widehat v^\T (L-H)Gw.
\end{equation}

We start with a bound on $\|\widehat{v}\|$.  By Lemma~\ref{lem:bdphi},
\[
    \|H\| \le R_{\max}\max_i|\phi_i(z)|  \le \frac{3R_{\max}}{2|z|} \le \frac{|z|}{4},
\]
where the last inequality follows from $|z|^2\ge 36R_{\max}$. Hence
\begin{equation}\label{eq:bdhatv}
    \|\widehat v\| \le \|(zI-H)^{-1}\| \le \frac{1}{|z|-\|H\|} \le \frac{4}{3|z|}.
\end{equation}

With this $\|\widehat{v}\|$ bound, the control of the first error term in \eqref{eq:local-law-two-terms} reduces to the following proposition, which we prove in Appendix \ref{app:local-part}. The proof follows the cumulant expansion strategy of He, Knowles and Rosenthal~\cite{HKR18}, which compares the random resolvent with its deterministic self-consistent approximation. Our application is simpler in two ways. First, we work only in the large-$z$ outlier regime, where resolvent derivatives have elementary bounds. Second, we need only a second-order cumulant expansion with a controlled third-order remainder. The main challenge is handling highly inhomogeneous and sparse variance profiles.
\begin{prop}
\label{prop:bdDvw}
Fix $D>0$. Let $v,w \in\mathbb R^n$ be deterministic unit vectors, and let
$z\in\mathbb C_+$ satisfy $|z|\ge 6\sqrt{R_{\max}}$. Under
Assumptions~\ref{assump:noise} and \ref{assump:extra}, with probability
at least $1-n^{-D-4}$,
\begin{equation}\label{eq:bdDvw}
    \left| v^\T(E-L)G w \right|
    \le \frac{\mathfrak m_D}{|z|} + Cn^{-500}.
\end{equation}
Here, the $Cn^{-500}$ term appears only in the general sub-Gaussian regime; in
the bounded sparse-entry regime, it may be omitted.
\end{prop}

The second error term in \eqref{eq:local-law-two-terms} is controlled using  Proposition \ref{prop:bdDvw} and a deterministic stability argument for the
diagonal self-consistent equation.  Its proof is given in Appendix \ref{app:local-part2}.
\begin{prop}\label{prop:local-part2}
Fix $D>0$. Under the assumption of Propositon \ref{prop:bdDvw}, with probability at least $1-n^{-D-4}$:   
\begin{align*}
    \left|{v}^\T(L-H)Gw \right| \le C\left(\frac{\mathfrak m_D}{|z|}+ n^{-500}\right).
\end{align*}
In the bounded sparse-entry regime, the $n^{-500}$ term may be omitted.
\end{prop}

\begin{proof}[Proof of Theorem~\ref{thm:local}]
We use the decomposition \eqref{eq:local-law-two-terms}. For the first error term in \eqref{eq:local-law-two-terms},
Proposition~\ref{prop:bdDvw} and \eqref{eq:bdhatv} yield
\begin{equation}\label{eq:first-term-local}
    \left| \widehat v^\T(E-L)Gw \right|
    \le
    \|\widehat v\|
    \left( \frac{\mathfrak m_D}{|z|} + Cn^{-500}\right)\le 
    C\frac{\mathfrak m_D}{|z|^2}+ C\frac{n^{-500}}{|z|}.
\end{equation}

Similarly, for the second error term in \eqref{eq:local-law-two-terms}, from Proposition~\ref{prop:local-part2} and \eqref{eq:bdhatv}, we get
\begin{equation}\label{eq:second-term-local}
    \left| \widehat v^\T (L-H)Gw \right|  \le
    C\frac{\mathfrak m_D}{|z|^2} + C\frac{n^{-500}}{|z|}.
\end{equation}
Combining \eqref{eq:first-term-local} and \eqref{eq:second-term-local} with
\eqref{eq:local-law-two-terms}, we obtain
\[
    \left| v^\T (G-\Phi)w \right| \le  C\frac{\mathfrak m_D}{|z|^2} + C\frac{n^{-500}}{|z|}.
\]
Since $R_{\max} \le n^{100}$ and $\sigma_{\max}^{-1} \le n^{100}$,  we have $\mathfrak m_D \ge c n^{-100}$ for an absolute constant $c>0$. Absorbing the $n^{-500}$-term, we get
\[
    \left| v^\T(G-\Phi)w \right| \le \frac{\mathfrak m_D}{|z|^2}.
\]
This completes the proof.
\end{proof}


\section{Proofs of the outlier eigenvalue locations: Theorems \ref{thm:outlier-ev} and \ref{thm:outlier-cluster}}\label{sec:outlier-ev}
\subsection{Proof of Theorem \ref{thm:outlier-ev}}
For notational convenience throughout this section, we replace the target index $k$ with $j$. We first consider the case that $\lambda_j\ge 6\sqrt{R_{\max}}$ for $j \in [r_+]$ is a simple eigenvalue. In this section, we use the normalized version
\[
    \alpha_j(z):=\lambda_j^{-1}-u_j^\T \Phi(z)u_j,
\]
so that $\widehat\alpha_j(z)=\lambda_j\alpha_j(z)$. The deterministic
root $\vartheta_j$ defined in Proposition \ref{prop:vartheta} is equivalently characterized by $\alpha_j(\vartheta_j)=0$.

We still work on the event $\mathcal E_D$. However, we only need a weaker signal assumption $\lambda_j \ge 6 \sqrt{R_{\max}} + 100\rho_j$ since the local law is applied on deterministic contours around $\lambda_j$, whereas in the eigenspace proof, it is applied at the random points $z=\widetilde\lambda_s.$

By Proposition \ref{prop:vartheta}, there exists a unique deterministic location $\vartheta_j$ satisfying \[1-\lambda_j u_j^\T \Phi(\vartheta_j)u_j=0\]
with \[\left| \vartheta_j - \lambda_j - \frac{\mathcal V_{jj}}{\lambda_j}\right|\le 145 \frac{R_{\max}^2}{\lambda_j^3}.\]
In particular, since $\mathcal V_{jj}\ge 0$, we have 
\begin{align}\label{eq:vartheta-upper}
    \vartheta_j\le \lambda_j + \frac{\mathcal V_{jj}}{\lambda_j}+145 \frac{R_{\max}^2}{\lambda_j^3} \le \lambda_j + \frac{1}{36}\lambda_j + \frac{145}{36^2}\lambda_j< 1.16 \lambda_j
\end{align}
and directly by Proposition \ref{prop:vartheta},
\begin{align}\label{eq:vartheta-lower}
\vartheta_j \ge \lambda_j
\end{align}

We will show that on the event $\mathcal E_D$, there exists exactly one eigenvalue $\widetilde{\lambda}_j$ lying inside the disk \[\mathtt{D}(\vartheta_j, \rho):=\{z\in \C : |z-\vartheta_j|\le \rho \}\] for $\rho \equiv \rho_j$ defined in \eqref{def:radius}. 
Denote 
\begin{align*}
    K(z):=\Lambda^{-1} -U^\T \Phi(z) U \quad \text{and} \quad B(z) := U^\T(G(z)-\Phi(z))U.
\end{align*}
Set
\begin{align}
    g(z):=\det(K(z)) \quad\text{and} \quad f(z)= \det(\Lambda^{-1} -U^\T G(z) U)=\det(K(z)-B(z)).
\end{align}
By Lemma 21 from \cite{OVW22}, the eigenvalues of $\widetilde{A}$ outside the $[-\|E\|,\|E\|]$ are the zeros of $f(z)$. Besides, the algebraic multiplicity of each eigenvalue matches the corresponding multiplicity of each zero. Next, we show that the zeros of $f(z)$ are close to the zeros of $g(z)$. 

The following result is a special case of the generalized Rouch\'e's theorem for operator functions (see \cite{gohberg1971operator}). We include a short proof for completeness. 
\begin{lemma}\label{lem:rouche}
    Let $\Gamma \subseteq \C$ be a simple closed contour. Assume $K(z)$ and $B(z)$ are analytic in a neighborhood of $\Gamma$ and its interior. If 
    \[ \sup_{z\in \Gamma} \|K^{-1}(z) B(z)\| <1,\]
    then $f$ and $g$ have the same number of zeros inside $\Gamma$, counting multiplicity. 
\end{lemma}
\begin{proof}
    For $t\in [0,1]$, define $F_t(z):=\det(K(z)-t B(z))$. For $z\in \Gamma$,
    \[K(z)-tB(z) = K(z) (I-t K^{-1}(z) B(z)).\]
    Since $ \|tK^{-1}(z) B(z)\|<1$, $I-t K^{-1}(z) B(z)$ is invertible. Hence, $F_t(z)\neq 0$ on $\Gamma$ for all $t\in [0,1].$ By the argument principle (\cite[Section 5.2]{ahlfors1979complex}), the number of zeros of $F_t$ inside $\Gamma$ is a constant in $t$. Taking $t=0$ and $t=1$ proves the claim. 
\end{proof}

We will take a circle \[\Gamma_j(\rho):=\{z\in \C:|z-\vartheta_j| = \rho \}\] with the radius 
\begin{align}\label{def:radius}
    \rho \equiv\rho_j := 10 M + 15\frac{\osc(\mathcal{R})}{\lambda_j}.
\end{align}
To apply Lemma \ref{lem:rouche} on $\Gamma_j(\rho)$, we bound \[\sup_{z\in \Gamma_j(\rho)} \|K^{-1}(z) B(z)\| \le \sup_{z\in \Gamma_j(\rho)} \|K^{-1}(z)\| \cdot \sup_{z\in \Gamma_j(\rho)} \| B(z)\|.\]
To bound $\sup_{z\in \Gamma_j(\rho)} \|B(z)\|$, first note that $\Gamma_j(\rho) \subseteq \mathcal S_{\rm out}$. To see this, since $\lambda_j\ge6\sqrt{R_{\max}}+100\rho_j$ and 
$\rho_j\le\lambda_j/100$, for every
$z\in\Gamma_j(\rho_j)$, 
\begin{align*}
    &|z| \ge \vartheta_j-\rho_j \ge  \lambda_j-\rho_j  \ge 6\sqrt{R_{\max}},\\
    &|z| \le \vartheta_j+\rho_j  \le 1.16\lambda_j+\frac{\lambda_j}{100}  \le 2R_{\max}^3
\end{align*}
Hence, on the event $\mathcal E_D$, we get
\begin{align}\label{eq:est-B-circle}
    \sup_{z\in \Gamma_j(\rho)} \|B(z)\| \le \left(\frac{100}{99} \right)^2 \frac{M}{\vartheta_j^2}.
\end{align}

To bound $\sup_{z\in \Gamma_j(\rho)} \|K^{-1}(z)\|$, we rewrite $K(z)$ for $z\in \Gamma_j(\rho)$ in block form. Without loss of generality, by relabeling the columns of $U$ and the diagonal entries of $\Lambda$, we may assume that the index $j$ is the first one. Namely, write
\[ U=(u_j,U_{-j}) \quad \text{and}\quad \Lambda=\begin{pmatrix}
    \lambda_j &0 \\
    0&\Lambda_{-j}
\end{pmatrix}.\]
Thus,
\[K(z)= \Lambda^{-1} -U^\T \Phi(z) U = \begin{pmatrix}
\alpha_j(z) & -b_j(z)^\T\\
-b_j(z) & D_j(z)
\end{pmatrix},
\]
where
\[\alpha_j(z):=\lambda_j ^{-1}-u_j^\T \Phi(z)u_j, \quad b_j(z):= U_{-j}^\T \Phi(z) u_j, \quad D_j(z) = \Lambda_{-j}^{-1} -U_{-j}^\T \Phi(z) U_{-j}.\]
The following lemma provides the estimates of $|\alpha_j(z)|, \|b_j(z)\|,$ and the smallest singular value $s_{\min}(D_j(z))$ on the circle $\Gamma_j(\rho)$. Its proof is deferred to Section \ref{sec:est-circle}.
\begin{lemma}\label{lem:est-circle} Fix $j\in [r_+]$. Assume $\lambda_j \ge 6\sqrt{R_{\max}} + 100\rho$ and $\lambda_j$ is a simple eigenvalue of $A$. Assume the following gap condition holds:
\begin{align}\label{def:gap-con1}
    \min_{l <j } (\vartheta_l -\vartheta_j) &\ge 50 \rho\quad \text{and} \quad \min_{l\in\mathcal{O}_+,\, l >j } (\vartheta_j -\vartheta_l) \ge 50 \rho.
\end{align}
where $\mathcal{O}_{+}:=\{l \in [r_+] : \lambda_l \ge 6\sqrt{R_{\max}}\}$ is the index set for the supercritical eigenvalues. 
Then for $z\in \Gamma_j(\rho)$, we have
\begin{align}
    &\frac{\rho}{4\vartheta_j^2}\le |\alpha_j(z)| \le \frac{3\rho}{\vartheta_j^2},\label{eq:est-alpha}\\
    &\|b_j(z)\|\le \frac{1}{10}\frac{\rho}{\vartheta_j^2}, \label{eq:est-bj}\\
    &s_{\min}(D_j(z)) \ge  3\frac{\rho}{\vartheta_j^2}\label{eq:est-sv-Dinverse}.
\end{align}
\end{lemma}

In particular, \eqref{eq:est-sv-Dinverse} suggests that $D_j(z)$ is invertible on $\Gamma_j(z)$ and 
\begin{align}\label{eq:bound-D-inverse}
  \|D^{-1}_j(z)\| \le \frac{1}{3} \frac{\vartheta_j^2}{\rho}.  
\end{align}
Define \[h_j(z):= \alpha_j(z) - b_j(z)^\T D_j^{-1}(z) b_j(z)\] for $z\in \Gamma_j(\rho)$. By Lemma \ref{lem:est-circle}, $h_j(z)$ is well-defined on $\Gamma_j(\rho)$. Besides, 
\[ |h_j(z)| \ge |\alpha_j(z)| - \|b_j(z)\|^2\|D_j^{-1}(z)\| \ge \frac{\rho}{4\vartheta_j^2} -\frac{1}{300} \frac{\rho}{\vartheta_j^2}=\frac{37}{150}\frac{\rho}{\vartheta_j^2}>0\]
and 
\begin{align}\label{eq:bd-h-inverse}
    |h_j^{-1}(z)| \le \frac{150}{37} \frac{\vartheta_j^2}{\rho}.
\end{align}
Now we are ready to estimate $\|K^{-1}(z)\|$ on $\Gamma_j(\rho)$. By Schur's complement,
\begin{align*}
K(z)^{-1}=\begin{pmatrix}
h_j^{-1}(z) & h_j^{-1}(z) b_j(z)^\T D_j^{-1}(z)\\
h_j^{-1}(z) D_j^{-1}(z) b_j(z)& D_j^{-1}(z)+ h_j^{-1}(z) D_j^{-1}(z)b_j(z) b_j(z)^\T D_j^{-1}(z)
\end{pmatrix}
\end{align*}
We decompose $K(z)^{-1}$ as
\begin{align*}
K(z)^{-1}= \begin{pmatrix}
    0 & 0\\
    0 & D_j^{-1}(z)
\end{pmatrix}
+ h_j^{-1}(z)\begin{pmatrix}
1 &  b_j(z)^\T D_j^{-1}(z)\\
D_j^{-1}(z) b_j(z)&  D_j^{-1}(z)b_j(z) b_j(z)^\T D_j^{-1}(z)
\end{pmatrix}.
\end{align*}
Note that by setting $\mathtt{v}:=b_j(z)^\T D_j^{-1}(z)$, we have
\[\mathtt{H}:=\begin{pmatrix}
1 &  b_j(z)^\T D_j^{-1}(z)\\
D_j^{-1}(z) b_j(z)&  D_j^{-1}(z)b_j(z) b_j(z)^\T D_j^{-1}(z)
\end{pmatrix} = \begin{pmatrix}
    1 & \mathtt{v}^\T\\
    \mathtt{v} & \mathtt{v}\mathtt{v}^\T
\end{pmatrix} =\begin{pmatrix}
    1\\
    \mathtt{v}
\end{pmatrix} \begin{pmatrix}
    1 & \mathtt{v}^\T
\end{pmatrix}.\]
Hence,
\[\|\mathtt{H}\|= 1+ \|\mathtt{v}\|^2 = 1 + \|b_j(z)^\T D_j^{-1}(z)\|^2 \le 1+ \|b_j(z)\|^2\cdot \|D_j^{-1}(z)\|^2.\]
It follows that
\begin{align}\label{eq:est-K-Bound}
    \|K(z)^{-1}\| &\le \|D_j^{-1}(z)\| + |h_j^{-1}(z)| \cdot\|\mathtt{H}\|\\
    &\le \|D_j^{-1}(z)\| + |h_j^{-1}(z)| \left( 1+ \|b_j(z)\|^2\cdot \|D_j^{-1}(z)\|^2\right).
\end{align}
Plugging the bounds in Lemma \ref{lem:est-circle}, \eqref{eq:bound-D-inverse}, and \eqref{eq:bd-h-inverse}, we get
\begin{align}\label{eq:bound-Ki-circle}
    \|K(z)^{-1}\| \le \frac{1}{3} \frac{\vartheta_j^2}{\rho} + \frac{150}{37} \frac{\vartheta_j^2}{\rho}\left(1+ \frac{1}{900} \right)< \frac{9}{2} \frac{\vartheta_j^2}{\rho}.
\end{align}
Together with \eqref{eq:est-B-circle}, we arrive at 
\begin{align*}
    \sup_{z\in \Gamma_j(\rho)} \|K^{-1}(z) B(z)\| \le \sup_{z\in \Gamma_j(\rho)} \|K^{-1}(z)\| \cdot \sup_{z\in \Gamma_j(\rho)} \| B(z)\| \le \frac{9}{2} \frac{\vartheta_j^2}{\rho} \cdot \left(\frac{100}{99} \right)^2 \frac{M}{\vartheta_j^2} <1,
\end{align*}
where we used $\rho \ge 5 M$ in the last inequality. Hence, by
Lemma~\ref{lem:rouche}, the functions $f(z)=\det(K(z)-B(z))$ and  $g(z)=\det (K(z))$
have the same number of zeros inside $\Gamma_j(\rho)$, counting
multiplicities. This comparison alone, however, does not yet identify that
number, since $K(z)$ is not diagonal.  We therefore perform the zero
counting by comparing $K(z)$ with its diagonal approximation.

Recall that $\mathtt{D}(\vartheta_l, \rho):=\{z\in \C : |z-\vartheta_l|\le\rho \}$ denotes the disk centered as $\vartheta_l$ with radius $\rho$. Set 
\begin{align}\label{def:D-pm}
    D_j^+ := \bigcup_{l<j}\, {\mathtt{D}(\vartheta_l,\rho)},
\qquad
D_j^- := D_j^+ \cup {\mathtt{D}(\vartheta_j,\rho)}, 
\end{align}
and let $\mathcal C_j^\pm:=\partial D_j^\pm$ be positively oriented on each connected component. By our gap assumption \eqref{def:gap-con1}, the disk ${\mathtt{D}(\vartheta_j,\rho)}$ is disjoint from $D_j^+$. Moreover, the
disks corresponding to roots $\vartheta_l$ with $l<j$ lie strictly to the
right of $ \mathtt D( \vartheta_j, \rho)$, while the supercritical roots
$ \vartheta_l$ with $ l>j $ lie strictly to its left.

Define the diagonal matrix of $K(z)$ as
\begin{align}\label{def:A0}
    A_0(z):=\diag(\alpha_1(z),\dots,\alpha_r(z)) \quad\text{with} \quad \alpha_l(z):=\lambda_l^{-1} - u_l^\T \Phi(z) u_l.
\end{align}
We show the next result whose proof is deferred to Section \ref{sec:compareAK}
\begin{lemma}\label{lem:compareAK} Under the assumption of Lemma \ref{lem:est-circle}, we have
    \[
\sup_{z\in \mathcal C_j^\pm}\|A_0(z)^{-1}(K(z)-A_0(z))\|<1,
\qquad
\sup_{z\in \mathcal C_j^\pm}\|K(z)^{-1}B(z)\|<1.
\]
\end{lemma}
Applying Lemma~\ref{lem:rouche} twice on the contours
$\mathcal C_j^\pm $, first to $A_0 (z)$ and $K(z)$, and then to
$K(z)$ and $ K(z)-B(z)$, shows that
\[
    \det (A_0(z)),\qquad g(z)=\det(K(z)),\qquad f(z)=\det(K(z)-B(z))
\]
have the same number of zeros inside each of the regions $D_j^+$ and
$ D_j^-$, counting multiplicities.

Now
\[\det (A_0(z))=\prod_{l=1}^r \alpha_l(z),\]
and among these factors on the right-hand side, only $\alpha_l(z)$ with $l\in\mathcal O_+$ have zeros in the supercritical regime, namely at $z=\vartheta_l$ (see Proposition \ref{prop:vartheta}). Since
$\lambda_j$ is supercritical, all $\lambda_l $ with $ l<j $ are also
supercritical. Hence, by construction, $\det (A_0(z))$ has exactly $ j-1$ zeros
inside $ D_j^+ $ and exactly $ j$ zeros inside $ D_j^-$. Consequently,
$ f(z)$ has exactly $ j-1$ zeros inside $ D_j^+$ and exactly $ j$ zeros
inside $D_j^- $.

Since $D_j^-\setminus D_j^+=\mathtt{D}(\vartheta_j,\rho)$ and $\Gamma_j(\rho)=\partial\, \mathtt{D}(\vartheta_j,\rho)$, it follows that $f(z)=\det(K(z)-B(z))$ has exactly one zero in $\mathtt{D}(\vartheta_j,\rho)$. This zero is the same zero enclosed by
$\Gamma_j(\rho)=\partial\mathtt D(\vartheta_j,\rho)$. Denote it by
$\widehat\lambda$. As noted above, the zeros of $f$ outside
$[-\|E\|,\|E\|] $ coincide with the eigenvalues of $ \widetilde A = A+E$, with
matching algebraic multiplicities. Moreover, by \eqref{eq:vartheta-lower} and
$ \rho\le \lambda_j/100$, the disk
$ \mathtt D(\vartheta_j,\rho)$ lies strictly to the right of the noise bulk
$ [-\|E\|,\|E\|]$. Therefore $ \widehat\lambda$ is an
eigenvalue of $\widetilde A$.

It remains only to identify its index. The zero count above shows that exactly
$ j-1$ zeros of $f$ lie in $ D_j^+$, which is strictly to the right of
$ \mathtt D(\vartheta_j,\rho)$. The gap condition also places all remaining
supercritical roots $ \vartheta_l$ with $ l>j $ strictly to the left of
$ \mathtt D(\vartheta_j,\rho)$; applying the same comparison on the
corresponding separating contours rules out additional zeros of $f$ in the
gaps between these root neighborhoods. Finally, all non-outlier eigenvalues of
$ \widetilde A $ lie in $[-\|E\|,\|E\|]$, which is also strictly to the left
of $ \mathtt D(\vartheta_j,\rho) $. Hence there are exactly $ j-1$ eigenvalues
of $\widetilde A$ larger than $\widehat\lambda$. Hence,
\[
\widehat\lambda=\widetilde\lambda_j.
\]

\subsubsection{Proof of Lemma \ref{lem:compareAK}}\label{sec:compareAK}
 We estimate $\|A_0(z)^{-1}\|$, $\|K(z)-A_0(z)\|$, and $\|K^{-1}(z)\|$ respectively for $z\in \mathcal{C}_j^{\pm}$. Since
\[ K(z)-A_0(z) = -\off(U^\T \Phi(z) U),
\]
we have
\[ \|K(z)-A_0(z)\| \le \|\off(U^\T \Phi(z) U)\|.\]
\noindent{{\bf Claim.}} For $z\in \mathcal{C}_j^{\pm}$, the following estimates hold:
\begin{align}
    & |\alpha_l(z)| \ge \frac{1}{2}\frac{\rho}{|z|^2},\label{eq:0414-alpha}\\
    &\|b_j(z)\| \le \frac{2}{25}\frac{\rho}{|z|^2},\label{eq:0414-bj}\\
    &\|\off(U^\T \Phi(z) U)\| \le \frac{1}{5}\frac{\rho}{|z|^2}, \label{eq:0414-Phi}\\
     &\|\off(U^\T \varepsilon(z) U)\| \le \frac{1}{10}\frac{\rho}{|z|^2}\label{eq:0414-epsilon}.
\end{align}
We defer the proof of this claim to the end of this section. Using these bounds, we first have
\begin{align*}
    \|A_0(z)^{-1}\| =\max_{l\in [r]} |\alpha_l^{-1}(z) = \frac{1}{|\min_{l\in [r]}\alpha_l(z)|}  \le 2\frac{|z|^2}{\rho}
\end{align*}
Hence, 
\begin{align*}
    \sup_{z\in \mathcal C_j^\pm}\|A_0(z)^{-1}(K(z)-A_0(z))\| &\le \sup_{z\in \mathcal C_j^\pm}\|A_0(z)^{-1}\| \sup_{z\in \mathcal C_j^\pm}\|K(z)-A_0(z))\|\\
    & \le 2\frac{|z|^2}{\rho} \cdot \frac{1}{5}\frac{\rho}{|z|^2}<1.
\end{align*}

Next, we estimate $\|K(z)^{-1}\|$ for $z\in \mathcal{C}_j^{\pm}$, which is similar to the previous estimates on $\Gamma_j(\rho)$. Since $\Big\|\off\big(U_{-j}^\T \Phi(z) U_{-j}\big)\Big\|\le \Big\|\off\big(U^\T \Phi(z) U\big)\Big\|$, from \eqref{eq:smin-D1} and by the \textbf{Claim}, we get
\begin{align}
             s_{\min}(D_j(z))&\ge \min_{l \in [r]\setminus \{j\}} |\alpha_l(z)| - \Big\|\off\big(U_{-j}^\T \Phi(z) U_{-j}\big)\Big\|\\
             &\ge \frac{1}{2}\frac{\rho}{|z|^2} - \frac{1}{5}\frac{\rho}{|z|^2} = \frac{3}{10} \frac{\rho}{|z|^2}.
\end{align}
and thus
\begin{align*}
    \|D_j^{-1}(z)\|\le \frac{10}{3} \frac{|z|^2}{\rho}.
\end{align*}

Recall that $h_j(z)= \alpha_j(z) - b_j(z)^\T D_j^{-1}(z) b_j(z)$. By the claim,
\[
|h_j(z)| \ge |\alpha_j(z)| - \|b_j(z)\|^2\|D_j^{-1}(z)\| \ge \frac{1}{2}\frac{\rho}{|z|^2} - \left(\frac{2}{25}\frac{\rho}{|z|^2}\right)^2\frac{3}{10} \frac{|z|^2}{\rho} =\frac{359}{750}\frac{\rho}{|z|^2}.
\]
Hence, $|h_j^{-1}(z)| \le \frac{750}{359}\frac{|z|^2}{\rho}$ and by \eqref{eq:est-K-Bound}, 
\begin{align*}
     \|K(z)^{-1}\| 
    &\le \|D_j^{-1}(z)\| + |h_j^{-1}(z)| \left( 1+ \|b_j(z)\|^2\cdot \|D_j^{-1}(z)\|^2\right)\\
    &\le \frac{10}{3} \frac{|z|^2}{\rho} + \frac{750}{359}\frac{|z|^2}{\rho}\left( 1+  \left(\frac{2}{25}\frac{\rho}{|z|^2}\right)^2 \left(\frac{3}{10} \frac{|z|^2}{\rho} \right)^2\right)\\
    &< 5.5 \frac{|z|^2}{\rho}.
\end{align*}
Finally, 
\begin{align*}
    \sup_{z\in \mathcal{C}_j^{\pm}} \|K^{-1}(z) B(z)\| \le \sup_{z\in \mathcal{C}_j^{\pm}} \|K^{-1}(z)\| \cdot \sup_{z\in \mathcal{C}_j^{\pm}} \| B(z)\| \le 5.5 \frac{|z|^2}{\rho} \frac{M}{|z|^2}<1,
\end{align*}
where we used $\rho \ge 10 M$ in the last inequality.

\medskip
To finish the proof, it remains to prove the claim. 
\paragraph{Proof of \eqref{eq:0414-alpha}} We bound $|\alpha_l(z)|$ for $z\in \mathcal{C}_j^{\pm}$. The proof is similar to that of \eqref{eq:0415-previous-alpha}. We split the discussion into three cases: $l \in \mathcal{O}_+,$ $l\in  [r_+]\setminus\mathcal{O}_+$, and $l \in [r]\setminus[r_+]$.

\medskip
\noindent{\emph{Case 1:} $l \in \mathcal{O}_{+}$.} We have $\lambda_l \ge 6\sqrt{R_{\max}}$. By Proposition \ref{prop:vartheta}, define $\vartheta_l$ as the unique solution to $\alpha_l(z)=0$. Also, $\vartheta_l \ge 5 \sqrt{R_{\max}}$.

For $z=x+\sqrt{-1} y \in \mathcal{C}_j^{\pm}$, we have $|y|\le \rho$ and $x\ge \vartheta_j - \rho \ge \frac{99}{100}\vartheta_j > 4\sqrt{R_{\max}}$. Note that using the same estimates as \eqref{eq:04112026-1} and \eqref{eq:04112026-2}, we have
\begin{align}\label{eq:0415-diff-alpha}
    \frac{4}{5} \frac{1}{x^2} \le \alpha'_l(x)=-u_l^\T\Phi'(x)u_l \le \frac{3}{x^2}.
\end{align}
Indeed, we refined the constant for the lower bound using $x\ge 4\sqrt{R_{\max}}$:
\begin{align*}
   \alpha_l' (x) = -u_l^\T \Phi'(x) u_l \ge \frac{1}{x^2} - \frac{243}{5} \frac{R_{\max}^2}{x^6} \ge \frac{1}{x^2}  \left( 1- \frac{243}{5\cdot 256}\right) > \frac{4}{5x^2}.
\end{align*}
Since $\alpha_l(\vartheta_l)=0$, from $|\alpha_l(x)| = |\int_{\vartheta_l}^x \alpha'_l(t)\,dt|$, we get
\begin{align}\label{eq:0415-alpha}
 \frac{4}{5} \frac{|x-\vartheta_l|}{x \vartheta_l}\le   |\alpha_l(x)| \le 3\frac{|x-\vartheta_l|}{x \vartheta_l}.
\end{align}
Now the Taylor expansion of $\alpha_l(z)$ at $x$ yields:
\[ \alpha_l(z) = \alpha_l(x) + \sqrt{-1} y \alpha_l'(x) + \frac{1}{2}\alpha_l''(\zeta) (\sqrt{-1} y)^2, \]
where $\zeta$ lies in a vertical segment connecting $x$ and $z$. In particular, $|\zeta| \ge x \ge 4\sqrt{R_{\max}}$. It follows from part (i) of Lemma~\ref{lem:bdphi} that
\[|\alpha_l''(\zeta)|=|u_l^\T \Phi''(\zeta) u_l|\le \|\Phi''(\zeta) \| \le \frac{28}{|\zeta|^3} \le \frac{28}{x^3}.\]
Next, plugging in \eqref{eq:0415-diff-alpha} and \eqref{eq:0415-alpha} yields
\begin{align}
    | \alpha_l(x) + \sqrt{-1} y \alpha_l'(x) | &= \sqrt{\alpha_l(x)^2 + y^2 \alpha_l'(x)^2}\\
    & \ge \sqrt{ \frac{16}{25} \frac{(x-\vartheta_l)^2}{x^2 \vartheta_l^2} + \frac{16}{25} \frac{y^2}{x^4} } = \frac{4}{5} \frac{1}{x^2} \sqrt{ (x-\vartheta_l)^2 \frac{x^2}{\vartheta_l^2} + y^2 }.
\end{align}
If $\vartheta_l \le \frac{100}{99}x$, then $x/\vartheta_l \ge 0.99$ and 
\begin{align*}
    | \alpha_l(x) + \sqrt{-1} y \alpha_l'(x) | \ge \frac{4}{5 x^2} 0.99\sqrt{(x-\vartheta_l)^2 + y^2} = 0.792\frac{|z-\vartheta_l|}{x^2}\ge 0.792\frac{\rho}{x^2}.
\end{align*}
If $\vartheta_l > \frac{100}{99}x$, we use the fact that the function $f(t)=\frac{t-x}{xt}$ is increasing for $t>x$. Its minimum value occurs at the boundary $t=\frac{100}{99}x$. By \eqref{eq:0415-alpha},
\begin{align*}
    | \alpha_l(x) + \sqrt{-1} y \alpha_l'(x) | \ge | \alpha_l(x)| \ge \frac{4}{5} \frac{\vartheta_l-x}{x \vartheta_l} \ge \frac{4}{5} \frac{\frac{100}{99}x-x}{x \cdot\frac{100}{99}x} = \frac{4x}{500 x^2}.
\end{align*}
Since $x \ge 99\rho$, we have $\frac{4x}{500x^2} \ge \frac{4 \cdot 99\rho}{500x^2} = 0.792 \frac{\rho}{x^2}$.

Hence, we obtain
\begin{align*}
    |\alpha_l(z)| &\ge |\alpha_l(x) + \sqrt{-1} y \alpha_l'(x)| - 14\frac{y^2}{x^3}\\
    & \ge 0.792\frac{\rho}{x^2} - \frac{14}{99} \frac{\rho}{x^2}> \frac{\rho}{2 x^2} \ge \frac{\rho}{2 |z|^2} ,
\end{align*}
where we used $|z|\ge x\ge 99\rho \ge 99|y|$.


\medskip
\noindent{\emph{Case 2:} $l \in [r_+]\setminus \mathcal{O}_{+}$.} The lower bound on $|\alpha_l(z)|$ follows the same approach as \eqref{eq:alpha-l-2}. We skip the proof. In this case, we have
\[|\alpha_l(z)| \ge 70 \frac{\rho}{\vartheta_j^2}\ge 70 \left(\frac{99}{100}\right)^2 \frac{\rho}{|z|^2}> 60\frac{\rho}{|z|^2}.\]

\medskip
\noindent{\emph{Case 3:} $l\in [r] \setminus [r_+]$.} We have $\lambda_l <0$. For $z = x+\sqrt{-1}y \in \mathcal{C}_j^\pm$, we have $x \ge \vartheta_j - \rho \ge \frac{99}{100}\vartheta_j  4\sqrt{R_{\max}}$. From \eqref{eq:Phi-expansion}, we expand
\begin{align*}
    \alpha_l(z) = \frac{1}{\lambda_l} - u_l^\T \Phi(z) u_l = \frac{1}{\lambda_l} - \frac{1}{z} - \frac{1}{z^3} u_l^\T \mathcal{R} u_l - u_l^\T \varepsilon(z) u_l.
\end{align*}
Note that 
\[\left| \frac{1}{\lambda_l} - \frac{1}{z} \right| \ge \left| \Re\left( \frac{1}{\lambda_l} - \frac{1}{z} \right) \right| = \frac{1}{|\lambda_l|} + \frac{x}{|z|^2} > \frac{x}{|z|^2}.\]
Furthermore, we have
\begin{align*}
    \left|\frac{1}{z^3} u_l^\T \mathcal{R} u_l + u_l^\T \varepsilon(z) u_l\right| \le \frac{R_{\max}}{|z|^3} + \frac{45}{8}\frac{R_{\max}^2}{|z|^5} \le \frac{173}{128}\frac{R_{\max}}{|z|^3} \le \frac{173}{128\cdot 16} \frac{x}{|z|^2}
\end{align*}
using the bound $|z|\ge x \ge 4 \sqrt{R_{\max}}$. Hence, 
\begin{align*}
    |\alpha_l(z)| \ge \frac{x}{|z|^2} -\frac{173}{128\cdot 16} \frac{x}{|z|^2} > 90 \frac{\rho}{|z|^2}.
\end{align*}
Combining the above estimates, we have proved \eqref{eq:0414-alpha}.

\paragraph{Proof of \eqref{eq:0414-bj}} The proof is almost identical to that of \eqref{eq:est-bj}. We briefly it here. Following the same approach as \eqref{eq:est-bj}, we first get
\begin{align*}
    \|b_j(z)\| \le \frac{\osc(\mathcal R)}{2|z|^3} + \frac{15 R_{\max}}{2|z|^5} \osc(\mathcal R)< \frac{\osc(\mathcal R)}{|z|^3}.
\end{align*}
Since \(|z|\ge \vartheta_j-\rho \ge \frac{99}{100}\vartheta_j\) and
\[\rho\ge 13\frac{\osc(\mathcal R)}{\vartheta_j},\]
we further get
\begin{align}\label{eq:0416-ratio}
    \frac{\osc(\mathcal R)}{|z|^3}
\le \frac{100}{99}\frac{\osc(\mathcal R)}{\vartheta_j}\frac{1}{|z|^2}
\le \frac{100}{99\cdot 13}\frac{\rho}{|z|^2},
\end{align}
and hence
\[
\|b_j(z)\|\le \frac{2}{25}\frac{\rho}{|z|^2}.
\]

\paragraph{Proofs of \eqref{eq:0414-Phi} and \eqref{eq:0414-epsilon}} The bound for $\off\big(U^\T \Phi(z) U\big)$ are almost identical to the bound for $\off\big(U_{-j}^\T \Phi(z) U_{-j}\big)$ derived  in \eqref{eq:est-sv-Dinverse}. We briefly sketch it here. 
It follows from  the same steps that
\[ 
\|\off(U_{-j}^\T \varepsilon(z) U_{-j})\| \le \frac{15 R_{\max}}{\sqrt{2}|z|^5} \osc(\mathcal R)\le \frac{15}{9\sqrt{2}}\frac{\osc(\mathcal R)}{|z|^3} \le \frac{1}{10}\frac{\rho}{|z|^2}
\]
and
\begin{align*}
    \big\|\off\big(U_{-j}^\T \Phi(z) U_{-j}\big)\big\| &\le \frac{\osc(\mathcal R)}{2|z|^3} + \frac{15 R_{\max}}{\sqrt{2}|z|^5} \osc(\mathcal R)\\
    &\le \frac{\osc(\mathcal R)}{|z|^3}\left(\frac 12 + \frac{15}{9\sqrt{2}} \right)\le  \frac 15\frac{\rho}{|z|^2} 
\end{align*}
where we used \eqref{eq:0416-ratio}.

\subsubsection{Proof of Lemma \ref{lem:est-circle}}\label{sec:est-circle}
Recall from \eqref{def:radius} that \[\rho = 5 M + 15\frac{\osc(\mathcal{R})}{\lambda_j}.\]  By our assumption, $\rho \le \lambda_j/100$. From the estimates in \eqref{eq:vartheta-lower}, we also have \[\rho\ge 15\cdot0.87\frac{\osc(\mathcal{R})}{\vartheta_j}>13\frac{\osc(\mathcal{R})}{\vartheta_j}.\]
For $z\in \Gamma_j(\rho)$ satisfying $|z-\vartheta_j|=\rho$, we have 
\[ \frac{99}{100}\vartheta_j \le \vartheta_j-\rho\le |z| \le \vartheta_j-\rho \le \frac{101}{100}\vartheta_j.\]

\paragraph{Proof of \eqref{eq:est-alpha}.} We first estimate $|\alpha_j(z)|$. Recall that $\alpha_j (\vartheta_j)=0$. By Taylor's theorem,
\[\alpha_j(z)= \alpha_j' (\vartheta_j)(z-\vartheta_j) + \mathtt{E}_2,\] where \[|\mathtt{E}_2| \le \frac{1}{2}|z-\vartheta_j|^2 \max_{|\zeta-\vartheta_j| \le \rho} |\alpha_j''(\zeta)|= \frac{1}{2}\rho^2 \max_{|\zeta-\vartheta_j| \le \rho} |\alpha_j''(\zeta)|.\]
By part (ii) of Lemma~\ref{lem:bdphi} (see also \eqref{eq:diff-f}), 
\begin{align}\label{eq:04112026-1}
   \alpha_j' (\vartheta_j) = -u_j^\T \Phi'(\vartheta_j) u_j \ge \frac{1}{\vartheta_j^2} - \frac{243}{5} \frac{R_{\max}^2}{\vartheta_j^6} \ge \frac{1}{\vartheta_j^2}  \left( 1- \frac{243}{5\cdot 81}\right) = \frac{2}{5\vartheta_j^2}.
\end{align}
For $\zeta$ satisfying $|\zeta-\vartheta_j| \le \rho$, we have $|\zeta| \ge \vartheta_j -\rho \ge \frac{99}{100}\vartheta_j > 4 \sqrt{R_{\max}}$. It follows from part (i) of Lemma~\ref{lem:bdphi} that
\[|\alpha_j''(\zeta)|=|u_j^\T \Phi''(\zeta) u_j|\le \|\Phi''(\zeta) \| \le \frac{28}{|\zeta|^3} \le 28\left(\frac{100}{99}\right)^3  \frac{1}{\vartheta_j^3}< \frac{29}{\vartheta_j^3}.\]
Thus, \[|\mathtt{E}_2| \le \frac{29}{2} \frac{\rho^2}{\vartheta_j^3} \le \frac{29}{200} \frac{\rho}{\vartheta_j^2}\] and
\[ |\alpha_j(z)| \ge \frac{2\rho}{5\vartheta_j^2}-\frac{29}{200} \frac{\rho}{\vartheta_j^2} > \frac{\rho}{4\vartheta_j^2}.\]
For the upper bound, observe that
\[ |\alpha_j(z)|= |\alpha_j(z)-\alpha_j(\vartheta_j)| \le \max_{|\zeta-\vartheta_j| \le \rho} |\alpha_j'(\zeta)| \cdot \rho.\]
By part (ii) of Lemma~\ref{lem:bdphi} (see also \eqref{eq:diff-f}), we have 
\begin{align}\label{eq:04112026-2}
|\alpha_j'(\zeta)| \le \frac{1}{|\zeta|^2} + \frac{3\mathcal V_{jj}}{|\zeta|^4} + \frac{243}{5} \frac{R_{\max}^2}{|\zeta|^6} \le \frac{2}{|\zeta|^2}\le \frac{3}{\vartheta_j^2}
\end{align}
using $|\zeta| \ge \frac{99}{100}\vartheta_j$ and $\vartheta_j > 5\sqrt{R_{\max}}$. Hence,
\[|\alpha_j(z)| \le \frac{3\rho}{\vartheta_j^2}. \]

\paragraph{Proof of \eqref{eq:est-bj}.} Next, we estimate $\|b_j(z)\| = \|U_{-j}^\T \Phi(z) u_j\|$. Since $U_{-j}^\T u_j=0$, using the expansion of $\Phi(z)$ in \eqref{eq:Phi-expansion}, we immediately get
\begin{align*}
    b_j(z) = \frac{1}{z^3} U_{-j}^\T \mathcal{R} u_j + U_{-j}^\T \varepsilon(z) u_j.
\end{align*}
Thus,
\begin{align*}
    \|b_j(z)\| \le \frac{\|U_{-j}^\T \mathcal{R} u_j\|}{|z|^3} + \|U_{-j}^\T \varepsilon(z) u_j\|.
\end{align*}
Since $U_{-j}^\T u_j=0$, 
\[ \|U_{-j}^\T \mathcal{R} u_j\| = \|U_{-j}^\T (\mathcal{R}-cI) u_j\| \le \|\mathcal{R}-cI\| \le \frac{1}{2}\osc(\mathcal R)\]
by selecting $c= \frac{1}{2}(R_{\max} + R_{\min})$. Similarly, we also have 
\[\|U_{-j}^\T \varepsilon(z) u_j\| \le \frac{1}{2}\osc(\varepsilon(z))\le \frac{15 R_{\max}}{2|z|^5} \osc(\mathcal R),\]
where the last inequality follows from Lemma \ref{lem:osc}. Thus
\begin{align*}
    \|b_j(z)\| \le \frac{\osc(\mathcal R)}{2|z|^3} + \frac{15 R_{\max}}{2|z|^5} \osc(\mathcal R)< \frac{\osc(\mathcal R)}{|z|^3}\le \left(\frac{100}{99}\right)^2\frac{\osc(\mathcal R)}{\vartheta_j^3} <\frac{1}{10}\frac{\rho}{\vartheta_j^2} 
\end{align*}
by our assumption $\rho \ge 15\frac{\osc(\mathcal R)}{\lambda_j}$.

\paragraph{Proof of \eqref{eq:est-sv-Dinverse}.} Finally, we show the lower bound on $s_{\min}(D_j(z)).$ For each $l\in [r]$, set 
\[\alpha_l(z):=\lambda_l ^{-1}-u_l^\T \Phi(z)u_l.\] We rewrite
\begin{align*}
    D_j(z) &= \Lambda_{-j}^{-1} -U_{-j}^\T \Phi(z) U_{-j} \\
    &=\Lambda_{-j}^{-1} -\diag(u_l^\T\Phi(z)u_l)_{l \in [r]\setminus \{j\}} -\off\big(U_{-j}^\T \Phi(z) U_{-j}\big)\\
    &= \diag(\alpha_l(z))_{l \in [r]\setminus \{j\}} -\off\big(U_{-j}^\T \Phi(z) U_{-j}\big).
\end{align*}
Thus,
\begin{align}\label{eq:smin-D1}
    s_{\min}(D_j(z)) \ge \min_{l \in [r]\setminus \{j\}}|\alpha_l(z)| - \Big\|\off\big(U_{-j}^\T \Phi(z) U_{-j}\big)\Big\|.
\end{align}
We first bound $\off\big(U_{-j}^\T \Phi(z) U_{-j}\big)$. By \eqref{eq:Phi-expansion}, 
\begin{align}\label{eq:bd-off-Phi1}
    \off\big(U_{-j}^\T \Phi(z) U_{-j}\big) &= \off\big(\frac{1}{z}I + \frac{1}{z^3}U_{-j}^\T \mathcal{R} U_{-j} + U_{-j}^\T \varepsilon(z) U_{-j}\big)\nonumber\\
    &=\off\big(\frac{1}{z^3}U_{-j}^\T \mathcal{R} U_{-j}\big) + \off\big(U_{-j}^\T \varepsilon(z) U_{-j}\big).
\end{align}
Note that \[\off(U_{-j}^\T (\mathcal{R}-c_1I) U_{-j})=\off(U_{-j}^\T \mathcal{R} U_{-j}).\]
Choosing $c_1=\frac{1}{2}(R_{\max} + R_{\min})$, we get
\[ \Big\|\off\big(U_{-j}^\T \mathcal{R} U_{-j}\big)\Big\| \le \|\mathcal{R}-c_1I\| \le \frac{1}{2}\osc(\mathcal R).\]
Similarly, we also have \[\|\off(U_{-j}^\T \varepsilon(z) U_{-j})\|=\|\off(U_{-j}^\T (\varepsilon(z)-c_2I) U_{-j})\| \le \|\varepsilon(z)-c_2 I\|.\] 
Decompose $\varepsilon(z) =\diag(\varepsilon_k(z))= \diag(x_k) + \sqrt{-1}\diag(y_k)$. Set $c_2 = c_x+ \sqrt{-1}c_y=\frac{1}{2}(\max_k x_k +\min_k x_k) + \sqrt{-1} (\max_k y_k +\min_k y_k).$ Then 
\[ |\varepsilon_k(z) - c_2 |^2 = |x_k-c_x|^2 + |y_k-c_y|^2| \le \frac{1}{4}\osc(x)^2 + \frac{1}{4}\osc(y)^2\le \frac{1}{2}\osc(\varepsilon(z))^2.\]
It follows that
\[ \|\off(U_{-j}^\T \varepsilon(z) U_{-j})\| \le \frac{1}{\sqrt{2}}\osc(\varepsilon(z)) \le \frac{15 R_{\max}}{\sqrt{2}|z|^5} \osc(\mathcal R).\]
The last inequality is by Lemma \ref{lem:osc}. Combining these estimates with \eqref{eq:bd-off-Phi1} yields
\begin{align}\label{eq:2nd-term-Di}
    \big\|\off\big(U_{-j}^\T \Phi(z) U_{-j}\big)\big\| &\le \frac{\osc(\mathcal R)}{2|z|^3} + \frac{15 R_{\max}}{\sqrt{2}|z|^5} \osc(\mathcal R)\\
    &\le \frac{\osc(\mathcal R)}{|z|^3}\left(\frac 12 + \frac{15}{9\sqrt{2}} \right)\\
    &< 1.8\frac{\osc(\mathcal R)}{\vartheta_j^3} \le \frac{1.8}{13} \frac{\rho}{\vartheta_j^2}
\end{align}
where we used $|z|\ge \frac{99}{100}\vartheta_j$ and $\rho> 13\frac{\osc(\mathcal R)}{\vartheta_j}$.

It remains to bound $\min_{l \in [r]\setminus \{j\}}|\alpha_l(z)|$ where \[\alpha_l(z)=\lambda_l^{-1} - u_l^\T \Phi(z) u_l.\] Define \[\mathcal{O}_{+}:=\{l \in [r_+] : \lambda_l \ge 6\sqrt{R_{\max}}\}.\] We split the discussion into three cases: $l \in \mathcal{O}_{+}$, $l \in [r_+]\setminus \mathcal{O}_{+}$ and $l\in [r] \setminus [r_+]$.

\medskip
\noindent{\emph{Case 1:} $l \in \mathcal{O}_{+}$.} We have $\lambda_l \ge 6\sqrt{R_{\max}}$. By Proposition \ref{prop:vartheta}, define $\vartheta_l$ as the unique solution to $\alpha_l(z)=0$. Also, $\vartheta_l \ge 5 \sqrt{R_{\max}}$. Using similar estimates as \eqref{eq:04112026-1}, we get, for any real $t\ge 3\sqrt{R_{\max}}$,
\[\alpha_l'(t)=- u_l^\T \Phi'(t) u_l \ge \frac{2}{5t^2}.\]
Since $\alpha_l(\vartheta_l)=0$, 
\begin{align*}
    |\alpha_l(\vartheta_j)| =\left|\int_{\vartheta_l}^{\vartheta_j} \alpha_l'(t)\,dt  \right| \ge \left|\int_{\vartheta_l}^{\vartheta_j} \frac{2}{5t^2}\,dt \right|=\frac{2}{5} \frac{|\vartheta_l -\vartheta_j|}{\vartheta_l \vartheta_j} \ge \frac{2}{5} \frac{|\vartheta_l -\vartheta_j|}{\vartheta_j(\vartheta_j + |\vartheta_l -\vartheta_j|)}.
\end{align*}
The last step is from $\vartheta_l \le \vartheta_j + |\vartheta_l -\vartheta_j|.$  
Since $|\vartheta_l -\vartheta_j| \ge 50\rho$ by our assumption \eqref{def:gap-con1}, together with $\rho\le \vartheta_j /100$, we get
\[ |\alpha_l(\vartheta_j)|  \ge \frac{40}{3} \frac{\rho}{\vartheta_j^2}.\]

Next, using the same estimates as in \eqref{eq:04112026-2}, we also have 
\[|\alpha_l(\vartheta_j)-\alpha_l(z)| \le \max_{|\zeta-\vartheta_j| \le \rho} |\alpha_l'(\zeta)| \cdot \rho \le 3\frac{\rho}{\vartheta_j^3}.\]
Hence, 
\begin{align}\label{eq:alpha-l-1}
   |\alpha_l(z)| \ge  |\alpha_l(\vartheta_j)| - |\alpha_l(\vartheta_j)-\alpha_l(z)| \ge \left(\frac{40}{3} -3 \right)\frac{\rho}{\vartheta_j^2} = \frac{10}{3}\frac{\rho}{\vartheta_j^2}.
\end{align}

\medskip
\noindent{\emph{Case 2:} $l \in [r_+]\setminus \mathcal{O}_{+}$.} We have $0<\lambda_l< 6\sqrt{R_{\max}}\le \lambda_j$. From
\[\alpha_l(z) - \alpha_j(z) = \left(\lambda_l^{-1} - \lambda_j^{-1}\right) - \left(u_l^\T \Phi(z) u_l - u_j^\T \Phi(z) u_j \right),\]
observe that 
\begin{align}\label{eq:compare-alpha-lj}
    |\alpha_l(z) - \alpha_j(z)| &\ge \frac{\lambda_j-\lambda_l}{\lambda_j\lambda_l} -\left| u_l^\T \Phi(z) u_l - u_j^\T \Phi(z) u_j \right|.
\end{align}
Furthermore, by \eqref{eq:Phi-expansion},
\begin{align*}
    \left| u_l^\T \Phi(z) u_l - u_j^\T \Phi(z) u_j \right| &= \left| \frac{1}{z^3}u_l^\T \mathcal{R} u_l - \frac{1}{z^3} u_j^\T \mathcal{R} u_j -(u_l^\T \varepsilon(z) u_l-u_j^\T \varepsilon(z) u_j)\right|\\
    &\le \frac{1}{|z|^3} |u_l^\T \mathcal{R} u_l - u_j^\T \mathcal{R} u_j| + |u_l^\T \varepsilon(z) u_l-u_j^\T \varepsilon(z) u_j|\\
    &\le \frac{\osc(\mathcal{R})}{|z|^3} + \osc(\varepsilon(z)).
\end{align*}
For the last inequality above, since $u_l,u_j$ are unit vectors, their quadratic forms represent convex combinations of the diagonal entries. Therefore, the difference between them is strictly bounded by the maximum pairwise distance in the complex plane, which is the oscillation.

By Lemma \ref{lem:osc}, 
\[\left| u_l^\T \Phi(z) u_l - u_j^\T \Phi(z) u_j \right| \le \frac{\osc(\mathcal{R})}{|z|^3}\left(1 + \frac{15 R_{\max}}{|z|^2} \right) \le \frac{8}{3}\frac{\osc(\mathcal{R})}{|z|^3}. \]
Continuing from \eqref{eq:compare-alpha-lj}, together with $|\alpha_j(z)|\le \frac{3\rho}{\vartheta_j^2}$ from \eqref{eq:est-alpha} and the triangle inequality, we obtain
\begin{align*}
    |\alpha_l(z)| \ge |\alpha_l(z) - \alpha_j(z)| - |\alpha_j(z)| \ge \frac{\lambda_j-\lambda_l}{\lambda_j^2} - \frac{8}{3}\frac{\osc(\mathcal{R})}{|z|^3} - \frac{3\rho}{\vartheta_j^2}.
\end{align*}
Since $\lambda_l< 6\sqrt{R_{\max}}$ and $\lambda_j \ge 6\sqrt{R_{\max}} + 100\rho$, we have $\lambda_j-\lambda_l \ge 100\rho$. Using our assumptions, $|z|\ge \frac{99}{100}\vartheta_j$,  $0.87\lambda_j \le \vartheta_j \le 1.16 \lambda_j$ from \eqref{eq:vartheta-lower}, and $\rho> 13 \frac{\osc(\mathcal{R})}{\vartheta_j}$. Thus we  get 
\begin{align}\label{eq:alpha-l-2}
    |\alpha_l(z)| \ge (100\cdot 0.87^2-3) \frac{\rho}{\vartheta_j^2} -\frac{8}{3} \left(\frac{100}{99}\right)^3\frac{\osc(\mathcal{R})}{\vartheta_j^3} > 70 \frac{\rho}{\vartheta_j^2}.
\end{align}

\medskip
\noindent{\emph{Case 3:} $l\in [r] \setminus [r_+]$.} We have $\lambda_l <0$. Using the expansion \eqref{eq:Phi-expansion}, we get
\begin{align*}
    u_l^\T \Phi(\vartheta_j) u_l  &= \frac{1}{\vartheta_j} + \frac{\mathcal{V}_{ll}}{\vartheta_j^3} + \epsilon_l(\vartheta_j) \\
    &\ge \frac{1}{\vartheta_j} - \|\epsilon_l(\vartheta_j)\| \ge \frac{1}{\vartheta_j} -\frac{45}{8}\frac{R_{\max}^2}{\vartheta_j^5} \\
    &\ge \frac{1}{\vartheta_j} - \frac{5}{72} \frac{1}{\vartheta_j} = \frac{67}{72} \frac{1}{\vartheta_j}.
\end{align*}
Thus, $\alpha_l(\vartheta_j)= \lambda_l^{-1} - u_l^\T \Phi(\vartheta_j) u_l<0$ and 
\[|\alpha_l(\vartheta_j)| \ge \frac{67}{72} \frac{1}{\vartheta_j} \ge \frac{6700}{72}\frac{\rho}{\vartheta_j^2}. \]
Using similar estimates as \eqref{eq:04112026-1}, we have $|\alpha_j'(\zeta)| \le \frac{3}{\vartheta_j^2}$ for any $\zeta$ satisfying $|\zeta-\vartheta_j|\le \rho$. Consequently, 
\begin{align*}
    |\alpha_l(z) -\alpha_l(\vartheta_j)| \le \rho \cdot \max_{|\zeta-\vartheta_j|\le \rho}|\alpha_j'(\zeta)| \le 3\frac{\rho}{\vartheta_j^2} 
\end{align*}
and 
\begin{align}\label{eq:alpha-l-3}
    |\alpha_l(z)| \ge|\alpha_l(\vartheta_j)|-|\alpha_l(z) -\alpha_l(\vartheta_j)| \ge \left( \frac{6700}{72} -3\right)\frac{\rho}{\vartheta_j^2}> 30\frac{\rho}{\vartheta_j^2}.
\end{align}

Finally, combining \eqref{eq:alpha-l-1}, \eqref{eq:alpha-l-2}, and \eqref{eq:alpha-l-3}, we obtain
\begin{align}\label{eq:0415-previous-alpha}
    \min_{l\in [r]\setminus\{j\}} |\alpha_l(z)| \ge \frac{10}{3}\frac{\rho}{\vartheta_j^2}.
\end{align}
Continuing from \eqref{eq:smin-D1}, together with \eqref{eq:2nd-term-Di}, we arrive at 
\begin{align*}
   s_{\min}(D_j(z)) \ge \min_{l \in [r]\setminus \{j\}}|\alpha_l(z)| - \Big\|\off\big(U_{-j}^\T \Phi(z) U_{-j}\big)\Big\|  \ge \left(\frac{10}{3}- \frac{1.8}{13}\right)\frac{\rho}{\vartheta_j^2}> 3\frac{\rho}{\vartheta_j^2}. 
\end{align*}
This completes the proof.

\subsection{Proof of Theorem \ref{thm:outlier-cluster}}
Theorem \ref{thm:outlier-cluster} will be proved using similar techniques and estimates as Theorem \ref{thm:outlier-ev}. We briefly sketch it here. Recall that $\rho\equiv \rho_j = 10 M+ 15 \frac{\osc(\mathcal R)}{\lambda_j}$. 

To establish the lower bound for the top-$j$ outlier cluster, we define the block contour enclosing the deterministic roots $\vartheta_1, \dots, \vartheta_j$. Following our previous notation \eqref{def:D-pm}, let
\begin{align*}
D_j^- := \bigcup_{l=1}^j \mathtt{D}(\vartheta_l,\rho_j),
\end{align*}
and let $\mathcal{C}_j^- := \partial D_j^-$ be its positively oriented boundary. By our gap assumption $\min_{l>k} (\underline{\vartheta}_k - \vartheta_l) \ge 50\rho_k$, every supercritical root excluded from $D_k^-$ lies at a distance $49\rho$ from $\mathcal{C}_k^-$, ensuring the contour is strictly isolated. 

Recall from \eqref{def:A0} that 
\[
A_0(z):=\diag(\alpha_1(z),\dots,\alpha_r(z)) \quad\text{with} \quad \alpha_l(z):=\lambda_l^{-1} - u_l^\T \Phi(z) u_l.
\]
The main goal is to establish an analogous result to Lemma \ref{lem:compareAK}:
\begin{lemma}\label{lem:compareAK-cluster} Under the assumption of Theorem \ref{thm:outlier-cluster}, we have
    \[
\sup_{z\in \mathcal C_j^-}\|A_0(z)^{-1}(K(z)-A_0(z))\|<1,
\qquad
\sup_{z\in \mathcal C_j^-}\|K(z)^{-1}B(z)\|<1.
\]
\end{lemma}
Because the boundary $\mathcal{C}_j^-$ maintains a distance of at least $\rho_j$ from any enclosed root $\vartheta_l$, the uniform estimates established in the \textbf{Claim} of Section \ref{sec:compareAK} hold identically for all $z \in \mathcal{C}_j^-$. Consequently, Lemma \ref{lem:compareAK-cluster} holds.

Hence, by Lemma \ref{lem:compareAK-cluster} and Lemma \ref{lem:rouche}, $\det (A_0(z))$, $g(z)=\det (K(z))$, and $f(z)=\det(K(z)-B(z))$ have the same number of zeros inside the region $D_j^-$.

Since $\det(A_0(z)) = \prod_{l=1}^r \alpha_l(z)$, its zeros in the supercritical regime occur at $\vartheta_l$'s. By construction, $D_j^-$ contains exactly $j$ such roots: $\vartheta_1, \dots, \vartheta_j$. Therefore, $f(z)$ has exactly $j$ zeros inside $D_j^-$. As noted previously, the zeros of $f(z)$ outside $[-\|E\|,\|E\|]$ correspond exactly to the eigenvalues of the perturbed matrix $\widetilde{A}$. Since the eigenvalues are strictly ordered (i.e., $\widetilde{\lambda}_1 \ge \dots \ge \widetilde{\lambda}_r$), these $j$ enclosed roots inside $D_j^-$ must be exactly the top $j$ perturbed eigenvalues $\widetilde{\lambda}_1, \dots, \widetilde{\lambda}_j$. Because all $j$ of these eigenvalues lie strictly inside $D_j^-$, their locations are bounded from below by the leftmost boundary of the region. 

Recalling that $\underline{\vartheta}_j := \min_{1 \le l \le j} \vartheta_l$, the real part of any $z \in D_k^-$ satisfies $\Re(z) \ge \underline{\vartheta}_k - \rho_k$. Thus, we conclude that for all $s \in [k]$,
\[
\widetilde{\lambda}_s \ge \underline{\vartheta}_k - \rho_k.
\]
The proof is complete.


\appendix
\section{Proofs of Lemmas \ref{lem:bdphi}, Lemma \ref{lem:osc}, Proposition \ref{prop:vartheta}, and Lemma \ref{lem:bdE}}
\subsection{Proof of Lemma \ref{lem:bdphi}}\label{sec:proof-lem:bdphi}
Fix $z$ with $|z|^2\ge 6 R_{\max}$. 

\medskip
\noindent{\emph{\underline{Proof of Part (i)}:}} Define a mapping $F:\C^n \to \C^n$ componentwise as 
\[\left( F(y) \right)_i= \frac{1}{z-\sum_j \sigma_{ij}^2 y_j}=\frac{1}{z-(\Sigma y)_i}.\] 
Consider the metric $\|y\|_\infty=\max_i |y_i|$ on $\C^n$. Note that $\phi$ is a fixed point of $F$, i.e., $F(\phi)=\phi$. Define a closed set
\[\mathcal B:=\left\{y\in \C^n : \|y\|_\infty \le \frac{3}{2|z|} \right\}.\]
We first show $F(\mathcal B)\subseteq \mathcal B$. For $y\in \mathcal B$ and  every $i\in [n]$,
\[|(\Sigma y)_i|=\big|\sum_j \sigma_{ij}^2 y_j \big| \le \|y\|_\infty R_{\max} \le \frac{2}{|z|} R_{\max}  \le \frac{|z|}{3}\]
by our supposition $|z|^2\ge 6 R_{\max}$. Hence,
\begin{align}\label{eq:lembdphi-1}
   \big|z- (\Sigma y)_i \big|   \ge |z|-\big|(\Sigma y)_i \big| \ge \frac{2}{3}|z|
\end{align}
and 
\[ \big|\left( F(y) \right)_i \big| = \frac{1}{|z-(\Sigma y)_i|} \le \frac{3}{2|z|}.\]
This shows that $F(y)\in \mathcal B.$

Next, we show that $F:\mathcal B\to \mathcal B$ is a contraction. Take $y,y'\in \mathcal B$. For each $i\in [n]$, using \eqref{eq:lembdphi-1}, we have 
\begin{align*}
   \big|\left( F(y) \right)_i-\left( F(y') \right)_i  \big| &=\left|\frac{\sum_j \sigma_{ij}^2 (y_j-y_j')}{(z- (\Sigma y)_i)(z- (\Sigma y')_i)} \right|\\
   &\le \frac{R_{\max} \|y-y'\|_{\infty}}{(2|z|/3)^2} = \frac{9}{4}\frac{R_{\max}}{|z|^2}\|y-y'\|_{\infty} \le \frac{3}{8}\|y-y'\|_{\infty}.
\end{align*}
Since $(\mathcal B, \|\cdot\|_\infty)$ is a complete metric space, by Banach's fixed point theorem, $F$ has a unique fixed point in $y^*\in \mathcal B$. Since the solution $\phi$ to the QVE \eqref{def:QVE} is unique and $\phi(z)\to 0$ as $|z|\to\infty$. Therefore, $\phi=y^* \in \mathcal B$. In particular, 
\begin{align}\label{eq:phibdcrude}
\max_i |\phi_i(z)| \le \frac{3}{2|z|}.
\end{align}
We have 
\begin{align}\label{eq:sumphi}
\Big|\big(\Sigma\phi(z)\big)_i\Big|=\Big|\sum_j \sigma_{ij}^2 \phi_j(z) \Big| \le \frac{3}{2|z|} R_{\max}.
\end{align}
Hence, $\Big|\sum_j \sigma_{ij}^2 \phi_j(z) \Big| \le\frac{1}{4}|z| $ and 
 \[\frac{3}{4}|z| \le \big|z-\sum_j \sigma_{ij}^2 \phi_j(z)\big| \le \frac{5}{4}|z|.\]
To show the bounds on $\phi'(z)$, differentiating the equation
$\phi_i(z)^{-1}=z-\sum_{j=1}^n \sigma_{ij}^2\phi_j(z)$
gives
\[-\phi_i(z)^{-2}\phi_i'(z)=1-\sum_{j=1}^n \sigma_{ij}^2\phi_j'(z),\]
hence
\[\phi_i'(z)=-\phi_i(z)^2\Bigl(1-\sum_{j=1}^n \sigma_{ij}^2\phi_j'(z)\Bigr)\]
and
\[|\phi_i'(z)|\le|\phi_i(z)|^2\Bigl(1+\sum_{j=1}^n \sigma_{ij}^2|\phi_j'(z)|\Bigr).\]
Let \(M_1:=\max_i |\phi_i'(z)|\). By \eqref{eq:phibdcrude}, 
\[M_1 \le \frac{9}{4|z|^2}(1+ R_{\max}M_1).\]
Since $R_{\max}/|z|^2\le 1/6$, we solve $M_1 \le \frac{9}{4|z|^2} + \frac{3}{8}M_1$ to obtain 
\[ M_1=\max_i |\phi_i'(z)| \le \frac{18}{5|z|^2}. \]
To bound $\phi''(z)$, similarly, differentiating the equation one more time yields 
\[ \phi_i''(z)=-2\phi_i(z)\phi_i'(z)\Bigl(1-\sum_{j=1}^n \sigma_{ij}^2\phi_j'(z)\Bigr)
+\phi_i(z)^2 \sum_{j=1}^n \sigma_{ij}^2\phi_j''(z).\]
Hence,
\[|\phi_i''(z)|\le
2|\phi_i(z)|\,|\phi_i'(z)|
\Bigl(1+\sum_{j=1}^n \sigma_{ij}^2|\phi_j'(z)|\Bigr)+
|\phi_i(z)|^2 \sum_{j=1}^n \sigma_{ij}^2|\phi_j''(z)|.\]
Let $M_2:=\max_i |\phi_i''(z)|$. Using the previous bounds on $M_1$ and \eqref{eq:phibdcrude}, we have \[M_2\le
2\cdot \frac{3}{2|z|}\cdot \frac{18}{5|z|^2}\cdot \Bigl(1+\frac35\Bigr)+
\frac{9R_{\max}}{4|z|^2}M_2.\]
Thus,
\[ M_2=\max_i |\phi_i''(z)|\le \frac{3456}{125|z|^3} \le \frac{28}{|z|^3}.\]

\medskip
\noindent{\emph{\underline{Proof of Part (ii)}:}} In the following, we prove \eqref{eq:phi-expand}. From $\phi_i(z) = \frac{1}{z-(\Sigma \phi)_i}$ and the exact identity $\frac{1}{z-x} = \frac{1}{z} + \frac{x}{z(z-x)}$, treating $x=(\Sigma \phi)_i$, we get
\begin{align}\label{eq:phi-1stexp}
    \phi_i(z) = \frac{1}{z}+ \frac{(\Sigma \phi)_i}{z (z-(\Sigma \phi)_i)} = \frac{1}{z}+  \frac{(\Sigma \phi)_i}{z}\phi_i(z):= \frac{1}{z} + \delta_i(z),
\end{align}
where 
\[|\delta_i(z)| \le \frac{|(\Sigma \phi)_i|}{|z|}|\phi_i(z)| \le \frac{9}{4} \frac{R_{\max}}{|z|^3}\]
by \eqref{eq:phibdcrude} and \eqref{eq:sumphi}. Next, from the algebraic identity $\frac{1}{z-x} = \frac{1}{z} + \frac{x}{z^2} + \frac{x^2}{z^2(z-x)}$, we also have 
\begin{align}\label{eq:phi-2ndexp}
    \phi_i(z) = \frac{1}{z}+ \frac{(\Sigma \phi)_i}{z^2} + \frac{(\Sigma \phi)_i^2}{z^2}\phi_i(z) :=\frac{1}{z}+ \frac{(\Sigma \phi)_i}{z^2} + \tilde{\delta}_i(z),
\end{align}
where, by \eqref{eq:phibdcrude} and \eqref{eq:sumphi}, we have the estimate 
\[ |\tilde{\delta}_i(z)| \le \frac{|(\Sigma \phi)_i|^2}{|z|^2}\big|\phi_i(z)\big| \le \frac{27}{8} \frac{R_{\max}^2}{|z|^5}.\]
Now we plug \eqref{eq:phi-1stexp} into the term $\frac{(\Sigma \phi)_i}{z^2}$ in \eqref{eq:phi-2ndexp}: 
\[\frac{(\Sigma \phi)_i}{z^2} = \frac{1}{z^2}\sum_j \sigma_{ij}^2\left(\frac{1}{z} + \delta_j(z) \right) = \frac{R_i}{z^3} + \frac{\sum_j \sigma_{ij}^2 \delta_j(z)}{z^2}.\]
Observe that \[\left|\frac{\sum_j \sigma_{ij}^2 \delta_j(z)}{z^2} \right| \le \frac{R_{\max} }{|z|^2}\max_j |\delta_j(z)| \le \frac{9}{4} \frac{R_{\max}^2}{|z|^5}.\]
Combining the above estimates, we arrive at the final expansion:
\[\phi_i(z) = \frac{1}{z} + \frac{R_i}{z^3} + \left(\tilde{\delta}_i(z) + \frac{\sum_j \sigma_{ij}^2 \delta_j(z)}{z^2}\right):= \frac{1}{z} + \frac{R_i}{z^3} + \varepsilon_i(z)\]
with the error term $\varepsilon_i(z)$ satisfying 
\[|\varepsilon_i(z)| \le |\tilde{\delta}_i(z)| +\left|\frac{\sum_j \sigma_{ij}^2 \delta_j(z)}{z^2} \right| \le \frac{27}{8} \frac{R_{\max}^2}{|z|^5} + \frac{9}{4} \frac{R_{\max}^2}{|z|^5}= \frac{45}{8}\frac{R_{\max}^2}{|z|^5}. \]
This proves \eqref{eq:phi-expand}. 

To prove \eqref{eq:Phi-prime-expansion}, we start with 
\[\phi_i(z) = \left(\frac{1}{z} + \frac{R_i}{z^3} +\varepsilon_i(z)\right)' = -\frac{1}{z^2}-\frac{3R_i}{z^4}+\varepsilon_i'(z).\] Recall that 
\[\varepsilon_i(z)=\tilde\delta_i(z)+\frac{1}{z^2}\sum_j \sigma_{ij}^2\delta_j(z),\] where \[\delta_i(z)=\frac{(\Sigma\phi)_i}{z}\phi_i(z)\quad \text{ and }\quad \tilde\delta_i(z)=\frac{(\Sigma\phi)_i^2}{z^2}\phi_i(z).\] It remains to bound $|\varepsilon_i'(z)|$. Note that 
\begin{align*}
&\delta_i'(z)=-\frac{(\Sigma\phi)_i}{z^2}\phi_i(z)+\frac{(\Sigma\phi')_i}{z}\phi_i(z)
+\frac{(\Sigma\phi)_i}{z}\phi_i'(z),\\
&\tilde\delta_i'(z)=-2\frac{(\Sigma\phi)_i^2}{z^3}\phi_i(z)+
2\frac{(\Sigma\phi)_i(\Sigma\phi')_i}{z^2}\phi_i(z)
+\frac{(\Sigma\phi)_i^2}{z^2}\phi_i'(z).
\end{align*}
Plugging in the previous estimates 
\[
|\phi_i(z)|\le \frac{3}{2|z|},
\quad
|\phi_i'(z)|\le \frac{18}{5|z|^2},
\quad
|(\Sigma\phi)_i|\le \frac{3R_{\max}}{2|z|},
\quad
|(\Sigma\phi')_i|\le \frac{18R_{\max}}{5|z|^2}
\]
yields
\[ |\delta_i'(z)| \le 13.05\frac{R_{\max}}{|z|^4} \quad \text{and} \quad |\tilde\delta_i'(z)| \le 31.05 \frac{R_{\max}^2}{|z|^6}.\]
Together with $|\delta_i(z)|\le \frac{9}{4} \frac{R_{\max}}{|z|^3}$, we obtain that
\begin{align*}
    |\varepsilon_i'(z)|\le |\tilde\delta_i'(z)| + \left|-\frac{2}{z^3} \sum_{j}\sigma_{ij}^2 \delta_j(z) + \frac{1}{z^2}\sum_j \sigma_{ij}^2\delta_j'(z) \right| \le 48.6 \frac{R_{\max}^2}{|z|^6}.
\end{align*}
We have proved \eqref{eq:Phi-prime-expansion}.


\subsection{Proof of Lemma \ref{lem:osc}}\label{app:proof-lem-osc}

Fix $z\in \C$ with $|z|\ge\sqrt{6R_{\max}}>0$. Let $\sqrt{\cdot}$ denote the principal branch of the complex square root on $\C\setminus(-\infty,0]$. Recall that in \eqref{def:QVE}, $\phi_i(z)= \frac{1}{z-\sum_{j}\sigma_{ij}^2 \phi_j(z)}$. For comparison, for each $i\in[k]$, if $R_i>0$, define an auxiliary function
\begin{align}\label{def:mi}
    m_i(z)= \frac{z-\sqrt{z^2-4R_i}}{2R_i} = \frac{z}{2R_i} \left( 1- \sqrt{1-\frac{4R_i}{z^2}}\right) 
\end{align}
and if $R_i=0$, define $m_i(z) = \frac{1}{z}$. Then $m_i(z)$ satisfies the quadratic equation
\begin{align}\label{eq:mi-eq}
    m_i(z) = \frac{1}{z - R_i m_i(z)}.
\end{align}

Set \[d_i(z) = \phi_i(z)-m_i(z) \quad\text{and}\quad d(z)=(d_1(z),\dots,d_n(z)).\]  
Then
\begin{align*}
    \Phi(z)=\diag(m_1(z),\dots,m_n(z)) + d(z).
\end{align*}
We can rewrite 
\begin{align*}
    \varepsilon_i(z)= \phi_i(z) - \frac{1}{z} - \frac{R_i}{z^3} = m_i(z) - \frac{1}{z} - \frac{R_i}{z^3} + d_i(z),
\end{align*}
which yields
\begin{align}\label{eq:osc-split}
    \osc(\varepsilon_i(z)) &\le \osc\left(m_i(z) - \frac{1}{z} - \frac{R_i}{z^3} \right) + 2\|d(z)\|_{\infty}.
\end{align}
We now bound each term on the right-hand side separately. 

We start with the first term on the right-hand side of \eqref{eq:osc-split}. Since $|\frac{4R_i}{z^2}|\le \frac{2}{3}<1$,  the binomial series for $\sqrt{1-w}$ is absolutely convergent at $w=\frac{4R_i}{z^2}$. Hence, 
\[ \sqrt{1-w} = 1-2 \sum_{k=0}^{\infty} \frac{1}{k+1}\binom{2k}{k} \left(\frac{w}{4}\right)^{k+1} = 1 - 2 \sum_{k=0}^{\infty} C_k \left(\frac{w}{4}\right)^{k+1},\]
where $C_k = \frac{1}{k+1}\binom{2k}{k}$ denotes the Catalan number. Applying this expansion \[m_i(z)= \frac{z}{2R_i} \left( 1- \sqrt{1-\frac{4R_i}{z^2}}\right), \]
we get
\begin{align}\label{eq:expansion-mi}
    m_i(z) &= \frac{z}{2R_i}\left[1- \left( 1- 2\sum_{k=0}^{\infty} C_k \left( \frac{4R_i/z^2}{4}\right)^{k+1} \right) \right]= \frac{z}{R_i}\left(\sum_{k=0}^{\infty}  C_k \frac{R_i^{k+1}}{z^{2k+2}}\right)\nonumber\\
    &= \frac{1}{z} + \frac{R_i}{z^3} + \sum_{k=2}^{\infty} C_k \frac{R_i^{k}}{z^{2k+1}}:= \frac{1}{z} + \frac{R_i}{z^3} + v_i(z)
\end{align} 
and thus
\begin{align}\label{eq:decompose-Phi}
    \Phi(z)= \frac{1}{z}I + \frac{1}{z^3}\mathcal{R} + \diag(v_i(z))_{1\le i \le n} + d(z). 
\end{align}
It follows that 
\begin{align*}
    \osc\left(m_i(z) - \frac{1}{z} - \frac{R_i}{z^3} \right) &=\max_{i,j}|v_i(z) - v_j(z)|\\
    & = \max_{i,j} \left| \sum_{k=2}^{\infty} C_k \frac{R_i^{k}-R_j^{k}}{z^{2k+1}}\right| \\
    &\le \left(\sum_{k=2}^{\infty} C_k\frac{k R_{\max}^{k-1}}{|z|^{2k+1}}\right)\osc(\mathcal R),
\end{align*}
where we used \[|R_i^{k}-R_j^{k}|=|(R_i-R_j)(R_i^{k-1}+R_i^{k-2}R_j+\dots+R_j^{k-1})| \le \osc(\mathcal{R}) \cdot k R_{\max}^{k-1}.\]
Since $R_{\max}/|z|^2 \le 1/6$ by our assumption, we continue to have 
\begin{align}\label{eq:04102026}
     \osc\left(m_i(z) - \frac{1}{z} - \frac{R_i}{z^3} \right) 
    &\le \frac{R_{\max}}{|z|^5}\left(\sum_{k=2}^{\infty} kC_k\frac{ R_{\max}^{k-2}}{|z|^{2(k-2)}}\right)\osc(\mathcal R)\nonumber\\
    &\le \frac{R_{\max}}{|z|^5}\osc(\mathcal R) \left(\sum_{k=2}^{\infty} \frac{kC_k}{6^{k-2}} \right)\nonumber\\
    &\le \frac{11R_{\max}}{|z|^5}\osc(\mathcal R) .
\end{align}
In the last step, we used $\sum_{k=2}^{\infty} \frac{kC_k}{6^{k-2}}\le 11$. This can be derived by taking the generating function of the Catalan numbers 
\begin{align}\label{def:generating}
c(x) = \sum_{k=0}^\infty C_k x^k = \frac{1 - \sqrt{1-4x}}{2x}.
\end{align}
The sum is $S = \sum_{k=2}^\infty k C_k x^{k-2}$ evaluated at $x = 1/6$. Differentiating the generating function gives $c'(x) = \sum_{k=1}^\infty k C_k x^{k-1} = 1 + xS$. Therefore, the sum is exactly given by $S = \frac{c'(x) - 1}{x}$. Taking the derivative yields $c'(x) = \frac{1 - 2x - \sqrt{1-4x}}{2x^2\sqrt{1-4x}}$. Evaluating this at $x = 1/6$ gives $c'(1/6) = 12\sqrt{3} - 18$. Substituting this back into the sum formula yields the exact value:$$\sum_{k=2}^\infty  \frac{kC_k}{6^{k-2}} = 6(12\sqrt{3} - 19) \le 11.$$
Furthermore, we obtain the following bound for later estimates:
\begin{align}\label{eq:mi-diff-bd}
    |m_i(z)-m_j(z)| &\le  \frac{|R_i-R_j|}{|z|^3} + \osc\left(m_i(z) - \frac{1}{z} - \frac{R_i}{z^3} \right) \nonumber\\
    &\le \frac{\osc{(\mathcal R)}}{|z|^3} + \frac{11R_{\max}}{|z|^5}\osc(\mathcal R)\nonumber\\
    &\le \left( 1 + \frac{1}{6} \right)\frac{\osc{(\mathcal R)}}{|z|^3}= \frac{7}{6}\frac{\osc{(\mathcal R)}}{|z|^3}.
\end{align}
In the last inequality, we applied $R_{\max}/|z|^2 \le 1/6.$

Next, we turn to the second term on the right-hand side of \eqref{eq:osc-split}. For each $i$, we have 
\begin{align}\label{eq:04092026}
    |d_i(z)| = |\phi_i(z)-m_i(z)|  &= \left| \frac{1}{z-\sum_j \sigma_{ij}^2 \phi_j(z)} - \frac{1}{z-R_i m_i(z)}\right|\nonumber\\
    &= \frac{|R_i m_i(z) - \sum_j \sigma_{ij}^2 \phi_j(z)|}{|z-\sum_j \sigma_{ij}^2 \phi_j(z)| \cdot |z-R_i m_i(z)|}\nonumber\\
    &\le \frac{2}{|z|^2} \Big|R_i m_i(z) - \sum_j \sigma_{ij}^2 \phi_j(z)\Big|.
\end{align}
For the last step above, observe that $|z-\sum_j \sigma_{ij}^2 \phi_j(z)|\ge \frac{3}{4}|z|$ by Lemma \ref{lem:bdphi} (i). In addition, we claim that \[|z-R_i m_i(z)|\ge \frac{2|z|}{3}.\] To see this, from \eqref{eq:mi-eq}, we have  $|z-R_i m_i(z)|=\frac{1}{|m_i(z)|}.$
Using the series expansion of $m_i(z)$ in \eqref{eq:expansion-mi} and $R_{\max}/|z|^2 \le 1/6$, we get
\begin{align*}
    |m_i(z)|\le \frac{1}{|z|} \sum_{k=0}^{\infty} C_k\Bigr( \frac{R_i}{|z|^2}\Bigl)^k \le \frac{1}{|z|}\sum_{k=0}^{\infty}\frac{C_k}{6^k}=\frac{c(1/6)}{|z|}<\frac{1.5}{|z|},
\end{align*}
where $c(x)$ is defined in \eqref{def:generating}. The claim is proved.

It remains to estimate $|R_i m_i(z) - \sum_j \sigma_{ij}^2 \phi_j(z)|$. We write 
\begin{align*}
    \Big|R_i m_i(z) - \sum_j \sigma_{ij}^2 \phi_j(z)\Big| &= |\sum_j \sigma_{ij}^2 (m_i(z)-m_j(z)) + \sum_j \sigma_{ij}^2 (m_j(z)-\phi_j(z))|\\
    &\le \sum_j \sigma_{ij}^2 |m_i(z)-m_j(z)| + R_{\max} \|d(z)\|_{\infty}.
\end{align*}
Plugging in \eqref{eq:mi-diff-bd}, we further have
\begin{align*}
    |R_i m_i(z) - \sum_j \sigma_{ij}^2 \phi_j(z)|  \le \frac{7}{6}\frac{R_{\max}}{|z|^3}\osc{(\mathcal R)} + R_{\max} \|d(z)\|_{\infty}
\end{align*}
and thus from \eqref{eq:04092026},
\begin{align}
    \|d(z)\|_{\infty} = \max_i |d_i(z)| &\le \frac{2}{|z|^2}\left( \frac{7}{6}\frac{R_{\max}}{|z|^3}\osc{(\mathcal R)} + R_{\max} \|d(z)\|_{\infty}\right)\\
    &\le \frac{7}{3}\frac{R_{\max}}{|z|^5}\osc{(\mathcal R)} +\frac{1}{3} \|d(z)\|_{\infty}.
\end{align}
Rearranging terms, we arrive at 
\begin{align}\label{eq:infbd-d}
    \|d(z)\|_{\infty}\le \frac{7 R_{\max}}{2|z|^5}\osc{(\mathcal R)}.
\end{align}

Finally, we conclude that 
\begin{align*}
     \osc(\varepsilon_i(z)) \le \osc\left(m_i(z) - \frac{1}{z} - \frac{R_i}{z^3} \right) + 2\|d(z)\|_{\infty} \le 14.5 \frac{R_{\max}}{z^5}\osc{(\mathcal R)}.
\end{align*}

\subsection{Proof of Proposition \ref{prop:vartheta}}\label{app:proof-prop-vartheta}
By part (ii) of Lemma~\ref{lem:bdphi},
\begin{align}\label{eq:04082026}
    u_j^\T \Phi(z)u_j=\frac{1}{z}+\frac{\mathcal V_{jj}}{z^3}+\epsilon_j(z),
\qquad
|\epsilon_j(z)|=|u_j^\T \varepsilon(z) u_j|\le \frac{45}{8}\frac{R_{\max}^2}{z^5},
\end{align}
for all real $z\ge 3\sqrt{R_{\max}}$. For brevity, denote $f(z)=u_j^\T \Phi(z) u_j$ and thus $\widehat\alpha_j(z)=1-\lambda_j f(z)$.

We first show that $\widehat\alpha_j$ has a unique zero. Note that $\mathcal V_{jj}\ge0$. At $z_0:=3\sqrt{R_{\max}}$,
\[f(z_0)=u_j^\T \Phi(z_0)u_j
\ge \frac{1}{z_0}-\frac{45}{8}\frac{R_{\max}^2}{z_0^5}
= \left(\frac13-\frac{45}{8\cdot 3^5}\right)\frac1{\sqrt{R_{\max}}}
= \frac{67}{216}\frac1{\sqrt{R_{\max}}},\]
hence
\[\widehat\alpha_j(z_0)
\le 1-\frac{67}{216}\frac{\lambda_j}{\sqrt{R_{\max}}}< 0\]
since $\lambda_j\ge 6\sqrt{R_{\max}}$. On the other hand, $\widehat\alpha_j(z)\to 1$ as $z\to\infty$. Next, by part (ii) of Lemma~\ref{lem:bdphi},
\begin{align}\label{eq:diff-f}
    f'(z)=u_j^\T \Phi'(z)u_j=
-\frac1{z^2}-\frac{3\mathcal V_{jj}}{z^4}+\epsilon_j'(z),
\qquad
|\epsilon_j'(z)|=|u_j^\T \varepsilon'(z) u_j|\le \frac{243}{5}\frac{R_{\max}^2}{z^6}.
\end{align}
Therefore, on the interval $[3\sqrt{R_{\max}},\infty)$,
\[\widehat\alpha_j'(z)=-f'(z) \ge \lambda_j\left(\frac1{z^2}-\frac{243}{5}\frac{R_{\max}^2}{z^6}\right)>0\]
since \[\frac{243}{5}\frac{R_{\max}^2}{z^6} \le \frac{243}{5}\frac{1}{81z^2}= \frac{243}{405}\frac{1}{z^2} < \frac{1}{z^2}.\] Therefore, $\widehat\alpha_j$ is strictly increasing, and has a unique and simple real zero $\vartheta_j \in [3\sqrt{R_{\max}},\infty)$.

Observe that at $z=\lambda_j/2 \ge 3\sqrt{R_{\max}}$ and $z=2\lambda_j$, we have
\begin{align*}
    \widehat\alpha_j(\lambda_j/2) = 1- \lambda_ju_j^\T \Phi(\lambda_j/2)u_j &\le 1- \lambda_j \left( \frac{2}{\lambda_j}- \frac{45}{8}\frac{2^5 R_{\max}^2}{\lambda_j^5} \right)\\
    &\le -1 +\frac{5}{36}<0
\end{align*}
and
\begin{align*}
    \widehat\alpha_j(2\lambda_j)=1- \lambda_ju_j^\T \Phi(2\lambda_j)u_j &\ge 1- \lambda_j \left( \frac{1}{2\lambda_j}+\frac{45}{8}\frac{R_{\max}^2}{2^5 \lambda_j^5} \right)\\
    &\ge \frac{1}{2} - \frac{45}{9216}>0.
\end{align*}
Since $\widehat\alpha_j$ is strictly increasing, we have 
\[ \frac{\lambda_j}{2}\le \vartheta_j \le 2\lambda_j.\]

Now we locate $\vartheta_j$. Note that $f(\vartheta_j)= \lambda_j^{-1}$. We introduce a scalar comparison function 
\begin{align}\label{eq:mv-approx}
    m_{\mathcal{V},j}(z)= \frac{z-\sqrt{z^2-4\mathcal{V}_{jj}}}{2\mathcal{V}_{jj}}
\end{align}
which satisfies the quadratic equation 
\begin{align}\label{eq:mv-approx-quad}
    m(z)= \frac{1}{z-\mathcal{V}_{jj} m(z)}.
\end{align}
This function approximates $f(z)=u_j^\T \Phi(z) u_j$ by noting $\mathcal{V}_{jj}= \sum_i u_j(i)^2 R_i$ and the definition of $\phi_i(z)$ in \eqref{def:QVE}. 

Denote \[x_j:= \lambda_j + \frac{\mathcal V_{jj}}{\lambda_j}.\] Then $x_j$ is a zero of $1-\lambda_j m_{\mathcal{V},j}(z)=0$. Indeed, plugging $m_{\mathcal{V},j}(z) = \frac{1}{\lambda_j}$ into \eqref{eq:mv-approx-quad}, we get 
$\lambda_j = z- \mathcal{V}_{jj} m_{\mathcal{V},j}(z) = z- \frac{\mathcal{V}_{jj}}{\lambda_j}$ and thus $z= \lambda_j + \frac{\mathcal{V}_{jj}}{\lambda_j}.$ 

Next, we expand $m_{\mathcal{V},j}(z)$ as a Catalan series the same as \eqref{eq:expansion-mi} at $z=x_j$:
\begin{align}\label{eq:expand-mv}
    \frac{1}{\lambda_j} = m_{\mathcal{V},j}(x_j) = \frac{1}{x_j} + \frac{\mathcal{V}_{jj}}{x_j^3} + \sum_{k=2}^{\infty} C_k \frac{\mathcal{V}_{jj}^k}{x_j^{2k+1}}.
\end{align}

To estimate $|\vartheta_j - x_j|$, we first show that $f(x_j) = u_j^\T \Phi(x_j) u_j \approx \frac{1}{\lambda_j}=f(\vartheta_j)$. Then we apply the mean value theorem to bound $|\vartheta_j - x_j|$.

In the proof of Lemma \ref{lem:osc}, we have the expansion of $\Phi(z)$ given in \eqref{eq:decompose-Phi}. Hence,
\begin{align*}
    f(x_j) = u_j^\T \Phi(x_j) u_j = \frac{1}{x_j} + \frac{\mathcal{V}_{jj}}{x_j^3} + \sum_{k=2}^{\infty} C_k \frac{\sum_i u_j(i)^2 R_i^k}{x_j^{2k+1}} + u_j^\T d(x_j) u_j.
\end{align*}
Subtracting with \eqref{eq:expand-mv}, we obtain
\begin{align*}
    |f(x_j)-f(\vartheta_j)|=|f(x_j)-\frac{1}{\lambda_j}|\le \sum_{k=2}^{\infty} \frac{C_k}{x_j^{2k+1}} \left|\mathcal{V}_{jj}^k - \sum_i u_j(i)^2 R_i^k \right| + \|d(x_j)\|_{\infty}.
\end{align*}
Note that 
\begin{align*}
\left|\mathcal{V}_{jj}^k - \sum_i u_j(i)^2 R_i^k \right| \le \sum_{i} u_j(i)^2\left|\mathcal{V}_{jj}^k - R_i^k \right|&\le k R_{\max}^k \sum_{i} u_j(i)^2 \left|\mathcal{V}_{jj} - R_i \right|\\
&\le k R_{\max}^k \osc(\mathcal R).
\end{align*}
Following the same step as \eqref{eq:04102026} and combining the bound on $\|d\|_{\infty}$ in \eqref{eq:infbd-d}, we arrive at 
\begin{align*}
    |f(x_j)-f(\vartheta_j)| \le 11\frac{R_{\max}}{x_j^5} \osc(\mathcal R) + \frac{7}{2}\frac{R_{\max}}{x_j^5} \osc(\mathcal R)\le14.5\frac{R_{\max}}{\lambda_j^5} \osc(\mathcal R).
\end{align*}
Finally, by mean value theorem, $|f(x_j)-f(\vartheta_j)| = |x_j-\vartheta_j|\cdot |f'(w)|$ for some $w$ between $x_j$ and $\vartheta_j$. Since $\lambda_j/2 \le \vartheta_j \le 2\lambda_j$ and $x_j \le 2\lambda_j$,  we have $w\le 2\lambda_j$. By \eqref{eq:diff-f}, we get
\[f'(w) \ge \frac{1}{w^2} - \frac{243}{5} \frac{R_{\max}^2}{w^6} \ge \frac{1}{w^2} \left(1- \frac{243}{5\cdot81}\right) = \frac{2}{5w^2}\ge \frac{1}{10 \lambda_j^2}.\]
Rearranging yields that 
\[ \left|\lambda_j + \frac{\mathcal V_{jj}}{\lambda_j} - \vartheta_j \right|= |x_j-\vartheta_j| \le 145\frac{R_{\max}}{\lambda_j^3} \osc(\mathcal R).\]
The proof is now complete. 

\subsection{Proof of Lemma \ref{lem:bdE}}\label{app:lem-bdE}
We apply Corollary 3.12 from \cite{BvH16} with a standard truncation technique.  Denote $$L_n:=K\sigma_{\max}((D+4)\log n)^{1/2}.$$ By a union bound, we have 
\begin{align*}
   \Prob( \max_{i,j} |E_{ij}|>L_n) \le n^2\cdot 2\exp(-(L_n/K)^2) = 2n^{-(D+2)}.
\end{align*}
We split $E$ into big and small parts according to $L_n$. Namely, we write $E=Z+M + B$, where we define these symmetric matrices via
\[Z_{ij}= E_{ij}\mathbf{1}_{\{|E_{ij}|\le L_n \}} - \E(E_{ij}\mathbf{1}_{\{|E_{ij}|\le L_n \}}),\,  M_{ij} = \E(E_{ij}\mathbf{1}_{\{|E_{ij}|\le L_n \}}),\, B_{ij}=E_{ij}\mathbf{1}_{\{|E_{ij}|> L_n \}}.\]
Note that if $\max_{i,j} |E_{ij}|\le L_n$, then $B=0$. We work on the event $\mathcal{E}=\{ \max_{i,j} |E_{ij}|\le L_n\}$ below and previous calculations suggest that $\Prob(\mathcal{E})\ge 1-2n^{-(D+2)}.$

Next, we bound $\|M\|$. Note that since $E_{ij}$ has mean 0, $M_{ij} = -\E(E_{ij}\mathbf{1}_{\{|E_{ij}|> L_n \}})$. Therefore,
\[ |M_{ij}| \le \E(|E_{ij}|\mathbf{1}_{\{|E_{ij}|> L_n \}}) = \int_{L_n}^{\infty} \Prob(|E_{ij}|\ge t)\,dt\le 2 L_n  n^{-(D+4)}\] since $\Prob(|E_{ij}|> L_n) \le 2n^{-(D+4)}.$ Then we have
\[\|M\| \le \max_i \sum_j |M_{ij}|\le 2 L_n n^{-(D+3)}. \]

Finally, for the bound on $\|Z\|$, we apply Corollary 3.12 from \cite{BvH16}. Note that $Z_{ij}$'s have mean 0 and $|Z_{ij}|\le 2L_n$. Also, $\E(Z_{ij}^2) \le \E(E_{ij}^2)=\sigma_{ij}^2$. Applying of \cite[Corollary 3.12 ]{BvH16}, together with a standard symmetrization (as in the proof of \cite[Corollary 3.3 ]{BvH16}), yields that for $t>0,$
\[ \Prob\left(\|Z\|\ge 2\sqrt{2}(1+\epsilon)\tilde{\sigma} +t \right)\le n \exp(-t^2/(c_{\epsilon}L_n)^2).\]
We choose $\epsilon$ such that $2\sqrt{2}(1+\epsilon)=2.9$ and denote $c_{\epsilon}=c$. The parameter $\tilde{\sigma}  =\max_i \sqrt{\sum_j \E(Z_{ij}^2)} =\sqrt{\max_i  \sigma_{ij}^2} = \sqrt{R_{\max}}$ by our supposition. Choosing $t=cL_n\sqrt{(D+3)\log n}$, we obtain that 
\[\Prob\left( \|Z\| \ge 2.9\sqrt{R_{\max}} + cL_n \sqrt{(D+3)\log n} \right) \le n^{-(D+2)}.\]

Combining the above estimates, we have that with probability at least $1-3n^{-(D+2)}$,
\begin{align*}
    \|E\| \le \|Z\|+\|M\|+\|B\| &\le 2L_n n^{-(D+3)} + 1.5\sqrt{R_{\max}} + cL_n \sqrt{(D+3)\log n}\\
    &\le 2.9\sqrt{R_{\max}} + c K \sigma_{\max}(D+4)\log n.
\end{align*}

\section{Proofs of Proposition \ref{prop:spacepert} and Proposition \ref{prop:detbd-vector}}\label{app:perturb-prepare}
\subsection{Proof of Proposition \ref{prop:spacepert}} Let $Q=I-UU^\T$ be the orthogonal projection onto the null space of $A$. Set $P_J = U_J U_J^\T$ where $J=[r]\setminus[k]$ and $U_J = (u_{k+1},\dots,u_r)$. Then \[P_{U_k^{\perp}} = P_J + Q.\] 

We first prove \eqref{eq:sin-det} in part (i). Starting from $\|\sin\angle(U_k, \widetilde{U}_k)\| = \|P_{U_k^{\perp}} P_{\widetilde{U}_k} \|$ (see, for instance, \cite[Exercises VII. 1. 9-1.11]{Bhatia}), we have
\begin{align*}
    \|\sin\angle(U_k, \widetilde{U}_k)\| = \|(P_J + Q) P_{\widetilde{U}_k} \| &\le \|P_J  P_{\widetilde{U}_k} \| + \|Q P_{\widetilde{U}_k} \|\\
    &\le \|U_J^\T \widetilde{U}_k\|_F + \|QP_{\widetilde{U}_k} \|.
\end{align*}
It remains to bound $\|QP_{\widetilde{U}_k} \| = \|Q \widetilde{U}_k\|$. From the spectral decomposition of $\widetilde{A}$, we have $(A+E) \widetilde{U}_k = \widetilde{U}_k \tilde{\Lambda}_k$. Multiplying by $Q$ on the left and noting that $QA=0$, we obtain $Q E \widetilde{U}_k = Q \widetilde{U}_k \widetilde{\Lambda}_k$. Since $\widetilde{\Lambda}_k$ is invertible, we have 
\[ Q \widetilde{U}_k  = Q E \widetilde{U}_k \widetilde{\Lambda}_k^{-1}.\] To see that $\widetilde{\Lambda}_k$ is indeed invertible, note that by Weyl's inequality, for every $1\le s \le k$, 
\begin{align}\label{eq:weylbd-eig}
 \frac{3}{2}\lambda_s \ge \lambda_s +\|E\|   \ge \widetilde{\lambda}_s \ge \lambda_s - \|E\|\ge \frac{1}{2}\lambda_s,
\end{align}
where the last inequality follows from our assumption $\lambda_s\ge\lambda_k \ge 2\|E\|$. Therefore, 
\begin{align}\label{eq:0603-proj-null}
\|QP_{\widetilde{U}_k} \| = \|Q E \widetilde{U}_k \widetilde{\Lambda}_k^{-1}\| \le \|Q E \widetilde{U}_k \| \cdot\|\widetilde{\Lambda}_k^{-1}\| \le \frac{\|E\|}{\widetilde{\lambda}_k} \le 2\frac{\|E\|}{{\lambda}_k}. 
\end{align}

Now we turn to part (ii). Note that \eqref{eq:row-det} is given in Eq. (10) of \cite{Wang24}. It remains to prove \eqref{eq:row-proj-det}. It follows from the decomposition $\widetilde{u}_s = P_{U_k}\widetilde{u}_s + P_J \widetilde{u}_s + Q \widetilde{u}_s $ for each $1\le s \le k$ that 
\[ \widetilde{U}_k = P_{U_k}\widetilde{U}_k + P_J \widetilde{U}_k + Q \widetilde{U}_k \]
and thus
\begin{align*}
    \|\widetilde{U}_k -P_{U_k}\widetilde{U}_k\|_{2,\infty} =  \|P_J \widetilde{U}_k + Q \widetilde{U}_k \|_{2,\infty} &\le \|U_J U_J^\T \widetilde{U}_k\|_{2,\infty} + \|Q\widetilde{U}_k \|_{2,\infty}\\
    &\le \|U\|_{2,\infty} \cdot \|U_J^\T \widetilde{U}_k\|_F + \|Q\widetilde{U}_k \|_{2,\infty}.
\end{align*}
To finish the proof, we show that 
\begin{align*}
    \|Q\widetilde{U}_k \|_{2,\infty} \le & \frac{13}{2}\sqrt{k} \|U\|_{2,\infty} \frac{R_{\max}}{\widetilde{\lambda}_k^2}+ 2 \sqrt{\sum_{s=1}^k \widetilde{\lambda}_s^2  \max_{1\le i\le n}\left\|  e_i^\T Q \Xi(z) U \right\|^2}.
\end{align*}
For $1\le s \le k$, from the $(A+E)\widetilde{u}_s = \widetilde{\lambda}_s \widetilde{u}_s$, we have $A \widetilde{u}_s =  (\widetilde{\lambda}_s -E) \widetilde{u}_s$. Therefore,
\begin{align}\label{eq:solve-pertb-u}
    \widetilde{u}_s = G(\widetilde{\lambda}_s) A \widetilde{u}_s = \Phi(\widetilde{\lambda}_s) A\widetilde{u}_s + \Xi(\widetilde{\lambda}_s) A \widetilde{u}_s
\end{align}
and for any canonical vector $e_l$, we obtain
\begin{align*}
e_l^\T Q \widetilde{U}_k = \big(e_l^\T Q \Phi(\widetilde{\lambda}_s) A\widetilde{u}_s\big)_{s\in[k]} + \big(e_l^\T Q \Xi(\widetilde{\lambda}_s) A\widetilde{u}_s\big)_{s\in[k]}:=\mathtt{q}_l^{(1)} + \mathtt{q}_l^{(2)}.
\end{align*}
It follows that \[\|Q\widetilde{U}_k \|_{2,\infty} = \max_{1\le l \le n} \| e_l^\T Q \widetilde{U}_k\|\le  \max_{1\le l \le n} \big(\|\mathtt{q}_l^{(1)}\| + \|\mathtt{q}_l^{(2)}\|\big).\]
We bound $\|\mathtt{q}_l^{(1)}\|$ and $ \|\mathtt{q}_l^{(2)}\|$ on the right hand side. We start from \[\|\mathtt{q}_l^{(1)}\|= \sqrt{\sum_{s=1}^k\big(e_l^\T Q \Phi(\widetilde{\lambda}_s) A\widetilde{u}_s\big)^2}.\]
For each $l\in [n]$, observe that 
\begin{align}\label{eq:20260313-1}
 |e_l^\T Q \Phi(\widetilde{\lambda}_s) A\widetilde{u}_s| = |e_l^\T Q \Phi(\widetilde{\lambda}_s)U \cdot \Lambda U^\T \widetilde{u}_s| \le \|e_l^\T Q \Phi(\widetilde{\lambda}_s)U\| \cdot \|\Lambda U^\T \widetilde{u}_s\|.   
\end{align}
Note that \[\|\Lambda U^\T \widetilde{u}_s\|\le 2 \widetilde{\lambda}_s.\] This is because from $(A+E) \widetilde{u}_s = \widetilde{\lambda}_s \widetilde{u}_s$, we have  $A\widetilde{u}_s= U\Lambda U^\T \widetilde{u}_s = (\widetilde{\lambda}_s I - E)\widetilde{u}_s$. Multiplying $U^\T$ on each side, we further get $\Lambda U^\T \widetilde{u}_s = U^\T (\widetilde{\lambda}_s - E)\widetilde{u}_s$ and 
\[ \|\Lambda U^\T \widetilde{u}_s\| \le (\widetilde{\lambda}_s +\|E\|)\le 2 \widetilde{\lambda}_s \] by \eqref{eq:weylbd-eig} and our assumption $\lambda_s\ge \lambda_k \ge 2\|E\|$. 

Next, we estimate $\|e_l^\T Q \Phi(\widetilde{\lambda}_s)U\|$. Since $\widetilde{\lambda}_s \ge \frac{1}{2}\lambda_s \ge 3\sqrt{R_{\max}}$, 
Note that $\Phi(\widetilde{\lambda}_s)$ is a diagonal matrix. Indeed, for a general real diagonal matrix $D = \diag(d_1,\cdots,d_n)$,  
\begin{align*}
     e_l^\T Q D U = e_l^\T (I - U U^\T) DU &= e_l^\T D U - e_l^\T U(U^\T DU)\\
     &= e_l^\T U (d_l I - U^\T D U) = e_l^\T U \cdot U^\T (d_l I - D )U
\end{align*}
 due to $e_l^\T D =d_l e_l^\T$ and hence,
\[ \|e_l^\T Q D U  \| \le \|e_l^\T U\| \max_i |d_i -d_l| \le \|U\|_{2,\infty} \osc(D).\]
We apply this to $D=\Phi(\widetilde{\lambda}_s)$. Since $\Phi(z) = \frac{1}{z}I_n + \frac{1}{z^3} \mathcal{R} + \varepsilon(z)$ from \eqref{eq:Phi-expansion}, by Lemma \ref{lem:osc}, we get
\[ \osc(\Phi(\widetilde{\lambda}_s)) \le \left( 1+ 15\frac{R_{\max}}{\widetilde{\lambda}_s^2}\right) \frac{\osc(\mathcal{R})}{\widetilde{\lambda}_s^3} \le \frac{8}{3} \frac{\osc(\mathcal{R})}{\widetilde{\lambda}_s^3},\]
where we used $\widetilde{\lambda}_s \ge \frac{1}{2}\lambda_s \ge 3 \sqrt{R_{\max}}$ from \eqref{eq:weylbd-eig} and our assumption $\lambda_k \ge 6\sqrt{R_{\max}}.$ Thus,
\[\|e_l^\T Q \Phi(\widetilde{\lambda}_s) U\| \le  \|U\|_{2,\infty}\osc(\Phi(\widetilde{\lambda}_s)) \le \frac{8}{3} \|U\|_{2,\infty} \frac{\osc(\mathcal{R})}{\widetilde{\lambda}_s^3}.\]
 Combining the above estimates, we continue from \eqref{eq:20260313-1} to achieve
\[  |e_l^\T Q \Phi(\widetilde{\lambda}_s) A\widetilde{u}_s| \le \frac{8}{3} \|U\|_{2,\infty} \frac{\osc(\mathcal{R})}{\widetilde{\lambda}_s^2}\le \frac{8}{3} \|U\|_{2,\infty} \frac{\osc(\mathcal{R})}{\widetilde{\lambda}_k^2}.\]
Thus 
\[\|\mathtt{q}_l^{(1)}\|= \sqrt{\sum_{s=1}^k\big(e_l^\T Q \Phi(\widetilde{\lambda}_s) A\widetilde{u}_s\big)^2} \le \frac{8}{3}\sqrt{k}  \|U\|_{2,\infty} \frac{\osc(\mathcal{R})}{\widetilde{\lambda}_k^2}. \]

The estimate for \[ \|\mathtt{q}_l^{(2)}\|= \sqrt{\sum_{s=1}^k\big(e_l^\T Q \Xi(\widetilde{\lambda}_s) A\widetilde{u}_s\big)^2}\] is simpler. Using an estimate similar to \eqref{eq:20260313-1}, we have
\begin{align*}
    |e_l^\T Q \Xi(\widetilde{\lambda}_s) A\widetilde{u}_s|=|e_l^\T Q \Xi(\widetilde{\lambda}_s)U \cdot \Lambda U^\T\widetilde{u}_s|& \le \|e_l^\T Q \Xi(\widetilde{\lambda}_s) U\| \cdot \|\Lambda U^\T\widetilde{u}_s\| \\
    &\le 2\widetilde{\lambda}_s \|e_l^\T Q \Xi(\widetilde{\lambda}_s) U\|.
\end{align*}
Thus
\[\|\mathtt{q}_l^{(2)}\| \le 4 \sqrt{\sum_{s=1}^k \widetilde{\lambda}_s^2  \max_{1\le i\le n}\left\|  e_i^\T Q \Xi(\widetilde{\lambda}_s) U \right\|^2}.\]
Combining all the above estimates, we complete the proof.

\subsection{Proof of Proposition \ref{prop:detbd-vector}}
Fix $s\in [k]$. Recall that $J=[r]\setminus[k]$. We partition $J$ into two parts \[\mathcal{J}:=\{k+1,\dots,r_{+}\}\quad \text{and}\quad \mathcal{N}:=\{r_{+}+1,\dots,r\}.\] Then
\[ \|U_J^\T \widetilde{u}_s\|^2 = \|U_{\mathcal{J}}^\T \widetilde{u}_s\|^2 + \|U_\mathcal{N}^\T \widetilde{u}_s\|^2.\]
We estimate each term separately. 

\medskip
\noindent{\bf Step 1. Bound on $\|U_{\mathcal{J}}^\T \widetilde{u}_s\|$.} Set $\mathcal{I}=[k]\cup\{r_{+}+1,\dots,r\}=[r]\setminus \mathcal{J}$. We then have the decomposition 
\[ A = U \Lambda U^\T=U_{\mathcal{J}} \Lambda_{\mathcal J} U_{\mathcal J}^\T + U_{\mathcal{I}} \Lambda_{\mathcal I} U_{\mathcal I}^\T.\]
Multiplying both sides of \eqref{eq:solve-pertb-u} by $U_{\mathcal J}^\T$ from the left yields 
\begin{align}
    U_{\mathcal J}^\T \widetilde{u}_s &= U_{\mathcal J}^\T \Phi(\widetilde{\lambda}_s) A \widetilde{u}_s + U_{\mathcal J}^\T \Xi(\widetilde{\lambda}_s) A \widetilde{u}_s\\
    &= U_{\mathcal J}^\T \Phi(\widetilde{\lambda}_s) U_{\mathcal J} \Lambda_{\mathcal J} \cdot U_{\mathcal J}^\T \widetilde{u}_s + U_{\mathcal J}^\T \Phi(\widetilde{\lambda}_s)U_{\mathcal{I}}\cdot \Lambda_{\mathcal I} U_{\mathcal I}^\T \widetilde{u}_s + U_{\mathcal J}^\T \Xi(\widetilde{\lambda}_s)U\cdot \Lambda U^\T\widetilde{u}_s,
\end{align}
where we have used the decomposition of $A$. Isolating the term $U_{\mathcal J}^\T \widetilde{u}_s$ gives
\begin{align}\label{eq:0325-eq1}
    \left( I - U_{\mathcal J}^\T \Phi(\widetilde{\lambda}_s) U_{\mathcal J} \Lambda_{\mathcal J} \right) U_{\mathcal J}^\T \widetilde{u}_s = U_{\mathcal J}^\T \Phi(\widetilde{\lambda}_s)U_{\mathcal{I}}\cdot \Lambda_{\mathcal I} U_{\mathcal I}^\T \widetilde{u}_s + U_{\mathcal J}^\T \Xi(\widetilde{\lambda}_s)U\cdot \Lambda U^\T\widetilde{u}_s.
\end{align}
Denote \[L_{\mathcal J}(\widetilde{\lambda}_s):= I - U_{\mathcal J}^\T \Phi(\widetilde{\lambda}_s) U_{\mathcal J} \Lambda_{\mathcal J}.\] We first bound the smallest singular value $s_{\min}(L_{\mathcal J}(\widetilde{\lambda}_s))$. We split 
\begin{align*}
  U_{\mathcal J}^\T \Phi(\widetilde{\lambda}_s) U_{\mathcal J} &= \diag(u_j^\T \Phi(\widetilde{\lambda}_s)u_j)_{j\in \mathcal{J}}  + \off(U_{\mathcal J}^\T \Phi(\widetilde{\lambda}_s) U_{\mathcal J})  \\
  &:=\diag(u_j^\T \Phi(\widetilde{\lambda}_s)u_j)_{j\in \mathcal{J}}  + Q_{\mathcal J}(\widetilde{\lambda}_s) 
\end{align*}
and rewrite
\begin{align*}
    L_{\mathcal J}(\widetilde{\lambda}_s) &= I -\diag(u_j^\T \Phi(\widetilde{\lambda}_s)u_j)_{j\in \mathcal{J}} \Lambda_{\mathcal J}- Q_{\mathcal J}(\widetilde{\lambda}_s)\Lambda_{\mathcal J}\\
    &=\diag\left( 1- \lambda_j u_j^\T \Phi(\widetilde{\lambda}_s)u_j \right)_{j\in \mathcal{J}}- Q_{\mathcal J}(\widetilde{\lambda}_s)\Lambda_{\mathcal J}.
\end{align*}
Then 
\begin{align}\label{eq:0325-bd-ss}
    s_{\min}(L_{\mathcal J}(\widetilde{\lambda}_s)) \ge \frac{\Delta_{\mathcal J}^{\Phi}(\widetilde{\lambda}_s)}{\widetilde{\lambda}_s} - \|Q_{\mathcal J}(\widetilde{\lambda}_s)\Lambda_{\mathcal J}\|
\end{align}
with $\Delta_{\mathcal J}^{\Phi}(\widetilde{\lambda}_s)=\widetilde{\lambda}_s\min_{j\in \mathcal{J}}|1-\lambda_j u_j^\T\Phi(\widetilde{\lambda}_s)u_j|$ defined in \eqref{eq:gap-det-J}. 

Now we bound $\|Q_{\mathcal J}(\widetilde{\lambda}_s)\Lambda_{\mathcal J}\|$ on the right-hand side of \eqref{eq:0325-bd-ss}. Since $\widetilde{{\lambda}_s}\ge 3\sqrt{R_{\max}}$, by the second-order expansion of $\Phi$ in \eqref{eq:phi-2ndexp}, we have
\begin{align*}
    Q_{\mathcal J}(\widetilde{\lambda}_s)=\off(U_{\mathcal J}^\T \Phi(\widetilde{\lambda}_s) U_{\mathcal J}) = \off\left( \frac{\mathcal V_{\mathcal{J}\mathcal{J}}}{\widetilde{\lambda}_s^3} \right) + \off\left( U_{\mathcal J}^\T \varepsilon(\widetilde{\lambda}_s)U_{\mathcal J} \right).
\end{align*}
Note, by orthogonality, that
\begin{align*}
    \left\|\off\left( U_{\mathcal J}^\T \varepsilon(\widetilde{\lambda}_s)U_{\mathcal J} \right)\right\| =\left\| \off\left( U_{\mathcal J}^\T \big(\varepsilon(\widetilde{\lambda}_s)-c^* I\big)U_{\mathcal J} \right)\right\|\le \|\varepsilon(\widetilde{\lambda}_s)-c^* I\| = \frac{1}{2}\osc(\varepsilon(\widetilde{\lambda}_s)),
\end{align*}
where we select $c^*= \frac{1}{2}\left(\max_i \varepsilon_i(\widetilde{\lambda}_s) +\min_i \varepsilon_i(\widetilde{\lambda}_s) \right)$.
It follows, by applying Proposition \ref{lem:osc}, that 
\begin{align*}
    \|Q_{\mathcal J}(\widetilde{\lambda}_s)\Lambda_{\mathcal L}\| \le \frac{\lambda_{k+1}}{\widetilde{\lambda}_s^3}\|\off(\mathcal V_{\mathcal{J}\mathcal{J}})\| + \frac{13}{2}\frac{R_{\max}\lambda_{k+1}}{\widetilde{\lambda}_s^5} \osc(\mathcal R).
\end{align*}
Our assumption implies that 
\[ \frac{\Delta_{\mathcal J}^{\Phi}(\widetilde{\lambda}_s)}{\widetilde{\lambda}_s} \ge 2 \|Q_{\mathcal J}(\widetilde{\lambda}_s)\Lambda_{\mathcal L}\|\]
and hence, from \eqref{eq:0325-bd-ss}, we get
\[s_{\min}(L_{\mathcal J}(\widetilde{\lambda}_s)) \ge \frac{1}{2}\frac{\Delta_{\mathcal J}^{\Phi}(\widetilde{\lambda}_s)}{\widetilde{\lambda}_s}.\]

Next, we turn to the term on the right-hand side of \eqref{eq:0325-eq1}. Applying the triangle inequality yields 
\begin{align}
    &\|U_{\mathcal J}^\T \Phi(\widetilde{\lambda}_s)U_{\mathcal{I}}\cdot \Lambda_{\mathcal I} U_{\mathcal I}^\T \widetilde{u}_s + U_{\mathcal J}^\T \Xi(\widetilde{\lambda}_s)U\cdot \Lambda U^\T\widetilde{u}_s\| \\
    &\le \|U_{\mathcal J}^\T \Phi(\widetilde{\lambda}_s)U_{\mathcal{I}}\|\cdot \|\Lambda_{\mathcal I} U_{\mathcal I}^\T \widetilde{u}_s\| + \|U_{\mathcal J}^\T \Xi(\widetilde{\lambda}_s)U\|\cdot \|\Lambda U^\T\widetilde{u}_s\|.
\end{align}
For the bound on $\|\Lambda U^\T\widetilde{u}_s\|$, starting from $(A+E)\widetilde{u}_s = \widetilde{\lambda}_s \widetilde{u}_s $, we have $U \Lambda U^\T \widetilde{u}_s = (\widetilde{\lambda}_s I - E) \widetilde{u}_s$. Multiplying $U^\T$ from the left on both sides and taking the norms yield
\[\|\Lambda U^\T \widetilde{u}_s \| = \|U^\T (\widetilde{\lambda}_s I - E) \widetilde{u}_s\| \le \widetilde{\lambda}_s + \|E\| \le  2\widetilde{\lambda}_s.\]
The last step follows from \eqref{eq:weylbd-eig} and our assumption $\lambda_s \ge 2\|E\|$. Similarly, we also have 
\[ \|\Lambda_{\mathcal I} U_{\mathcal I}^\T \widetilde{u}_s\| \le 2\widetilde{\lambda}_s.\]
The right-hand side of \eqref{eq:0325-eq1} is bounded by 
\begin{align*}
    2\widetilde{\lambda}_s \left( \|U_{\mathcal J}^\T \Phi(\widetilde{\lambda}_s)U_{\mathcal{I}}\| + \|U^\T \Xi(\widetilde{\lambda}_s)U\| \right).
\end{align*}
It remains to estimate  $\|U_{\mathcal J}^\T \Phi(\widetilde{\lambda}_s)U_{\mathcal{I}}\|$. We plug in the second-order expansion of $\Phi(\widetilde{\lambda}_s)$ in \eqref{eq:phi-2ndexp}. By the orthogonality of eigenvectors, we get
\begin{align*}
    U_{\mathcal J}^\T \Phi(\widetilde{\lambda}_s)U_{\mathcal{I}} &= \frac{U_{\mathcal J}^\T \mathcal{R} U_{\mathcal I}}{\widetilde{\lambda}_s^3} + U_{\mathcal J}^\T \varepsilon(\widetilde{\lambda}_s) U_{\mathcal I}\\
    &= \frac{\mathcal{V}_{\mathcal{J}\mathcal{I}}}{\widetilde{\lambda}_s^3} +U_{\mathcal J}^\T \left(\varepsilon(\widetilde{\lambda}_s) -c^* I \right)U_{\mathcal I},
\end{align*}
where $c^*= \frac{1}{2}\left(\max_i \varepsilon_i(\widetilde{\lambda}_s) +\min_i \varepsilon_i(\widetilde{\lambda}_s) \right).$ By Proposition \ref{lem:osc}, we obtain
\begin{align*}
    \|U_{\mathcal J}^\T \Phi(\widetilde{\lambda}_s)U_{\mathcal{I}}\| &\le \frac{\|\mathcal{V}_{\mathcal{J}\mathcal{I}}\|}{\widetilde{\lambda}_s^3} + \|\varepsilon(\widetilde{\lambda}_s) -c^* I\|\\
    &\le \frac{\|\mathcal{V}_{\mathcal{J}\mathcal{I}}\|}{\widetilde{\lambda}_s^3} + \frac{1}{2} \osc(\varepsilon(\widetilde{\lambda}_s))
    \le \frac{\|\mathcal{V}_{\mathcal{J}\mathcal{I}}\|}{\widetilde{\lambda}_s^3} + \frac{13}{2}\frac{R_{\max}}{\widetilde{\lambda}_s^5} \osc(\mathcal R).
\end{align*}
Therefore, the right-hand side of can be bounded by 
\[2\widetilde{\lambda}_s \left( \frac{\|\mathcal{V}_{\mathcal{J}\mathcal{I}}\|}{\widetilde{\lambda}_s^3} + \frac{13}{2}\frac{R_{\max}}{\widetilde{\lambda}_s^5} \osc(\mathcal R) + \|U^\T \Xi(\widetilde{\lambda}_s)U\| \right).\]

Finally, combining the above estimates, we continue from \eqref{eq:0325-eq1} to achieve the bound 
\begin{align*}
  \|U_{\mathcal J}^\T \widetilde{u}_s\| &\le  \frac{2\widetilde{\lambda}_s \left( \frac{\|\mathcal{V}_{\mathcal{J}\mathcal{I}}\|}{\widetilde{\lambda}_s^3} + \frac{13}{2}\frac{R_{\max}}{\widetilde{\lambda}_s^5} \osc(\mathcal R) + \|U^\T \Xi(\widetilde{\lambda}_s)U\| \right)}{s_{\min}(L_{\mathcal J}(\widetilde{\lambda}_s))} \\
  &\le \frac{4 }{\Delta_{\mathcal J}^{\Phi}(\widetilde{\lambda}_s)}\left(\frac{\|\mathcal{V}_{\mathcal{J}\mathcal{I}}\|}{\widetilde{\lambda}_s} + \frac{13}{2}\frac{R_{\max}}{\widetilde{\lambda}_s^3} \osc(\mathcal R) + \widetilde{\lambda}_s^2 \|U^\T \Xi(\widetilde{\lambda}_s)U\| \right).
\end{align*}

\medskip
\noindent{\bf Step 2. Bound on $\|U_{\mathcal{N}}^\T \widetilde{u}_s\|$.} Denote the index set $\mathcal{K}:=[r]\setminus \mathcal{N}$. We decompose
\[ A = U_{\mathcal{N}} \Lambda_{\mathcal N} U_{\mathcal N}^\T + U_{\mathcal{K}} \Lambda_{\mathcal K} U_{\mathcal K}^\T.\]
Following the same derivation of \eqref{eq:0325-eq1}, we have 
\begin{align}\label{eq:0325-eq-neg}
    \left( I - U_{\mathcal N}^\T \Phi(\widetilde{\lambda}_s) U_{\mathcal N} \Lambda_{\mathcal N} \right) U_{\mathcal N}^\T \widetilde{u}_s = U_{\mathcal N}^\T \Phi(\widetilde{\lambda}_s)U_{\mathcal{K}}\cdot \Lambda_{\mathcal K} U_{\mathcal K}^\T \widetilde{u}_s + U_{\mathcal N}^\T \Xi(\widetilde{\lambda}_s)U\cdot \Lambda U^\T\widetilde{u}_s.
\end{align}
Denote \[\widehat{L}_{\mathcal N}(\widetilde{\lambda}_s):= I - U_{\mathcal N}^\T \Phi(\widetilde{\lambda}_s) U_{\mathcal N} \Lambda_{\mathcal N}.\] We first bound $s_{\min}\big(\widehat{L}_{\mathcal N}(\widetilde{\lambda}_s)\big)$. Since $\Lambda_{\mathcal N}$ contains only negative eigenvalues, set $D_{\mathcal N}:=\diag(|\lambda_l|)_{l\in \mathcal{N}}=-\Lambda_{\mathcal N} \succ 0$. Denote $Q_{\mathcal N}:= U_{\mathcal N}^\T \Phi(\widetilde{\lambda}_s) U_{\mathcal N}$ for simplicity. Since $\phi_i (\widetilde{\lambda}_s)>0$, $Q_{\mathcal N} \succ 0$. Rewrite  
\begin{align*}
   \widehat{L}_{\mathcal N}(\widetilde{\lambda}_s) = I + Q_{\mathcal N} D_{\mathcal N} = Q_{\mathcal N}^{1/2} \left( I + Q_{\mathcal N}^{1/2} D_{\mathcal N} Q_{\mathcal N}^{-1/2}\right)Q_{\mathcal N}^{1/2}.
\end{align*}
Observe that $I + Q_{\mathcal N}^{1/2} D_{\mathcal N} Q_{\mathcal N}^{-1/2} \succeq I$. We get
\begin{align*}
    s_{\min}\big(\widehat{L}_{\mathcal N}(\widetilde{\lambda}_s)\big) \ge \sqrt{\frac{\lambda_{\min}(Q_{\mathcal N})}{\lambda_{\max}(Q_{\mathcal N})}}.
\end{align*}
To bound the eigenvalues of $Q_{\mathcal N}$, note that for any unit vector $v$, $v^\T Q_{\mathcal N} v = (U_{\mathcal N}v)^\T \Phi(\widetilde{\lambda}_s) (U_{\mathcal N} v).$ Combining Lemma \ref{lem:bdphi} (i),  we get 
\[ \frac{1}{2\widetilde{\lambda}_s} \le \min_i \phi(\widetilde{\lambda}_s) \le \lambda_{\min}(Q_{\mathcal{N}}) \le \lambda_{\max}(Q_{\mathcal N}) \le \max_i \phi(\widetilde{\lambda}_s)\le \frac{3}{2\widetilde{\lambda}_s}\]
and thus
\[s_{\min}\big(\widehat{L}_{\mathcal N}(\widetilde{\lambda}_s)\big) \ge \sqrt{\frac{1}{3}}.\]

By the same estimation as in {\bf Step 1}, we bound the right-hand side of \eqref{eq:0325-eq-neg} by 
\begin{align*}
    2\widetilde{\lambda}_s \left( \frac{\|\mathcal{V}_{\mathcal{N}\mathcal{K}}\|}{\widetilde{\lambda}_s^3} + \frac{13}{2}\frac{R_{\max}}{\widetilde{\lambda}_s^5} \osc(\mathcal R) + \|U^\T \Xi(\widetilde{\lambda}_s)U\| \right).
\end{align*}

Hence, 
\begin{align*}
   \|U_{\mathcal N}^\T \widetilde{u}_s\|  &\le 2\sqrt{3} \left( \frac{\|\mathcal{V}_{\mathcal{N}\mathcal{K}}\|}{\widetilde{\lambda}_s^2} + \frac{13}{2}\frac{R_{\max}}{\widetilde{\lambda}_s^4} \osc(\mathcal R) + \widetilde{\lambda}_s\|U^\T \Xi(\widetilde{\lambda}_s)U\|\right).
\end{align*}
This completes the proof.
    

\section{Proof of Proposition \ref{prop:bdDvw}}\label{app:local-part}
Before we prove Proposition \ref{prop:bdDvw}, we introduce the necessary notation and a cumulant expansion formula to be used in the proof (adapted from \cite{HKR18}). 

For any matrix $F\in\R^{n\times n}$ and vectors $v,w\in\R^n$,  we write
\[
F_{vw}:=v^{\T}Fw,\qquad
F_{iw}:=e_i^{\T}Fw,\qquad
F_{vj}:=v^{\T}Fe_j,\qquad
F_{ij}:=e_i^{\T}Fe_j,
\]
where $\{e_1,\dots,e_n\}$ denotes the standard basis of $\R^n$.
Moreover, for $1\le i,j\le n$ define
\begin{equation}\label{def:Deltaij}
\Delta^{ij}:=(1+\delta_{ij})^{-1}\bigl(e_i e_j^\T+e_j e_i^\T\bigr),
\end{equation}
and for a differentiable function $f=f(H)$ of a symmetric matrix $H=(H_{ij})$ (with $\{H_{ij}:i\le j\}$ treated as independent variables) set
\[
\partial_{ij} f \;:=\; \frac{\partial}{\partial H_{ij}} f(H)
\;=\; \left.\frac{\mathrm{d}}{\mathrm{d}t}\right|_{t=0}\, f\bigl(H+t\,\Delta_{ij}\bigr).
\]
Since there is potential ambiguity between $\bigl|\partial_{ij}^{\,r}(G_{vw})\bigr|$ and $\bigl|(\partial_{ij}^{\,r}G)_{vw}\bigr|$, we explicitly note that
\[
\partial_{ij}^{\,r}\bigl(v^{\T}Gw\bigr)\;=\;v^{\T}\bigl(\partial_{ij}^{\,r}G\bigr)w,
\]
and throughout the proof we therefore write $\bigl|\partial_{ij}^{\,r}G_{vw}\bigr|$ with no confusion.

\begin{lemma}[Cumulant expansion, Lemma~2.4 from \cite{HKR18}]
\label{lem:cumulant-expansion}
Let $h$ be a real random variable with finite moments of all orders. For $k\ge 1$ define the $k$-th cumulant by
\[
\mathcal{C}_k(h) := (-\mathrm{i})^{k}\left.\frac{\mathrm{d}^k}{\mathrm{d}t^k}\log \mathbb{E}\bigl[e^{\mathrm{i}th}\bigr]\right|_{t=0}.
\]
Let $f:\mathbb{R}\to\mathbb{C}$ be smooth, and write $f^{(k)}$ for its $k$-th derivative. Fix $\ell\in\mathbb{N}$. Then
\begin{equation}\label{eq:cumulant-expansion}
\mathbb{E}\bigl[hf(h)\bigr] = \sum_{k=0}^{\ell}\frac{1}{k!}\,\mathcal{C}_{k+1}(h)\,\mathbb{E}\bigl[f^{(k)}(h)\bigr] + \mathtt R_{\ell+1},
\end{equation}
where the remainder $R_{\ell+1}$ satisfies, for any $t>0$,
\begin{align}\label{eq:cumulant-remainder}
|\mathtt R_{\ell+1}|
\le&\ 2^{\ell}(\ell+2) \cdot \left(
\mathbb{E}\Bigl[\sup_{|x|\le |h|}\bigl|f^{(\ell+1)}(x)\bigr|^{2}\Bigr]\cdot
\mathbb{E}\bigl[|h|^{2\ell+4}\mathbf{1}(|h|>t)\bigr]
\right)^{1/2} \nonumber \\
&+ 2^{\ell}(\ell+2) \cdot \mathbb{E}|h|^{\ell+2}\cdot
\sup_{|x|\le t}\bigl|f^{(\ell+1)}(x)\bigr|.
\end{align}
\end{lemma}

Now we provide the proof. Define
\begin{equation}\label{eq:defL}
L:=\mathrm{diag}(l_1,\dots,l_n) \quad \text{with} \quad l_k := \sum_{s=1}^n \sigma_{k s}^2 G_{ss}.
\end{equation}
The main technical task is to establish:
\begin{prop}
\label{prop:bdDvw-app}
Fix $D>0$. Let $v,w\in\mathbb R^n$ be deterministic unit vectors, and let
$z\in\mathbb C$ satisfy $|z|\ge 6\sqrt{R_{\max}}$.

\textup{(i)} Under the general sub-Gaussian regime, with probability at least
$1-6n^{-D-4}$,
\begin{equation*}
   \left| v^\T (E-L) G w \right| \le C (D+6)^{3/2} \frac{\beta K\sigma_{\max}\log^2 n}{|z|}+C n^{-500},
\end{equation*}
where $C\ge 1$ is an absolute constant.

\textup{(ii)} Under the bounded sparse-entry regime, namely if $\max_{i,j}|E_{ij}|\le c_0 K\sigma_{\max}$ almost surely, then, with probability at least $1-6n^{-D-4}$,
\begin{equation*}
   \left| v^\T (E-L) G w \right| \le C(D+6)c_0 
    \frac{K\sigma_{\max}\log n}{|z|},
\end{equation*}
where $C>0 $ is an absolute constant.
\end{prop}

The key to proving Proposition \ref{prop:bdDvw-app} is the following moment estimates, whose proof is deferred to Section \ref{sec:proof-moments}. For brevity, we denote
\begin{equation}\label{eq:defD}
\mathcal D := (E-L)G.
\end{equation}
and bound the moments of $$\mathcal{D}_{vw} = v^\T (E-L)G w.$$
We define a regime-dependent truncation parameter 
\begin{equation}\label{def:Ln}
\mathtt L_n:=
\begin{cases}
100\sqrt{D+6} K\sigma_{\max}\log n, & \text{in the general sub-Gaussian regime},\\
c_0K\sigma_{\max}, & \text{in the bounded sparse-entry regime}.
\end{cases}
\end{equation}
Define 
\[
\Omega_n^0 :=\left\{
\|E\|\le 3\sqrt{R_{\max}}, \quad \max_{a,b}|E_{ab}|\le \mathtt L_n \right\}
\]
and a slightly larger event
\begin{align}\label{def:entrybound}
\Omega_n:=
\left\{\|E\|\le 3.5\sqrt{R_{\max}},\quad \max_{a,b}|E_{ab}|\le 2\mathtt L_n \right\}.
\end{align}
Let $\chi=\chi(E) \in [0,1]$ be a smooth cutoff such that
\[
\chi(E)=1\quad\text{on }\Omega_n^0, \qquad \chi(E)=0\quad\text{outside }\Omega_n.
\]
Our main technical estimate is the following bound:
\begin{lemma} \label{lem:proof-moments}For any integer $p\ge 3$, 
\begin{align*}
\E(|\mathcal{D}_{vw}|^{2p} \chi(E)) \le& CK \beta \sum_{m=1}^{3} p^{m-1} \frac{\sigma_{\max}^m }{|z|^m}\E\left[ \left( |\mathcal D_{vw}|+ 22\frac{\mathtt L_n}{|z|} \right)^{2p-m} \chi(E)\right] \\
& + C K^2 \beta^2 \sum_{m=1}^4 p^{m-1}\frac{\sigma_{\max}^{m+1}}{|z|^{m+1}}\E\left[ \left( |\mathcal D_{vw}|+ 22\frac{\mathtt L_n}{|z|} \right)^{2p-m} \chi(E)\right] + \varepsilon_p 
\end{align*}
for an absolute constant $C\ge 1$. Here, we define
\begin{align}\label{def:varepsilon_p}
    \varepsilon_p:=
    C^p p^3
    \exp \bigl(-5000(D+6)(\log n)^2\bigr)
    \dfrac{(\beta K\sigma_{\max})^2}{|z|^2}
\end{align}
in the general sub-Gaussian regime and $\varepsilon_p:=0$ in the bounded sparse-entry regime.
\end{lemma}
Using this lemma, we first obtain a bound on 
\[ \mathtt{Z}_p:= \left(\E(|\mathcal{D}_{vw}|^{2p} \chi(E))\right)^{1/2p}.\]
By the H\"older inequality, for any $q\le 2p$, 
\[
\E[ |\mathcal{D}_{vw}|^{q} \chi(E) ] \le \left( \E |\mathcal{D}_{vw}|^{2p} \chi(E)\right)^{\frac{q}{2p}} = \mathtt{Z}_p^q.
\]
The same argument gives
\[
\E\left[ \left( |\mathcal D_{vw}|+ 22\frac{\mathtt{L}_n}{|z|} \right)^{q} \chi(E)\right] \le \left( \E \left[ \left( |\mathcal D_{vw}|+ 22\frac{\mathtt{L_n}}{|z|} \right)^{2p} \chi(E)\right]\right)^{\frac{q}{2p}}.
\]
Note that by the Minkowski inequality, 
\[\left( \E \left[ \left( |\mathcal D_{vw}|+ 22\frac{\mathtt{L_n}}{|z|} \right)^{2p} \chi(E)\right]\right)^{\frac{1}{2p}} \le \left( \E |\mathcal{D}_{vw}|^{2p} \chi(E)\right)^{\frac{1}{2p}}  + 22\frac{\mathtt{L_n}}{|z|}= \mathtt{Z}_p + 22\frac{\mathtt{L_n}}{|z|}.\]
We get
\[
\E\left[ \left( |\mathcal D_{vw}|+ 22\frac{\mathtt{L_n}}{|z|} \right)^{q} \chi(E)\right] \le \left( \mathtt{Z}_p + 22\frac{\mathtt{L}_n}{|z|}\right)^q.
\]
Therefore, Lemma~\ref{lem:proof-moments} implies 
\begin{align}\label{eq:0526-Zp}
\mathtt{Z}_p^{2p} \le& CK \beta \sum_{m=1}^{3} p^{m-1} \frac{\sigma_{\max}^m }{|z|^m} \left( \mathtt{Z}_p + 22\frac{\mathtt{L}_n}{|z|}\right)^{2p-m} \nonumber\\
& + C K^2 \beta^2 \sum_{m=1}^4 p^{m-1}\frac{\sigma_{\max}^{m+1}}{|z|^{m+1}}\left( \mathtt{Z}_p + 22\frac{\mathtt{L}_n}{|z|}\right)^{2p-m} + \varepsilon_p,
\end{align}
where $\varepsilon_p$ is defined in \eqref{def:varepsilon_p}.

We claim that
\begin{equation}
\label{eq:Zp-bound}
\mathtt Z_p \le
C\left(
p\frac{\mathtt{L_n}}{|z|}+p\frac{\beta K\sigma_{\max}}{|z|}+\left(\frac{\beta K\sigma_{\max}}{|z|}\right)^2+\varepsilon_p^{1/(2p)}
\right).
\end{equation}
To see this, consider the case where $\mathtt Z_p \le 88p\frac{\mathtt{L_n}}{|z|} + 2\varepsilon_p^{1/(2p)}$. In this regime, \eqref{eq:Zp-bound} holds trivially. Now, let us assume $\mathtt Z_p > 88p\frac{\mathtt{L_n}}{|z|}$ and $\mathtt Z_p > 2\varepsilon_p^{1/(2p)}$. Thus, $\mathtt Z_p+22\frac{\mathtt{L}_n}{|z|} \le (1+\frac{1}{4p})\mathtt Z_p$ and for $1\le m\le 4$,
\[\left(\mathtt Z_p+22\frac{\mathtt{L_n}}{|z|}\right)^{2p-m} \le \left( 1+ \frac{1}{4p}\right)^{2p}\mathtt Z_p^{2p-m}\le e^{1/2}\mathtt Z_p^{2p-m}.\]
Moreover, the last term in \eqref{eq:0526-Zp} satisfies $\varepsilon_p \le 2^{-2p}\mathtt Z_p^{2p}$ and can be absorbed into the left-hand side. Hence, \eqref{eq:0526-Zp} simplifies to 
\begin{align}\label{eq:0526-Zp-simple}
\mathtt Z_p^{2p} \le  \sum_{m=1}^{3}   C K\beta p^{m-1} \frac{\sigma_{\max}^m}{|z|^m} \mathtt Z_p^{2p-m} +  \sum_{m=1}^{4} C K^2\beta^2 p^{m-1} \frac{\sigma_{\max}^{m+1}}{|z|^{m+1}}
\mathtt Z_p^{2p-m}
\end{align}
for an absolute constant $C\ge1$. Next, we apply Young's inequality for $\mathtt{B} \cdot \mathtt Z_p^{2p-m}$. Let $\delta>0$ be a small constant to be specified later. Since
\begin{align*}
\mathtt{B} \cdot \mathtt Z_p^{2p-m} = (\delta\mathtt Z_p^{2p-m})\cdot  \frac{\mathtt{B}}{\delta}& \le \frac{(\delta\mathtt Z_p^{2p-m})^{2p/(2p-m)}}{2p/(2p-m)} + \frac{(\mathtt{B}/\delta)^{2p/m}}{2p/m}\\
&= \frac{\delta^{2p/(2p-m)}}{2p/(2p-m)}\mathtt Z_p^{2p} + \frac{(1/\delta)^{2p/m}}{2p/m} \mathtt{B}^{\frac{2p}{m}}.
\end{align*}
Choose $\delta=1/100$. Then $\frac{\delta^{2p/(2p-m)}}{2p/(2p-m)} < 1/10$ for $1\le m\le 4$. We absorb the $\mathtt Z_p^{2p}$ terms to the left-hand side of \eqref{eq:0526-Zp-simple}. Continuing from \eqref{eq:0526-Zp-simple}, we have
\begin{align*}
\mathtt Z_p^{2p} \le  \sum_{m=1}^{3} \left( C K\beta p^{m-1} \frac{\sigma_{\max}^m}{|z|^m}\right)^{\frac{2p}{m}} + \sum_{m=1}^{4} \left( C K^2 \beta^2 p^{m-1} \frac{\sigma_{\max}^{m+1}}{|z|^{m+1}}\right)^{\frac{2p}{m}}.
\end{align*}
Taking the $2p$-th roots on both sides yields
\begin{align*}
\mathtt Z_p \le C\max_{1\le m\le3}
\left( K\beta p^{m-1} \frac{\sigma_{\max}^m}{|z|^m} \right)^{\frac 1m} + C\max_{1\le m\le4} \left( K^2\beta^2 p^{m-1} \frac{\sigma_{\max}^{m+1}}{|z|^{m+1}} \right)^{\frac 1m}.
\end{align*}
Since $K\beta \ge 1$, $p\ge 3$, and $|z|\ge \sigma_{\max}$, we have for $1\le m \le 3$,
\[\left( K\beta p^{m-1} \frac{\sigma_{\max}^m}{|z|^m} \right)^{\frac 1m} = (K\beta)^{\frac 1m} p^{\frac{m-1}{m}} \frac{\sigma_{\max}}{|z|}\le p\frac{\beta K\sigma_{\max}}{|z|}\]
and for $2\le m \le 4$,
\begin{align*}
    (K\beta)^{\frac 2m} p^{\frac{m-1}{m}} \left(\frac{\sigma_{\max}}{|z|}\right)^{1 + \frac 1m} &= \left(\frac{\beta K\sigma_{\max}}{|z|}\right) \cdot (K\beta)^{\frac{2-m}{m}} \cdot p^{\frac{m-1}{m}} \cdot \left(\frac{\sigma_{\max}}{|z|}\right)^{\frac 1m}\\
    &\le p\frac{\beta K\sigma_{\max}}{|z|}.
\end{align*}
The above inequality simplifies to 
\begin{align*}
\mathtt Z_p\le C\left(p\frac{\beta K\sigma_{\max}}{|z|}+\left(\frac{\beta K\sigma_{\max}}{|z|}\right)^2\right).
\end{align*}
Combining this bound with the exceptional case proves \eqref{eq:Zp-bound}.

\medskip

Now we prove Proposition~\ref{prop:bdDvw}. Let us first consider the general sub-Gaussian regime.  By the union bound and our sub-Gaussian tail assumption,
\begin{align*}
    \mathbb{P}\left(\max_{a,b}|E_{ab}|> 2\mathtt{L}_n \right) \le \sum_{a,b} \mathbb{P}(|E_{ab}| > \mathtt{L}_n)\le 2n^2 \exp\left(-10^4(D+6)(\log n)^2\right) \le 2n^{-D-4}.
\end{align*}
Moreover,  by Lemma \ref{lem:bdE}, we have $$\mathbb{P}(\|E\|> 3\sqrt{R_{\max}}) \le 3n^{-D-4}$$
by selecting the constant $C_0$ in Assumption~\ref{assump:extra} sufficiently large. 
Hence, 
\begin{align*}
\mathbb{P} \big((\Omega_n^0)^c \big) \le 5n^{-D-4}.
\end{align*}

From the decomposition $ |\mathcal{D}_{vw}| = |\mathcal{D}_{vw}|\mathbf{1}_{\Omega_n^0} + |\mathcal{D}_{vw}|\mathbf{1}_{(\Omega_n^0)^c} $ and $\chi(E)=1$ on $\Omega_n^0$, we have  $|\mathcal{D}_{vw}|\mathbf{1}_{\Omega_n^0} \le |\mathcal{D}_{vw}|\chi(E)$. Thus
$$ |\mathcal{D}_{vw}| \le |\mathcal{D}_{vw}|\chi(E) + |\mathcal{D}_{vw}|\mathbf{1}_{(\Omega_n^0)^c}.$$
It follows that 
$$ \mathbb{P}(|\mathcal{D}_{vw}| \ge t) \le \mathbb{P}(|\mathcal{D}_{vw}|\chi(E) \ge t) + \mathbb{P}\big((\Omega_n^0)^c\big).$$
Then by Markov's inequality and $0\le \chi \le 1$, for any integer $p\ge 3$,
\begin{align*}
\mathbb{P}( |\mathcal{D}_{vw}| \ge t) &\le  \frac{\E(|\mathcal{D}_{vw}|^{2p} \chi(E))}{t^{2p}} + 5n^{-D-4}\\
&= \left(  \frac{\mathtt Z_p}{t} \right)^{2p}+ 5n^{-D-4}.
\end{align*}
We choose $p = \lceil(D+4)\log n\rceil$ and in view of \eqref{eq:Zp-bound},
\[
t = e \cdot C\left( p\frac{\mathtt{L}_n}{|z|} + p\frac{\beta K\sigma_{\max}}{|z|} + \left(\frac{\beta K\sigma_{\max}}{|z|}\right)^2 + \varepsilon_p^{1/(2p)} \right).
\]
With this choice of $t$, 
\[
\left( \frac{\mathtt Z_p}{t} \right)^{2p} \le e^{-2p} \le \exp\big(-2(D+4)\log n\big) = n^{-2D-8} \le n^{-D-4}.
\]
It follows that 
\[
\mathbb{P}( |\mathcal{D}_{vw}| \ge t) \le n^{-D-4} + 5n^{-D-4} = 6n^{-D-4}.
\]

Finally, let us simplify the terms in $t$ under the choice of $p$. For the term $\varepsilon_p^{1/(2p)}$, note that by Assumption~\ref{assump:extra},  $R_{\max}\le n\sigma_{\max}^2$ and $ \sqrt{R_{\max}}\ge C_0 (D+8) K\sigma_{\max}\log n.$ We have
\[ K\le \frac{\sqrt n}{C_0(D+8)\log n}.\]
As discussed in Remark~\ref{rem:beta}, $\beta \le C\sqrt{\log(eK)} \le C \sqrt{\log n}$. Under our assumption $|z| \ge 6\sqrt{R_{\max}}$, $\frac{\beta K\sigma_{\max}}{|z|} \le C'$. Thus
\begin{align*}
\varepsilon_p^{1/(2p)} 
&= C^{1/2} p^{3/(2p)} \exp\left( - \frac{10(D+6)\log^2 n}{2p} \right) \left(\frac{\beta K\sigma_{\max}}{|z|}\right)^{1/p} \\
&\le C\exp\left( - \frac{5000(D+6)\log^2 n}{2(D+5)\log n} \right)\le C \exp(-500 \log n) = C n^{-500}.
\end{align*}

Next, substituting $p \le (D+5)\log n$ and $\mathtt{L}_n = 5\sqrt{D+6}K\sigma_{\max}\log n$, we also have
\[
p\frac{\mathtt{L}_n}{|z|} \le 5(D+6)^{3/2} \frac{K\sigma_{\max}\log^2 n}{|z|} 
\quad \text{and} \quad 
p\frac{\beta K\sigma_{\max}}{|z|} \le (D+5) \frac{\beta K\sigma_{\max}\log n}{|z|}.
\]
Note that 
\begin{align}\label{eq:0604-absorb}
    p\frac{\beta K\sigma_{\max}}{|z|} \ge \left(\frac{\beta K\sigma_{\max}}{|z|}\right)^2
\end{align}
by selecting the constant $C_0$ sufficiently large in Assumption~\ref{assump:extra}. Thus, we conclude that with probability at least $1 - 6n^{-D-4}$,
\begin{align*}
|\mathcal{D}_{vw}| &\le C \left( 5(D+6)^{3/2} \frac{K\sigma_{\max}\log^2 n}{|z|} +  (D+5) \frac{\beta K\sigma_{\max}\log n}{|z|} + n^{-500}\right)\\
&\le C (D+6)^{3/2} \frac{K\sigma_{\max}\log^2 n}{|z|} + C n^{-500}.
\end{align*}

We now prove the bounded sparse-entry case. In this regime,
$\mathtt L_n=c_0K\sigma_{\max}$ and  $\varepsilon_p=0$.
Thus from \eqref{eq:Zp-bound}, together with $\eqref{eq:0604-absorb}$, we obtain
\[
    \mathtt Z_p
    \le
    C
    \left(
        p\frac{c_0K\sigma_{\max}}{|z|}
        +
        p\frac{\beta K\sigma_{\max}}{|z|}
    \right).
\]
Moreover, the entrywise truncation event is deterministic, i.e., $\max_{i,j}|E_{ij}|\le \mathtt L_n$ almost surely. Also, $\beta\le c_0$.
Taking
$p=\lceil(D+4)\log n\rceil$ and applying the same Markov argument as above, together with the operator norm
event, yields
\[
    |\mathcal D_{vw}|
    \le C (D+6) \frac{c_0K\sigma_{\max}\log n}{|z|}
\]
with probability at least $1-6n^{-D-4}$. This proves part \textup{(ii)} and completes the proof of Proposition \ref{prop:bdDvw-app}.

\subsection{Proof of Lemma \ref{lem:proof-moments}}\label{sec:proof-moments}
We prove Lemma \ref{lem:proof-moments} in this section. Let
\[
Q:=\mathcal D_{vw}^{p-1}\overline{\mathcal D}_{vw}^{p}.
\]
We rewrite
\begin{equation}
\label{eq:E|Dvw|^2p}
    \E[|\mathcal{D}_{vw}|^{2p}\chi(E)] = \E\left[ (EG)_{vw} Q \chi(E)\right] - \E\left[ \left(LG\right)_{vw} Q \chi(E)\right].
\end{equation}
We analyze the first term on the right-hand side of \eqref{eq:E|Dvw|^2p} using the cumulant expansion. For $1\le i,j\le n$, define the auxiliary functions
\[
h_{ij}(E):=G_{jw}Q, \qquad f_{ij}^{\chi}(E):=h_{ij}(E)\chi(E).
\] 
Thus
\[\E\left[(EG)_{vw}Q\,\chi(E)\right] = \sum_{i,j}v_i\,\E\left[E_{ij}f_{ij}^{\chi}(E)\right].\]
Let $\chi\equiv \chi(E)$. For $l \le 3$, the product rule
gives
\[ \partial_{ij}^l f_{ij}^{\chi} = \chi \cdot\partial_{ij}^l h_{ij} +
\sum_{a=1}^l \binom{l}{a} (\partial_{ij}^{l-a}h_{ij})(\partial_{ij}^a\chi).
\]
The first term on the right-hand side is the main term estimated below. The
remaining terms contain at least one derivative of the cutoff function.  We choose
\(\chi\) so that, for \(1\le a\le3\),
\[
|\partial_{ij}^a\chi(E)|\le C_a \mathtt L_n^{-a}.
\]
These derivative terms are supported in
\(\Omega_n\setminus\Omega_n^0\). Using the derivative bounds proved below,
together with Cauchy-Schwarz inequality and the high-probability estimate for
\((\Omega_n^0)^c\), the total contribution of the terms containing
\(\partial_{ij}^a\chi\), \(a\ge1\), is absorbed into an error
\(O(n^{-D-5})\), after increasing \(C_0\) if necessary.

In the estimates below, we write only the terms in which no derivative falls
on $\chi$. The terms involving derivatives of $\chi$ are supported on
$\Omega_n\setminus\Omega_n^0$. By the same derivative estimates used below
and the probability bound for this transition region, their total contribution
is negligible compared with the final high-probability bound. We therefore
record only the main terms in the moment estimate. The remainder
$\varepsilon_p$ below refers only to the tail part of the cumulant remainder;
in the bounded sparse-entry regime this tail part is zero. Therefore, to simplify notation, we suppress the cutoff $\chi(E)$ from the main estimates. From this point on, we simply work with
\[
f_{ij}(E):=G_{jw}Q
\]

Since $E$ is symmetric,  the independent random variables are $E_{ij}$ ($1\le i\le j \le n)$. For $i<j$, we write $\partial_{ij}$ for differentiation with respect to $E_{ij} = E_{ji}$, and we set $\partial_{ji}:=\partial_{ij}$. For the diagonal entries, 
$\partial_{ii}$ denotes differentiation with respect to $E_{ii}$. With this convention, for $x,y\in\mathbb{R}^n$,
\begin{align}\label{eq:partialG}
    \partial_{ij}G_{xy} = (1+\delta_{ij})^{-1}(G_{xi}G_{jy} + G_{xj}G_{iy}),
\end{align}

Now, from 
\[ (EG)_{vw}Q =\sum_i v_i E_{ii}G_{iw}Q + \sum_{i<j}E_{ij}\bigl(v_iG_{jw}+v_jG_{iw}\bigr)Q =\sum_{i,j}v_iE_{ij}f_{ij}(E),
\]
we apply the cumulant expansion Lemma \ref{lem:cumulant-expansion} to each independent entry $E_{ij}$ ($1\le i\le j \le n)$ with $l=2$ to get
\begin{align}\label{cumulant_expansion}
\E\!\left[(EG)_{vw}Q\right]:= X_1+X_2+\mathtt{R}_{3},
\end{align}
where for $k=1,2$, 
\begin{equation}\label{def:X1X2}
X_k := \sum_{i, j} v_i\frac{1}{k!} \mathcal{C}_{k+1}(E_{ij}) \mathbb{E}\left[\partial_{ij}^k f_{ij}(E)\right].
\end{equation} 
Here the $k=0$ cumulant term vanishes because the entries of $E$ are
centered. The remainder is written as
\begin{equation}\label{def:R3}
\mathtt{R}_{3}
:= \sum_{i,j} v_i \mathtt{R}_3^{ij}
\end{equation}
where 
\begin{equation*}
\mathtt R_3^{ij}:=\E\bigl[E_{ij}f_{ij}(E)\bigr]-\sum_{k=1}^{2}\frac{1}{k!}\mathcal C_{k+1}(E_{ij})\E\bigl[\partial_{ij}^{\,k}f_{ij}(E)\bigr].
\end{equation*}

With these notations, we have
\begin{equation}
\label{eq:moments-Dvw}
    \E[|\mathcal{D}_{vw}|^{2p}] = \left(X_1 - \E\left[ \left(LG\right)_{vw} Q \right]\right) + X_2 + \mathtt{R}_{3}.
\end{equation}
The next lemma bounds each term on the right-hand side of \eqref{eq:moments-Dvw}.
\begin{lemma}\label{Key_Lemma}
\begin{enumerate}
    \item[(i)] For $X_1$,
    $$\left| X_1 -  \E\left[ \left(LG\right)_{vw}Q \right] \right|
    \le \frac{2}{3}\sigma_{\max}\frac{\E\left[|\mathcal{D}_{vw}|^{2p-1}\right]}{|z|} 
    + 36 p\sigma_{\max}^2 \frac{\E[|\mathcal{D}_{vw}|^{2p-2}]}{|z|^2}.$$
    
    \item[(ii)] For $X_2$,
    \[ 
    |X_2| \le \frac{4}{9}\beta K\sigma_{\max} \frac{\E\left[|\mathcal{D}_{vw}|^{2p-1}\right]}{|z|} +178 p \beta K\sigma_{\max}^2 \frac{\E[|\mathcal{D}_{vw}|^{2p-2}]}{|z|^2} + 700 p^2 \beta K\sigma_{\max}^3 \frac{\E[|\mathcal{D}_{vw}|^{2p-3}]}{|z|^3}.
    \]
    
    \item[(iii)] For $\mathtt R_3$, recalling $\mathtt{L}_n$ from \eqref{def:Ln}, we have
    \begin{align*}
    |\mathtt R_3| \le &CK^2 \beta^2 \sum_{m=1}^4 p^{m-1}\frac{\sigma_{\max}^{m+1}}{|z|^{m+1}}
\E\left[ \left( |\mathcal D_{vw}|+ 22\frac{\mathtt{L}_n}{|z|} \right)^{2p-m} \right]\\
& + C^p p^3 \exp\big(-10(D+6)\log^2 n\big)\frac{(\beta K \sigma_{\max})^2}{|z|^2}
    \end{align*}
    for an absolute constant $C\ge 1$.
\end{enumerate}
\end{lemma}
Again, we emphasize that in Lemma \ref{Key_Lemma}, we suppress $\chi(E)$ from the notation. More precisely, every expectation in the bounds for $X_1$, $X_2$, and
$\mathtt R_3$ is understood to include the factor $\chi(E)$, except for
the explicitly separated cutoff-derivative terms.  

With this lemma, the conclusion of Lemma \ref{lem:proof-moments} follows from  \eqref{eq:moments-Dvw}:
\begin{align*}
\E\big[|\mathcal{D}_{vw}|^{2p}\chi(E) \big] &\le\left| X_1 - \E\left[ \left(LG\right)_{vw} Q \right]\right| + |X_2| + |\mathtt{R}_{3} |\\
& \le CK \beta \sum_{m=1}^{3} p^{m-1} \frac{\sigma_{\max}^m }{|z|^m}\E\left[ |\mathcal{D}_{vw}|^{2p-m}\chi(E) \right] \\
&\quad + C K^2 \beta^2 \sum_{m=1}^4 p^{m-1}\frac{\sigma_{\max}^{m+1}}{|z|^{m+1}}\E\left[ \left( |\mathcal D_{vw}|+ 22\frac{\mathtt{L}_n}{|z|} \right)^{2p-m}\chi(E) \right] \\
&\quad + C^p p^3 \exp\big(-10(D+6)\log^2 n\big) \frac{(\beta K\sigma_{\max} )^2}{|z|^2}.
\end{align*}

In the remainder of this section, we prove Lemma \ref{Key_Lemma}. 
\paragraph{Proof of Lemma \ref{Key_Lemma} (i).}
From the definition of $L=\diag(\sum_{j=1}^n \sigma_{k j}^2 G_{jj})_{1\le k \le n}$ in \eqref{eq:defL}, we have
\begin{equation}
\label{eq:LGD^2p-1}
    \E [(LG)_{vw} Q]=\E [\sum_{i,j} v_i L_{ij} G_{jw} Q]=\sum_{i,j} v_i \sigma_{ij}^2\E [  G_{jj} G_{iw} Q].
\end{equation}

From \eqref{def:X1X2}, by product rule, we get
\begin{align*}
\label{eq:X1}
    X_1 &= \sum_{i,j} v_i  \mathcal{C}_{2}(E_{ij}) \E[\partial_{ij} f_{ij}(E) ] \\
    &= \sum_{i,j} v_i  \sigma_{ij}^2 \E[\partial_{ij} (G_{jw} Q)] \\
    &=\sum_{i,j} v_i \sigma_{ij}^2 \E\left[ (\partial_{ij}G_{jw}) Q \right] + \sum_{i,j} v_i \sigma_{ij}^2 \E\left[ G_{jw} (\partial_{ij}Q) \right].
\end{align*}
Plugging in $\partial_{ij}G_{jw}=(1+\delta_{ij})^{-1}\bigl(G_{ji}G_{jw}+G_{jj}G_{iw}\bigr)$ by \eqref{eq:partialG}, we further obtain
\begin{align*}
X_1=  \sum_{i,j} v_i \sigma_{ij}^2 \E\left[G_{jj}G_{iw}Q \right] +\sum_{i\neq j} v_i \sigma_{ij}^2 \E\left[G_{ji}G_{jw} Q \right]+ \sum_{i,j} v_i \sigma_{ij}^2 \E\left[ G_{jw} (\partial_{ij}Q) \right].
\end{align*}

Comparing with \eqref{eq:LGD^2p-1}, we arrive at
\begin{align}
    X_1 - \E\left[ \left(LG\right)_{vw} Q \right] 
    = T^{(1)} + T^{(2)},
\end{align}
where
\begin{align*}
T^{(1)}:=\E\Big[\sum_{i\neq j} v_i \sigma_{ij}^2 G_{ji}G_{jw} Q \Big], \qquad T^{(2)}:=\E\Big[ \sum_{i,j} v_i \sigma_{ij}^2 G_{jw} (\partial_{ij}Q) \Big].
\end{align*}

We now provide upper bounds for $T^{(1)}$ and $T^{(2)}$. 

\medskip
\noindent{\bf Step 1.} \emph{Estimates for $T^{(1)}$.} Applying Cauchy-Schwarz, we have 
\begin{align*}
\Big|\sum_{i\neq j} v_i \sigma_{ij}^2 G_{ji}G_{jw}\Big| \le \Big(\sum_{i \neq j} v_i^2 \sigma_{ij}^2 |G_{ij}|^2\Big)^{1/2} \Big(\sum_{i \neq j} \sigma_{ij}^2 |G_{jw}|^2\Big)^{1/2}.
\end{align*}
Note that 
\[ \sum_{i \neq j} v_i^2 \sigma_{ij}^2 |G_{ij}|^2 \le \sum_i v_i^2 \big(\sum_j\sigma_{ij}^2|G_{ij}|^2\big) \le \sigma_{\max}^2 \sum_i v_i^2 \|Ge_i\|^2\le \sigma_{\max}^2 \|G\|^2\]
and 
\[\sum_{i \neq j} \sigma_{ij}^2 |G_{jw}|^2 \le \sum_j |G_{jw}|^2\big(\sum_i \sigma_{ij}^2\big) \le R_{\max} \|Gw\|^2 \le R_{\max} \|G\|^2.  \]
We get 
\begin{align*}
\Big|\sum_{i,j} v_i \sigma_{ij}^2 G_{ji}G_{jw}\Big| 
\le \sigma_{\max} \sqrt{R_{\max}}\|G\|^2 
\le \frac{4\,\sigma_{\max}\sqrt{R_{\max}}}{|z|^2} \le \frac{2}{3}\frac{\sigma_{\max}}{|z|}.
\end{align*}
by applying \eqref{eq:crudebdG} and the assumption $|z|\ge 6\sqrt{R_{\max}}$. 
Consequently, 
\[
|T^{(1)}| \le \frac{2}{3}\frac{\sigma_{\max}}{|z|} \E\left[|\mathcal{D}_{vw}|^{2p-1}\right]. \] 

\medskip
\noindent{\bf Step 2.} \emph{Estimates for $T^{(2)}$.}  Recall that
\begin{equation*}
    T^{(2)} := \sum_{i,j} v_i \sigma_{ij}^2 \E\left[ G_{jw} \partial_{ij}(\mathcal{D}_{vw}^{p-1}\overline{\mathcal{D}}_{vw}^{p}) \right].
\end{equation*}
For $p\ge 2$, by the product rule,
 \[\partial_{ij}(\mathcal{D}_{vw}^{p-1}\overline{\mathcal{D}}_{vw}^{p}) = (p-1)(\partial_{ij} \mathcal{D}_{vw})\mathcal{D}_{vw}^{p-2}\overline{\mathcal{D}}_{vw}^{p} + p(\partial_{ij} \overline{\mathcal{D}}_{vw})|\mathcal{D}_{vw}|^{2(p-1)}.\]
For $p=1$, the first term is absent, and the same final bound follows. Hence, it suffices to bound the two quantities
\[ \mathtt{A}:= \sum_{i,j}v_i \sigma_{ij}^2 G_{jw}(\partial_{ij} \mathcal{D}_{vw}), \qquad \bar{\mathtt{A}} := \sum_{i,j}v_i \sigma_{ij}^2 G_{jw}(\partial_{ij} \overline{\mathcal{D}}_{vw}).\]

We first estimate $\mathtt{A}.$ From the identity $(zI-E)G=I$, we have $EG=zG-I$. Plugging into $\mathcal{D}=(E-L)G=EG-LG$ yields
\[
    \mathcal{D}=(zI-L)G-I.
\]
Thus $\partial_{ij}\mathcal D=(zI-L)\partial_{ij}G-(\partial_{ij}L)G$. Recall $\Delta_{ij}$ from \eqref{def:Deltaij}. We also have 
\begin{align*}
    \partial_{ij} \mathcal{D}_{vw} = v^\T(zI-L) G\Delta^{ij}G w - v^\T (\partial_{ij} L)G w.
\end{align*}
Denote $\widetilde{v}^\T:=v^\T(zI-L)$. 
From the definition of $L$ in \eqref{eq:defL} and the partial derivative formula \eqref{eq:partialG}, we compute \[\partial_{ij} l_s =\sum_{k} \sigma_{sk}^2 \partial_{ij} G_{kk} = 2(1+\delta_{ij})^{-1}\sum_{k}\sigma_{sk}^2 G_{ki}G_{kj}.\]
Hence, 
\[ v^\T (\partial_{ij} L)G w = 2(1+\delta_{ij})^{-1} \sum_{s,k} v_s \sigma_{s,k}^2 G_{ki} G_{kj} G_{sw}.\]
It follows that 
\begin{align}\label{eq:firstpartialD}
    \partial_{ij} \mathcal{D}_{vw} = (1+\delta_{ij})^{-1}\Big(G_{\widetilde{v} i} G_{jw} + G_{\widetilde{v} j} G_{iw} - 2\sum_{s,k} v_s \sigma_{sk}^2 G_{ki} G_{kj} G_{sw}\Big). 
\end{align}
We are ready to bound 
\begin{align*}
 |\mathtt{A}|=  \Big|\sum_{i,j}v_i \sigma_{ij}^2 G_{jw}(\partial_{ij} \mathcal{D}_{vw})\Big| \le &  \Big|\sum_{i,j}v_i \sigma_{ij}^2 G_{jw}^2G_{\widetilde{v} i}\Big| + \Big|\sum_{i,j}v_i \sigma_{ij}^2 G_{jw}G_{\widetilde{v} j} G_{iw}\Big|\\
   &+ 2\Big|\sum_{i,j}v_i \sigma_{ij}^2 G_{jw}\sum_{s,k} v_s \sigma_{s,k}^2 G_{ki} G_{kj} G_{sw}\Big| \\
   &:=T_1^{(2)} + T_2^{(2)} +T_3^{(2)}.
\end{align*}
Note that \[\|L\| \le R_{\max}\|G\| \le 2\frac{R_{\max}}{|z|} \le \frac{1}{18}|z|\] by \eqref{eq:crudebdG} and our supposition $|z| > 6 \sqrt{R_{\max}}$. Thus
$\|\widetilde{v}\| \le \|zI-L\| \le |z| + \|L\| \le \frac{19}{18}|z|$. We will repeatedly use $\sum_i |G_{iw}|^2 = \|Gw\|^2 \le \|G\|^2.$

Now we estimate $T_1^{(2)},T_2^{(2)}$ and $T_3^{(2)}$ respectively. For $T_1^{(2)}$, since $\max_{i,j} \sigma_{ij}^2 \le \sigma_{\max}^2$, 
\begin{align*}
  T_1^{(2)} &\le \sigma_{\max}^2 \sum_{i} |v_i| |G_{\widetilde{v} i}|\cdot\sum_j|G_{jw}|^2 \le  \sigma_{\max}^2 \|Gw\|^2 \sqrt{\sum_i |v_i|^2 \sum_i |G_{\widetilde{v} i}|^2 }\\
  &\le  \sigma_{\max}^2\|G\|^2 \|\widetilde{v}^\T G\| \le  \sigma_{\max}^2 \|G\|^3 \|\widetilde{v}\| \le \frac{76\sigma_{\max}^2}{9|z|^2}. 
\end{align*}
Similarly, 
\begin{align*}
  T_2^{(2)} &\le \sigma_{\max}^2 \sum_{i} |v_i| |G_{iw}|\cdot\sum_j |G_{jw}| |G_{\widetilde{v} j}| \le  \sigma_{\max}^2 \sqrt{\sum_i |v_i|^2 \sum_i |G_{iw}|^2 } \sqrt{\sum_j |G_{jw}|^2 \sum_j |G_{\widetilde{v} j}|^2}\\
  &\le  \sigma_{\max}^2 \|G\|^3 \|\widetilde{v}\| \le \frac{76\sigma_{\max}^2}{9|z|^2}. 
\end{align*}
For $T_3^{(2)}$, observe that
\[
\sum_s |v_s\sigma_{sk}^2| |G_{sw}| \le \sigma_{\max}^2 \sum_s |v_s| |G_{sw}| \le \sigma_{\max}^2 \|G\|\] and \[\sum_k |G_{ki} G_{kj}| \le \sqrt{\sum_k |G_{ki}|^2\sum_k |G_{kj}|^2} \le \|G\|^2
\]  
by Cauchy-Schwarz inequality. We further have
\begin{align*}
  T_3^{(2)} &\le 2\sigma_{\max}^2 \|G\|^3 \sum_{i,j} |v_i\sigma_{ij}| \cdot |\sigma_{ij} G_{jw}| \le 2\sigma_{\max}^2 \|G\|^3 \sqrt{\sum_{i,j} |v_i|^2\sigma_{ij}^2} \sqrt{\sum_{i,j} \sigma_{ij}^2 |G_{jw}|^2}\\
  &\le 2\sigma_{\max}^2 \|G\|^3 \sqrt{R_{\max}} \sqrt{R_{\max} \|Gw\|^2} \le 2\sigma_{\max}^2 R_{\max} \|G\|^4 \le \frac{8\sigma_{\max}^2}{9|z|^2}.
\end{align*}
Combining the above estimates, we obtain
\[ |\mathtt{A}| \le \frac{76\sigma_{\max}^2}{9|z|^2}+\frac{76\sigma_{\max}^2}{9|z|^2}+\frac{8\sigma_{\max}^2}{9|z|^2} = \frac{160\sigma_{\max}^2}{9|z|^2}\le \frac{18\sigma_{\max}^2}{|z|^2}.\]

The same estimate holds for $\bar{\mathtt A}$.  Indeed, since $\partial_{ij}\overline{\mathcal D}_{vw}=\overline{\partial_{ij}\mathcal D_{vw}},$ the preceding argument applies with $G,L,z$ replaced by
$\overline G,\overline L,\overline z.$ Hence, 
\[ |\bar{\mathtt A}| \le \frac{18\sigma_{\max}^2}{|z|^2}.\]

Finally, we bound 
\begin{align}
 |T^{(2)}| &\le (p-1) |\E \big[ |\mathtt{A}| \cdot \mathcal{D}_{vw}^{p-2}\overline{\mathcal{D}}_{vw}^{p} \big] | +  p |\E \big[ |\bar{\mathtt A}| \cdot |\mathcal{D}_{vw}|^{2(p-1)} \big] | \nonumber\\
 &\le \frac{36 p\sigma_{\max}^2}{|z|^2} \E[|\mathcal{D}_{vw}|^{2p-2}].
\end{align}
This completes the proof of Lemma \ref{Key_Lemma} (i).

\paragraph{Proof of Lemma \ref{Key_Lemma} (ii).} For brevity, we suppress the harmless distinction between $\mathcal D_{vw}^{p-1}\overline{\mathcal D}_{vw}^{p}$ and $\mathcal D_{vw}^{2p-1}$. The argument with complex conjugates is identical, using $\partial_{ij}\overline{\mathcal D}_{vw} = \overline{\partial_{ij}\mathcal D_{vw}}.$

Recall that
\begin{align}
\label{df:Xk}
X_2 &= 
\sum_{i,j} v_i\, \frac{1}{2}\, \mathcal{C}_{3}(E_{ij})\, \E\! \left[\partial_{ij}^{\,2}\! \left(G_{jw}\, \mathcal{D}_{vw}^{2p-1}\right)\right].
\end{align}
Using the Leibniz rule and our assumption $| \mathcal{C}_{3}(E_{ij})| =|\E(E_{ij}^3)| \le \beta K\sigma_{\max} \sigma_{ij}^2$, we get
\begin{align}\label{def:X2-bd}
    |X_2| \le \beta K\sigma_{\max} \sum_{s=0}^2 \E\Big[ \sum_{i,j}   |v_i|\sigma_{ij}^2 |\partial_{ij}^{2-s} (G_{jw})|\cdot |\partial_{ij}^{s} (\mathcal{D}_{vw}^{2p-1} )|\Big].
\end{align}
We estimate the three values $s=0,1,2$ separately. 

\medskip

\noindent $\bullet$ \textbf{Case 1: $s=0$. } We estimate
\[\sum_{i,j}   |v_i|\sigma_{ij}^2 |\partial_{ij}^{2} (G_{jw})| .\]
Using 
\[\partial_{ij}^2 G =2!G(\Delta^{ij}G)^2,\]
one obtains
\begin{align*}
\sum_{i,j}   |v_i|\sigma_{ij}^2 |\partial_{ij}^{2} (G_{jw})| \le 2 \sum_{i,j} |v_i| \sigma_{ij}^2 |e_j^\T G(\Delta^{ij}G)^2 w| = 2 \sum_{j} \big(\sum_i |v_i| \sigma_{ij}^2\big) \cdot |e_j^\T G(\Delta^{ij}G)^2 w|.
\end{align*}
Let $\mathtt{v} = (|v_i|)$ and $\mathtt{v} $ is still a unit vector. By Cauchy-Schwartz, together with $\|\Delta^{ij}\| \le 1$, we get
\begin{align*}
\sum_{i,j}   |v_i|\sigma_{ij}^2 |\partial_{ij}^{2} (G_{jw})| &\le 2 \sqrt{\sum_j \big(\sum_i  \sigma_{ij}^2 \mathtt{v}_i\big)^2} \sqrt{\sum_j |e_j^\T G(\Delta^{ij}G)^2 w|^2}\\
&= 2  \|\Sigma \mathtt{v} \| \cdot \|G(\Delta^{ij}G)^2 w\| \le 2\|\Sigma\| \cdot \|G\|^3.
\end{align*}
Since $\|\Sigma\| \le \sqrt{\|\Sigma\|_1 \cdot \|\Sigma\|_{\infty}} = R_{\max}$, $\|G\| \le 2/|z|$ and $|z|\ge 6\sqrt{R_{\max}}$, this gives
\begin{align*}
\sum_{i,j}   |v_i|\sigma_{ij}^2 |\partial_{ij}^{2} (G_{jw})| \le 16\frac{R_{\max}}{|z|^3}\le \frac{4}{9} \frac{1}{|z|}.
\end{align*}
Therefore, the contribution from the $s=0$ term in \eqref{def:X2-bd} is bounded by 
\begin{align*}
 \frac{4}{9}\beta K\sigma_{\max} \frac{\E|\mathcal D_{vw}|^{2p-1}}{|z|}.
\end{align*}

\medskip

\noindent $\bullet$ \textbf{Case 2: $s=1$. } We estimate 
\[
\E\Big[ \sum_{i,j}   |v_i|\sigma_{ij}^2 |\partial_{ij} G_{jw} |\cdot |\partial_{ij}(\mathcal{D}_{vw}^{2p-1} )|\Big].
\]  
From \eqref{eq:partialG}, we first have
\[ \partial_{ij}G_{jw} = (1+\delta_{ij})^{-1}(G_{ij}G_{jw} + G_{jj}G_{iw}).\]
By the chain rule,
\[ |\partial_{ij}(\mathcal{D}_{vw}^{2p-1} )| \le 2p |\mathcal D_{vw}|^{2p-2}
|\partial_{ij}\mathcal D_{vw}|.\]
Recall from \eqref{eq:firstpartialD} that
\begin{align*}
    \partial_{ij} \mathcal{D}_{vw} = (1+\delta_{ij})^{-1}\Big(G_{\widetilde{v} i} G_{jw} + G_{\widetilde{v} j} G_{iw} - \mathcal{K}\Big),
\end{align*}
where for brevity we define \[\mathcal{K}:=2\sum_{s,k} v_s \sigma_{sk}^2 G_{ki} G_{kj} G_{sw}.\] We also recall that $\widetilde{v}^\T:=v^\T(zI-L)$ and $\|\widetilde{v}\| \le \frac{19}{18}|z|$. 

For vectors $\mathsf{v}, \mathsf{w}$, we will repeatedly use the following estimates:
\begin{align}\label{eq:ward}
  \sum_i |G_{i\mathsf{w}}|^2 = \|G \mathsf{w}\|^2 \le \|G\|^2\|\mathsf{w}\|^2 
\end{align}
Consequently, by the Cauchy-Schwarz inequality, 
\begin{align}\label{eq:wardCS}
    \sum_i |\mathsf{v}_i G_{i\mathsf{w}}| \le \|G\| \|\mathsf{v}\| \|\mathsf{w}\| \quad\text{and} \quad \sum_{j} |G_{\mathsf{v} j} G_{j\mathsf{w}}| \le \|G \mathsf{v}\| \|G\mathsf{w}\| \le \|G\|^2\|\mathsf{v}\| \|\mathsf{w}\|. 
\end{align}
It follows that
\begin{align}\label{eq:bdRvw}
    |\mathcal{K}| \le 2\sigma_{\max}^2 \big|\sum_s v_s G_{sw}\big| \cdot \big|\sum_k G_{ik} G_{kj} \big| 
    \le 2 \sigma_{\max}^2 \|G\|^3.
\end{align}

Now we are ready to estimate
\begin{align}\label{eq:mathcalT1}
\sum_{i,j} \big|v_i \sigma_{ij}^2 (\partial_{ij}G_{jw}) (\partial_{ij} \mathcal{D}_{vw})\big| 
\le& \, \sum_{i,j} \big|v_i \sigma_{ij}^2 (G_{ij}G_{jw} + G_{jj}G_{iw})(G_{\widetilde{v} i} G_{jw} + G_{\widetilde{v} j} G_{iw}) \big| \nonumber\\
&+ \sum_{i,j} |v_i \sigma_{ij}^2 (G_{ij}G_{jw} + G_{jj}G_{iw}) \mathcal{K}|.
\end{align}

For the first term on the right-hand side of \eqref{eq:mathcalT1}, expanding the product yields four terms. We bound each respectively using \eqref{eq:ward} and \eqref{eq:wardCS}:
\begin{align*}
    \sum_{i,j} |v_i \sigma_{ij}^2| | G_{ij}G_{jw} G_{\widetilde{v} i} G_{jw}| 
    \le \sigma_{\max}^2 \|G\| \sum_i|v_i G_{\widetilde{v} i}| \sum_j |G_{jw}|^2 
    \le \sigma_{\max}^2 \|G\|^4 \|\widetilde{v}\|.
\end{align*}
\begin{align*}
    \sum_{i,j} |v_i \sigma_{ij}^2| | G_{ij}G_{jw} G_{\widetilde{v} j} G_{iw}| 
    \le \sigma_{\max}^2 \|G\| \sum_i|v_i G_{iw}| \sum_j |G_{jw}||G_{\widetilde{v}j}| 
    \le \sigma_{\max}^2 \|G\|^4 \|\widetilde{v}\|.
\end{align*}
\begin{align*}
    \sum_{i,j} |v_i \sigma_{ij}^2| | G_{jj}G_{iw} G_{\widetilde{v} i} G_{jw}| 
    \le \sigma_{\max} \|G\|^2 \sum_i|v_i G_{\widetilde{v} i}| \sum_j |\sigma_{ij}G_{jw}| 
    \le \sigma_{\max} \sqrt{R_{\max}}\, \|G\|^4 \|\widetilde{v}\|.
\end{align*}
\begin{align*}
    \sum_{i,j} |v_i\sigma_{ij}^2| | G_{jj}G_{iw} G_{\widetilde{v} j} G_{iw}| 
    \le \sigma_{\max} \|G\|^2 \sum_i|v_i G_{iw}| \sum_j |\sigma_{ij} G_{\widetilde{v}j}| 
    \le \sigma_{\max} \sqrt{R_{\max}}\, \|G\|^4 \|\widetilde{v}\|.
\end{align*}
Summing these estimates and using the fact that $\sigma_{\max}^2 \le \sigma_{\max}\sqrt{R_{\max}}$, we obtain
\[ \sum_{i,j} \big| v_i \sigma_{ij}^2 (G_{ij}G_{jw} + G_{jj}G_{iw})(G_{\widetilde{v} i} G_{jw} + G_{\widetilde{v} j} G_{iw}) \big| \le 4\sigma_{\max} \sqrt{R_{\max}}\,\|G\|^4 \|\widetilde{v}\|.\]

For the second term on the right-hand side of \eqref{eq:mathcalT1}, note that by \eqref{eq:bdRvw}
\begin{align*}
\sum_{i,j} |v_i \sigma_{ij}^2 (G_{ij}G_{jw} + G_{jj}G_{iw}) \mathcal{K}| 
&\le 2\sigma_{\max}^2 \|G\|^3 \sum_{i,j} |v_i \sigma_{ij}^2| (|G_{ij} G_{jw}| + |G_{jj}G_{iw}| ) \\
&\le 4 \sigma_{\max}^2 R_{\max} \|G\|^5,
\end{align*}
where in the last inequality, we apply \eqref{eq:wardCS} to bound  
\[ \sum_{i,j} |v_i \sigma_{ij}^2| |G_{ij}| |G_{jw}| 
\le \sqrt{\sum_{i,j} v_i^2 \sigma_{ij}^2 |G_{ij}|^2} \sqrt{\sum_{i,j} \sigma_{ij}^2 |G_{jw}|^2} 
\le \sigma_{\max} \sqrt{R_{\max}} \|G\|^2\]
and 
\[ \sum_{i,j} |v_i \sigma_{ij}^2| |G_{jj}G_{iw}|
\le \|G\| \sum_i |v_i| |G_{iw}| \sum_j \sigma_{ij}^2 
\le R_{\max} \|G\|^2.\]

Combining the above estimates, we have  
\begin{align*}
\sum_{i,j} |v_i \sigma_{ij}^2 (\partial_{ij}G_{jw}) (\partial_{ij} \mathcal{D}_{vw})| 
&\le 4\sigma_{\max} \sqrt{R_{\max}}\|G\|^4 \|\widetilde{v}\| + 4 \sigma_{\max}^2 R_{\max} \|G\|^5.
\end{align*}
Using \eqref{eq:crudebdG} and our assumption $|z| \ge 6\sqrt{R_{\max}}$, we further have
\begin{align*}
\sum_{i,j} \big|v_i \sigma_{ij}^2 (\partial_{ij}G_{jw}) (\partial_{ij} \mathcal{D}_{vw})\big| \le 4\sigma_{\max} \frac{|z|}{6}\frac{2^4}{|z|^4} \frac{19 |z|}{18} + 4\sigma_{\max} \frac{|z|^2}{36} \frac{2^5}{|z|^5} 
\le \frac{64 \sigma_{\max}}{|z|^2}.
\end{align*}
Consequently, the contribution from the $s=1$ term in \eqref{def:X2-bd} is bounded by 
\begin{align*}
\beta K\sigma_{\max} \E\Big[ \sum_{i,j}   |v_i|\sigma_{ij}^2 |\partial_{ij} G_{jw} |\cdot |\partial_{ij}(\mathcal{D}_{vw}^{2p-1} )|\Big] \le 128 p\beta K\sigma_{\max}^2\frac{\E|\mathcal D_{vw}|^{2p-2}}{|z|^2}.
\end{align*}

\medskip

\noindent $\bullet$ \textbf{Case 3: $s=2$. } We estimate 
\[\E\Big[ \sum_{i,j}   |v_i|\sigma_{ij}^2 |G_{jw}| |\partial_{ij}^{2}  (\mathcal{D}_{vw}^{2p-1} )|\Big].\]
The second derivative gives
\[
|\partial_{ij}^2\left(\mathcal D_{vw}^{2p-1}\right)|
\le  2 p |\mathcal D_{vw}|^{2p-2} |\partial_{ij}^2\mathcal D_{vw}| + 4 p^2 |\mathcal D_{vw}|^{2p-3} |\partial_{ij}\mathcal D_{vw}|^2 .
\]
Hence, 
\begin{align}\label{eq:0518-3rd}
&\E\Big[ \sum_{i,j}   |v_i|\sigma_{ij}^2 |G_{jw}| |\partial_{ij}^{2} (\mathcal{D}_{vw}^{2p-1} )|\Big] \nonumber\\
&\le  2p \E \Big[\sum_{i,j}   |v_i|\sigma_{ij}^2 |G_{jw}| |\partial_{ij}^2\mathcal D_{vw}| |\mathcal D_{vw}|^{2p-2} \Big] + 4p^2 \E \Big[\sum_{i,j}   |v_i|\sigma_{ij}^2 |G_{jw}| |\partial_{ij}\mathcal D_{vw}|^2 |\mathcal D_{vw}|^{2p-3} | \Big] .
\end{align}
We first estimate 
\[ \sum_{i,j}   |v_i|\sigma_{ij}^2 |G_{jw}| |\partial_{ij}^2\mathcal D_{vw}|.\]
From \eqref{eq:firstpartialD}, we have
\[\partial_{ij}^2 \mathcal{D}_{vw} = (1+\delta_{ij})^{-1} \Big[ \partial_{ij}(G_{\widetilde{v} i} G_{jw}) + \partial_{ij}(G_{\widetilde{v} j} G_{iw}) - \partial_{ij}\mathcal{K}  \Big]\]
and thus 
\begin{align}\label{eq:case3part}
    \sum_{i,j} |v_i \sigma_{ij}^2 G_{jw} (\partial_{ij}^2 \mathcal{D}_{vw})| \le & \sum_{i,j}|v_i \sigma_{ij}^2 G_{jw} \partial_{ij}(G_{\widetilde{v} i} G_{jw})| + \sum_{i,j}|v_i \sigma_{ij}^2 G_{jw} \partial_{ij}(G_{\widetilde{v} j} G_{iw})|\nonumber\\
    &+ \sum_{i,j}|v_i \sigma_{ij}^2 G_{jw} \partial_{ij} \mathcal{K}|.
\end{align}
The first two terms on the right-hand side of \eqref{eq:case3part} can be estimated similarly.  We provide the calculation for the first term.  By \eqref{eq:partialG},
\begin{align*}
    &\sum_{i,j}|v_i \sigma_{ij}^2 G_{jw} \partial_{ij}(G_{\widetilde{v} i} G_{jw})| \\
    &\le \sum_{i,j}\Big[  \big|v_i \sigma_{ij}^2 G_{jw}^2 (G_{\widetilde{v}i} G_{ij} + G_{\widetilde{v} j} G_{ii}) \big| + \big|v_i \sigma_{ij}^2 G_{jw}G_{\widetilde{v} i} (G_{ij}G_{jw} + G_{jj}G_{iw}) \big| \Big].
\end{align*}
Note that each term on the right-hand side is bounded by $(\sigma_{\max}^2 + \sigma_{\max} \sqrt{R_{\max}}\,)\|G\|^4 \|\widetilde{v}\|$ by applying \eqref{eq:ward} and \eqref{eq:wardCS}. For instance,
\begin{align*}
  &\sum_{i,j}  \big|v_i \sigma_{ij}^2 G_{jw}^2 (G_{\widetilde{v}i} G_{ij} + G_{\widetilde{v} j} G_{ii}) \big| \\
  \le &\,\sigma_{\max}^2 \|G\| \sum_{i} |v_i G_{\widetilde{v} i}| \sum_j |G_{jw}|^2 + \sigma_{\max} \|G\|^2 \|\widetilde{v}\| \sum_j |G_{jw}|^2 \sum_{i} |v_i \sigma_{ij}| \\
  \le &\, (\sigma_{\max}^2 + \sigma_{\max} \sqrt{R_{\max}}\,)\|G\|^4 \|\widetilde{v}\|.
\end{align*}

Continuing from \eqref{eq:case3part}, we have
\begin{align*}
    \sum_{i,j} |v_i \sigma_{ij}^2 G_{jw} (\partial_{ij}^2 \mathcal{D}_{vw})| \le 8\sigma_{\max} \sqrt{R_{\max}} \,\|G\|^4 \|\widetilde{v}\| +  \sum_{i,j}|v_i G_{jw} \partial_{ij} \mathcal{K}|.
\end{align*}
It remains to estimate $ \sum_{i,j} |v_i G_{jw} \partial_{ij} \mathcal{K} |$, where 
\[ \partial_{ij} \mathcal{K} = 2 \sum_{s,k} v_s \sigma_{sk}^2 \partial_{ij}(G_{ki} G_{kj} G_{sw}).\]
We first estimate $|\partial_{ij} \mathcal{K}|$. Expanding $\partial_{ij}(G_{ki} G_{kj} G_{sw})$ using the product rule and \eqref{eq:partialG} yields several terms, each of which admits a similar treatment. We present the calculation for one representative term. For instance, 
\begin{align*}
    \sum_{s,k} |v_s \sigma_{sk}^2 (\partial_{ij} G_{ki}) G_{kj} G_{sw}| 
    &\le \sigma_{\max}^2 \sum_{s,k} |v_s (G_{ki}G_{jj} + G_{kj} G_{ij})G_{kj} G_{sw}| \\
    &\le \sigma_{\max}^2 \|G\| \sum_{s} |v_s G_{sw}| \sum_k (|G_{ki}| + |G_{kj}|) |G_{kj}|\\
    &\le 2 \sigma_{\max}^2 \|G\|^4
\end{align*}
where we applied \eqref{eq:ward} and \eqref{eq:wardCS} in the last inequality. Hence, 
\[|\partial_{ij} \mathcal{K}| \le 12\, \sigma_{\max}^2 \|G\|^4 \]
and by Cauchy-Schwarz and \eqref{eq:ward},
\begin{align*}
    \sum_{i,j} |v_i \sigma_{ij}^2 G_{jw} \partial_{ij} \mathcal{K}|
    \le 12 \sigma_{\max}^2 \|G\|^4 \sum_{i,j} |v_i \sigma_{ij}^2 G_{jw}| 
    \le 12 \sigma_{\max}^2 \|G\|^4 \cdot R_{\max} \|G\|
    = 12 \sigma_{\max}^2 R_{\max} \|G\|^5. 
\end{align*}
Finally, using $|z|\ge 6\sqrt{R_{\max}}$ and $\|G\|\le 2/|z|$, we have 
\begin{align}\label{eq:0520-previous1}
    \sum_{i,j} |v_i \sigma_{ij}^2 G_{jw} (\partial_{ij}^2 \mathcal{D}_{vw})| 
    &\le 8\sigma_{\max} \sqrt{R_{\max}} \|G\|^4 \|\widetilde{v}\| + 12 \sigma_{\max}^2 R_{\max} \|G\|^5 \nonumber\\
    &\le 8\sigma_{\max} \frac{|z|}{{6}} \frac{2^4}{|z|^4} \frac{19 |z|}{18} + 12\sigma_{\max} \frac{|z|^3}{{6}^3} \frac{2^5}{|z|^5} < 25\frac{\sigma_{\max}}{|z|^2}.
\end{align}


Next, we bound
\[ \sum_{i,j}   |v_i|\sigma_{ij}^2 |G_{jw}| |\partial_{ij}\mathcal D_{vw}|^2 .\]
Continuing from \eqref{eq:firstpartialD}, we get
\begin{align*}
(\partial_{ij} \mathcal{D}_{vw})^2 = (1+\delta_{ij})^{-2} \Big[(G_{\widetilde{v} i} G_{jw} + G_{\widetilde{v} j} G_{iw})^2 + \mathcal{K}^2 -2\mathcal{K}(G_{\widetilde{v} i} G_{jw} + G_{\widetilde{v} j} G_{iw}) \Big],
\end{align*}
where \[\mathcal{K}=2\sum_{s,k} v_s \sigma_{sk}^2 G_{ki} G_{kj} G_{sw}.\]
Plugging this expansion into $\sum_{i,j}|v_i  G_{jw}  (\partial_{ij} \mathcal{D}_{vw})^2|$, we first observe from \eqref{eq:bdRvw} that 
\[ \sum_{i,j}|v_i \sigma_{ij}^2 G_{jw}| |\mathcal K|^2 
\le (2\sigma_{\max}^2 \|G\|^3)^2 \sum_{i,j}|v_i \sigma_{ij}^2| |G_{jw}| 
\le 4\sigma_{\max}^4 \|G\|^6 \cdot R_{\max} \|G\|.\]
Next, using \eqref{eq:ward} and \eqref{eq:wardCS}, we also have
\begin{align*}
    &2\sum_{i,j}|v_i \sigma_{ij}^2 G_{jw}(G_{\widetilde{v} i} G_{jw} + G_{\widetilde{v} j} G_{iw})| |\mathcal{K}| \\
    \le \,&4 \sigma_{\max}^4 \|G\|^3 \left( \sum_{i,j}|v_i G_{\widetilde{v} i}| \cdot|G_{jw}|^2 + \sum_{i}|v_i G_{iw}| \sum_j |G_{jw} G_{\widetilde{v} j}|  \right)\\
    \le \,&8 \sigma_{\max}^4 \|G\|^6 \|\widetilde{v}\|.
\end{align*}
Finally, we claim that  
\[\sum_{i,j}|v_i \sigma_{ij}^2 G_{jw} (G_{\widetilde{v} i} G_{jw} + G_{\widetilde{v} j} G_{iw})^2| \le 4 \sigma_{\max}^2 \|G\|^5 \|\widetilde{v}\|^2.\]
To see this, we expand the square in the summation and show that each term is bounded by $\|G\|^5\|\widetilde{v}\|^2$, using a similar calculation with \eqref{eq:ward} and \eqref{eq:wardCS}. For instance, 
\[\sum_{i,j}|v_i  G_{jw} G_{\widetilde{v} i}^2 G_{jw}^2| \le \|G\|^2 \|\widetilde{v}\| \sum_i |v_i G_{\widetilde{v} i}| \sum_j |G_{jw}|^2 \le \|G\|^5\|\widetilde{v}\|^2,\]
and the remaining terms are estimated in the same manner.
\vspace{0.3em}

Hence, together with $\|\widetilde{v}\| \le \frac{19}{18}|z|$, we get
\begin{align}\label{eq:0520-previous2}
\sum_{i,j}|v_i \sigma_{ij}^2 G_{jw}  (\partial_{ij} \mathcal{D}_{vw})^2| &\le 4\sigma_{\max}^4 R_{\max} \|G\|^7 + 8\sigma_{\max}^4 \|G\|^6 \|\widetilde{v}\| + 4\sigma_{\max}^2 \|G\|^5\|\widetilde{v}\|^2 \nonumber\\
&\le \frac{2^3 \sigma_{\max}^2}{|z|^3} (4\frac{|z|^4}{6^4} \frac{2^4}{|z|^4} +  8\frac{|z|^2}{6^2} \frac{2^3}{|z|^3} \frac{19|z|}{18} + 4\frac{2^2}{|z|^2} \frac{19^2|z|^2}{18^2}) < 175\frac{\sigma_{\max}^2}{|z|^3}.
\end{align}

We continue from \eqref{eq:0518-3rd} to get
\begin{align*}
\E\Big[ \sum_{i,j}   |v_i|\sigma_{ij}^2 |G_{jw}| |\partial_{ij}^{2} (\mathcal{D}_{vw}^{2p-1} )|\Big] \le  50p \sigma_{\max}\frac{\E |\mathcal D_{vw}|^{2p-2} }{|z|^2}  + 700p^2 \sigma_{\max}^2 \frac{\E |\mathcal D_{vw}|^{2p-3}}{|z|^3}.
\end{align*}
Thus, the contribution from the $s=2$ term in \eqref{def:X2-bd} is bounded by 
\begin{align*}
&\beta K\sigma_{\max} \E\Big[ \sum_{i,j}   |v_i|\sigma_{ij}^2 |G_{jw}| |\partial_{ij}^{2} (\mathcal{D}_{vw}^{2p-1} )|\Big] \\
&\le  50p \beta K\sigma_{\max}^2\frac{\E |D_{vw}|^{2p-2} }{|z|^2}  + 700p^2 \beta K\sigma_{\max}^3 \frac{\E |\mathcal D_{vw}|^{2p-3}}{|z|^3}.
\end{align*}
Combining all the above estimates completes the proof of Lemma \ref{Key_Lemma} (ii).

\paragraph{A third-order bound for the remainder term.}
Before turning to the remainder term, we record an additional derivative estimate required by the cumulant remainder formula in Lemma~\ref{lem:cumulant-expansion}. This estimate serves as a third-order analogue of those used in bounding $X_2$ above, corresponding to the third-order derivatives of
\[
G_{jw}\mathcal D_{vw}^{2p-1}.
\]
The argument follows the structure of the $s=0,1,2$ cases in the proof of Lemma~\ref{Key_Lemma}~(ii). We isolate this estimate here for use in the remainder analysis.

Recall from the truncation parameter $\mathtt L_n$ from \eqref{def:Ln}.  For $1\le i,j\le n$ and $|x|\le \mathtt{L}_n$, let
\[
E^{ij(x)}:=E+(x-E_{ij})\Delta^{ij},
\qquad
G^{ij(x)}(z):=(zI-E^{ij(x)})^{-1}.
\]
Define
\[
L^{ij(x)}:=\operatorname{diag}(l_1^{ij(x)},\dots,l_n^{ij(x)})
\qquad\text{with}\quad
l_s^{ij(x)}:=\sum_{k=1}^n \sigma_{sk}^2G_{kk}^{ij(x)},
\]
and
\[ \mathcal D^{ij(x)}:=(E^{ij(x)}-L^{ij(x)})G^{ij(x)}(z). \]
We use the notation $\mathcal D^{ij(x)}_{vw} = v^\T \mathcal D^{ij(x)} w$ and $G_{vw}^{ij(x)} = v^\T G^{ij(x)} w$. 
Set 
\begin{align}\label{def:fij-x}
f_{ij}(E^{ij(x)}):=G_{jw}^{ij(x)} \left(\mathcal D^{ij(x)}_{vw}\right)^{2p-1}.
\end{align}
Recall from \eqref{def:entrybound} the event 
\[
\Omega_n:=
\left\{\|E\|\le 3.5\sqrt{R_{\max}},\quad \max_{a,b}|E_{ab}|\le 2\mathtt{L}_n \right\}.
\]
The goal is to prove the following estimate on the event $\Omega_n$ (the support of $\chi(E)$):
\begin{align}\label{eq:0520-bd-3rd}
\sum_{i,j}|v_i|\sigma_{ij}^2 \sup_{|x|\le \mathtt L_n} \left| \partial_{ij}^3 f_{ij}(E^{ij(x)})\right| \le C \sum_{m=1}^4 p^{m-1}\frac{\sigma_{\max}^{m-1}}{|z|^{m+1}}
\left( |\mathcal D_{vw}|+ 22\frac{\mathtt{L}_n}{|z|} \right)^{2p-m}
\end{align}
for an absolute constant $C>0$.

We start with some preliminary estimates.  On the event $\Omega_n$, we have, for $|x|\le \mathtt{L}_n$,
\begin{align}\label{eq:0520-Eijx}
\|E^{ij(x)}\| \le \|E\| + (|x| + |E_{ij}|)\|\Delta_{ij}\| \le 3.5\sqrt{R_{\max}}+3\mathtt{L}_n \le 4\sqrt{R_{\max}}
\end{align}
provided the constant $C_0$ in Assumption~\ref{assump:extra} is chosen sufficiently large. Hence, for $|z|\ge6\sqrt{R_{\max}}$, by the Neumann expansion, we have 
\begin{align}\label{eq:bdGijx}
 \|G^{ij(x)}(z)\| \le \frac{3}{|z|}.
\end{align}
Next, we show that 
\begin{align}\label{eq:bdDvwij-diff}
\sup_{|x|\le \mathtt{L}_n}|\mathcal D^{ij(x)}_{vw} -\mathcal D_{vw}| \le 22\frac{\mathtt{L}_n}{|z|}.
\end{align}
To see this, for $|x|\le \mathtt L_n$, we bound
\begin{align*}
&|\mathcal D^{ij(x)}_{vw} -\mathcal D_{vw}| = |v^\T (\mathcal D^{ij(x)} -\mathcal D)w|  \le \|\mathcal D^{ij(x)} -\mathcal D\|.
\end{align*}
From the definitions of $\mathcal D^{ij(x)}$ and $\mathcal D$, we further have
\begin{align*}
|\mathcal D^{ij(x)}_{vw} -\mathcal D_{vw}|&\le \|(E^{ij(x)}-L^{ij(x)})G^{ij(x)}(z) - (E-L) G(z)\|\\
&= \| [(E^{ij(x)}- E) -(L^{ij(x)} - L)] G^{ij(x)}(z) + (E-L)(G^{ij(x)}(z) -G(z)) \|\\
&\le (\|E^{ij(x)}- E\| + \|L^{ij(x)} - L\|) \|G^{ij(x)}(z)\| + (\|E\|+\|L\|)\|G^{ij(x)}(z) -G(z)\|.
\end{align*}
For simplicity, denote $P^{ij}:=E^{ij(x)}- E = (x-E_{ij})\Delta^{ij}$. Similar to \eqref{eq:0520-Eijx}, we have
\[ \|E^{ij(x)}- E \| = \|P^{ij}\| \le (|x| + |E_{ij}|)\|\Delta_{ij}\| \le 3 \mathtt{L}_n.\]
Next, by the resolvent identity, we have
\begin{align*}
G^{ij(x)}(z) -G(z) = G^{ij(x)}(z) P^{ij} G(z)
\end{align*}
Thus, by \eqref{eq:bdGijx} and \eqref{eq:crudebdG}, 
\begin{align}
\|G^{ij(x)}(z) -G(z) \|= \|G^{ij(x)}(z)\| \cdot \|P^{ij}\| \cdot\| G(z)\| \le 18 \frac{\mathtt{L}_n}{|z|^2}.
\end{align}
Then 
\begin{align*}
\|L^{ij(x)} - L\| = \max_s |l_s^{ij(x)} -l_s |&= \max_s \sum_{k=1}^n \sigma_{sk}^2 |G_{kk}^{ij(x)}(z) - G_{kk}(z)|\\
&\le R_{\max} \|G^{ij(x)}(z) -G(z) \| \le 18 \mathtt{L}_n \frac{R_{\max}}{|z|^2} \le \frac{1}{2} \mathtt{L}_n
\end{align*}
by the assumption $|z|\ge 6\sqrt{R_{\max}}$.
Similarly, 
\begin{align*}
\|L\| =\max_s |l_s | \le \max_s \sum_{k=1}^n \sigma_{sk}^2 |G_{kk}(z)|\le R_{\max} \|G(z)\| \le 2\frac{R_{\max}}{|z|}.
\end{align*}
It follows that
\[ \|E\|+\|L\| \le 3.5\sqrt{R_{\max}} + 2\frac{R_{\max}}{|z|} \le \left(\frac{7}{12} + \frac{1}{18} \right)|z|.\]
Combining the above estimates, we conclude that 
\begin{align*}
|\mathcal D^{ij(x)}_{vw} -\mathcal D_{vw}|\le 10.5\frac{\mathtt{L}_n}{|z|} + 18\left(\frac{7}{12} + \frac{1}{18}\right)\frac{\mathtt{L}_n}{|z|} = 22\frac{\mathtt{L}_n}{|z|}.
\end{align*}
Thus
\begin{align}\label{eq:bdDvwij}
\sup_{|x|\le \mathtt{L}_n}|\mathcal D^{ij}_{vw}(x) | \le |\mathcal D_{vw}| + 22\frac{\mathtt{L}_n}{|z|}.
\end{align}

Now, we proceed to the proof of \eqref{eq:0520-bd-3rd}. By the Leibniz rule, 
\[\partial_{ij}^3 f_{ij}(E^{ij(x)}) = \sum_{s=0}^3 \binom{3}{s} \left(\partial_{ij}^{3-s}G_{jw}^{ij(x)}\right) \left(\partial_{ij}^s (\mathcal D^{ij(x)}_{vw})^{2p-1}\right).\]
We estimate the four cases $s=0,1,2,3$ separately. The cases $s=0,1,2$ follow similar arguments to those used for $X_2$ in \eqref{def:X2-bd}.

First, for $s=0$, using \[ \partial_{ij}^3 G^{ij(x)}= 3!G^{ij(x)} (\Delta^{ij}G^{ij(x)})^3 \] and the same estimate used for the $s=0$ term in the $X_2$ bound, we get
\begin{align*}
\sum_{i,j}|v_i|\sigma_{ij}^2 |\partial_{ij}^3G_{jw}^{ij(x)}| \le 6R_{\max} \|G^{ij(x)}\|^4 \le 6R_{\max} \frac{3^4}{|z|^4}\le \frac{27}{2|z|^2}.
\end{align*}
Therefore,
\[
\sum_{i,j}|v_i|\sigma_{ij}^2 |\partial_{ij}^3G_{jw}^{ij(x)}| \cdot |(\mathcal D^{ij(x)}_{vw})^{2p-1}| \le  \frac{27}{2|z|^2} | \mathcal D^{ij(x)}_{vw}|^{2p-1}.
\]

Next, for $s=1$, since
\[ |\partial_{ij}(\mathcal D^{ij(x)}_{vw})^{2p-1} | \le 2p |\mathcal D^{ij(x)}_{vw}|^{2p-2} |\partial_{ij}\mathcal D^{ij(x)}_{vw}|,\]
we aim to bound 
\[ \sum_{i,j}|v_i|\sigma_{ij}^2
|\partial_{ij}^2G_{jw}^{ij}(x)|
|\partial_{ij}\mathcal D_{vw}^{ij(x)}|.
\]
This is proved exactly as the $s=1$ estimate for $X_2$, but with an extra factor of $G$. We sketch the key steps here. Expand
\[
\partial_{ij}\mathcal D_{vw}^{ij(x)} = (1+\delta_{ij})^{-1}
\left( G_{\widetilde v_x i}^{ij(x)}G_{jw}^{ij(x)} +
G_{\widetilde v_x j}^{ij(x)}G_{iw}^{ij(x)} - \mathcal K_x \right),
\]
where $\widetilde v_x^\T:=v^\T(zI-L^{ij(x)})$ and \[
\mathcal K_x := 2\sum_{s,k}v_s\sigma_{sk}^2 G_{ki}^{ij(x)}G_{kj}^{ij(x)}G_{sw}^{ij(x)}.
\]
Since $$\|L^{ij(x)}\| \le \max_s \sum_{k=1}^n \sigma_{sk}^2 |G^{ij(x)}(z)|\le R_{\max} \|G^{ij(x)}(z)\| \le  3 \frac{R_{\max}}{|z|}\le \frac{1}{12}|z|,$$
we have $\|\widetilde v_x^\T\| \le \frac{13}{12}|z|.$ Similar to \eqref{eq:bdRvw}, we obtain
\[ |\mathcal K_x| \le 2 \sigma_{\max}^2 \|G^{ij(x)}\|^3.\]
For the second derivative of the resolvent
\begin{align*}
\partial_{ij}^2 G_{jw}^{ij(x)} = 2(1+\delta_{ij})^{-2} \Big[ \big((G^{ij(x)}_{ij})^2 + G^{ij(x)}_{ii}G^{ij(x)}_{jj} \big)G^{ij(x)}_{jw} + 2G^{ij(x)}_{ij}G^{ij(x)}_{jj} G^{ij(x)}_{iw} \Big],
\end{align*}
we apply the crude upper bound
\begin{align*}
|\partial_{ij}^2 G_{jw}^{ij(x)}| \le 4\|G^{ij(x)}\|^2 \big( |G_{jw}^{ij(x)}| + |G_{iw}^{ij(x)}| \big).
\end{align*}
Multiplying this by the expansion of $\partial_{ij}\mathcal D_{vw}^{ij(x)}$ and applying the identical Cauchy-Schwarz estimates used for $X_2$, we obtain 
\begin{align*}
\sum_{i,j}|v_i|\sigma_{ij}^2 |\partial_{ij}^2G_{jw}^{ij(x)}| |\partial_{ij}\mathcal D_{vw}^{ij(x)}|
&\le 16\sigma_{\max}\sqrt{R_{\max}} \|G^{ij(x)}\|^5 \|\widetilde v_x\| + 16\sigma_{\max}^2 R_{\max}\|G^{ij(x)}\|^6 \\
&\le 16\sigma_{\max} \left(\frac{|z|}{6}\right) \frac{3^5}{|z|^5} \left(\frac{13|z|}{12}\right) + 16\sigma_{\max} \left(\frac{|z|^3}{216}\right) \frac{3^6}{|z|^6} \\
&\le 702 \frac{\sigma_{\max}}{|z|^3} + 54 \frac{\sigma_{\max}}{|z|^3} = 756 \frac{\sigma_{\max}}{|z|^3}.
\end{align*}
Hence the $s=1$ contribution is bounded by
\begin{align*}
1512 p \frac{\sigma_{\max}}{|z|^3}| \mathcal D^{ij(x)}_{vw}|^{2p-2}.
\end{align*}

For $s=2$, we use
\[
|\partial_{ij}^2\left((\mathcal D^{ij(x)}_{vw})^{2p-1}\right)|
\le  2 p |\mathcal D^{ij(x)}_{vw}|^{2p-2} |\partial_{ij}^2\mathcal D^{ij(x)}_{vw}| + 4 p^2 |\mathcal D^{ij(x)}_{vw}|^{2p-3} |\partial_{ij}\mathcal D^{ij(x)}_{vw}|^2 .
\]
The first estimate needed is 
\[
\sum_{i,j}|v_i|\sigma_{ij}^2 |\partial_{ij}G_{jw}^{ij}(x)|
|\partial_{ij}^2\mathcal D_{vw}^{ij}(x)|. 
\]
Using $|\partial_{ij} G_{jw}^{ij(x)}| \le 2\|G^{ij(x)}\| (|G_{jw}^{ij(x)}| + |G_{iw}^{ij(x)}|)$ and multiplying by the bounds established in \eqref{eq:0520-previous1}, we pull out an extra factor of $2\|G^{ij(x)}\|$ and double the sum to account for $|G_{iw}|$:
\begin{align*}
\sum_{i,j}|v_i|\sigma_{ij}^2 |\partial_{ij}G_{jw}^{ij(x)}| |\partial_{ij}^2\mathcal D_{vw}^{ij(x)}|
&\le 32\sigma_{\max} \sqrt{R_{\max}} \|G^{ij(x)}\|^5 \|\widetilde{v}_x\| + 48 \sigma_{\max}^2 R_{\max} \|G^{ij(x)}\|^6 \\
&\le 1404 \frac{\sigma_{\max}}{|z|^3} + 162 \frac{\sigma_{\max}}{|z|^3} 
= 1566 \frac{\sigma_{\max}}{|z|^3}.
\end{align*}
We also require the estimate for 
$$ \sum_{i,j}|v_i|\sigma_{ij}^2 |\partial_{ij}G_{jw}^{ij}(x)| |\partial_{ij}\mathcal D_{vw}^{ij(x)}|^2. $$
Again, using $|\partial_{ij} G_{jw}^{ij(x)}|\le 2\|G^{ij(x)}\| (|G_{jw}^{ij(x)}| + |G_{iw}^{ij(x)}|)$ to the estimates from \eqref{eq:0520-previous2} yields:
\begin{align*}
\sum_{i,j}|v_i|\sigma_{ij}^2 |\partial_{ij}G_{jw}^{ij(x)}| |\partial_{ij}\mathcal D_{vw}^{ij(x)}|^2
&\le 24 \sigma_{\max}^2 \|G^{ij(x)}\|^6 \|\widetilde{v}_x\|^2 + 16 \sigma_{\max}^4 R_{\max} \|G^{ij(x)}\|^8 \\
&\le 24 \sigma_{\max}^2 \frac{3^6}{|z|^6} \left(\frac{169|z|^2}{144}\right) + 16 \sigma_{\max}^2 \left(\frac{|z|^4}{1296}\right) \frac{3^8}{|z|^8} \\
&\le C \frac{\sigma_{\max}^2}{|z|^4}.
\end{align*}

Therefore, the total $s=2$ contribution is bounded by
\[
C p \frac{\sigma_{\max}}{|z|^3}|\mathcal D_{vw}^{ij(x)}|^{2p-2}
+
C p^2\frac{\sigma_{\max}^2}{|z|^4}|\mathcal D_{vw}^{ij(x)}|^{2p-3}.
\]

\medskip
Finally, for $s=3$, applying the generalized chain rule (Faà di Bruno's formula \cite[p.~137]{CL12}) to the third derivative of $(\mathcal{D}_{vw}^{ij(x)})^{2p-1}$ yields:
\begin{align*}
\Big|\partial_{ij}^3\Big((\mathcal D_{vw}^{ij(x)})^{2p-1}\Big) \Big|
&\le 2p |\mathcal D_{vw}^{ij(x)}|^{2p-2}|\partial_{ij}^3\mathcal D_{vw}^{ij(x)}| \\
&\quad+ 12p^2 |\mathcal D_{vw}^{ij(x)}|^{2p-3}|\partial_{ij}\mathcal D_{vw}^{ij(x)}||\partial_{ij}^2\mathcal D_{vw}^{ij(x)}|  \\
&\quad+  8 p^3 |\mathcal D_{vw}^{ij(x)}|^{2p-4}|\partial_{ij}\mathcal D_{vw}^{ij(x)}|^3.
\end{align*}
Bounding this requires three separate estimates:
\begin{align}
\sum_{i,j}|v_i| \sigma_{ij}^2 |G_{jw}^{ij(x)}| |\partial_{ij}^3\mathcal D_{vw}^{ij(x)}| &\le C_1 \frac{\sigma_{\max}}{|z|^3},\label{eq:D3-contraction}\\
\sum_{i,j}|v_i| \sigma_{ij}^2 |G_{jw}^{ij(x)}| |\partial_{ij}\mathcal D_{vw}^{ij(x)}| |\partial_{ij}^2\mathcal D_{vw}^{ij(x)}| &\le C_2 \frac{\sigma_{\max}^2}{|z|^4},\label{eq:D1D2-contraction}\\
\sum_{i,j}|v_i| \sigma_{ij}^2 |G_{jw}^{ij(x)}| |\partial_{ij}\mathcal D_{vw}^{ij(x)}|^3 &\le C_3 \frac{\sigma_{\max}^3}{|z|^5},\label{eq:D1cubed-contraction}
\end{align}
where $C_1, C_2, C_3 > 0$ are absolute constants (which are all bounded $10^5$). 

We briefly sketch the proof of these estimates. For \eqref{eq:D3-contraction}, from the identity $\mathcal D=(zI-L)G-I$, we have$$ \partial_{ij}^3\mathcal D = (zI-L)\partial_{ij}^3G - 3(\partial_{ij}L)\partial_{ij}^2G - 3(\partial_{ij}^2L)\partial_{ij}G - (\partial_{ij}^3L)G. $$ Then we multiply this expansion with $|v_i| \sigma_{ij}^2 |G_{jw}^{ij(x)}|$ and sum over $i,j$. The first term is bounded by expanding $\partial_{ij}^3G=3!G(\Delta^{ij}G)^3$. The factor $|G_{jw}^{ij(x)}|$ combined with the variance weight generates one factor of $\sigma_{\max}$, while the remaining Green functions are controlled by \eqref{eq:ward} and \eqref{eq:wardCS}. This yields the $C_1 \sigma_{\max}|z|^{-3}$ bound.
The terms containing $\partial_{ij}^a L$ are handled by expanding
$$ \partial_{ij}^a l_s = \sum_k\sigma_{sk}^2\partial_{ij}^a G_{kk}^{ij(x)}, $$
together with
$$ |\partial_{ij}^aG_{kk}^{ij(x)}| \le \widetilde{C}_a \|G^{ij(x)}\|^{a-1} \sum_{s\in\{i,j\}} |G_{ks}^{ij(x)}|^2. $$
and the Cauchy-Schwarz inequality. This proves \eqref{eq:D3-contraction}.  

The estimates \eqref{eq:D1D2-contraction} and \eqref{eq:D1cubed-contraction} follow by
combining the already proved bounds for $\partial_{ij}\mathcal D_{vw}^{ij(x)}$ and $\partial_{ij}^2\mathcal D_{vw}^{ij(x)}$, and the factor $|G_{jw}^{ij(x)}|$. Note that each additional $\partial_{ij}\mathcal D$ contributes one factor
$\sigma_{\max}/|z|$ in the estimate, exactly as in the $s=2$ case of the $X_2$ estimate.

Consequently, the \(s=3\) contribution is bounded by
$$ 2p\,C_1\frac{\sigma_{\max}}{|z|^3}|\mathcal D_{vw}^{ij(x)}|^{2p-2} + 12p^2\,C_2\frac{\sigma_{\max}^2}{|z|^4}|\mathcal D_{vw}^{ij(x)}|^{2p-3} + 8p^3\,C_3\frac{\sigma_{\max}^3}{|z|^5}|\mathcal D_{vw}^{ij(x)}|^{2p-4}. $$

Combining the four cases $s=0,1,2,3$ yields
\begin{align}\label{eq:0522-3rd-order}
\sum_{i,j}|v_i|\sigma_{ij}^2 \sup_{|x|\le \mathtt L_n} \left| \partial_{ij}^3 \left( G_{jw}^{ij(x)}(\mathcal D_{vw}^{ij(x)})^{2p-1} \right) \right| \le C \sum_{m=1}^4 p^{m-1}\frac{\sigma_{\max}^{m-1}}{|z|^{m+1}} |\mathcal D_{vw}^{ij(x)}|^{2p-m}. 
\end{align}
Using $|\mathcal D_{vw}^{ij(x)}| \le |\mathcal D_{vw}| + 22 \mathtt{L}_n/|z|$ from \eqref{eq:bdDvwij}, we finally obtain 
\begin{align*}
\sum_{i,j}|v_i|\sigma_{ij}^2 \sup_{|x|\le \mathtt{L}_n} \left| \partial_{ij}^3 \left( G_{jw}^{ij(x)}(\mathcal D_{vw}^{ij(x)})^{2p-1} \right) \right| \le C \sum_{m=1}^4 p^{m-1}\frac{\sigma_{\max}^{m-1}}{|z|^{m+1}} \left( |\mathcal D_{vw}|+ 22\frac{\mathtt{L}_n}{|z|} \right)^{2p-m}. 
\end{align*}
This proves \eqref{eq:0520-bd-3rd}.

\paragraph{Proof of Lemma \ref{Key_Lemma} (iii).} We estimate the remainder term in the cumulant expansion with $l=2$. Recall from \eqref{def:R3} that 
\[ \mathtt{R}_{3}
= \sum_{i,j} v_i \mathtt{R}_3^{ij}.\]
By Lemma \ref{lem:cumulant-expansion}, 
\begin{align*}
|\mathtt{R}_{3}|&\le \sum_{i,j} | v_i \mathtt{R}_{3}^{ij}|\\
&\le  16 \sum_{i,j} | v_i| \E\left[|E_{ij}|^{4}\right] \E\left[\sup_{|x|\le t}|\partial_{ij}^{3}f_{ij}(E^{ij}+x\Delta^{ij})|\right] \quad &(\text{Term 1}) \\
&+  16\sum_{i,j} | v_i| \left( \E\left[\sup_{|x|\le|E_{ij}|} |\partial_{ij}^{3}f_{ij}(E^{ij(x)})|^2\right] \right)^{1/2} \left( \E\left[|E_{ij}|^{8} \mathbf{1}_{|E_{ij}|>t}\right] \right)^{1/2}\quad &(\text{Term 2})
\end{align*}
We choose $t=\mathtt{L}_n$. In the general sub-Gaussian regime, $\mathtt{L}_n=5\sqrt{D+6} K\sigma_{\max}\log n,$ , while in the bounded sparse-entry regime, it is equal to zero. Here, as before, the cutoff $\chi$ is suppressed from the notation. Also, as explained in the beginning of Section \ref{sec:proof-moments}, terms with derivatives applied to $\chi$ are negligible. 

Let us consider the general sub-Gaussian regime. We first estimate Term 1. Set 
\[\theta:=\beta K\sigma_{\max}.\]
By Assumption~\ref{assump:noise},
\[ \E |E_{ij}|^4 \le \beta^2K^2\sigma_{\max}^2\sigma_{ij}^2 = \theta^2\sigma_{ij}^2.\]
Hence, by \eqref{eq:0520-bd-3rd}, the contribution from Term 1 is bounded by
\begin{align*}
&16  \sum_{i,j} | v_i| \E\left[|E_{ij}|^{4}\right] \E\left[\sup_{|x|\le \mathtt{L}_n}|\partial_{ij}^{3}f_{ij}(E^{ij}+x\Delta^{ij})|\right] \\
&\qquad \le 16 \theta^2 \sum_{i,j}|v_i|\sigma_{ij}^2 \E \left[ \sup_{|x|\le \mathtt L_n} \left| \partial_{ij}^3 f_{ij}(E^{ij(x)})\right| \right]\\
& \qquad \le C\theta^2 \sum_{m=1}^4 p^{m-1}\frac{\sigma_{\max}^{m-1}}{|z|^{m+1}}
\E\left[ \left( |\mathcal D_{vw}|+ 22\frac{\mathtt{L}_n}{|z|} \right)^{2p-m} \right].
\end{align*}

It remains to estimate Term 2. By the sub-Gaussian tail assumption, we have 
\begin{equation}\label{eq:eighth-tail-bound}
\left(\E\left[ |E_{ij}|^8\mathbf 1_{\{|E_{ij}|>\mathtt{L}_n\}} \right] \right)^{1/2}\le 
\theta^2 \sigma_{ij}^2 \exp(-10(D+6)\log^2 n).
\end{equation}
To see this, note that the case $\sigma_{ij}=0$ is trivial. For $\sigma_{ij}>0$, write $E_{ij}=\sigma_{ij}\xi_{ij}$ with $\|\xi_{ij}\|_{\psi_2}\le K$. A standard sub-Gaussian tail integration yields 
\begin{align}\label{eq:0526-revise}
\E\left[
|E_{ij}|^8\mathbf 1_{\{|E_{ij}|>\mathtt{L}_n\}} \right] &\le C(K\sigma_{ij})^8
\left(1+\frac{\mathtt{L}_n}{K\sigma_{ij}}\right)^8 \exp\left( -\frac{\mathtt{L}_n^2}{K^2\sigma_{ij}^2} \right)\nonumber\\
&= C \left( K\sigma_{ij} + \mathtt{L}_n \right)^8 \exp\left( -\frac{\mathtt{L}_n^2}{K^2\sigma_{ij}^2} \right).
\end{align}
Because $\sigma_{ij} \le \sigma_{\max}$ and $\mathtt{L}_n=5\sqrt{D+6} K\sigma_{\max}\log n$, we bound \[\left( K\sigma_{ij} + \mathtt{L}_n \right)^8 \le C D^4 K^8 \sigma_{\max}^8 (\log n)^8.\] Moreover, we write 
\begin{align*}
\exp\left( -\frac{\mathtt L_n^2}{K^2\sigma_{ij}^2} \right)&= \frac{\sigma_{ij}^4}{\sigma_{\max}^4} \left[ \frac{\sigma_{\max}^4}{\sigma_{ij}^4} \exp\left( -25(D+6)\frac{\sigma_{\max}^2}{\sigma_{ij}^2}\log^2 n \right) \right]\\
&\le \frac{\sigma_{ij}^4}{\sigma_{\max}^4} \exp\big(-25 (D+6) \log^2 n\big),
\end{align*}
where we used the fact that the function $y \mapsto y^2 \exp(-cy)$ is strictly decreasing for $y \ge 1$ and large $c>0$, and then $y^2 \exp(-cy)\le \exp(-c)$.

Thus, the right-hand side of \eqref{eq:0526-revise} is bounded by $C D^4 K^8 \sigma_{\max}^4 \sigma_{ij}^4 (\log n)^8 \exp(-25 (D+6) \log^2 n).$
To get the desired bound \eqref{eq:eighth-tail-bound}, we only need this quantity 
\[
    D^4 K^8 \sigma_{\max}^4 \sigma_{ij}^4 (\log n)^8 \exp(-25 (D+6) \log^2 n) \le  \theta^4 \sigma_{ij}^4 \exp(-20(D+6)\log^2 n),
\]
where $\theta=\beta  K \sigma_{\max}$. Since $\beta \ge 1$, we require
\begin{align}\label{eq:0521-smallbd}
\exp(25 (D+6) \log^2 n) \ge D^4 K^4 \exp(20(D+6)\log^2 n).
\end{align}
We now use Assumption~\ref{assump:extra}. Since $R_{\max}\le n\sigma_{\max}^2$ and $ \sqrt{R_{\max}}\ge C_0 (D+8) K\sigma_{\max}\log n,$
we have
\[ K\le \frac{\sqrt n}{C_0 (D+8)\log n}.\]
Hence, \eqref{eq:0521-smallbd} holds and  \eqref{eq:eighth-tail-bound} is proved. 

Now, Term 2 is bounded by 
\begin{align*}
& 16\sum_{i,j} | v_i| \left( \E\left[\sup_{|x|\le|E_{ij}|} |\partial_{ij}^{3}f_{ij}(E^{ij(x)})|^2\right] \right)^{1/2} \left( \E\left[|E_{ij}|^{8} \mathbf{1}_{|E_{ij}|>t}\right] \right)^{1/2}\\
&\quad \le C \exp(-10(D+6)\log^2 n) \theta^2 \sum_{i,j}  \left( \E\left[\sup_{|x|\le|E_{ij}|} \left| |v_i| \sigma_{ij}^2 \partial_{ij}^{3}f_{ij}(E^{ij(x)}) \right|^2\right] \right)^{1/2}.
\end{align*}
Recall that the cutoff $\chi$ is suppressed from the notation. On the support $\Omega_n$ of $\chi$, \eqref{eq:0522-3rd-order} holds. Moreover, by \eqref{eq:0520-Eijx} and \eqref{eq:bdGijx}, we also have this crude bound
\[
\|\mathcal D^{ij(x)}\| =\|(E^{ij(x)}-L^{ij(x)})G^{ij(x)}\| \le C.
\]
Therefore,
\[
|\mathcal D_{vw}^{ij(x)}|^{2p-m}\le C^p, \qquad 1\le m\le 4.
\]
It follows from \eqref{eq:0522-3rd-order} that
\begin{align*}
\left( \E\left[\sup_{|x|\le|E_{ij}|} \left| |v_i| \sigma_{ij}^2 \partial_{ij}^{3}f_{ij}(E^{ij(x)}) \right|^2\right] \right)^{1/2} \le C^p \sum_{m=1}^4 p^{m-1}\frac{\sigma_{\max}^{m-1}}{|z|^{m+1}}.
\end{align*}
Since $\sigma_{\max}\le \sqrt{R_{\max}} \le |z|/6$, we simplify this bound to 
\begin{align*}
\left( \E\left[\sup_{|x|\le|E_{ij}|} \left| |v_i| \sigma_{ij}^2 \partial_{ij}^{3}f_{ij}(E^{ij(x)}) \right|^2\right] \right)^{1/2} \le C^p p^3\frac{1}{|z|^2}.
\end{align*}
Combining the above bound and \eqref{eq:eighth-tail-bound}, the total contribution of Term 2 is bounded by
\begin{align*}
& 16\sum_{i,j} | v_i| \left( \E\left[\sup_{|x|\le|E_{ij}|} |\partial_{ij}^{3}f_{ij}(E^{ij(x)})|^2\right] \right)^{1/2} \left( \E\left[|E_{ij}|^{8} \mathbf{1}_{|E_{ij}|>t}\right] \right)^{1/2}\\
&\quad \le C^p p^3 \exp\big(-10(D+6)\log^2 n \big) \frac{\theta^2}{|z|^2}.
\end{align*}
Finally, we arrive at
\begin{align*}
|\mathtt R_3| &\le C\theta^2 \sum_{m=1}^4 p^{m-1}\frac{\sigma_{\max}^{m-1}}{|z|^{m+1}}
\E\left[ \left( |\mathcal D_{vw}|+ 22\frac{\mathtt{L}_n}{|z|} \right)^{2p-m} \right] + C^p p^3 \exp\big(-10(D+6)\log^2 n \big) \frac{\theta^2}{|z|^2}\\
&= C(K\beta)^2 \sum_{m=1}^4 p^{m-1}\frac{\sigma_{\max}^{m+1}}{|z|^{m+1}}
\E\left[ \left( |\mathcal D_{vw}|+ 22\frac{\mathtt{L}_n}{|z|} \right)^{2p-m} \right] + C^p p^3 \exp\big(-10(D+6)\log^2 n \big) \frac{\theta^2}{|z|^2}.
\end{align*}
The bounded sparse-entry regime follows the same approach. Indeed, the Term 2 in the bound for $|\mathtt R_3|$ is equal to zero because
$\max_{i,j}|E_{ij}|\le \mathtt L_n$ almost surely. This simplifies the derivation. We skip the details. The proof of Lemma \ref{Key_Lemma} is now complete.

\section{Proof of Proposition \ref{prop:local-part2}}\label{app:local-part2}
Recall that
\begin{align*}
  &  H=\diag(h_1,\ldots,h_n), \qquad h_i=\sum_{j=1}^n\sigma_{ij}^2\phi_j;\\
&   L=\operatorname{diag}(l_1,\ldots,l_n), \qquad  l_i=\sum_{j=1}^n\sigma_{ij}^2G_{jj}.
\end{align*}
Set
\[
    \mathcal D:=(E-L)G, \qquad  \mathcal D_{ii}:=e_i^\T \mathcal D e_i.
\]
Note that $\max_i|\phi_i(z)|\le \frac{3}{2|z|}$ by Lemma~\ref{lem:bdphi}. Define
\[ d_i:=\phi_i-G_{ii}, \qquad d:=(d_1,\ldots,d_n)^\T \]
We first bound $\|d\|_\infty$. From $ (zI-E)G=I,$
we have $EG=zG-I.$
Therefore, taking the $(i,i)$-entry of  $ \mathcal D=(E-L)G $ yields
\[
    \mathcal D_{ii} =  (EG)_{ii}-l_iG_{ii}  = (z-l_i)G_{ii}-1.
\]
Equivalently,
\[
    (z-l_i)G_{ii}=1+\mathcal D_{ii}.
\]
On the other hand, we have from the QVE that $(z-h_i)\phi_i=1.$
Thus
\[
\begin{aligned}
    (z-h_i)(\phi_i-G_{ii})
    &=
    (z-h_i)\phi_i-(z-h_i)G_{ii}  \\
    &=
    1-(z-h_i)G_{ii} \\
    &=
    (z-l_i)G_{ii}-\mathcal D_{ii}-(z-h_i)G_{ii} \\
    &=
    (h_i-l_i)G_{ii}-\mathcal D_{ii}.
\end{aligned}
\]
Multiplying by \(\phi_i\), since \(\phi_i=(z-h_i)^{-1}\), yields
\begin{equation}\label{eq:diagonal-stability-equation-LH}
    d_i = \phi_iG_{ii}(h_i-l_i)-\phi_i\mathcal D_{ii}.
\end{equation}
Moreover,
\[
    h_i-l_i =\sum_{j=1}^n\sigma_{ij}^2(\phi_j-G_{jj}) = \sum_{j=1}^n\sigma_{ij}^2d_j.
\]
Substituting this into \eqref{eq:diagonal-stability-equation-LH}, we obtain
\[
    d_i = \phi_iG_{ii} \Big(\sum_{j=1}^n\sigma_{ij}^2 \Big) d_j - \phi_i\mathcal D_{ii}.
\]

Thus, by taking the $\ell_\infty$ norm, we get
\[
    \|d\|_\infty \le \max_i|\phi_iG_{ii}| \cdot R_{\max}\|d\|_\infty + \max_i|\phi_i|\cdot\max_i|\mathcal D_{ii}|.
\]
Since
\[
    \max_i|\phi_iG_{ii}| \le \frac{3}{2|z|}\cdot\frac{2}{|z|} =  \frac{3}{|z|^2},
\]
we bound
\[
    \|d\|_\infty  \le \frac{3R_{\max}}{|z|^2}\|d\|_\infty +  \frac{3}{2|z|}\max_i|\mathcal D_{ii}|.
\]
Since $|z|^2\ge36R_{\max}$ by our assumption, we further obtain
\begin{equation}\label{eq:d-infty-bound-LH}
    \|d\|_\infty \le \frac{2}{|z|}\max_i|\mathcal D_{ii}|.
\end{equation}

We now return to $v^\T (L-H)G w$. Since $v,w$ are unit vectors,
\[
    |v^\T (L-H)Gw| \le \|L-H\|\cdot\|G\|,
\]
where
\begin{align*}
    \|L-H\| =  \max_i \Big| \sum_{j=1}^n\sigma_{ij}^2(G_{jj}-\phi_j) \Big|   \le
    R_{\max}\|d\|_\infty.
\end{align*}
Combining this with \(\|G\|\le2/|z|\) and
\eqref{eq:d-infty-bound-LH}, we arrive at
\begin{equation}\label{eq:LH-reduced-to-Dii}
    |v^\T (L-H)Gw| \le \frac{4R_{\max}}{|z|^2} \max_i |\mathcal D_{ii}| \le \frac{1}{9}\max_i|\mathcal D_{ii}|.
\end{equation}

To bound $\mathcal D_{ii}$,  we apply
Proposition~\ref{prop:bdDvw} with $v=w=e_i$ and take a union bound to obtain 
\begin{align*}
    \max_{1\le i\le n}|\mathcal D_{ii}| \le C\left(  \frac{\mathfrak m_D}{|z|} +n^{-500}
    \right)
\end{align*}
with probability at least $1-n^{-D-4}$. In the bounded sparse-entry regime,
the analogous bound follows without the $n^{-500}$ term.

Combining the above estimate with \eqref{eq:LH-reduced-to-Dii} proves Proposition~\ref{prop:local-part2}.


\section{Proofs of Corollary~\ref{lem:local-law-uniform} and Corollary~\ref{coro:row-adaptive}}\label{app:local-law-uniform}
\subsection{Proof of Corollary~\ref{lem:local-law-uniform}}
Write 
\[\Xi(z)= G(z)-\Phi(z).\]
We start with the following consequence of the isotropic local law:
\begin{coro}[Compressed local law on the signal subspace]\label{coro:0603-locallaw-U} Fix $D>0$. Suppose the assumptions of Theorem~\ref{thm:local} hold. Then, for every fixed $z\in\C_+$ with $ |z|\ge 6\sqrt{R_{\max}},$
with probability at least $1-n^{-D-2}$,
    \[\left\|U^\T (G(z)-\Phi(z)) U \right\| \le C r\frac{\mathfrak m_{D}}{|z|^2}.\]
\end{coro}
\begin{proof}[Proof of Corollary \ref{coro:0603-locallaw-U}]
   Since 
\begin{align*}
   \left\|U^\T (G(z)-\Phi(z)) U \right\|  \le \left\|U^\T (G(z)-\Phi(z)) U \right\|_F&=\sqrt{\sum_{i,j=1}^r \big|{e_i}^\T U^\T (G(z)-\Phi(z)) U e_j\big|^2}\\
   &\le r\max_{1\le i,j \le r} \big|{e_i}^\T U^\T (G(z)-\Phi(z)) U e_j\big|,
\end{align*}
the conclusion of Corollary \ref{coro:0603-locallaw-U} follows directly from Theorem \ref{thm:local} and the union bounds.
\end{proof}

To extend the pointwise bound Corollary \ref{coro:0603-locallaw-U} to the spectral domain 
\[\mathcal S_{\rm out}:=\{ z\in\mathbb C: 6\sqrt{R_{\max}}\le |z|\le 2R_{\max}^3 \},\]
we use a standard $\varepsilon$-net argument. Set 
\[  \eta := \min\left\{1, \frac{r\mathfrak m_D}{20}\right\}.\]
Let $\mathcal N_{\eta}$ be a $\eta$-net of $\mathcal S_{\rm out}$. A simple volume argument (see for instance \cite[Lemma 3.3]{OW}) shows that 
\[
    |\mathcal N_{\eta}| \le \left(1+\frac{8R_{\max}^3}{\eta}\right)^2.
\]
Since $R_{\max} \le n^{100}$ and $\sigma_{\max}^{-1} \le n^{100}$, we have $\mathfrak m_D\ge c n^{-100}$ for an absolute constant $c>0$. Thus $R_{\max}^3 / \eta \le C' n^{400}$ and $|\mathcal N_{\eta}|  \le C n^{800}$.

Applying Corollary \ref{coro:0603-locallaw-U} with $D_1=D+810$ and taking a union bound over $z\in\mathcal N_{\eta}$, we obtain 
\[
  \Prob\left(  \max_{z\in\mathcal N_{\eta}} |z|^2\left\|U^\T\Xi(z)U\right\| \ge C r \mathfrak m_{D_1} \right)
    \le |\mathcal N_{\eta}|  n^{-D-814} < n^{-D-10}.
\]
Using $\mathfrak m_{D_1}= \mathfrak m_{D+810} \le C' \mathfrak m_{D}$, we establish that with probability at least $1-n^{-D-10}$,
\[ \max_{z\in\mathcal N_{\eta}} |z|^2\left\|U^\T\Xi(z)U\right\| \le C r \mathfrak m_{D}.\]

Now we extend the estimate from $\mathcal N_{\eta}$ to the whole domain $\mathcal S_{\rm out}$. Define
\[
    f(z):=z^2U^\T G(z)U,
    \qquad
    g(z):=z^2U^\T\Phi(z)U.
\]
We first show that $f$ and $g$ are Lipschitz on $\mathcal S$.

Let $z,w\in\mathcal S_{\rm out}$ and assume without loss of generality that $|z|\ge |w|\ge 6\sqrt{R_{\max}}$. By the resolvent identity,
\[
    G(z)-G(w) = (w-z)G(z)G(w).
\]
Using $\|G(z)\|\le 2/|z|$ and $\|G(w)\|\le 2/|w|$, we get
\begin{align*}
    \|f(z)-f(w)\|
    &\le |z||z-w|\|G(z)\| + |w||z-w|\|G(z)\| \\
    &\quad + |w|^2\|G(z)-G(w)\| \\
    &\le 2|z-w| + 2\frac{|w|}{|z|}|z-w| + 4\frac{|w|}{|z|}|z-w| \\
    &\le 8|z-w|.
\end{align*}
Hence $f$ is $8$-Lipschitz on $\mathcal S$.

Next we prove the Lipschitz bound for $g$. For each $i\in[n]$, the QVE gives
\[
    \frac{1}{\phi_i(z)}
    = z-\sum_j\sigma_{ij}^2\phi_j(z),
    \qquad
    \frac{1}{\phi_i(w)}
    = w-\sum_j\sigma_{ij}^2\phi_j(w).
\]
Subtracting the two identities and multiplying by $\phi_i(z)\phi_i(w)$, we obtain
\[
    \phi_i(z)-\phi_i(w)
    = \phi_i(z)\phi_i(w) \left[ w-z+\sum_j\sigma_{ij}^2(\phi_j(z)-\phi_j(w)) \right].
\]
Therefore, by Lemma \ref{lem:bdphi},
\begin{align*}
    \|\Phi(z)-\Phi(w)\|
    \le \frac{9}{4|z||w|} \left( |z-w|+R_{\max}\|\Phi(z)-\Phi(w)\| \right).
\end{align*}
Since $|z||w|\ge 36R_{\max}$, we have
\[
    \frac{9R_{\max}}{4|z||w|}
    \le \frac{1}{16}.
\]
Thus,
\[
    \|\Phi(z)-\Phi(w)\|
    \le \frac{12}{5}\frac{|z-w|}{|z||w|}.
\]
Using again Lemma \ref{lem:bdphi}, we have $\|\Phi(z)\|\le 3/(2|z|)$. Hence
\begin{align*}
    \|g(z)-g(w)\|
    &\le |z||z-w|\|\Phi(z)\| +|w||z-w|\|\Phi(z)\| \\
    &\quad +|w|^2\|\Phi(z)-\Phi(w)\| \\
    &\le \frac{3}{2}|z-w| + \frac{3}{2}\frac{|w|}{|z|}|z-w| + \frac{12}{5}\frac{|w|}{|z|}|z-w| \\
    &< 6|z-w|.
\end{align*}
Therefore $g$ is $6$-Lipschitz on $\mathcal S$.

Now take any $z\in\mathcal S_{\rm out}$. By the definition of the $\eta$-net, there exists $w\in\mathcal N_{\eta}$ such that $|z-w|\le \eta$. Then
\begin{align*}
    |z|^2\left\|U^\T\Xi(z)U\right\|
    &= \|f(z)-g(z)\| \\
    &\le \|f(z)-f(w)\|+\|f(w)-g(w)\|+\|g(w)-g(z)\| \\
    &\le 14\eta + |w|^2\left\|U^\T\Xi(w)U\right\|.
\end{align*}
Since $\eta \le \frac{r\mathfrak m_D}{20}$, we get 
\[
    \sup_{z\in\mathcal S_{\rm out}} |z|^2\left\|U^\T\Xi(z)U\right\|
    \le \frac{14 r\mathfrak m_D}{20} + C r\mathfrak m_D \le C' r\mathfrak m_D
\]
with probability at least $1-n^{-D-10}$. This completes the proof of \eqref{eq:uniform-compressed-local-law}.

The proof of \eqref{eq:uniform-row-signal-local-law} is identical to the proof
of \eqref{eq:uniform-compressed-local-law}. One applies Theorem \ref{thm:local} to the pairs $(e_i,u_a)$ for $ i\in[n]$ and $a\in [r]$, takes a union bound
over $i,a$, and uses the same net argument over $\mathcal S_{\rm out}$.

Finally,
\[
    e_i^\T Q\Xi(z)U  = e_i^\T\Xi(z)U - e_i^\T U\cdot U^\T\Xi(z)U,
\]
so \eqref{eq:uniform-row-signal-local-law-Q} follows from
\eqref{eq:uniform-row-signal-local-law} and
\eqref{eq:uniform-compressed-local-law}.

\subsection{Proof of Corollary~\ref{coro:row-adaptive}}\label{app:proof-coro:row-adaptive}
We first collect two preliminary estimates used in the proof.  We begin with the concentration estimate used in the leave-one-out
argument. It follows from the scalar Bernstein inequality
\cite[Theorem~2.8.4]{Vbook}, combined with a standard net
argument; in the sub-Gaussian regime, we first truncate the variables at
the appropriate high-probability level.
\begin{lemma}[Vector Bernstein inequality]
\label{lem:vector-bernstein}
Let $X_1,\ldots,X_n$ be independent centered real random variables and
$c_1,\ldots,c_n \in\R^r$ be deterministic. Set
$v^2:=\sum_{j=1}^n \E (X_j^2)\|c_j\|^2$ and  $c_\infty:=\max_{1\le j\le n }\|c_j\|.$
There exists an absolute constant $C>0$ such that the following hold.
\begin{enumerate}
\item[\rm(a)] If $|X_j|\le L$ almost surely for every $j$, then, for
$t\ge C r$,
\[
    \Big\|\sum_{j=1}^n X_jc_j\Big\|  \le   C\Big(\sqrt{t}\,v+tLc_\infty \Big)
\]
with probability at least $1-2e^{-ct}$.

\item[\rm(b)] Suppose $\|X_j\|_{\psi_2}\le K\tau_j$, where
$\tau_j^2=\E (X_j^2) \le\sigma_*^2$. Then, for $t\ge C(r +\log n)$,
\[
    \Big\|\sum_{j=1}^n X_jc_j \Big\|
    \le C \Big(  \sqrt{t}\,v  +tK\sigma_*\sqrt{t}\,c_\infty \Big)
\]
with probability at least $1-4e^{-ct}$.
\end{enumerate}
The same conclusions hold conditionally if the vectors $c_j$ are measurable with respect to a $\sigma$-field independent of $X_1,\ldots,X_n$.
\end{lemma}
In particular, taking $X_j=E_{ij}$ and $t=t_{D,r}=C(D+r)\log n$, in either regime
we have
\begin{equation} \label{eq:unified-vector-concentration}
    \Big\|\sum_j E_{ij}c_j \Big\| \le C\Big[   \sqrt{t_{D,r}} \big( \sum_j\sigma_{ij}^2\|c_j\|_2^2
        \big)^{1/2} + t_{D,r}L_{D,r}^{\rm ra} \max_j\|c_j\|_2 \Big]
\end{equation}
with probability at least $1-4e^{-ct_{D,r}}$. 

\medskip

Next, we collect several resolvent identities needed for the leave-one-out approach. For $i\in[n]$, let $E^{(i)}$ be the principal minor obtained from $E$ by
deleting its $i$th row and column, and write
\[
    h_i:=(E_{ji})_{j\ne i}\in\R^{n-1},  \qquad   G^{(i)}(z):=(zI-E^{(i)})^{-1}.
\]
We also write $B^{(i)}\in\R^{(n-1)\times r}$ for the matrix obtained from $B$ by deleting its $i$th row. We defer a short proof of these resolvent identities to the end of this section. 

\begin{lemma}\label{lem:schur-identities}
For every $i\in[n]$,
\begin{align} 
    G_{ii} &= \left(  z-E_{ii}-h_i^\T G^{(i)}h_i \right)^{-1}, \label{eq:schur-Gii}\\
    e_i^\T GB &= G_{ii}\left(  e_i^\T B+h_i^\T G^{(i)}B^{(i)} \right), \label{eq:schur-row}\\
    G^{(i)}B^{(i)} &= (GB)^{(i)} - G^{(i)} h_i(e_i^\T GB), \label{eq:minor-full}
\end{align}
where $(GB)^{(i)}$ denotes the matrix obtained from $GB$ by deleting its $i$th row.
\end{lemma}

We first prove the row-adaptive local law for a fixed spectral parameter. Recall the definition of $L_{D,r}$ from \eqref{eq:L-ra} and  from \eqref{eq:t-Dm} that \[t_{D,r}:=C(D+r)\log n.\] 

\begin{prop}
\label{prop:fixed-z}
Let $B\in\R^{n\times r}$ be deterministic and suppose $|z|\ge C' \sqrt{t_{D,r}R_{\max}}$ for a sufficiently large absolute constant $C'$. Then, with
probability at least $1-Cn^{-D-2},$
\begin{equation}
    |z|^2
    \|(G(z)-\Phi(z))B\|_{2,\infty}
    \le C\left[ \sqrt{t_{D,r}}\,\Gamma(B)
        + \left( t_{D,r}L_{D,r}^{\rm ra}+\mathfrak m_D \right) \|B\|_{2,\infty}\right].
    \label{eq:fixed-z-main}
\end{equation}
\end{prop}

\begin{proof}
Our goal is to control $\|(G-\Phi)B\|_{2,\infty}$ by a
leave-one-out argument. We first note that $\|\Phi\| = \max_i |\phi_i| \le \frac{3}{2|z|}$ by Lemma \ref{lem:bdphi} (i) and thus
\begin{equation}
    \|GB\|_{2,\infty} \le \|\Phi B\|_{2,\infty} +\|(G-\Phi)B\|_{2,\infty}
    \le \frac{3}{2|z|}\|B\|_{2,\infty} +\|(G-\Phi)B\|_{2,\infty}.
    \label{eq:GB-bootstrap}
\end{equation}
This inequality will be used to close the bootstrap below.

Fix $i\in[n]$ and condition on $E^{(i)}$ and $E_{ii}$. Then
$h_i$ is independent of the conditioning $\sigma$-field, whereas
$G^{(i)}B^{(i)}$ is measurable with respect to it. Hence, combining
\eqref{eq:unified-vector-concentration} with a union bound
over $i$ yields
\begin{align}\label{eq:conditional-row-bound}
    \|h_i^\T G^{(i)}B^{(i)}\|
    \le C\Big[ \sqrt{t_{D,r}} \big(  \sum_{j\neq i}\sigma_{ij}^2
            \|e_j^\T G^{(i)}B^{(i)}\|^2 \big)^{1/2}+ t_{D,r}L_{D,r}^{\rm ra} \|G^{(i)}B^{(i)}\|_{2,\infty}
    \Big]
\end{align}
for all $i$, with probability at least $1-4ne^{-ct_{D,r}}$.

We now estimate the two terms on the right-hand side of
\eqref{eq:conditional-row-bound}. By \eqref{eq:minor-full}, for
$j\ne i$,
\[
    e_j^\T G^{(i)}B^{(i)} =  e_j^\T GB - e_j^\T G^{(i)}h_i\cdot e_i^\T GB.
\]
Decomposing
\[
    e_j^\T GB = e_j^\T \Phi B + e_j^\T(G-\Phi)B= \phi_j e_j^\T B + e_j^\T(G-\Phi)B
\]
and applying Minkowski's inequality, we obtain
\begin{align*}
    \big( \sum_{j\ne i}  \sigma_{ij}^2 \|e_j^\T G^{(i)}B^{(i)}\|^2 \big)^{1/2}
    \le&
    \big(  \sum_{j\ne i} \sigma_{ij}^2 |\phi_j|^2 \|e_j^\T B\|^2 \big)^{1/2}
    + \sqrt{R_i}\,
    \|(G-\Phi)B\|_{2,\infty}\\
    & +
    \big(  \sum_{j\ne i}  \sigma_{ij}^2  |e_j^\T G^{(i)}h_i|^2 \big)^{1/2} \|GB\|_{2,\infty}.
\end{align*}
Using $\max_j|\phi_j|\le 3/(2|z|)$ and $\sigma_{ij}\le\sigma_{\max}$, we further get
\begin{equation}
    \big(  \sum_{j\ne i} \sigma_{ij}^2  \|e_j^\T G^{(i)}B^{(i)}\|^2 \big)^{1/2}
    \le \frac{3}{2|z|}\gamma_i(B) + \sqrt{R_i}\, \|(G-\Phi)B\|_{2,\infty} + \sigma_{\max}\|G^{(i)}h_i\|\cdot \|GB\|_{2,\infty}.
    \label{eq:weighted-minor-bound}
\end{equation}
Likewise, \eqref{eq:minor-full} yields
\begin{equation}
    \|G^{(i)}B^{(i)}\|_{2,\infty}
    \le  \bigl(1+\|G^{(i)}h_i\| \bigr)  \|GB\|_{2,\infty}.
    \label{eq:minor-row-bound}
\end{equation}
Combining \eqref{eq:conditional-row-bound},
\eqref{eq:weighted-minor-bound}, and
\eqref{eq:minor-row-bound}, and using $\sigma_{\max}\le L_{D,r}^{\rm ra}$ with
\[
    \|G^{(i)}h_i\|
    \le \|G^{(i)}\| \, \|E e_i\| \le \|G^{(i)}\| \, \|E\| \le \frac{6\sqrt{R_{\max}}}{|z|},
\]
we obtain
\begin{align*}
    &\|h_i^\T G^{(i)}B^{(i)}\|\le C\Bigg[\frac{\sqrt{t_{D,r}}}{|z|}\gamma_i(B) +
        \sqrt{t_{D,r}R_i}\,  \|(G-\Phi)B\|_{2,\infty} \\
        &\qquad\qquad\qquad\qquad +  t_{D,r}L_{D,r}^{\rm ra}
        \left(  1+\frac{\sqrt{R_{\max}}}{|z|} \right)\|GB\|_{2,\infty}
    \Bigg].
\end{align*}
Since $|z|\ge C'\sqrt{t_{D,r}R_{\max}}$ and by \eqref{eq:GB-bootstrap},
\[
    \|GB\|_{2,\infty} \le \frac{3}{2|z|}\|B\|_{2,\infty}  +\|(G-\Phi)B\|_{2,\infty},
\]
this simplifies to
\begin{align}
    \|h_i^\T G^{(i)}B^{(i)}\| &\le C\Bigg[
    \frac{\sqrt{t_{D,r}}}{|z|}\gamma_i(B) +
        \frac{t_{D,r}L_{D,r}^{\rm ra}}{|z|} \|B\|_{2,\infty}\\
        &\qquad\quad+
        \left( \sqrt{t_{D,r}R_i}+t_{D,r}L_{D,r}^{\rm ra}\right) \|(G-\Phi)B\|_{2,\infty}
    \Bigg].
    \label{eq:loo-row-final}
\end{align}

To proceed, note that, by \eqref{eq:schur-row},
\[
    e_i^\T(G-\Phi)B  =  (G_{ii}-\phi_i)e_i^\T B +  G_{ii}h_i^\T G^{(i)}B^{(i)}.
\]
Using \eqref{eq:loo-row-final} and taking the maximum over $i$, together with $\|G\| \le \frac{2}{|z|}$, we get
\begin{align}
    \|(G-\Phi)B\|_{2,\infty}
    \le{}&
    \max_i|G_{ii}-\phi_i|\,\|B\|_{2,\infty}  \notag\\
    &+ \frac{C}{|z|^2} \left[ \sqrt{t_{D,r}}\,\Gamma(B)  + t_{D,r}L_{D,r}^{\rm ra}\|B\|_{2,\infty}
    \right] \notag\\
    &+ C\frac{\sqrt{t_{D,r}R_{\max}} +t_{D,r}L_{D,r}^{\rm ra}}{|z|} \|(G-\Phi)B\|_{2,\infty}.
    \label{eq:bootstrap-inequality}
\end{align}
By Assumption~\ref{assump:row-adaptive},
\[
    t_{D,r}L_{D,r}^{\rm ra} \lesssim\sqrt{t_{D,r}R_{\max}} \le (1/C')|z|,
\]
so the last coefficient in \eqref{eq:bootstrap-inequality} is at most
$1/2$ when $C'$ is sufficiently large. Therefore,
\[
    \|(G-\Phi)B\|_{2,\infty}
    \le 2\max_i|G_{ii}-\phi_i|\,\|B\|_{2,\infty} + \frac{C}{|z|^2}
    \left[  \sqrt{t_{D,r}}\,\Gamma(B)  + t_{D,r}L_{D,r}^{\rm ra}\|B\|_{2,\infty} \right].
\]
Finally, it follows by applying Theorem \ref{thm:local} with $v=w=e_i$ and taking a union
bound over $i$ that
\[
    \max_i |G_{ii}-\phi_i|  \le \frac{\mathfrak m_D}{|z|^2},
\]
which proves \eqref{eq:fixed-z-main}.
\end{proof}

We now pass from the fixed-$z$ estimate to the uniform row-adaptive
local law. The argument follows the proof of
Corollary~\ref{lem:local-law-uniform}, so we only sketch it.
\begin{proof}[Proof of Corollary~\ref{coro:row-adaptive}]
It remains to make Proposition~\ref{prop:fixed-z} uniform over
$z\in\mathcal S_{\rm ra}$ defined in \eqref{eq:S-ra}. On the event $\{\|E\|\le3\sqrt{R_{\max}}\}$, the same resolvent calculation as in the proof of Corollary~\ref{lem:local-law-uniform}
shows that the maps rowwise$z\mapsto z^2G(z)B$ and $z\mapsto z^2\Phi(z)B$ are both
$8\|B\|$-Lipschitz.

We now apply the same $\varepsilon$-net argument as in the proof of
Corollary~\ref{lem:local-law-uniform}. By the polynomial bounds in
Assumption~\ref{assump:noise}, a suitable net of
$\mathcal S_{\rm ra}$ has cardinality polynomial in $n$. Applying
Proposition~\ref{prop:fixed-z} on this net and taking a union bound,
the choice
\[
    t_{D,r}=C(D+r)\log n
\]
absorbs the net cardinality and yields failure probability at most
$n^{-D}$. The Lipschitz bounds above extend the estimate from the net
to all $z\in\mathcal S_{\rm ra}$. This proves \eqref{eq:row-adaptive-main}.
\end{proof}

Finally, we provide the proof of Lemma \ref{lem:schur-identities}. 
\begin{proof}[Proof of Lemma \ref{lem:schur-identities}]
Without loss of generality, denote \[
    zI-E = \begin{pmatrix}
        z-E_{ii}&-h_i^\T\\
        -h_i&zI-E^{(i)}
    \end{pmatrix}.
\]
Write accordingly
\[
    G= \begin{pmatrix}
        G_{ii}&G_{i,-i}\\
        G_{-i,i}&G_{-i,-i}
    \end{pmatrix},
\]
where the subscript $-i$ means that the $i$th coordinate is
deleted.  The block inverse formula (the Schur complement) yields \eqref{eq:schur-Gii} and 
\[
    G_{i,-i}=G_{ii}h_i^\T G^{(i)},\qquad
    G_{-i,-i}  =  G^{(i)} +G^{(i)}h_iG_{ii}h_i^\T G^{(i)},
\]
where $-i$ denotes deletion of the $i$th coordinate.
Since
\[
    B= \begin{pmatrix}
        e_i^\T B\\ B^{(i)}
    \end{pmatrix},
\]
multiplying the first block row by $B$ leads to
\[
    e_i^\T GB =
    G_{ii}\left( e_i^\T B+h_i^\T G^{(i)}B^{(i)} \right),
\]
which is \eqref{eq:schur-row}. Similarly, the lower block row gives
\[
    (GB)^{(i)} = G^{(i)}B^{(i)} +G^{(i)}h_iG_{ii}  \left( e_i^\T B+h_i^\T G^{(i)}B^{(i)} \right)
    = G^{(i)}B^{(i)} +G^{(i)}h_i(e_i^\T GB),
\]
which is equivalent to \eqref{eq:minor-full}.
\end{proof}

\section{Extensions to general signal eigenspaces} \label{app:general-simple-facts}
Theorem~\ref{thm:top-k-eigenspace} is stated for the top positive eigenspace
only to simplify notation. The proof extends directly to general eigenspaces $U_{l:s}=(u_l,\dots,u_s)$ for $1\le l < s \le r$. We also provide a simple alternative approach using the bounds for top eigenspaces.
     
For negative population spikes, we apply our framework to $-A$ and $-(A+E)$.  In our ordering $\lambda_1 \ge \dots \ge \lambda_{r_+} > 0 > \lambda_{r_{+} +1}\ge... \ge \lambda_{r}$, the eigenspaces corresponding to the most negative eigenvalues become the top eigenspaces. This symmetry yields the same perturbation guarantees for the bottom eigenspaces without additional analysis.

For intermediate blocks $U_{l:s}$ and $\widetilde{U}_{l:s}$ where $1\le l < s \le r_+$. We use the decomposition $P_{U_{l:s}} = P_{U_{s}}-P_{U_{l-1}}$ to obtain
     \begin{align}\label{eq:0423-general-sine}
         &\|\sin\angle(U_{l:s}, \widetilde{U}_{l:s})\| \le   \|\sin\angle(U_{l-1}, \widetilde{U}_{l-1})\| +\|\sin\angle(U_s, \widetilde{U}_s)\|,
         \end{align}
     and
          \begin{align}\label{eq:0423-general-2inf}
         &\min_{O} \|\widetilde{U}_{l:s} - U_{l:s}O\|_{2,\infty} \le \|(I - P_{U_{l:s}}) \widetilde{U}_{l:s}\|_{2,\infty} + \|U_{l:s}\|_{2,\infty} \|\sin\angle(U_{l:s}, \widetilde{U}_{l:s})\|^2,
     \end{align}
     where 
     \[\|(I - P_{U_{l:s}}) \widetilde{U}_{l:s}\|_{2,\infty} \le \|(I - P_{U_{s}}) \widetilde{U}_{s}\|_{2,\infty} + \|U_{l-1}\|_{2,\infty} \|\sin\angle(U_{l-1}, \widetilde{U}_{l-1})\|. \]
Thus, the perturbation bounds for $U_{l:s}, \widetilde{U}_{l:s}$ follow by substituting the estimates for the top-$s$ and top-$(l-1)$ eigenspaces.  By symmetry, analogous decompositions hold for intermediate blocks of negative spikes. 

We provide the proofs of (\ref{eq:0423-general-sine}) and (\ref{eq:0423-general-2inf}) below.
\medskip

\noindent\underline{Proof of \eqref{eq:0423-general-sine}:}  First, by \cite[Exercises VII. 1. 9– 1.11]{Bhatia},
\[  
\|\sin\angle(U_{l:s}, \widetilde{U}_{l:s})\| = \|(I - P_{U_{l:s}}) \widetilde{U}_{l:s}\| = \|P_{U_{l:s}^\perp} \widetilde{U}_{l:s}\|.
\]
From the decomposition $P_{U_{l:s}} = P_{U_{s}} - P_{U_{l-1}}$, we can express
\[
P_{U_{l:s}^\perp} = I - P_{U_{l:s}} = (I - P_{U_{s}}) + P_{U_{l-1}} = P_{U_{s}^\perp} + P_{U_{l-1}}.
\]
It follows that
\begin{equation}\label{eq:sine-triangle}
\|P_{U_{l:s}^\perp} \widetilde{U}_{l:s}\| \le \|P_{U_{s}^\perp} \widetilde{U}_{l:s}\| + \|P_{U_{l-1}} \widetilde{U}_{l:s}\|.
\end{equation}
For the first term, observe that the columns of $\widetilde{U}_{l:s}$ are a subset of the columns of $\widetilde{U}_{s}$. Thus
\[ 
\|P_{U_{s}^\perp} \widetilde{U}_{l:s}\| \le \|P_{U_{s}^\perp} \widetilde{U}_{s}\| = \|\sin\angle(U_{s}, \widetilde{U}_{s})\|.
\]
For the second term in \eqref{eq:sine-triangle}, since $\widetilde{U}_{l:s}$ is  orthogonal to $\widetilde{U}_{l-1}$, we have $P_{\widetilde{U}_{l-1}^\perp} \widetilde{U}_{l:s} = \widetilde{U}_{l:s}$. Then
\[
P_{U_{l-1}} \widetilde{U}_{l:s} = P_{U_{l-1}} P_{\widetilde{U}_{l-1}^\perp} \widetilde{U}_{l:s},
\]
and therefore,
\[
 \|P_{U_{l-1}} P_{\widetilde{U}_{l-1}^\perp} \widetilde{U}_{l:s}\| \le \|P_{U_{l-1}} P_{\widetilde{U}_{l-1}^\perp}\|\cdot \|\widetilde{U}_{l:s}\| = \|\sin\angle(U_{l-1}, \widetilde{U}_{l-1})\|.
 \]
Substituting these two bounds back into \eqref{eq:sine-triangle} completes the proof of \eqref{eq:0423-general-sine}.
\medskip

\noindent\underline{Proof of \eqref{eq:0423-general-2inf}:} Let $U_{l:s}^\T \widetilde{U}_{l:s} = W D V^\T$ be the singular value decomposition. We choose $O = W V^\T$ and then decompose 
\begin{align*}
\widetilde{U}_{l:s} - U_{l:s}O 
&= (I - P_{U_{l:s}})\widetilde{U}_{l:s} + P_{U_{l:s}}(\widetilde{U}_{l:s} - U_{l:s}O) \\ 
&= (I - P_{U_{l:s}})\widetilde{U}_{l:s} + U_{l:s}(U_{l:s}^\T \widetilde{U}_{l:s} - O).
\end{align*}
Taking the $2 \to \infty$ norm and applying the triangle inequality yields
\begin{equation}\label{eq:2inf-decomp-main}
\|\widetilde{U}_{l:s} - U_{l:s}O \|_{2,\infty} \le \|(I - P_{U_{l:s}})\widetilde{U}_{l:s} \|_{2,\infty} + \|U_{l:s}(U_{l:s}^T \widetilde{U}_{l:s} - O)\|_{2,\infty}.
\end{equation}
By our choice of $O$, the second term in \eqref{eq:2inf-decomp-main} is bounded by
\[   
\|U_{l:s} W(D - I)V^\T \|_{2,\infty} \le \|U_{l:s}\|_{2,\infty} \|W(I - D)V^\T \|  = \|U_{l:s}\|_{2,\infty} \|I - D\|.
\]
Note that the diagonal entries of $D$ are the cosines of the principal angles $\theta_i$. The entries of $I-D^2$ are exactly $\sin^2\theta_i$. Because $D \succeq 0$, we have $I - D \preceq I - D^2$, which implies
\[
 \|I - D\| \le \|I - D^2\| = \|\sin\angle(U_{l:s}, \widetilde{U}_{l:s})\|^2.
\] 
This proves the main inequality in \eqref{eq:0423-general-2inf}. It remains to bound $\|(I - P_{U_{l:s}}) \widetilde{U}_{l:s}\|_{2,\infty}$. Using $I - P_{U_{l:s}} = (I - P_{U_{s}}) + P_{U_{l-1}}$, we expand
\[
(I - P_{U_{l:s}}) \widetilde{U}_{l:s} = (I - P_{U_{s}}) \widetilde{U}_{l:s} + P_{U_{l-1}} \widetilde{U}_{l:s}.
\]
Thus
\begin{equation}\label{eq:2inf-split}
\|(I - P_{U_{l:s}}) \widetilde{U}_{l:s}\|_{2,\infty} \le \|(I - P_{U_{s}}) \widetilde{U}_{l:s}\|_{2,\infty} + \|P_{U_{l-1}} \widetilde{U}_{l:s}\|_{2,\infty}.
\end{equation}
For the first term, since $(I - P_{U_{s}}) \widetilde{U}_{l:s}$ consists of a subset of the columns of $(I - P_{U_{s}}) \widetilde{U}_{s}$, we get
\[
\|(I - P_{U_{s}}) \widetilde{U}_{l:s}\|_{2,\infty} \le \|(I - P_{U_{s}}) \widetilde{U}_{s}\|_{2,\infty}.
\]
For the second term, plugging in $P_{U_{l-1}} = U_{l-1}U_{l-1}^\T$, we find
\[
 \|P_{U_{l-1}} \widetilde{U}_{l:s}\|_{2,\infty} = \|U_{l-1} (U_{l-1}^\T \widetilde{U}_{l:s})\|_{2,\infty} \le \|U_{l-1}\|_{2,\infty} \|U_{l-1}^\T \widetilde{U}_{l:s}\|.
\] 
Finally, because $\widetilde{U}_{l:s}$ is orthogonal to $\widetilde{U}_{l-1}$, we obtain
\[
\|U_{l-1}^\T \widetilde{U}_{l:s}\| = \|U_{l-1}^\T P_{\widetilde{U}_{l-1}^\perp} \widetilde{U}_{l:s}\| \le \|U_{l-1}^\T P_{\widetilde{U}_{l-1}^\perp}\| \|\widetilde{U}_{l:s}\| = \|\sin\angle(U_{l-1}, \widetilde{U}_{l-1})\|.
\]
Substituting these bounds into \eqref{eq:2inf-split} finishes the proof.

\section{Extension to rectangular matrices: singular subspace perturbations}\label{app:rectangular}
We briefly discuss the extension of current results to the rectangular model: 
\begin{align*}
    \widetilde{X} = X + E \in \mathbb{R}^{n_1 \times n_2},
\end{align*}
where $X = U \Sigma V^\T = \sum_{i=1}^r \sigma_i u_i v_i^\T$ is a rank $r$ matrix with singular values $\sigma_1 \ge \sigma_2 \ge \dots \ge \sigma_r > 0$. The random noise matrix $E$ has i.i.d. entries with $\E(E_{ij})=0$ and $\mathbb{E}(E_{ij}^2) = \sigma_{ij}^2. $
Denote $\Sigma = (\sigma_{ij}^2) \in \mathbb{R}^{n_1 \times n_2}.$ The entries of $E$ are sub-Gaussian with a variance profile: $\|E_{ij}\|_{\psi_2}\le  K \sigma_{ij} .$

Using the standard Hermitian dilation (or linearization), we embed the rectangular matrices into larger symmetric block matrices by considering the $(n_1 + n_2) \times (n_1 + n_2)$ augmented symmetric matrices:
$$\mathbf{A} = \begin{pmatrix} 0 & X \\ X^\T & 0 \end{pmatrix}, \qquad \mathbf{E} = \begin{pmatrix} 0 & E \\ E^\T & 0 \end{pmatrix},\quad \widetilde{\mathbf{A}} = \mathbf{A}+\mathbf{E}.$$
The nonzero eigenvalues of $\mathbf{A}$ are
\[
\lambda_1=\sigma_1\ge \dots\ge\lambda_r=\sigma_r>0>\lambda_{r+1}=-\sigma_1\ge\dots\ge\lambda_{2r}=-\sigma_r
\]
with corresponding eigenvectors
\[
w_j:=\frac1{\sqrt2}\binom{u_j}{v_j},
\qquad
w_{j+r}:=\frac1{\sqrt2}\binom{u_j}{-v_j},
\qquad 1\le j\le r.
\] 
Let $\mathbf{A} = WDW^\T$, where 
\[ W= \frac{1}{\sqrt 2}\begin{pmatrix}
    U & U\\
    V & -V
\end{pmatrix}.\]
For $I \subseteq [2r]$, we denote by $W_I$ the submatrix of $W$ formed by columns $\{w_i : i \in I\}$. All information about the top-$k$ singular subspaces $U_k$ and $V_k$ is encoded in the eigenspace $W_{I_k}$, where
\[ I_k := \{1,\dots,k\} \cup \{k+1,\dots,k+r\}.\]
Hence, we can directly translate our eigenspace perturbation results to singular subspaces. We now introduce the relevant parameters. Define
\[
R_i^{(L)}:=\sum_{j=1}^{n_1} \sigma_{ij}^2,\qquad
R_j^{(R)}:=\sum_{i=1}^{n_2} \sigma_{ij}^2,
\]
and let $\mathcal R_L:=\diag(R_i^{(L)})_{1\le i \le n_1}$, $\mathcal R_R:=\diag(R_i^{(R)})_{1\le i \le n_2}$,
\[
\mathcal R_{\rm aug}
:=
\diag(\mathcal R_L, \mathcal R_R),
\qquad
R_{\max}^{\rm aug}
:=
\max\Bigl\{\max_i R_i^{(L)},\,\max_j R_j^{(R)}\Bigr\}.
\]
Define $M_{\rm aug}\equiv M_{\rm aug, D}$ as in \eqref{def:M} with the dimension factor $n_1+n_2$ and set 
\[
\rho_k^{\rm aug}:=10M_{\rm aug}+15\frac{\osc(\mathcal R_{\rm aug})}{\sigma_k},
\]
and
\[
\delta_k:=\sigma_k-\sigma_{k+1}.
\]
Also, we denote $\mathcal V = W^\T \mathcal{R}_{\rm aug} W$. Note that 
\[\mathcal{V}_{\mathcal{J}\mathcal{I}} = W_{\mathcal{J}}^\T \mathcal{R}_{\rm aug} W_{\mathcal{I}} = \frac{1}{2} \left( U_{\mathcal{J}}^\T \mathcal{R}_L U_{\mathcal{I}} + V_{\mathcal{J}}^\T \mathcal{R}_R V_{\mathcal{I}} \right).\]
We define $\mathcal B_k^{\rm aug}$ similar to $\mathcal B_k$ \eqref{def:bias-k}:
\[\mathcal B_k^{\rm aug} :=\frac{10\sqrt{k}}{\delta_k \sigma_k}\left( \|\mathcal{V}_{\mathcal{J}\mathcal{I}}\|+8\frac{R_{\max}^{\rm aug}}{\sigma_k^2}\osc(\mathcal R_{\rm aug})\right) +4 \sqrt{k}\frac{\|\mathcal{V}_{\mathcal{N}\mathcal{K}}\|}{\sigma_k^2}.\]

We can then directly translate eigenspace perturbation results from the symmetric model $\widetilde{\mathbf{A}} = \mathbf{A}+\mathbf{E}$ to the rectangular model $\widetilde{X} = X + E$ using the revised parameters above and the following key relationships (see \cite{Wang24} for details):  for any unitarily invariant norm $\vvvert \cdot \vvvert,$
\begin{align*}
    \max\{ \vvvert \sin\angle(U_k, \widetilde{U}_k) \vvvert, \vvvert \sin\angle(V_k, \widetilde{V}_k) \vvvert\} \le \vvvert \sin\angle(W_I, \widetilde{W}_k) \vvvert
\end{align*}
and 
\begin{align*}
    \max\{ \|\widetilde{U}_k - P_{U_k} \widetilde{U}_k\|_{2,\infty},\|\widetilde{V}_k - P_{V_k} \widetilde{V}_k\|_{2,\infty}\} \le \|\widetilde{W}_I - P_{W_I} \widetilde{W}_I\|_{2,\infty}.
\end{align*}

Specifically, this implies that the following $2 \to \infty$ bound  holds with  probability at least $1-n^{-D}$:
\begin{align}\label{eq:symmetric-2inf}
&\min_{O\in\mathbb O(k)}
\|\widetilde U_k-U_kO\|_{2,\infty} \nonumber\\
&\qquad\le 4\max\{\|U\|_{2,\infty},\|V\|_{2,\infty}\}
\left(
{\mathcal B}_k^{\rm aug}
+
30\frac{R_{\max}^{\rm aug}}{\sigma_k^2}
\right)
+
120\sqrt{k}\frac{M_{\rm aug}}{\sigma_k}.
\end{align}
The same bound also holds for $\min_{O\in\mathbb O(k)}
\|\widetilde V_k-V_kO\|_{2,\infty}$.

We conclude by clarifying the nature of the coupled terms in \eqref{eq:symmetric-2inf}. While the dependence on both variance profiles ($\mathcal{R}_L$ and $\mathcal{R}_R$) reflects the intrinsic statistical coupling in the asymmetric noise matrix, the dependence on the maximum joint incoherence $\max\{\|U\|_{2,\infty},\|V\|_{2,\infty}\}$ and the total dimension $n_1+n_2$ is an artifact of the symmetric embedding framework. We plan to address this limitation in future work.

\section{Proofs of Lemma~\ref{lem:J-block-stability}, Proposition~\ref{prop:oracle-block-debias}, (\ref{eq:plugin-sin-short}) and (\ref{eq:plugin-2inf-short})}\label{app:de-bias}
\subsection{Proof of Lemma~\ref{lem:J-block-stability}}\label{app:J-block-stability}
For brevity, write
\[
L_{\mathcal J,s}:=I-U_{\mathcal J}^\T \Phi(\widetilde\lambda_s)U_{\mathcal J}\Lambda_{\mathcal J}.
\]
For $j\in \mathcal J$, set
\[
\alpha_j(z):=\lambda_j^{-1}-u_j^\T \Phi(z)u_j.
\]
Recall the quantity defined in Proposition \ref{prop:detbd-vector}
\[
\Delta_{\mathcal J}^\Phi(\widetilde\lambda_s)
:=
\widetilde\lambda_s \min_{j\in {\mathcal J}}\left|1-\lambda_j u_j^\T\Phi(\widetilde\lambda_s)u_j\right|
=\min_{j\in {\mathcal J}}\left| \widetilde\lambda_s \lambda_j\alpha_j(\widetilde\lambda_s)\right|.
\]

We first relate $s_{\min}(L_{{\mathcal J},s})$ to $\Delta_{\mathcal J}^\Phi(\widetilde\lambda_s)$.
Decompose
\[
U_{\mathcal J}^\T\Phi(\widetilde\lambda_s)U_{\mathcal J}
=
\diag\bigl(u_j^\T\Phi(\widetilde\lambda_s)u_j\bigr)_{j\in {\mathcal J}}
+
Q_{\mathcal J}(\widetilde\lambda_s),
\]
where
\[
Q_{\mathcal J}(\widetilde\lambda_s):=\off\bigl(U_{\mathcal J}^\T\Phi(\widetilde\lambda_s)U_{\mathcal J}\bigr).
\]
Then
\[
L_{{\mathcal J},s}
=
\diag\left(1-\lambda_j u_j^\T \Phi(\widetilde\lambda_s) u_j\right)_{j\in {\mathcal J}}
-
Q_{\mathcal J}(\widetilde\lambda_s)\Lambda_{\mathcal J}.
\]
Hence
\begin{equation}
\label{eq:J-stability-start}
s_{\min}(L_{{\mathcal J},s})
\ge
\frac{\Delta_{\mathcal J}^\Phi(\widetilde\lambda_s)}{\widetilde\lambda_s}
-
\|Q_{\mathcal J}(\widetilde\lambda_s)\Lambda_{\mathcal J}\|.
\end{equation}

We next estimate the off-diagonal term $\|Q_{\mathcal J}(\widetilde\lambda_s)\Lambda_{\mathcal J}\|$. Using the expansion
\[
\Phi(z)=\frac1z I+\frac1{z^3}R+\varepsilon(z),
\]
we have
\[
Q_{\mathcal J}(\widetilde\lambda_s)
=
\frac1{\widetilde\lambda_s^3}\off(\mathcal V_{{\mathcal J}{\mathcal J}})
+
\off\bigl(U_{\mathcal J}^\T \varepsilon (\widetilde\lambda_s)U_{\mathcal J}\bigr).
\]
Since the eigenvalues in ${\mathcal J}$ are positive and bounded above by
$\lambda_{k+1}$,
\[
\|Q_{\mathcal J}(\widetilde{\lambda}_s)\Lambda_{\mathcal J}\|
\le
\frac{\lambda_{k+1}}{\widetilde{\lambda}_s^3}\|\off(\mathcal V_{{\mathcal J} {\mathcal J}})\|
+
\lambda_{k+1}
\left\|
\off\bigl(U_{\mathcal J}^\T \varepsilon(\widetilde{\lambda}_s)U_{\mathcal J}\bigr)
\right\|.
\]
By Lemma~\ref{lem:osc},
\[
\left\|
\off\bigl(U_{\mathcal J}^\T\varepsilon(\widetilde{\lambda}_s)U_{\mathcal J}\bigr)
\right\|
\le
\frac12\osc(\varepsilon(\widetilde{\lambda}_s))
\le
\frac{15}{2}\frac{R_{\max}}{\widetilde{\lambda}_s^5}\osc(\mathcal R).
\]
 Therefore
\begin{equation}
\label{eq:QJLambda-bound}
\|Q_{\mathcal J}(\widetilde{\lambda}_s)\Lambda_{\mathcal J}\|
\le
\frac{\lambda_{k+1}}{\widetilde{\lambda}_s^3}\|\off(V_{{\mathcal J} {\mathcal J}})\|
+
\frac{15}{2}\frac{\lambda_{k+1}R_{\max}}{\widetilde{\lambda}_s^5}\osc(\mathcal R).
\end{equation}
We now verify that the off-diagonal term in \eqref{eq:J-stability-start} is
small compared with the diagonal part. 
Using $\widetilde{\lambda}_s\ge \frac{1}{2}\lambda_k$, $R_{\max}\le 36\lambda_k^2$, and
$
\|\off(\mathcal V_{{\mathcal J} {\mathcal J}})\|\le \osc(\mathcal R),
$
from \eqref{eq:QJLambda-bound}, we get
\[
\widetilde{\lambda}_s\|Q_{\mathcal J}(\widetilde{\lambda}_s)\Lambda_{\mathcal J}\|
\le
\left(8+\frac{5}{3} \right)\frac{\lambda_{k+1}}{\lambda_k^2}\osc(\mathcal R)
\le
\frac{1}{20}\left(8+\frac{5}{3} \right)\delta_k,
\]
where the last inequality again follows from the gap assumption. Since
$\Delta_{\mathcal J}^\Phi(\widetilde{\lambda}_s)\ge\delta_k/2$, we get
\[
\|Q_{\mathcal J}(\widetilde{\lambda}_s)\Lambda_{\mathcal J}\|
\le
\frac{1}{10}\left(8+\frac{5}{3} \right)\frac{\Delta_{\mathcal J}^\Phi(\widetilde{\lambda}_s)}{\widetilde{\lambda}_s}.
\]
Returning to \eqref{eq:J-stability-start}, we obtain
\[
s_{\min}(L_{{\mathcal J},s})
\ge
\frac{1}{30}\frac{\Delta_J^\Phi(\widetilde{\lambda}_s)}{\widetilde{\lambda}_s}
\ge
\frac{1}{60}\frac{\delta_k}{\widetilde{\lambda}_s}.
\]
Thus $\|L_{{\mathcal J},s}^{-1}\| \le 60 \frac{\widetilde{\lambda}_s}{\delta_k}$. This proves the desired bound.


\subsection{Proof of Proposition \ref{prop:oracle-block-debias}}\label{app:oracle-block-debias}
For $s\in [k]$, the eigen-equation gives
\[\widetilde u_s=G(\widetilde\lambda_s)A\widetilde u_s = \Phi(\widetilde\lambda_s)A\widetilde u_s + \Xi(\widetilde\lambda_s)A\widetilde u_s,\]
where $\Xi(z):=G(z)-\Phi(z).$ It follows that
\[U^\T\widetilde u_s =
U^\T\Phi(\widetilde\lambda_s)U\Lambda U^\T \widetilde u_s +
U^\T\Xi(\widetilde\lambda_s)U\Lambda U^\T \widetilde u_s.\]
Taking the $\mathcal J$-block yields
\begin{align*}
    U_{\mathcal J}^\T \widetilde u_s = U_{\mathcal J}^\T\Phi(\widetilde{\lambda}_s)U_{\mathcal I} \Lambda_{\mathcal I}\cdot U_{\mathcal I}^\T \widetilde{u}_s + U_{\mathcal J}^\T\Phi(\widetilde{\lambda}_s)U_{\mathcal J} \Lambda_{\mathcal J}\cdot U_{\mathcal J}^\T \widetilde{u}_s + U_{\mathcal J}^\T\Xi(\widetilde\lambda_s)U \Lambda U^\T \widetilde{u}_s.
\end{align*}
Collecting the $U_{\mathcal J}^\T \widetilde u_s$ terms gives
\[
(I- U_{\mathcal{J}}^\T \Phi(\widetilde{\lambda}_s) U_{\mathcal J}\Lambda_{\mathcal J}) U_{\mathcal J}^\T \widetilde u_s = U_{\mathcal J}^\T\Phi(\widetilde{\lambda}_s)U_{\mathcal I} \Lambda_{\mathcal I}\cdot U_{\mathcal I}^\T \widetilde{u}_s + U_{\mathcal J}^\T\Xi(\widetilde\lambda_s)U \Lambda U^\T \widetilde{u}_s.
\]
Recall $\mathcal T_s$ from \eqref{def:Ts}. We further get
\begin{align}\label{eq:0601-J-correct}
U_{\mathcal J}^\T \widetilde u_s - \mathcal T_s U_{\mathcal I}^\T \widetilde{u}_s = (I- U_{\mathcal{J}}^\T \Phi(\widetilde{\lambda}_s) U_{\mathcal J}\Lambda_{\mathcal J})^{-1} U_{\mathcal J}^\T\Xi(\widetilde\lambda_s)U \Lambda U^\T \widetilde{u}_s.
\end{align}
By definition and $\mathcal{N}=\{r_{+}+1,\dots,r\}$,
\[
P_U\widetilde u_s =
U_kU_k^\T\widetilde u_s +
U_{\mathcal J} U_{\mathcal J}^\T\widetilde u_s +
U_{\mathcal N}U_{\mathcal N}^\T\widetilde u_s.
\]
Combining this with \eqref{eq:0601-J-correct}, we get
\begin{align}\label{eq:w-orc-01}
   w_s^{\orc}&= P_U\widetilde u_s - U_{\mathcal J}{\mathcal T}_s U_{\mathcal I}^\T \widetilde u_s \nonumber\\
   &=U_kU_k^\T\widetilde u_s + U_{\mathcal N}U_{\mathcal N}^\T\widetilde u_s + U_{\mathcal J}(I- U_{\mathcal{J}}^\T \Phi(\widetilde{\lambda}_s) U_{\mathcal J}\Lambda_{\mathcal J})^{-1} U_{\mathcal J}^\T\Xi(\widetilde\lambda_s)U\cdot \Lambda U^\T \widetilde{u}_s.
\end{align}
Note that we have established in \eqref{eq:0601-lateruse} that
\begin{align}\label{eq:0601-UN-later}
\|U_{\mathcal N}^\T \widetilde{u}_s\| < 4 \frac{\|\mathcal{V}_{\mathcal{N}\mathcal{K}}\|}{\lambda_k^2} +24\frac{R_{\max}}{{\lambda}_k^4} \osc(\mathcal R)+ 4 \frac{M}{\lambda_k}. 
\end{align}

Now we estimate the last term on the right-hand side of \eqref{eq:w-orc-01}. We first have $$\|\Lambda U^\T \widetilde{u}_s\| \le 2|\widetilde\lambda_s|.$$
To see this, from $(A+E)\widetilde u_s=\widetilde\lambda_s\widetilde u_s,$ we have 
\[A\widetilde u_s=(\widetilde\lambda_s I-E)\widetilde u_s.\]Multiplying by $U^\T$ gives $\Lambda U^\T \widetilde u_s
=
U^\T A\widetilde u_s = U^\T(\widetilde\lambda_s I-E)\widetilde u_s.$ Hence, \[\|\Lambda U^\T \widetilde u_s\|
\le
|\widetilde\lambda_s|+\|E\|\le 2|\widetilde\lambda_s|\]
by our assumption $\lambda_k\ge 9\sqrt{R_{\max}}$ and Weyl's inequality. 

On the event $\mathcal E_D$, we have
\[\|U_{\mathcal J}^\T\Xi(\widetilde\lambda_s)U\| \le \frac{M}{\widetilde\lambda_s^2}. \]
Combining these estimates with Lemma \ref{lem:J-block-stability} yields
\begin{align}\label{eq:0601-error-bs}
&\|U_{\mathcal J}(I- U_{\mathcal{J}}^\T \Phi(\widetilde{\lambda}_s) U_{\mathcal J}\Lambda_{\mathcal J})^{-1} U_{\mathcal J}^\T\Xi(\widetilde\lambda_s)U\cdot \Lambda U^\T \widetilde{u}_s \|\nonumber\\
&\qquad\le \|(I- U_{\mathcal{J}}^\T \Phi(\widetilde{\lambda}_s) U_{\mathcal J}\Lambda_{\mathcal J})^{-1} U_{\mathcal J}^\T\Xi(\widetilde\lambda_s)U\cdot \Lambda U^\T \widetilde{u}_s \|\le 120\frac{M}{\delta_k}.
\end{align}

From \eqref{eq:w-orc-01} and the orthogonality of the signal subspaces, $w_s^{\orc} - U_kU_k^\T\widetilde u_s$ lies entirely in the span of $U_k^\perp$. Combining the above bounds, we have
\begin{align*}
    \|w_s^{\orc} - U_kU_k^\T\widetilde u_s\| & \le \|U_{\mathcal N}^\T \widetilde{u}_s\| + \|U_{\mathcal J}(I- U_{\mathcal{J}}^\T \Phi(\widetilde{\lambda}_s) U_{\mathcal J}\Lambda_{\mathcal J})^{-1} U_{\mathcal J}^\T\Xi(\widetilde\lambda_s)U\cdot \Lambda U^\T \widetilde{u}_s \|\\
    &\le 4 \frac{\|\mathcal{V}_{\mathcal{N}\mathcal{K}}\|}{\lambda_k^2} + 24\frac{R_{\max}}{{\lambda}_k^4} \osc(\mathcal R) + 124\frac{M}{\delta_k}.
\end{align*}
Denote $$\Delta^{\orc} := (w_1^{\orc} - U_kU_k^\T\widetilde u_1, \ldots, w_k^{\orc} - U_kU_k^\T\widetilde u_k).$$It follows immediately that 
\[ \|\Delta^{\orc}\|\le \|\Delta^{\orc}\|_F \le \sqrt{k}\left( 4 \frac{\|\mathcal{V}_{\mathcal{N}\mathcal{K}}\|}{\lambda_k^2} + 24\frac{R_{\max}}{{\lambda}_k^4} \osc(\mathcal R) + 124\frac{M}{\delta_k}\right)=\mathcal{E}_{k}^{\osc}. \]

We rewrite
\[W_k^{\orc}=(w_1^{\orc},\ldots,w_k^{\orc})= U_k U_k^\T \widetilde U_k + \Delta^{\orc}.\]
Note that 
\begin{align}\label{eq:Uk-orc}
 U_k^{\orc}=\orth(W_k^{\orc}) = W_k^{\orc} ((W_k^{\orc})^\T W_k^{\orc})^{-1/2}.   
\end{align}
We aim to bound 
\[\|\sin\angle(U_k, U_k^{\orc})\| = \|P_{U_k^\perp} U_k^{\orc}\| = \|P_{U_k^\perp}W_k^{\orc} ((W_k^{\orc})^\T W_k^{\orc})^{-1/2}\| \le \frac{\|P_{U_k^\perp}W_k^{\orc}\|}{s_{\min}(W_k^{\orc})}.\]
First,
\[\|P_{U_k^\perp}W_k^{\orc}\| = \|P_{U_k^\perp}(U_k U_k^\T \widetilde U_k + \Delta^{\orc})\| = \|P_{U_k^\perp} \Delta^{\orc}\| \le \|\Delta^{\orc}\|. \]
Next, since $W_k^{\orc} = U_k U_k^\T \widetilde U_k + \Delta^{\orc}$, and $\Delta^{\orc}$ is orthogonal to $U_k$, we have
\begin{align*}
\sigma_{\min}(W_k^{\orc}) &\ge \sigma_{\min}(U_k U_k^\T \widetilde U_k) = \sigma_{\min}(U_k^\T \widetilde U_k) = \sqrt{1-\|\sin\angle(U_k, \widetilde U_k)\|^2}.
\end{align*}
On the event that Theorem \ref{thm:top-k-eigenspace} holds, $\|\sin\angle(U_k, \widetilde U_k)\|\le 1/100$. Thus, \[\sigma_{\min}(W_k^{\orc}) \ge \frac{1}{2}.\]
We conclude that 
\[ \|\sin\angle(U_k, U_k^{\orc})\| \le 2 \|\Delta^{\orc}\| \le 2\mathcal{E}_{k}^{\osc}.
\]

For the $2\to \infty$ norm bound, we start with the standard decomposition
\begin{align*}
    \min_{O\in\mathbb O(k)}\|U_k^{\orc}-U_kO\|_{2,\infty} \le \|(I-P_{U_k})U_k^{\orc}\|_{2,\infty} + \|U_k\|_{2,\infty}\|\sin\angle(U_k,U_k^{\orc})\|^2.
\end{align*}
It remains to bound $\|(I-P_{U_k})U_k^{\orc}\|_{2,\infty}$. Since $P_{U_k^\perp} W_k^{\orc} = \Delta^{\orc}$, from \eqref{eq:Uk-orc}, we get
\[
(I-P_{U_k})U_k^{\orc} = \Delta^{\orc} ((W_k^{\orc})^\T W_k^{\orc})^{-1/2}.
\]
Thus
\begin{align}\label{eq:2toinf_projection_orc}
\|(I-P_{U_k})U_k^{\orc}\|_{2,\infty} &\le \|\Delta^{\orc}\|_{2,\infty} \|((W_k^{\orc})^\T W_k^{\orc})^{-1/2}\| \\
&\le \frac{\|\Delta^{\orc}\|_{2,\infty}}{\sigma_{\min}(W_k^{\orc})}\le 2\|\Delta^{\orc}\|_{2,\infty}.
\end{align}
We bound $\|\Delta^{\orc}\|_{2,\infty}$. From \eqref{eq:w-orc-01}, 
\begin{align*}
    \Delta^{\orc} &= (w_1^{\orc} - U_kU_k^\T\widetilde u_1, \ldots, w_k^{\orc} - U_kU_k^\T\widetilde u_k)\\
    &= U_{\mathcal N}\cdot U_{\mathcal N}^\T \widetilde U_k + U_{\mathcal J}\cdot \mathtt B
\end{align*}
where $\mathtt B:=[\mathtt{b_1,\cdots,b_k}]$ with columns 
\[\mathtt{b}_s=(I- U_{\mathcal{J}}^\T \Phi(\widetilde{\lambda}_s) U_{\mathcal J}\Lambda_{\mathcal J})^{-1} U_{\mathcal J}^\T\Xi(\widetilde\lambda_s)U\cdot \Lambda U^\T \widetilde{u}_s.\]
In \eqref{eq:0601-error-bs}, we have 
\[\|\mathtt b_s\| \le 120\frac{M}{\delta_k}\]
and thus 
\[ \|\mathtt B\|\le \|\mathtt B\|_F \le 120\sqrt{k}\frac{M}{\delta_k}.\]
Also, by \eqref{eq:0601-UN-later}, we get
\begin{align*}
\| U_{\mathcal N}^\T \widetilde U_k \| \le \| U_{\mathcal N}^\T \widetilde U_k \|_F &\le \sqrt{k}\max_{s\in k}\|U_{\mathcal N}^\T \widetilde{u}_s\|\\
&\le 4 \sqrt{k}\frac{\|\mathcal{V}_{\mathcal{N}\mathcal{K}}\|}{\lambda_k^2} +24\sqrt{k} \frac{R_{\max}}{{\lambda}_k^4} \osc(\mathcal R)+ 4 \sqrt{k}\frac{M}{\lambda_k}.
\end{align*}
Now we get
\begin{align*}
  \|(I-P_{U_k})U_k^{\orc}\|_{2,\infty} &\le 2\|\Delta^{\orc}\|_{2,\infty}=2\max_{1\le i \le n} \|e_i^\T(U_{\mathcal N}\cdot U_{\mathcal N}^\T \widetilde U_k + U_{\mathcal J}\cdot \mathtt B)\|\\
  &\le 2 \|U\|_{2,\infty} (\| U_{\mathcal N}^\T \widetilde U_k \| + \|\mathtt B\|)\\
  &\le \|U\|_{2,\infty} \left(8\sqrt{k} \frac{\|\mathcal{V}_{\mathcal{N}\mathcal{K}}\|}{\lambda_k^2} + 48\sqrt{k}\frac{R_{\max}}{{\lambda}_k^4} \osc(\mathcal R) + 248\sqrt{k}\frac{M}{\delta_k} \right)\\
  &= 2\|U\|_{2,\infty} \mathcal{E}_{k}^{\osc}.
\end{align*}
Finally, combining the bound on $\|\sin\angle(U_k,U_k^{\orc})\|$, we obtain
\begin{align*}
    \min_{O\in\mathbb O(k)}\|U_k^{\orc}-U_kO\|_{2,\infty} &\le \|(I-P_{U_k})U_k^{\orc}\|_{2,\infty} + \|U_k\|_{2,\infty}\|\sin\angle(U_k,U_k^{\orc})\|^2\\
    &\le 3\|U\|_{2,\infty} \mathcal{E}_{k}^{\osc}.
\end{align*}
The proof is complete. 

\subsection{Proof of (\ref{eq:plugin-sin-short}) and (\ref{eq:plugin-2inf-short})}\label{app:0603-db-bound}
By the definitions of $\widehat{w}_s$ and $w_s^{\orc}$, their difference is
\[
    \widehat w_s - w_s^{\orc} = P_{\widetilde U}\widetilde u_s - P_U\widetilde u_s - \left( \widetilde U_{\mathcal J} \widehat{\mathcal T}_s \widetilde U_{\mathcal I}^{\T}\widetilde u_s - U_{\mathcal J}\mathcal T_s U_{\mathcal I}^{\T}\widetilde u_s \right).
\]
Since $P_{\widetilde U}\widetilde u_s = \widetilde u_s$, we get
\[
    \widehat W_k - W_k^{\orc} = (I-P_U)\widetilde U_k - \mathfrak D_k.
\]
Let $\Delta := \|\widehat W_k - W_k^{\orc}\|$. Applying the triangle inequality, we bound \begin{equation}\label{eq:Delta_short_bound}
    \Delta \le \|(I-P_U)\widetilde U_k\| + \|\mathfrak D_k\| \le 2\frac{\|E\|}{\lambda_k} + \|\mathfrak D_k\|.
\end{equation}
The bound in the second inequality was established in \eqref{eq:0603-proj-null}.

Similarly, taking the $2\to\infty$ norm yields
\begin{equation}\label{eq:Delta_2inf_short_bound}
    \Delta_{2,\infty} := \|\widehat W_k - W_k^{\orc}\|_{2,\infty} \le \|(I-P_U)\widetilde U_k\|_{2,\infty} + \|\mathfrak D_k\|_{2,\infty}.
\end{equation}

Note that the assumption $\mathcal E_k^{\orc} + \frac{\|E\|}{\lambda_k} + \|\mathfrak D_k\| \le c$ guarantees that $\mathcal E_k^{\orc} + \Delta \le c \le 1/10$ for a sufficiently small absolute constant $c$.  

Recall from the proof of  Proposition \ref{prop:oracle-block-debias} that $W_k^{\orc} = U_k U_k^\T \widetilde U_k + \Delta^{\orc}$, where  $\|P_{U_k^\perp} W_k^{\orc}\| \le \|\Delta^{\orc}\| \le \mathcal{E}_k^{\orc}$ and  $\|P_{U_k^\perp} W_k^{\orc}\|_{2,\infty} \le \|U\|_{2,\infty}\mathcal{E}_k^{\orc}$. 

From the decomposition $\widehat W_k = W_k^{\orc} + (\widehat W_k - W_k^{\orc})$, we obtain
\begin{equation}\label{eq:W_hat_projection}
    \|P_{U_k^\perp} \widehat W_k\| \le \|P_{U_k^\perp} W_k^{\orc}\| + \|\widehat W_k - W_k^{\orc}\| \le \mathcal{E}_k^{\orc} + \Delta.
\end{equation}
Similarly, the $2\to\infty$ norm satisfies
\begin{equation}\label{eq:W_hat_projection_2inf}
    \|P_{U_k^\perp} \widehat W_k\|_{2,\infty} \le \|P_{U_k^\perp} W_k^{\orc}\|_{2,\infty} + \|\widehat W_k - W_k^{\orc}\|_{2,\infty} \le \|U\|_{2,\infty}\mathcal{E}_k^{\orc} + \Delta_{2,\infty}.
\end{equation}

Next, we lower-bound $s_{\min}(\widehat W_k)$. We established previously that $s_{\min}(W_k^{\orc}) \ge 1/2$. By Weyl's inequality, 
\[
    s_{\min}(\widehat W_k) \ge s_{\min}(W_k^{\orc}) - \|\widehat W_k - W_k^{\orc}\| \ge \frac{1}{2} - \Delta \ge \frac{1}{2} - \frac{1}{10} = \frac{2}{5}
\]
since $\Delta \le 1/10$.

Similar to the proof of  Proposition \ref{prop:oracle-block-debias}, we get
\[
    \|\sin\angle(U_k, \widehat U_k^{\db})\| = \|P_{U_k^\perp} \widehat U_k^{\db}\| \le \frac{\|P_{U_k^\perp} \widehat W_k\|}{\sigma_{\min}(\widehat W_k)} \le \frac{5}{2}(\mathcal{E}_k^{\orc} + \Delta).
\]

For the $2 \to \infty$ bound, we start with 
\[
    \min_{O\in\mathbb O(k)}\|\widehat U_k^{\db}-U_kO\|_{2,\infty} \le \|(I-P_{U_k})\widehat U_k^{\db}\|_{2,\infty} + \|U_k\|_{2,\infty}\|\sin\angle(U_k,\widehat U_k^{\db})\|^2.
\]
Since $(I-P_{U_k})\widehat U_k^{\db} = P_{U_k^\perp} \widehat W_k (\widehat W_k^\T \widehat W_k)^{-1/2}$, we bound the first term using \eqref{eq:W_hat_projection_2inf}:
\[
    \|(I-P_{U_k})\widehat U_k^{\db}\|_{2,\infty} \le \frac{\|P_{U_k^\perp} \widehat W_k\|_{2,\infty}}{\sigma_{\min}(\widehat W_k)} \le \frac{5}{2}\left(\|U\|_{2,\infty}\mathcal{E}_k^{\orc} + \Delta_{2,\infty}\right).
\]
Substituting this and the bound on $\|\sin\angle(U_k, \widehat U_k^{\db})\|^2$ back into the decomposition yields the desired bound.


\section{Proofs of Propositions~\ref{prop:refined-top-k} and \ref{prop:refined-weighted-embedding}}\label{app:proof-embedding-approx}
All proofs are carried out on the event $\mathcal E_D$, as mentioned at the beginning of Section~\ref{sec:top-k-eigenspace}.

\subsection{Proof of Proposition~\ref{prop:refined-top-k}}
We first record the projected form of Corollary~\ref{coro:row-adaptive}. Let $Q=I-UU^\T$. Since
\[
    Q(G-\Phi)U  = (G-\Phi)U - U\bigl[U^\T(G-\Phi)U\bigr],
\]
we have 
\[
    \|Q(G-\Phi)U\|_{2,\infty} \le \|(G-\Phi)U\|_{2,\infty} +  \|U\|_{2,\infty} \|U^\T(G-\Phi)U\|.
\]
Denote $\Xi(z)=G(z)-\Phi(z)$. Applying Corollary~\ref{coro:row-adaptive} with $B=U$ to the first
term and Corollary~\ref{lem:local-law-uniform} to the second term, we obtain that with probability at least
$1-n^{-D}$,
\begin{equation}
    \sup_{z\in\mathcal S_{\rm ra}} |z|^2\|Q\Xi(z) U\|_{2,\infty} \le C\mathfrak R_D(U).
    \label{eq:projected-row-adaptive}
\end{equation}
Here we used that $\mathcal S_{\rm ra}\subseteq\mathcal S_{\rm out}$ and $\mathfrak m_D \le M_D$.

On the event $\|E\|\le3\sqrt{R_{\max}}$, by Weyl's inequality and our assumption on $\lambda_k$,  we have 
$$ \widetilde\lambda_s  \ge \lambda_k-3\sqrt{R_{\max}}  \ge C_{\rm ra}\sqrt{(D+r)R_{\max}\log n}.$$
Thus,
$\widetilde\lambda_s \in \mathcal S_{\rm ra}$ for all $s\le k$. Hence \eqref{eq:projected-row-adaptive} applies at $z=\widetilde\lambda_s$ for $s\in [k]$.

We now follow the proof of Theorem \ref{thm:top-k-eigenspace}. The only estimate that changes is the stochastic null-space term
\begin{align*}
\sqrt{\sum_{s=1}^k \widetilde{\lambda}_s^2 \left\|Q \Xi(\widetilde{\lambda}_s) U \right\|_{2,\infty}^2}
\end{align*}
on the right-hand side of \eqref{eq:row-proj-det}. By \eqref{eq:projected-row-adaptive}, we bound it by 
\begin{align*}
 \sqrt{\sum_{s=1}^k \widetilde{\lambda}_s^2 \left\|Q \Xi(\widetilde{\lambda}_s) U \right\|_{2,\infty}^2} \le \sqrt{\sum_{s=1}^k \frac{C^2 \mathfrak R_D(U)^2}{\widetilde{\lambda}_s^2}} \le C\sqrt{k} \frac{\mathfrak R_D(U)}{\lambda_k},
\end{align*}
where the last inequality follows from $\widetilde\lambda_s\ge c\lambda_k$.
All remaining estimates are unchanged from the proof of Theorem~\ref{thm:top-k-eigenspace}, which proves \eqref{eq:refined-top-k-row}. 

 \subsection{Proof of Proposition~\ref{prop:refined-weighted-embedding}}
Recall that
\[ \widetilde{Y}_k := \widetilde{U}_k \widetilde{\Lambda}_k^{1/2}, \qquad Y_k := U_k \Lambda_k^{1/2}.\]
We start with a deterministic decomposition. By inserting the projection operators $P_{U_k} = U_k U_k^\T$ and $P_{U_k^\perp} = I - P_{U_k}$, we have for any $O\in \mathbb O(k)$,
\begin{align*}
 \widetilde{Y}_k O -Y_k = P_{U_k^\perp} \widetilde{Y}_k O + U_k (U_k^\T	 \widetilde{Y}_k O - \Lambda_k^{1/2}).
\end{align*}
Thus, 
\begin{align}\label{eq:0624-wt-2}
\| \widetilde U_k\widetilde\Lambda_k^{1/2} O- U_k\Lambda_k^{1/2} \|_{2,\infty} \le \|P_{U_k^\perp} \widetilde{Y}_k \|_{2,\infty} + \|U_k\|_{2,\infty} \| U_k^\T	\widetilde{Y}_k  O -\Lambda_k^{1/2}\|.
\end{align}

We first control the second term on the right-hand side of \eqref{eq:0624-wt-2} by choosing a specific $O \in \mathbb{O}(k)$.  For simplicity, denote $Z:=U_k^\T \widetilde U_k\widetilde\Lambda_k^{1/2}= U_k^\T \widetilde{Y}_k $. Consider the polar decomposition $Z= H R$, where $H = (Z Z^\T)^{1/2}$ is a symmetric positive semi-definite matrix and $R \in \mathbb{O}(k)$. We set $O=R^{\T}$. Thus,
\[
ZO = HR R^\T = H.
\]
Consequently, 
\[
\| U_k^\T	\widetilde{Y}_k  O -\Lambda_k^{1/2}\| =\|ZO-\Lambda_k^{1/2}\| = \|H-\Lambda_k^{1/2}\|.
\]
To bound this difference, noting that both $H$ and $\Lambda_k$ are positive semi-definite, we apply Theorem 6.2 from \cite{Higham} to get
\[
\|H-\Lambda_k^{1/2}\|
    \le \frac{\|H^2-\Lambda_k\|}{\lambda_{\min}(H) + \lambda_{\min}(\Lambda_k^{1/2})} 
    \le \frac{\|H^2-\Lambda_k\|}{\sqrt{\lambda_k}},
\]
where the final inequality follows from the fact that $\lambda_{\min}(H) \ge 0$. 

Now we turn to bound $\|H^2-\Lambda_k\|$. It follows from the definition of $H$ that
\[
 H^2 = U_k^\T \big(\widetilde U_k\widetilde\Lambda_k \widetilde U_k^\T\big) U_k. 
 \]
Since $\widetilde{A}_k = \widetilde{U}_k \widetilde{\Lambda}_k \widetilde{U}_k^\T$, we also have $H^2 = U_k^\T \widetilde{A}_k U_k$. 

Plugging in  $\Lambda_k = U_k^\T A U_k$, the difference can be written as
\begin{align*}
    H^2 - \Lambda_k = U_k^\T (\widetilde{A}_k - A) U_k.
\end{align*}
Let $\widetilde U_{>k}\in \mathbb R^{n\times(n-k)}$ be an
orthonormal basis of $\Span(\widetilde U_k)^\perp$,
formed by the remaining eigenvectors of $\widetilde A$, and let
$\widetilde\Lambda_{>k}$ be the corresponding diagonal matrix of
eigenvalues. Then
\[
\widetilde A
= \widetilde U_k\widetilde\Lambda_k\widetilde U_k^\T +\widetilde U_{>k}\widetilde\Lambda_{>k}\widetilde U_{>k}^\T:= \widetilde{A}_k + \widetilde{A}_{>k}.
\]
From $\widetilde{A} = A + E$ and the decomposition $\widetilde{A} = \widetilde{A}_k + \widetilde{A}_{>k}$, we substitute $\widetilde{A}_k - A = E - \widetilde{A}_{>k}$ to obtain
\begin{align*}
    H^2 - \Lambda_k = U_k^\T E U_k - U_k^\T \widetilde{A}_{>k} U_k.
\end{align*}
Applying the triangle inequality yields
\begin{align*}
     \|H^2 - \Lambda_k\|  \le \|U_k^\T E U_k\| + \|U_k^\T \widetilde{A}_{>k} U_k\|.
\end{align*}
For the second term, observe that $U_k^T \widetilde{A}_{>k} U_k = (U_k^\T \widetilde{U}_{>k}) \widetilde{\Lambda}_{>k} (\widetilde{U}_{>k}^\T U_k)$. Taking the operator norm, we get
\begin{align*}
    \|U_k^\T \widetilde{A}_{>k} U_k\| \le \|\widetilde{\Lambda}_{>k}\|\cdot \|U_k^\T \widetilde{U}_{>k}\|^2.
\end{align*}
By Weyl's inequality,  $\|\widetilde{\Lambda}_{>k}\|\le \|{\Lambda}_{>k}\| + \|E\|$. Furthermore, $\|U_k^\T \widetilde{U}_{>k}\| = \|\sin\angle(U_k, \widetilde{U}_k)\|$. Substituting these bounds back into the triangle inequality yields 
\begin{align}\label{eq:0625-bd-diagonal}
   \|H^2 - \Lambda_k\| \le \|U_k^\T E U_k\| + (\|{\Lambda}_{>k}\| + \|E\|) \|\sin\angle(U_k, \widetilde{U}_k)\|^2.
\end{align}
Finally, we obtain
\begin{align*}
\| U_k^\T	\widetilde U_k\widetilde\Lambda_k^{1/2} O -\Lambda_k^{1/2}\| \le \frac{1}{\sqrt{\lambda_k}}\left( \|U_k^\T E U_k\| + ( \|{\Lambda}_{>k}\| + \|E\|) \|\sin\angle(U_k, \widetilde{U}_k)\|^2 \right).
\end{align*}

Next, we bound $\|P_{U_k^\perp} \widetilde{Y}_k \|_{2,\infty}=\|P_{U_k^\perp} \widetilde{U}_k \widetilde{\Lambda}_k^{1/2} \|_{2,\infty}$. First, we obtain a deterministic bound by directly applying the component-wise estimates established in the proof of Proposition \ref{prop:spacepert}. Recall the decomposition $P_{U_k^\perp} = P_J + Q$, where $P_J = U_J U_J^\T$ with $U_J = (u_{k+1},\dots,u_r)$, and $Q = I - U U^\T$. By triangle inequality, we get
\begin{align*}
\| P_{U_k^\perp} \widetilde{Y}_k \|_{2,\infty} &\le \| P_J \widetilde{U}_k \widetilde{\Lambda}_k^{1/2}\|_{2,\infty} + \|Q \widetilde{U}_k \widetilde{\Lambda}_k^{1/2}\|_{2,\infty}\\
&\le \|U\|_{2,\infty} \cdot \|U_J^\T \widetilde{U}_k \widetilde{\Lambda}_k^{1/2}\|+ \|Q \widetilde{U}_k \widetilde{\Lambda}_k^{1/2}\|_{2,\infty}\\
&\le \|U\|_{2,\infty} \cdot \|U_J^\T \widetilde{U}_k \widetilde{\Lambda}_k^{1/2}\|_F+ \|Q \widetilde{U}_k \widetilde{\Lambda}_k^{1/2}\|_{2,\infty}.
\end{align*}
Now we evaluate the rows of $Q \widetilde{U}_k \widetilde{\Lambda}_k^{1/2}$. The squared Euclidean norm of the $l$-th row is given by $\sum_{s=1}^k \widetilde{\lambda}_s (e_l^\T Q \widetilde{u}_s)^2$. In the proof of Proposition \ref{prop:spacepert}, we established the resolvent decomposition $e_l^\top Q \widetilde{u}_s = \mathtt{q}_{ls}^{(1)} + \mathtt{q}_{ls}^{(2)}$, where the individual terms are deterministically bounded by
\begin{align*}
|\mathtt{q}_{ls}^{(1)}| \le \frac{8}{3}\|U\|_{2,\infty}\frac{\osc(\mathcal R)}{\widetilde\lambda_s^2} \qquad \text{and} \qquad |\mathtt{q}_{ls}^{(2)}| \le 2\widetilde{\lambda}_s \| e_l^\T Q \Xi(\widetilde{\lambda}_s) U\|_2.
\end{align*}
Next, we multiply these squared bounds by the weight $\widetilde{\lambda}_s$ and sum over $s \in [k]$ to get
\begin{align*}
\sqrt{\sum_{s=1}^k \widetilde{\lambda}_s |\mathtt{q}_{ls}^{(1)}|^2}
\le \frac{8}{3} \|U\|_{2,\infty} \osc(\mathcal R)
\sqrt{  \sum_{s=1}^k \frac{1}{\widetilde\lambda_s^3} } \le 3\sqrt{k}
\|U\|_{2,\infty} \frac{\osc(\mathcal R)}{\lambda_k^{3/2}}, 
\end{align*}
and
\begin{align*}
\sqrt{\sum_{s=1}^k \widetilde{\lambda}_s |\mathtt{q}_{ls}^{(2)}|^2}
\le 2 \sqrt{\sum_{s=1}^k \widetilde{\lambda}_s^3 \|e_l^\T Q \Xi(\widetilde{\lambda}_s) U \|_2^2}.
\end{align*}

Combining these estimates via the triangle inequality yields the final deterministic bound for the weighted orthogonal projection:
\begin{align}\label{eq:row-proj-det-weighted}
\| P_{U_k^\perp} \widetilde{Y}_k\|_{2,\infty} &\le \| U\|_{2,\infty} \|U_J^\T \widetilde{U}_k \widetilde{\Lambda}_k^{1/2}\|_F + 3\sqrt{k}
\|U\|_{2,\infty} \frac{\osc(\mathcal R)}{\lambda_k^{3/2}} \nonumber\\
&\quad + 2 \sqrt{\sum_{s=1}^k \widetilde{\lambda}_s^3 \max_{1\le i\le n} \| e_i^\T Q \Xi(\widetilde{\lambda}_s) U \|_2^2}.
\end{align}   

It remains to bound 
\[
\|U_J^\T \widetilde{U}_k \widetilde{\Lambda}_k^{1/2}\|_F = \sqrt{\sum_{s=1}^k \widetilde{\lambda}_s \|U_J^\T \widetilde{u}_s\|^2}.
\] 
For $J\neq \emptyset$, from the bound \eqref{eq:0506-individual}, we get
\begin{align*}
   \sqrt{\widetilde{\lambda}_s}\|U_J^\T \widetilde{u}_s\|  \le  &\frac{4}{\Delta_{\mathcal J}^{\Phi}(\widetilde{\lambda}_s)}\left(\frac{\|\mathcal{V}_{\mathcal{J}\mathcal{I}}\|}{\widetilde{\lambda}_s^{1/2}} + \frac{13}{2}\frac{R_{\max}}{\widetilde{\lambda}_s^{5/2}}\osc(\mathcal{R}) + \widetilde{\lambda}_s^{5/2} \|U^\T \Xi(\widetilde{\lambda}_s)U\| \right) \nonumber \\
   &+ 2\sqrt{3} \left( \frac{\|\mathcal{V}_{\mathcal{N}\mathcal{K}}\|}{\widetilde{\lambda}_s^{3/2}} + \frac{13}{2}\frac{R_{\max}}{\widetilde{\lambda}_s^{7/2}} \osc(\mathcal R) + \widetilde{\lambda}_s^{3/2}\|U^\T \Xi(\widetilde{\lambda}_s)U\|\right)
\end{align*} 
Substituting $\Delta_{\mathcal J}^{\Phi}(\widetilde{\lambda}_s) \ge \frac{1}{2}\delta_k$ and $\|U^\top \Xi(\widetilde{\lambda}_s)U\| \le M/\widetilde{\lambda}_s^2$ into the above bound, we obtain
\begin{align*}
    \sqrt{\widetilde{\lambda}_s}\|U_J^\top \widetilde{u}_s\|  
    &\le \frac{8}{\delta_k} \frac{\|\mathcal{V}_{\mathcal{J}\mathcal{I}}\|}{\widetilde{\lambda}_s^{1/2}} 
    + \frac{52}{\delta_k} \frac{R_{\max}}{\widetilde{\lambda}_s^{5/2}}\osc(\mathcal{R}) 
    + \frac{8 M \widetilde{\lambda}_s^{1/2}}{\delta_k} \\
    &\quad + 2\sqrt{3} \frac{\|\mathcal{V}_{\mathcal{N}\mathcal{K}}\|}{\widetilde{\lambda}_s^{3/2}} 
    + 13\sqrt{3} \frac{R_{\max}}{\widetilde{\lambda}_s^{7/2}} \osc(\mathcal R) 
    + 2\sqrt{3} \frac{M}{\widetilde{\lambda}_s^{1/2}}.
\end{align*}

To compute the Frobenius norm $\|U_J^\top \widetilde{U}_k \widetilde{\Lambda}_k^{1/2}\|_F = \big\| \big( \sqrt{\widetilde{\lambda}_s}\|U_J^\top \widetilde{u}_s\| \big)_{s=1}^k \big\|_2$, we  sum over all $k$ components and use $\frac{3}{2} \lambda_s \ge \widetilde{\lambda}_s \ge \frac{99}{100}\lambda_k$, which yields
\begin{align*}
  &  \|U_J^\top \widetilde{U}_k \widetilde{\Lambda}_k^{1/2}\|_F \\
    &\quad \le \frac{8.1\sqrt{k}}{\delta_k \lambda_k^{1/2}} \|\mathcal{V}_{\mathcal{J}\mathcal{I}}\| + \frac{3.6\sqrt{k}}{\lambda_k^{3/2}} \|\mathcal{V}_{\mathcal{N}\mathcal{K}}\|  + 78\sqrt{k} \frac{R_{\max} \osc(\mathcal{R})}{\delta_k \lambda_k^{5/2}} + \frac{10 M}{\delta_k} \sqrt{\sum_{s=1}^k {\lambda}_s} + \frac{3.5\sqrt{k}}{\lambda_k^{1/2}} M.
\end{align*}

By \eqref{eq:projected-row-adaptive}, for every $s\le k$,
\[
    \max_{1\le i\le n} \|e_i^\T Q\Xi(\widetilde\lambda_s)U\|_2 \le C\frac{\mathfrak R_D(U)}  {\widetilde\lambda_s^2}.
\]

Thus,
\begin{align*}
2 \sqrt{\sum_{s=1}^k \widetilde{\lambda}_s^3 \max_{1\le i\le n} \| e_i^\T Q \Xi(\widetilde{\lambda}_s) U \|_2^2} \le C\mathfrak R_D(U) \sqrt{\sum_{s=1}^k \frac{1}{\widetilde{\lambda}_s}} \le C\sqrt{k}\frac{\mathfrak R_D(U) }{\sqrt{\widetilde{\lambda}_k}} < C\sqrt{k}\frac{\mathfrak R_D(U) }{\sqrt{{\lambda}_k}}.
\end{align*}

Finally, combining these estimates into \eqref{eq:row-proj-det-weighted} and noting that $\frac{13}{2}\widetilde{\lambda}_k^{-3/2} \le 7\lambda_k^{-3/2}$, we obtain 
\begin{align*}
    \| P_{U_k^\perp} \widetilde{Y}_k\|_{2,\infty} 
    &\le \|U\|_{2,\infty} \left( \frac{8.1\sqrt{k}}{\delta_k \lambda_k^{1/2}} \|\mathcal{V}_{\mathcal{J}\mathcal{I}}\|  + \frac{3.6\sqrt{k}}{\lambda_k^{3/2}} \|\mathcal{V}_{\mathcal{N}\mathcal{K}}\|  + 78\sqrt{k} \frac{R_{\max} \osc(\mathcal{R})}{\delta_k \lambda_k^{5/2}} \right)\mathbf{1}_{\{ k<r\}} \\
    & + \|U\|_{2,\infty} \left( \frac{10 M}{\delta_k} \sqrt{\sum_{s=1}^k {\lambda}_s} + \frac{3.5\sqrt{k}}{\lambda_k^{1/2}} M \right)\mathbf{1}_{\{ k<r\}}  + 3\sqrt{k}
\|U\|_{2,\infty} \frac{\osc(\mathcal R)}{\lambda_k^{3/2}} + C\sqrt{k}\frac{\mathfrak R_D(U)} {\lambda_k^{1/2}}.
 \end{align*}   
Recall that 
\[\mathcal{B}_k:=\frac{10\sqrt{k}}{\delta_k \lambda_k}\left( \|\mathcal{V}_{\mathcal{J}\mathcal{I}}\|+8\frac{R_{\max}}{\lambda_k^2}\osc(\mathcal R)\right) +4 \sqrt{k}\frac{\|\mathcal{V}_{\mathcal{N}\mathcal{K}}\|}{\lambda_k^2} .\]
We simplify the bound to
\begin{align}\label{eq:0625-2nd-bd}    
   \| P_{U_k^\perp} \widetilde{Y}_k\|_{2,\infty}  \le & \|U\|_{2,\infty} \sqrt{\lambda_k } \mathcal{B}_k \mathbf{1}_{\{ k<r\}}+ \|U\|_{2,\infty}  \frac{10 M}{\delta_k} \sqrt{\sum_{s=1}^k {\lambda}_s} \mathbf{1}_{\{ k<r\}}\nonumber \\
   &+ C\sqrt{k}\frac{\mathfrak R_D(U)} {\lambda_k^{1/2}} +3\sqrt{k}
\|U\|_{2,\infty} \frac{\osc(\mathcal R)}{\lambda_k^{3/2}}.
\end{align}

Now, from \eqref{eq:0624-wt-2}, we get
\begin{align*}
&\| \widetilde U_k\widetilde\Lambda_k^{1/2} O- U_k\Lambda_k^{1/2} \|_{2,\infty} \\
&\qquad \le \|P_{U_k^\perp} \widetilde{Y}_k \|_{2,\infty} +  \frac{\|U_k\|_{2,\infty}}{\sqrt{\lambda_k}}\left[ \|U_k^\T E U_k\| + ( \|{\Lambda}_{>k}\| + \|E\|) \|\sin\angle(U_k, \widetilde{U}_k)\|^2 \right],
\end{align*}
where $\|P_{U_k^\perp} \widetilde{Y}_k \|_{2,\infty}$ is bounded by \eqref{eq:0625-2nd-bd}. Applying the bound for $\|\sin\angle(U_k, \widetilde{U}_k)\|$ from Theorem~\ref{thm:top-k-eigenspace} completes the proof.

\section{Plug-in de-biasing implementation and limitations}\label{app:plug-in}  To keep the discussion simple, assume in this paragraph that $r$ and $r_+$ are known and the variance profile $\Sigma$ is known, so that
$\Phi(z)$ is computable from the QVE.. 

Define
\[
\widetilde U_{\mathcal J}
:=
(\widetilde u_{k+1},\ldots,\widetilde u_{r_+}),
\qquad
\widetilde U_{\mathcal I}
:=
(\widetilde u_1,\ldots,\widetilde u_k,
\widetilde u_{r_++1},\ldots,\widetilde u_r),
\]
with corresponding diagonal eigenvalue matrices $\widetilde\Lambda_{\mathcal J}$
and $\widetilde\Lambda_{\mathcal I}$. Also write
\[
\widetilde U:=(\widetilde u_1,\ldots,\widetilde u_r),
\qquad
P_{\widetilde U}:=\widetilde U\widetilde U^\T.
\]
If $\mathcal J=\emptyset$, no correction is needed and we  set $\widehat U_k^{\db}:=\widetilde U_k$. Otherwise, for $s\in [k]$, define
\begin{align}
\label{def:That_s-knownS}
\widehat{\mathcal T}_s :=
\left( I- \widetilde U_{\mathcal J}^{\T}
\Phi(\widetilde\lambda_s)
\widetilde U_{\mathcal J}
\widetilde\Lambda_{\mathcal J}
\right)^{-1}
\widetilde U_{\mathcal J}^{\T}
\Phi(\widetilde\lambda_s)
\widetilde U_{\mathcal I}
\widetilde\Lambda_{\mathcal I},
\end{align}
whenever the inverse exists. The plug-in corrected vector is
\begin{equation}
\label{eq:def-w-plugin-knownS}
\widehat w_s :=
P_{\widetilde U}\widetilde u_s -
\widetilde U_{\mathcal J}
\widehat{\mathcal T}_s
\widetilde U_{\mathcal I}^{\T}
\widetilde u_s,
\qquad s=1,\ldots,k.
\end{equation}
Note that for $s\le k$, $P_{\widetilde U}\widetilde u_s=\widetilde u_s$. Finally, set
\[
\widehat W_k:=(\widehat w_1,\ldots,\widehat w_k),
\qquad
\widehat U_k^{\db}:=\orth(\widehat W_k).
\]
To compare this estimator with the oracle estimator, define the correction discrepancy
\[
\mathfrak D_k
:=
\left(
\widetilde U_{\mathcal J}
\widehat{\mathcal T}_s
\widetilde U_{\mathcal I}^{\T}\widetilde u_s
-
U_{\mathcal J}
\mathcal T_s
U_{\mathcal I}^{\T}\widetilde u_s
\right)_{s=1}^k.
\]
Then the difference between the plug-in and oracle matrices is given by
\begin{align}\label{eq:0603-Wk-diff}
\widehat W_k-W_k^{\orc}
=
(I-P_U)\widetilde U_k-\mathfrak D_k.
\end{align}
Consequently, the accuracy of the plug-in correction is governed by the norms of this difference. As shown in \eqref{eq:0603-proj-null} below, $\|(I-P_U)\widetilde U_k\| \le 2\|E\|/\lambda_k$. On the event that
\[
\mathcal E_k^{\orc}
+
\frac{\|E\|}{\lambda_k}
+
\|\mathfrak D_k\|
\le c
\]
for a sufficiently small absolute constant $c>0$, the same
argument as in Proposition~\ref{prop:oracle-block-debias}
gives
\begin{equation}
\label{eq:plugin-sin-short}
\|\sin\angle(\widehat U_k^{\db}, U_k)\|
\le
C\left(
\mathcal E_k^{\orc}
+
\frac{\|E\|}{\lambda_k}
+
\|\mathfrak D_k\|
\right).
\end{equation}
Moreover,
\begin{align}
\label{eq:plugin-2inf-short}
\min_{O\in\mathbb O(k)}
\|\widehat U_k^{\db}-U_kO\|_{2,\infty}
\le C\Big[
&
\|U\|_{2,\infty}\mathcal E_k^{\orc}
+
\|(I-P_U)\widetilde U_k\|_{2,\infty}
+
\|\mathfrak D_k\|_{2,\infty}
\nonumber\\
&+
\|U_k\|_{2,\infty}
\left(
\mathcal E_k^{\orc}
+
\frac{\|E\|}{\lambda_k}
+
\|\mathfrak D_k\|
\right)^2
\Big].
\end{align}
The proofs of \eqref{eq:plugin-sin-short} and \eqref{eq:plugin-2inf-short} are given in Appendix \ref{app:0603-db-bound} for completeness. 

\begin{remark}[Limitation and useful regimes]The plug-in correction $\widehat U_k^{\db}$ is not a universal replacement for the raw empirical eigenspace $\widetilde U_k$. The term $\mathfrak D_k$ measures the error in estimating the oracle correction from the empirical lower-positive $\mathcal J$-block. This term could be large in two situations. First, if the lower-positive block $\mathcal J$ contains weak or poorly
separated spikes, then the empirical block
$\widetilde U_{\mathcal J}$ may not estimate $U_{\mathcal J}$ accurately. Second, when eigenvalues are tightly clustered ($\delta_k$ is small),  the matrices $I- U_{\mathcal J}^\T \Phi(\widetilde\lambda_s) U_{\mathcal J}\Lambda_{\mathcal J}$ become ill-conditioned. Inverting them will make $\|\mathfrak D_k\|$ blow up.

Thus the plug-in estimator is useful only when 
\[
\|\mathfrak D_k\| \lesssim \mathcal E_k^{\orc}
\]
and the removed bias is larger than the residual terms
\[
\frac{\sqrt{k}}{\delta_k\lambda_k} \left(
\|\mathcal V_{\mathcal J\mathcal I}\| +
\frac{R_{\max}}{\lambda_k^2}\osc(\mathcal R) \right)
\gg \mathcal E_k^{\orc} +\frac{\|E\|}{\lambda_k}.
\]
This condition does not require the geometric bias to vanish. It may hold, for instance, when $\mathcal J$ contains strong, well-separated spikes that can be accurately estimated, even if the variance profile $\mathcal R$ has a large coupling between $U_{\mathcal J}$ and $U_{\mathcal I}$. A systematic analysis of optimal fully data-driven
bias correction is left for future work.
\end{remark}

\section{Proofs for the applications}
\label{app:network-proofs}

\subsection{Proofs for the DCSBM application}\label{app:dcsbm-proof}
We first record the population geometry. The calculation is standard; see, for example, \cite{LeiRinaldo2015,Jin15}. Let
\[
    S:=Z^\top\Theta^2Z=\diag(s_1,\ldots,s_r)
    \quad \text{and} \quad
    X_0:=\Theta ZS^{-1/2}.
\]
Note that $S$ has rank $r$ and $X_0^\top X_0=I_r$.

\begin{lemma}[Population geometry of the DCSBM]
\label{lem:dcsbm-population-geometry}
Write the eigen-decomposition $S^{1/2}BS^{1/2}=V\Lambda V^\top$ where $V\in \mathbb O(r)$. Then the eigenspace of $P$ may be represented as $U=X_0V$, and
\begin{equation}
    e_i^\top U = \mu_i e_{g(i)}^\top V,
    \qquad
    \frac{e_i^\top U}{\|{e_i^\top U}\|} =  e_{g(i)}^\top V.
    \label{eq:dcsbm-population-rows}
\end{equation}
Consequently,
\[
    \| e_i^\top U\|=\mu_i,
    \qquad
    \|U\|_{2,\infty}=\mu_+,
    \qquad
    \min_i \|e_i^\top U\|_2=\mu_-,
\]
and the normalized population centers of two distinct communities are
separated by distance $\sqrt2$.
\end{lemma}

\begin{proof}
Note that 
\[
    P =\Theta ZBZ^\top\Theta = X_0\bigl(S^{1/2}BS^{1/2}\bigr)X_0^\top
    = (X_0V)\Lambda(X_0V)^\top.
\]
Thus we may take $U=X_0V$. If $g(i)=a$, then
\[
    e_i^\top X_0 = \frac{\theta_i}{\sqrt{s_a}}e_a^\top  = \mu_i e_a^\top,
\]
which proves \eqref{eq:dcsbm-population-rows}. Since the rows of $V$ are
orthonormal, $\|e_a^\top V-e_b^\top V\|=\sqrt2$ for $a\ne b$.
\end{proof}

We use the following specialization of the refined full-signal perturbation bound
from the main results. It is obtained by applying Corollary~\ref{cor:full-signal} to the positive and negative eigenspaces separately (see Remark~\ref{rem:general-space}). Recall that $t_{D,r} = (D+r)\log n.$

\begin{prop}[Perturbation input]
\label{prop:dcsbm-perturbation-input}
Let $P=U\Lambda U^\top$ be a symmetric rank-$r$ probability matrix and let
$ \lambda_\star:=\min_{1\le j\le r} |\lambda_j(P)|.$
For centered Bernoulli noise, suppose $R_{\max}\ge C_Dt_{D,r}$ and $  \lambda_\star\ge C_D\sqrt{t_{D,r}R_{\max}}.$
Then, with probability at least $1-n^{-D}$,
\begin{align*}
    \min_{O\in\mathbb O(r)}\big\|{\widetilde U_rO-U}\big\|_{2,\infty}
    \le C_D\Bigg[ \sqrt r\,
        \frac{  \sqrt{t_{D,r}}\,\Gamma(U) +t_{D,r}\|{U}\|_{2,\infty} }{\lambda_\star}
        +
        (\sqrt r+1)\|{U}\|_{2,\infty}
        \frac{R_{\max}}{\lambda_\star^2} \Bigg].
\end{align*}
Here
$ \Gamma(U)  :=
    \max_i \big( \sum_jP_{ij}(1-P_{ij})\|e_j^\top U\|^2 \big)^{1/2}.
$
\end{prop}
Now we are in the position to prove Theorem~\ref{thm:dcsbm-general}.
\begin{proof}[Proof of Theorem~\ref{thm:dcsbm-general}]
By Lemma~\ref{lem:dcsbm-population-geometry},
$ \|{U}\|_{2,\infty}=\mu_+$ and $\Gamma(U)=\Gamma_{\mathrm{dc}}.$
As given in \eqref{eq:dcsbm-compare}, Assumption~\ref{assum:DCSBM} suggests
\begin{equation}
    \tau_n=\lambda_\star=\min_{1\le j\le r} |\lambda_j(P)| \ge \sigma_{\min}(B) s_-.
 \label{eq:dcsbm-lambda-lower}
\end{equation}
Since $P$ is symmetric and entrywise nonnegative,
\[
    \|P\| \le d_{\max}= \max_i \sum_j P_{ij}.
\]
Moreover, $p_{\max}=o(1)$ implies
\begin{equation}
    R_i
    =
    \sum_jP_{ij}(1-P_{ij})
    \asymp
    d_i
    \quad\text{uniformly in }i,
    \qquad
    R_{\max}\asymp d_{\max},
    \qquad
    \osc(\mathcal R)\le R_{\max}.
    \label{eq:dcsbm-general-variance-comparison}
\end{equation}

For the centered Bernoulli entries and fixed $D$, we have $L_{D,r}^{\rm ra}\lesssim1$ and $M_D\lesssim_Dt_{D,r}$, and therefore (recall $\mathfrak R_D(U)$ in \eqref{eq:signal-row-scale})
\begin{equation}
    \mathfrak R_D(U)
    \lesssim_D \sqrt{t_{D,r}}\,\Gamma_{\mathrm{dc}} +t_{D,r}\mu_+.
    \label{eq:dcsbm-general-RD}
\end{equation}
By  \eqref{eq:dcsbm-general-size}, \eqref{eq:dcsbm-lambda-lower} and
\eqref{eq:dcsbm-general-variance-comparison}, the conditions in Proposition~\ref{prop:dcsbm-perturbation-input} are satisfied. 

Applying Proposition~\ref{prop:dcsbm-perturbation-input} and using
\eqref{eq:dcsbm-lambda-lower}--\eqref{eq:dcsbm-general-RD} yields
\[
    \min_{O\in\mathbb O(r)}\|\widetilde U_rO-U\|_{2,\infty}
    \le C_{2,D}\cdot\varepsilon_{D,n}(\tau_n),
\]
which proves \eqref{eq:dcsbm-general-rowwise}.

It remains to prove exact recovery. Let $O_*$ attain the bound \eqref{eq:dcsbm-general-rowwise} and write
\[
    \widehat u_i:=e_i^\top\widetilde U_rO_*,
    \qquad
    u_i:=e_i^\top U.
\]
If \eqref{eq:dcsbm-general-recovery} holds with $c_D$ sufficiently small,
then every $\widehat u_i$ is nonzero and
\[
    \left\| \frac{\widehat u_i}{\|\widehat u_i\|} -
        \frac{u_i}{\|{u_i}\|}
    \right\|
    \le \frac{2\|\widehat u_i-u_i\|}{\|{u_i}\|} \le \delta_D <\sqrt2/4.
\]
Thus two normalized empirical rows from the same community are at distance at
most $2\delta_D$, whereas rows from distinct communities are at distance at
least $\sqrt2-2\delta_D$.

We now verify the farthest-first step. Suppose the selected centers represent a
proper subset of the communities. Every point from a represented community lies
within $2\delta_D$ of one of the selected centers, while every point from an
unrepresented community has distance at least $\sqrt2-2\delta_D>2\delta_D$
from all selected centers. Hence the next farthest point belongs to an
unrepresented community. Induction shows that the $r$ selected centers contain
exactly one representative from each community. The same distance comparison
shows that nearest-center assignment is exact. A common orthogonal transformation
preserves row norms and Euclidean distances, so the Procrustes alignment does not
affect the algorithm.
\end{proof}

\begin{proof}[Proof of Corollary~\ref{cor:dcsbm-unbalanced}]
Under the assumptions of the corollary,
\begin{equation}
    s_a\asymp n_a\theta_0^2,
    \qquad
    s_-\asymp d_-,
    \qquad
    \mu_+\asymp m_n^{-1/2},
    \qquad
    \mu_-\asymp N_n^{-1/2}
    \label{eq:dcsbm-unbalanced-mu}
\end{equation}
where the constants depend only on $c_\theta$ and $C_\theta$.

Because $r$ is fixed, $B$ is entrywise nonnegative, and its singular values are
bounded above and below, every row and column of $B$ has total mass bounded above
and below by positive constants. Consequently,
\begin{equation}
    d_{\max}\asymp N_n\theta_0^2=\kappa_nd_-.
    \label{eq:dcsbm-unbalanced-dmax}
\end{equation}
Also, $p_{\max}\lesssim\theta_0^2=o(1)$, so
$R_{\max}\asymp d_{\max}$.

We next evaluate the $\Gamma_{\rm dc}$. If $i$ belongs to community
$a$, then Lemma~\ref{lem:dcsbm-population-geometry} gives
\begin{align*}
    \sum_jP_{ij}(1-P_{ij})\|e_j^\top U\|_2^2 &\le
    \theta_i \sum_{b=1}^rB_{ab}  \frac{\sum_{j:g(j)=b}\theta_j^3}{s_b}\lesssim
    \theta_0^2\sum_{b=1}^rB_{ab}
    \lesssim \sigma_{\max}(B)\theta_0^2.
\end{align*}
Thus
\begin{equation}
    \Gamma_{\mathrm{dc}}\lesssim\theta_0\sqrt{\sigma_{\max}(B)}.
    \label{eq:dcsbm-unbalanced-Gamma}
\end{equation}

Since $r$ is fixed, $t_{D,r}\asymp_D\log n$. By \eqref{eq:dcsbm-unbalanced-mu},
\begin{equation}
    \tau_n\ge \sigma_{\min}(B) s_- \asymp \sigma_{\min}(B) d_-.
    \label{eq:dcsbm-unbalanced-tau}
\end{equation}
The condition \eqref{eq:dcsbm-unbalanced-degree} implies the two
conditions in Theorem~\ref{thm:dcsbm-general}. Indeed, using
$\sigma_{\max}(B)\ge \sigma_{\min}(B)$, $\kappa_n\ge1$, and the lower bound in
\eqref{eq:dcsbm-unbalanced-dmax},
\[
    d_{\max} \gtrsim \sigma_{\min}(B)\kappa_n d_- \gtrsim_D \log n.
\]
Using the upper bound in \eqref{eq:dcsbm-unbalanced-dmax},
\[
    \sqrt{d_{\max}\log n} \lesssim \sqrt{\sigma_{\max}(B) \kappa_n d_-\log n}\lesssim_D \sigma_{\min}(B) d_-,
\]
and \eqref{eq:dcsbm-unbalanced-degree} also gives
$\log n\lesssim_D \sigma_{\max}(B) d_-$. Hence
\[
    \tau_n \gtrsim \sigma_{\max}(B)d_- \gtrsim_D \sqrt{d_{\max}\log n}.
\]

Substituting \eqref{eq:dcsbm-unbalanced-mu}, \eqref{eq:dcsbm-unbalanced-dmax},\eqref{eq:dcsbm-unbalanced-Gamma}, and \eqref{eq:dcsbm-unbalanced-tau} into \eqref{eq:dcsbm-intrinsic-error-tau} yields
\begin{align*}
    \varepsilon_{D,n}(\tau_n)
    &\lesssim_D
    \frac{\theta_0\sqrt{\sigma_{\max}(B) \log n}}{\sigma_{\min}(B) d_-}+
    \frac{\log n}{\sigma_{\min}(B)\sqrt{m_n}\,d_-}+\frac{\sigma_{\max}(B) \kappa_n}{\sigma_{\min}(B)^2\sqrt{m_n}\,d_-} \\
    &\lesssim_D
    \frac{1}{\sqrt{N_n}}
    \left[ \frac{\sqrt{\sigma_{\max}(B)}}{\sigma_{\min}(B)}\sqrt{\frac{\kappa_n\log n}{d_-}}
        +\frac{1}{\sigma_{\min}(B)}\frac{\sqrt{\kappa_n}\log n}{d_-}
        +\frac{\sigma_{\max}(B)}{\sigma_{\min}(B)^2}\frac{\kappa_n^{3/2}}{d_-}
    \right].
\end{align*}
By our assumption on $d_-$, we have 
\[\frac{1}{\sigma_{\min}(B)}\frac{\sqrt{\kappa_n}\log n}{d_-} \lesssim \frac{\sqrt{\sigma_{\max}(B)}}{\sigma_{\min}(B)}\sqrt{\frac{\kappa_n\log n}{d_-}}.\]
 Then \eqref{eq:dcsbm-unbalanced-rowwise} follows from Theorem~\ref{thm:dcsbm-general}.

Finally, divide the last display by
$\mu_-\asymp N_n^{-1/2}$. The two terms on the right-hand side of \eqref{eq:dcsbm-unbalanced-rowwise} are made sufficiently small by \eqref{eq:dcsbm-unbalanced-degree}.  After increasing $C_{1,D}$, we therefore have
$\varepsilon_{D,n}(\tau_n) \le c_D\mu_-$. Exact recovery follows from Theorem~\ref{thm:dcsbm-general}.

\end{proof}

\subsection{Proofs for the GRDPG application}
\label{app:grdpg-proof}

Throughout this section, the rank $r=p+q$ is fixed, and all implicit
constants may depend on $D$, $r$, and the constants in
Assumption~\ref{ass:grdpg}, but not on $n$. 

We first record the estimation of related parameters under Assumption~\ref{ass:grdpg}. Recall $d_n = n\rho_n$.

\begin{lemma}
\label{lem:grdpg-scales}
Suppose Assumption~\ref{ass:grdpg} holds. Then $P=U\Lambda U^\T$ has rank $r$, with exactly
$p$ positive and $q$ negative nonzero eigenvalues. Denote $\lambda_{\rm sig} := \min_{1\le a\le r}|\lambda_a(P)|$. 
Then
\begin{equation}
    \lambda_{\rm sig}\asymp d_n,
    \qquad
    \|P\|\asymp d_n,
    \qquad
    \|U\|_{2,\infty}\lesssim n^{-1/2}.
    \label{eq:grdpg-proof-signal-scales}
\end{equation}
The row-variance profile satisfies
\begin{equation}
    R_i\asymp d_n
    \quad\text{for all }i\in [n],
    \qquad
    R_{\max}\asymp d_n,
    \qquad
    \osc(\mathcal R)\lesssim d_n.
    \label{eq:grdpg-proof-variance-scales}
\end{equation}
In addition,
\begin{equation}
    \Gamma(U)^2 := \max_{i\in[n]}\sum_{j=1}^n P_{ij}(1-P_{ij})  \|e_j^\T U\|_2^2 \lesssim \frac{d_n}{n}.
    \label{eq:grdpg-proof-gamma}
\end{equation}
Finally, the matrix $Q_0$ in
\eqref{eq:grdpg-factor-equivalence} may be chosen so that
$\| Q_0\|+\|Q_0^{-1}\|\lesssim 1.$
\end{lemma}

\begin{proof}
For simplicity, set $S:=Y^\T Y=n\Delta_n.$ Since $\Delta_n$ is positive definite,
\[
    P = (YS^{-1/2}) \bigl(\rho_nS^{1/2}JS^{1/2}\bigr) (YS^{-1/2})^\T,
\]
and $YS^{-1/2}$ has orthonormal columns. Hence the nonzero
eigenvalues of $P$ are the eigenvalues of
\[
    \rho_nS^{1/2}JS^{1/2}
    =
    d_n\Delta_n^{1/2}J\Delta_n^{1/2}.
\]

By Sylvester's law of inertia, $\Delta_n^{1/2}J\Delta_n^{1/2}$, and therefore $P$, has exactly
$p$ positive and $q$ negative eigenvalues. Moreover,
\[
    \|\Delta_n^{1/2}J\Delta_n^{1/2}\| \lesssim 1,
    \qquad
    \|(\Delta_n^{1/2}J\Delta_n^{1/2})^{-1}\|
    = \| \Delta_n^{-1/2}J\Delta_n^{-1/2}\| \lesssim 1.
\]
Since $\Delta_n^{1/2}J\Delta_n^{1/2}$ is symmetric, all its eigenvalues have magnitudes bounded above and below by positive constants. This proves the first two relations in \eqref{eq:grdpg-proof-signal-scales}. For the last relation in \eqref{eq:grdpg-proof-signal-scales}, note that if $V\in\mathbb O(r)$ diagonalizes
$d_n\Delta_n^{1/2}J\Delta_n^{1/2}$, then  $U=YS^{-1/2}V.$
Consequently,
\[
    \|e_i^\T U\| \le \|y_i\| \cdot \|S^{-1/2}\| \lesssim n^{-1/2}.
\]

Now we prove \eqref{eq:grdpg-proof-variance-scales}. The expected degree of vertex $i$ is $\sum_{j=1}^n P_{ij} = d_n y_i^\T J\bar y_n  \asymp d_n.$ Since $0\le P_{ij}\le 1-c$,
\[
    c\sum_{j=1}^nP_{ij}
    \le R_i  = \sum_{j=1}^nP_{ij}(1-P_{ij})
    \le \sum_{j=1}^nP_{ij}.
\]
It also gives \eqref{eq:grdpg-proof-gamma}: $\Gamma(U)^2  \le \|{U}\|_{2,\infty}^2R_{\max}\lesssim \frac{d_n}{n}.$

Finally, since  $X$ and
$\sqrt{\rho_n}Y$ have full column rank and the same column space, 
there is an invertible matrix $Q_0$ such that
\[
    XQ_0=\sqrt{\rho_n}Y.
\]
Substituting this identity into $XJX^\T = \rho_nYJY^\T$
and using the full column rank of $X$, we get $Q_0JQ_0^\T=J$,
or equivalently $Q_0^\T JQ_0=J$. Thus $Q_0\in\mathcal O(p,q)$. Furthermore,
\begin{equation*}
    Q_0^\T|\Lambda|Q_0 =\rho_nY^\T Y = d_n\Delta_n.
\end{equation*}
Both $|\Lambda|$ and $d_n\Delta_n$ have eigenvalues comparable to
$d_n$.  This suggests that $\| Q_0\|+\|Q_0^{-1}\|\lesssim 1.$
\end{proof}

We now prove the full-embedding Theorem~\ref{thm:grdpg-full}. 
Let $E:=\widetilde A-P. $ Recall $\mathfrak R_D(U)$ from \eqref{eq:signal-row-scale}. 
In bounded-Bernoulli regime, we have, with probability at least $1-n^{-D}$, that 
\begin{equation}
    \|E\| \lesssim_D\sqrt{d_n},
    \qquad
    M_D\lesssim_D\log n,
    \qquad
    \mathfrak R_D(U)
    \lesssim_D
    \sqrt{\frac{d_n\log n}{n}}
    +
    \frac{\log n}{\sqrt n}.
    \label{eq:grdpg-proof-stochastic-scales}
\end{equation}
Here the last relation follows from
\eqref{eq:grdpg-proof-gamma} and
$\|U\|_{2,\infty}\lesssim n^{-1/2}$.

\begin{proof}[Proof of Theorem~\ref{thm:grdpg-full}]
By Lemma~\ref{lem:grdpg-scales}, every nonzero signal eigenvalue has
magnitude comparable to $d_n$, whereas
$R_{\max}\asymp d_n$. Hence, after choosing $C_{1,D}$ sufficiently
large, the condition $d_n\ge C_{1,D}\log n$ verifies the signal strength condition in Proposition~\ref{prop:refined-weighted-embedding}.

We apply Proposition~\ref{prop:refined-weighted-embedding} separately to the complete positive signal block of $P$ and to the complete negative signal block through $-P$. Combining the two orthogonal alignments yields
\begin{align}
    \min_{W\in\blockO{p}{q}}
    \|\widetilde X_rW-X\|_{2,\infty}
    \lesssim_D \frac{\mathfrak R_D(U)}{\sqrt{d_n}}
    + \|U\|_{2,\infty} \frac{R_{\max}+\osc(\mathcal R)}{d_n^{3/2}}+
    \frac{\|U\|_{2,\infty}}{\sqrt{d_n}} \eta_{\rm full},
    \label{eq:grdpg-proof-full-blackbox}
\end{align}
where $\eta_{\rm full}:=\eta_++\eta_-$ with
\[
    \eta_\pm := \|U_\pm^\T E U_\pm\|+
    \bigl(\|\Lambda_\mp\|+\|E\|\bigr)  \|\sin\angle(U_\pm,\widetilde U_\pm)\|^2.
\]

By Davis-Kahan theorem, 
\[
       \| \sin\angle(U_\pm,\widetilde U_\pm)\|  \lesssim_D \frac{\|E\|}{d_n} \lesssim_D d_n^{-1/2}.
\]
Consequently,
\[
    \eta_{\rm full} \lesssim_D\|E\|+ d_n \max_{\varsigma\in\{+,-\}}
        \|\sin\angle(U_\varsigma,\widetilde U_\varsigma)\|^2
    \lesssim_D \sqrt{d_n}.
\]
Substituting this estimate together with
\eqref{eq:grdpg-proof-signal-scales} and \eqref{eq:grdpg-proof-stochastic-scales}
into \eqref{eq:grdpg-proof-full-blackbox} yields
\[
    \min_{W\in\blockO{p}{q}} \|\widetilde X_rW-X\|_{2,\infty}
    \lesssim_D
    \sqrt{\frac{\log n}{n}} + \frac{\log n}{\sqrt{nd_n}} + \frac1{\sqrt n}
    +\frac1{\sqrt{nd_n}}
    \lesssim_D \sqrt{\frac{\log n}{n}},
\]
where we used $d_n\ge C_{1,D}\log n$ to get the last inequality. This proves
\eqref{eq:grdpg-full-spectral}.

Moreover, by Weyl's inequality, the $r$ eigenvalues of $\widetilde A = P+E$ of
largest magnitude are the signal outliers, with $p$ positive and $q$
negative eigenvalues. Indeed, the $r$ nonzero eigenvalues of $P$ have magnitudes of order $d_n$, whereas $\|E\|\lesssim_D\sqrt{d_n}$. 

It remains to pass from  $X$ to the latent positions $Y$. Let $W\in\blockO{p}{q}$ attaining \eqref{eq:grdpg-full-spectral}, and set $Q:=WQ_0.$ Both $W$ and $Q_0$ preserve $J$, so $Q\in\mathbb O(p,q)$. By $XQ_0=\sqrt{\rho_n}\,Y$ from \eqref{eq:grdpg-factor-equivalence},
\[
    \widetilde X_rQ-\sqrt{\rho_n}\,Y= (\widetilde X_rW-X)Q_0.
\]
Hence, by \eqref{eq:grdpg-full-spectral} and Lemma~\ref{lem:grdpg-scales},
\[
\| \widetilde X_rQ-\sqrt{\rho_n}\,Y\|_{2,\infty}\lesssim_D \sqrt{\frac{\log n}{n}}.
\]
Dividing by $\sqrt{\rho_n}=\sqrt{\frac{d_n}{n}}$ proves \eqref{eq:grdpg-full-latent}. The consistency statement follows immediately.
\end{proof}

Finally, we provide the proof of Proposition~\ref{prop:grdpg-truncated}.
\begin{proof}[Proof of Proposition~\ref{prop:grdpg-truncated}]
 By Lemma~\ref{lem:grdpg-scales} and Weyl's inequality, for $\varsigma \in \{+,-\}$ and every $a\le k_{\varsigma}$, we have $\widetilde \lambda_a^\varsigma \asymp \sqrt{d_n}$ and thus,
\[
    \tau_\varsigma :=
    \big(  \sum_{a\le k_\varsigma} \widetilde\lambda_a^\varsigma \big)^{1/2}
    \asymp \sqrt{d_n}.
\]
The assumptions $d_n\ge C_{1,D}\log n$ and $\delta_n\ge C_{1,D}\log n$ imply the signal strength and gap conditions in the corresponding signed versions of Proposition~\ref{prop:refined-weighted-embedding}.

We directly apply Proposition~\ref{prop:refined-weighted-embedding} to the retained positive block of $P$ and to the retained negative block through $-P$, and combine the two blockwise alignments. Note that the deterministic terms are precisely
$\GBias_{k_+,k_-}$ defined in \eqref{eq:grdpg-GBias}. Thus
\begin{align}
    \min_{W\in\blockO{k_+}{k_-}}
    \| \widetilde X_kW-X_k\|_{2,\infty}
    \le{}&
    \GBias_{k_+,k_-}
    +
    C_D\Bigg[ \frac{\mathfrak R_D(U)}{\sqrt{d_n}}+
        \sqrt{\frac{d_n}{n}}\, \log n
        \left(   \frac1{\delta_{+,n}} + \frac1{\delta_{-,n}} \right)
        \notag\\
        &\hspace{25mm}
        +
        \frac1{\sqrt{nd_n}}
        +
        \sum_{\varsigma\in\{+,-\}}
        \frac{
            \|U_{\varsigma,k_\varsigma}\|_{2,\infty} }{\sqrt{d_n}}
        \eta_\varsigma
    \Bigg],
    \label{eq:grdpg-proof-truncated-intermediate}
\end{align}
where, for each $\varsigma\in\{+,-\}$,
\begin{equation}
    \eta_\varsigma \lesssim_D
    \|U_{\varsigma,k_\varsigma}^\T E U_{\varsigma,k_\varsigma}\|
    + d_n \left( \mathcal B_{k_\varsigma}^{\varsigma}
        + \frac{\sqrt{k_\varsigma}\log n}{\delta_{\varsigma,n}}
        + \frac1{\sqrt{d_n}}\right)^2.
    \label{eq:grdpg-eta-signed}
\end{equation}
It remains to bound the bias quantities and the last term in
\eqref{eq:grdpg-proof-truncated-intermediate}.

We first note that if $A$ and $B$ have orthonormal columns and
$A^\T B=0$, then
\[
    \|A^\T\mathcal R B\| = \|A^\T(\mathcal R-c_0 I)B\|
    \le \frac12\osc(\mathcal R),
    \qquad
    c_0:=\frac{R_{\max}+R_{\min}}2.
\]
Hence, by \eqref{eq:grdpg-proof-variance-scales}, all variance-profile couplings appearing in
$\mathcal B_{k_\varsigma}^{\varsigma}$ (defined in \eqref{def:bias-k}) are $O(d_n)$ and thus 
\begin{equation}
    \mathcal B_{k_\varsigma}^{\varsigma}
    \lesssim \frac1{\delta_{\varsigma,n}} +  \frac1{d_n} \le \frac{1}{\delta_n}+  \frac1{d_n}.
    \label{eq:grdpg-proof-bias-scale}
\end{equation}
Since $r$ is fixed and $\delta_{\varsigma,n}\ge\delta_n$, \eqref{eq:grdpg-eta-signed} and \eqref{eq:grdpg-proof-bias-scale} give 
\begin{equation*} 
\eta_\varsigma \lesssim_D 
\|U_{\varsigma,k_\varsigma}^\T E U_{\varsigma,k_\varsigma}\| + d_n\left( \frac{\log n}{\delta_n} + \frac1{\sqrt{d_n}} \right)^2 \le \|E\| + d_n\left( \frac{\log n}{\delta_n} + \frac1{\sqrt{d_n}} \right)^2. 
\end{equation*} 
Here we used $d_n\gtrsim_D\log n$ to absorb the $d_n^{-1}$ term. Moreover, using $\|E\| \lesssim_D\sqrt{d_n} $ and $\|U_{\varsigma,k_\varsigma}\|_{2,\infty}\lesssim n^{-1/2}$, we get 
\begin{align*} 
\frac{\|U_{\varsigma,k_\varsigma}\|_{2,\infty}} {\sqrt{d_n}}\eta_\varsigma 
\lesssim_D \frac1{\sqrt n} + \sqrt{\frac{d_n}{n}} \left( \frac{\log n}{\delta_n} + \frac1{\sqrt{d_n}} \right)^2 
\lesssim_D \sqrt{\frac{\log n}{n}} + \sqrt{\frac{d_n}{n}}\, \frac{\log n}{\delta_n}.
\end{align*} 
Indeed, the assumption $\delta_n\ge C_{1,D}\log n$ implies $\log n/\delta_n\lesssim_D1$, while $d_n\gtrsim_D\log n$. Also, note that  $ \frac1{\delta_{+,n}}+\frac1{\delta_{-,n}} \le \frac{2}{\delta_n}$ and, by \eqref{eq:grdpg-proof-stochastic-scales}, 
\[ 
\frac{\mathfrak R_D(U)}{\sqrt{d_n}} + \frac1{\sqrt{nd_n}} \lesssim_D \sqrt{\frac{\log n}{n}}. 
\] 
Substituting these estimates into \eqref{eq:grdpg-proof-truncated-intermediate} proves \eqref{eq:grdpg-truncated-main}. 

\end{proof}

We wrap up this section with a short derivation of \eqref{eq:grdpg-truncated-macroscopic}.
\begin{proof}[Proof of \eqref{eq:grdpg-truncated-macroscopic}] Suppose $\delta_n\ge c d_n$. Then \eqref{eq:grdpg-proof-bias-scale} gives $ \mathcal B_{k_\varsigma}^{\varsigma} \lesssim \frac1{d_n}. $ Since $\lambda_{\varsigma,*}\asymp d_n$ and $\|U\|_{2,\infty}\lesssim n^{-1/2}$, \eqref{eq:grdpg-GBias} yields 
\begin{equation*} 
\GBias_{k_+,k_-} \lesssim \frac1{\sqrt{nd_n}}. 
\end{equation*} 
Moreover, \[ \sqrt{\frac{d_n}{n}}\, \frac{\log n}{\delta_n} \lesssim \frac{\log n}{\sqrt{nd_n}} \lesssim_D \sqrt{\frac{\log n}{n}}, \] where the last inequality is from $d_n\gtrsim_D\log n$. Substituting these estimates into \eqref{eq:grdpg-truncated-main} proves \eqref{eq:grdpg-truncated-macroscopic}. 
\end{proof}

\paragraph{Acknowledgment.} Ke Wang is supported by Hong Kong RGC 16304222 and 16305526.

\bibliography{NewPerturbation}
\bibliographystyle{abbrv}

\end{document}